\documentclass[11pt]{article}
\usepackage{pifont}
\usepackage{amsfonts,amsmath,amssymb,amsthm}
\usepackage{mathrsfs, amscd}
\usepackage{amsbsy}
\usepackage{xy}
\xyoption{all}
\usepackage{color}

\newcommand{\mput}{\multiput}
\newcommand{\bcen}{\begin{center}}
\newcommand{\ecen}{\end{center}}

\def\dz{\delta}
\def\lz{\lambda}

\def\A{{\cal A}}

\def\C{{\cal C}}

\def\F{{\cal F}}
\def\K{{\cal K}}
\def\U{{\cal U}}
\def\V{{\cal V}}
\def\M{{\cal M}}

\def\N{{\cal N}}
\def\T{{{\cal T}}}

\def\f{\textsf{f}}
\def\tl{\mbox l}
\def\tr{\mbox r}
\def\tm{\mbox m}
\def\td{\mbox d}
\def\te{\mbox e}

\def\l{\underline{l}}
\def\n{{\underline{n}}}
\def\m{{\underline{m}}}

\def\tt{{t\times t}}

\def\rad{\mbox{rad}}

\def\wt{\widetilde}

\def\Hom{\mbox{Hom}}

\def\Hom{\mbox{Hom}}

\def\dim{\mbox{dim}}
\def\mod{\mbox{mod}}

\def\Im{\mbox{Im}}

\def\rank{\mbox{rank}}
\def\mf{\mathfrak}

\def\mf{\mathfrak}

\def\IM{\hbox{\rm I \hskip -4.9pt M}}
\def\IE{\hbox{\rm I\hskip -2pt E}}
\def\IF{\hbox{\rm I\hskip -2pt F}}
\def\id{\hbox{\rm 1\hskip -2.2pt l}}

\begin{document}
\begin{center}{\huge\bf Algebras with homogeneous module category are tame} \end{center}

\begin{center} Zhang Yingbo$^1$ and Xu Yunge$^2$ \end{center}

\vspace{-3mm}\subitem {\footnotesize {1} School of Mathematis, Beijing Normal University, 100875 P.R.China, zhangyb@bnu.edu.cn }
\subitem {\footnotesize  {2} Faculty of Mathematics,  Hubei University, 430062 P.R.China, xuy@hubu.edu.cn}

\begin{abstract}

 The celebrated  Drozd's theorem asserts that
a finite-dimensional basic algebra  $\Lambda$  over an
algebraically closed field $k$   is either tame or wild,  whereas the Crawley-Boevey's theorem
states that   given a  tame algebra $\Lambda$  and a
dimension $d$,   all but finitely many isomorphism classes of indecomposable
$\Lambda$-modules of dimension $d$ are isomorphic to their
Auslander-Reiten translations and hence belong to homogeneous tubes.
In this paper, we prove the inverse of Crawley-Boevey's theorem, which
gives an internal description of tameness  in terms of
Auslander-Reiten quivers.

\end{abstract}

\renewcommand\thefootnote{}
\footnote{{\it 2010 Mathematics Subject Classification: 15A21, 16G20, 16G60,
16G70 }}
\footnote{The authors are supported by the National Natural Science Foundation of China,
No. 11271318, and 11371186.}
\vspace{-3mm}

\centerline{\Large{\bf Contents}}

Introduction

1. Matrix bimodule problems

\quad 1.1 Definition of matrix bimodule problems

\quad 1.2 Bi-comodule problems and bocses

\quad 1.3 Representation categories of matrix bimodule problems

\quad 1.4 Formal products and formal equations

2. Reductions on matrix bimodule problems

\quad 2.1 Admissible bimodules and induced matrix bimodule problems

\quad 2.2 Eight Reductions

\quad 2.3 Canonical forms

\quad 2.4 Defining systems

3. Classification of minimal wild bocses

\quad 3.1 An exact structure on representation categories of bocses

\quad 3.2 Almost split conflations in the process of reductions

\quad 3.3 Minimal wild bocses

\quad 3.4 Non-homogeneity in the cases of MW1-MW4

4. One-sided pairs

\quad 4.1 Definition of one-sided pairs

\quad 4.2 Differentials in one-sided pairs

\quad 4.3 Reduction sequences of one-sided pairs

\quad 4.4 Major pairs

\quad 4.5 Further reductions

\quad 4.6 Regularizations on non-effective $a$-class and all $b$-class arrows

5. Non-homogeneity of bipartite matrix bimodule problems of wild type

\quad 5.1 An inspiring example

\quad 5.2 Bordered matrices in bipartite case

\quad 5.3 Non-homogeneity in the case of MW5 and classification (I)

\quad 5.4 Bordered matrices in one-sided case

\quad 5.5 Non-homogeneity in the case of MW5 and classification (II)

\quad 5.6 Proof of Main theorem

References

\bigskip
\bigskip
\bigskip
\centerline{\Large\bf Introduction}
\bigskip

Throughout the paper,   we always assume that $k$ is an algebraically
closed field, and that  all  rings or algebras contain identities. We write
our maps  either on  the left or  on the right,
but always compose them as if they were written on the right.

We start with the following important definition of ``tame" and ``wild":

\medskip

{\bf Definition 1} \cite{D1, CB1, DS}\,  A finite-dimensional
$k$-algebra $\Lambda$ is of tame representation type, if for any
positive integer $d$, there are a finite number of localizations
$R_i=k[x,\phi_i(x)^{-1}]$ of $k[x]$ and $\Lambda$-$R_i$-bimodules
$T_i$ which are free  as right $R_i$-modules, such that all but finitely many iso-classes of indecomposable
$\Lambda$-modules of dimension at most $d$ are isomorphic to
$$T_i\otimes_{R_i}R_i/(x-\lz)^{m},$$
for some  $\lz\in k$ with $\phi_i(\lz)\ne 0$, and some positive integer $m$.

A finite-dimensional
$k$-algebra $\Lambda$ is  of wild representation type if there is a
finitely generated $\Lambda$-$k\langle x,y\rangle$-bimodule $T$,
which is free as a right $k\langle x,y\rangle$-module, such that the
functor
$$T\otimes_{k\langle x,y\rangle}-: k\langle
x,y\rangle\mbox{-mod}\rightarrow \Lambda\mbox{-mod}$$ preserves
indecomposability and isomorphism classes.

\medskip

Several authors worked on equivalent definitions of ``time" and ``wild",
for example in terms of generic modules \cite{CB3}.

In 1977 Drozd \cite{D1}  showed  that  a finite-dimensional algebra
over an algebraically closed field  is either of tame representation
type or of wild representation type.  This result  is known as Drozd's Tame-Wild Theorem, and
has been one of
the most fundamental results in the  representation theory of finite dimensional algebras.
On the other hand, however,
the proof of Drozd's Theorem is highly indirect.   Indeed,  the
argument relies on the notion of a bocs (the abbreviation for ``bimodule
of coalgebra structure"), introduced first by Rojter in \cite{Ro}.
In 1988, Crawley-Boevey \cite{CB1} formalized the theory of bocses and showed
that for a tame algebra $\Lambda$, and for each dimension $d$, all
but finitely many isomorphism classes of indecomposable
$\Lambda$-modules of dimension $d$ are isomorphic to their
Auslander-Reiten translations and hence belong to homogeneous tubes.

After the work \cite{CB1}, many authors   tried to prove the
converse of Crawley-Boevey's theorem,  aiming  to find
infinitely many non-isomorphic indecomposable representations
$\{M_i\mid i\in I\}$ of the same dimension in the representation
category of a layered bocs, such that $M_i\not\simeq DTr(M_i)$.   Somewhat  surprisingly,
 in 2000  the authors   \cite{BCLZ}  constructed
a  strongly homogeneous wild layered bocs
$\mf{B}$  for which  each representation
 is homogeneous (i.e., $DTr(M)\simeq M$),  and showed  that the
converse of the Crawley-Boevey's theorem does not hold true for general
layered bocses. Later on, C.M.Ringel proposed a concept of controlled wild, Y.Han described some classes of controlled
wild algebras \cite{H}, and H.Nagasy proved that a $\tau$-wild algebra is wild \cite{N}.
However, the converse of the Crawley-Boevey's theorem remains
open in the case of finite dimensional
$k$-algebras. Our main result in this paper, Theorem 3 below,  gives a full answer to
this problem.  To state our result, we need the following definition:

\medskip

{\bf Definition 2} \cite[2.1 Definitions]{BCLZ}
Let $\Lambda$ be  a finite-dimensional algebra over an algebraically
closed field. An indecomposable $\Lambda$-module $M$ is called
{\it homogeneous}, if $DTr(M)\cong M$. The category mod-$\Lambda$ is
said to be {\it homogeneous}, if for each dimension $d$
all but finitely many isomorphism classes of indecomposable
$\Lambda$-modules of dimension $d$ are homogeneous.

\medskip

We will prove the following main theorem throughout this whole article.

\medskip

{\bf Theorem 3} Let $\Lambda$ be  a finite-dimensional
algebra over an algebraically closed field. Then $\Lambda$ is of
tame representation type if and only if mod-$\Lambda$ is
homogeneous.\hfill$\Box$

\medskip

The necessity of Theorem 3  was previously proved by Crawley-Boevey \cite{CB1}.
We only need to prove the sufficiency part. The proof is
divided into five sections as shown in the contents.
Our proof relies on the notions of matrix bimodule
problems, their associated bocses, and  reduction
techniques. Since the matrix bimodule problems
associated to finite-dimensional algebras are bipartite, the key of our
argument is to find a full subcategory of
representation category of a bipartite matrix bimodule problem
which admits infinitely many isomorphism classes of non-homogeneous representations of dimension $d$.

\bigskip
\bigskip
\bigskip

\centerline{\Large\bf 1 Matrix bimodule problems}

\bigskip

In this section, a notion of matrix bimodule problems over a
minimal algebra is introduced, which is a generalization of bimodule problems
over a field $k$ defined by \cite{CB2} in terms of matrix. Then the associated bi-comodule problems
and bocses of matrix bimodule problems are discussed. Finally, a nice connection between a matrix bimodule problem
and its associated bocs is builded via the formal products of two structures.

\bigskip
\bigskip

\noindent{\large\bf 1.1 Definition of matrix bimodule problems}

\bigskip

The purpose of this subsection is two folds: 1) construct a $k$-algebra $\Delta$ based on
a minimal algebra $R$; 2) define matrix bimodule problems over
$\Delta$. The concepts and the results are proposed by S. Liu.

\medskip

Let ${\cal T}=\T_0\dot\cup\T_1$ be a vertex set, where
the subset $\T_0$ consists of trivial vertices, such that
$\forall\, X\in\T_0$, there is a $k$-algebra $R_X\simeq k$ with the identity $1_X$;
and $\T_1$ consists of non-trivial vertices, such that
$\forall\, X\in\T_1$, there is an algebra $R_X\simeq k[x, \phi_{_X}(x)^{-1}]$ with the identity $1_X$, the
finite localization of the polynomial ring $k[x]$ given by a non-zero polynomial
$\phi_{_X}(x)\in k[x]$, and $x$ is said to be the {\it
parameter associated to $X\in {\cal T}_1$}. Now we call the
$k$-algebra $R=\Pi_{X\in {\cal T}}R_X$ a {\it minimal algebra} over
$\cal T$ with a set of orthogonal primitive idempotents $\{1_{_X}\mid X\in\T\}$.

\medskip

We define a tensor product of $p\geqslant 1$ copies of $R$ over $k$ as follows:
$$\begin{array}{c}
R^{\otimes p}=R\otimes_k\cdots\otimes_k R=\sum_{(X_1,\cdots,X_p)
\in\T\times\cdots\times\T}R_{X_1}
\otimes_k\cdots\otimes_kR_{X_{p}}.\end{array}$$
There exists a natural left and right $R$-module structure on $R^{\otimes p}$:
$$\begin{array}{c}s\otimes_{_R}\alpha=(s\otimes_{_R}r_1)\otimes_kr_2\otimes_k\cdots\otimes_kr_p;\\
\alpha\otimes_{_R}s=r_1\otimes_k\cdots\otimes_kr_{p-1}\otimes_k(r_p\otimes_{_R}s),\end{array}\eqno{(1.1\mbox{-}1)}$$
for any $\alpha=r_1\otimes_k r_2\otimes_k\cdots\otimes_k
r_{p-1}\otimes_kr_p\in R^{\otimes p},s\in R$. If $r_i\in R_{X_i}, s\in R_Y$, then $s\otimes_{_{R}}\alpha=0$ for $Y\ne X_1$ and
$\alpha\otimes_{_R}s=0$ for $Y\ne X_p$. Thus $R^{\otimes p}$ can be viewed as an
{\it $R$-$R$-bimodule}, or simply an $R^{\otimes 2}$-module, with the module action for any $r,s\in R$:
$$\begin{array}{c}(r\otimes_ks)\otimes_{_{R^{\otimes 2}}}\alpha=r\otimes_{_R}\alpha\otimes_{_R}s=
\alpha\otimes_{_{R^{\otimes 2}}}(r\otimes_ks).
\end{array}\eqno{(1.1\mbox{-}2)}$$
Note that $\big(\alpha\otimes_{_{R^{\otimes 2}}}(r\otimes_{k}s)\big)\otimes_{_R}s'
=\alpha\otimes_{_{R^{\otimes 2}}}(r\otimes_kss'), \forall\, s'\in R$, and
$r'\otimes_{_R}\big(\alpha\otimes_{_{R^{\otimes 2}}}(r\otimes_ks)\big)=\alpha\otimes_{_{R^{\otimes 2}}}(r'r\otimes_ks),
\forall\, r'\in R$. The direct sum of $R^{\otimes p}$ for $p=1,2,\cdots$, is still an $R^{\otimes 2}$-module:
$$\begin{array}{c}\Delta=\oplus_{p=1}^\infty R^{\otimes p},\quad \mbox{let}\,\,
\bar\Delta=\oplus_{p=2}^\infty R^{\otimes p},\, \Delta=R\oplus\bar\Delta.\end{array}\eqno{(1.1\mbox{-}3)}$$
We define a multiplication on $R^{\otimes 2}$-module $\Delta$, given by $\Delta\times\Delta
\rightarrow\Delta\otimes_{_R}\Delta\subseteq\Delta$:
$$\begin{array}{c}\Delta^{\otimes p}\otimes_{_R}\Delta^{\otimes q}
\subseteq\Delta^{\otimes (p+q-1)},\quad\alpha\otimes_{_R}\beta=r_1\otimes_k
\cdots\otimes_k(r_ps_1)\otimes_ks_2\cdots\otimes_ks_q,
\end{array}\eqno{(1.1\mbox{-}4)}$$
where $\beta=s_1\otimes_k\cdots\otimes_ks_q$. And if $r_i\in R_{X_i},s_j\in R_{Y_j}$,
$\alpha\otimes_{_R}\beta=0$ for $X_p\ne Y_1$. Thus we obtain {\it an associative non-commutative
$k$-algebra $(\Delta,\otimes_{_R},1_R)$} with the set of orthogonal primitive idempotents $\{1_X\mid X\in\T\}$.
Moreover,  $\Delta\otimes_{_R}\Delta$ can be viewed as an $R^{\otimes 3}$-module:
for any $\alpha,\beta\in\Delta, r,s,w\in R$,
$$\begin{array}{c}(r\otimes_ks\otimes_kw)
\otimes_{_{R^{\otimes 3}}}(\alpha\otimes_{_R}\beta)=(\alpha\otimes_k\beta)
\otimes_{_{R^{\otimes 3}}}(r\otimes_ks\otimes_kw)\\
=r\otimes_{_R}\alpha\otimes_{_R}s\otimes_{_R}\beta\otimes_{_R}w.
\end{array}\eqno{(1.1\mbox{-}5)}$$

Denote by $\IM_{m\times n}(\Delta)$ the set of
{\it matrices over $\Delta$} of size $m\times n$; and
by $\mathbb T_n(\Delta), \mathbb N_n(\Delta)$, $\mathbb D_n(\Delta)$ the sets  of
upper triangular, strictly upper triangular, and diagonal $\Delta$-matrices of size $n\times n$ respectively.
The product of two $\Delta$-matrices is the {\it usual matrix product}.
If $H=(h_{ij})\in$IM$_{m\times n}(R)$, $U=(u_{ij})\in\,\IM_{m\times n}(R\otimes_kR)$, $\alpha\in\Delta$,
define
$$\begin{array}{c} H\otimes_{_R}\alpha=(h_{ij}\otimes_{_R}\alpha)\in\IM_{m\times n}(\Delta),\\
\alpha\otimes_{_R}H=(\alpha\otimes_{_R} h_{ij})\in\IM_{m\times n}(\Delta);\\[1.5mm]
U\otimes_{_{R^{\otimes 2}}}\alpha=\big(\alpha\otimes_{_{R^{\otimes 2}}}u_{ij}
\big)=\alpha\otimes_{_{R^{\otimes 2}}}U\in\IM_{m\times n}(\Delta).
\end{array}\eqno{(1.1\mbox{-}6)}$$
The first two are based on Formula (1.1-1). For the last one,
note that $(\alpha\otimes_{_{R^{\otimes 2}}}U)\otimes_{_R}H=\alpha\otimes_{_{R^{\otimes 2}}}(UH),
H\otimes_{_R}(\alpha\otimes_{_{R^{\otimes 2}}}U)=\alpha\otimes_{_{R^{\otimes 2}}}(HU)$ by the note stated under Formula (1.1-2).
Let $U=(u_{ij})\in$IM$_{m\times n}(R),V=(v_{jl})\in$IM$_{n\times r}(R\otimes_kR)$
and $\alpha,\beta\in\Delta$.  Then
$$\begin{array}{c}\begin{array}{ll}(\alpha\otimes_{_{R^{\otimes 2}}}U)(\beta\otimes_{_{R^{\otimes 2}}}V)
&=(\sum_{j=1}^n\big((\alpha\otimes_{_{R^{\otimes 2}}}u_{ij})\otimes_{_R}
(\beta\otimes_{_{R^{\otimes 2}}}v_{jl})\big)_{i,l}\\
&=\big(\sum_{j=1}^n(\alpha\otimes_{_R}\beta)\otimes_{_{R^{\otimes 3}}}
(u_{ij}\otimes_{_R}v_{jl})\big)_{i,l}\\
&=(\alpha\otimes_{_R}\beta)\otimes_{_{R^{\otimes 3}}}(UV).\end{array}\\
\end{array}\eqno{(1.1\mbox{-}7)}$$

An $R$-$R$-bimodule
${\mathcal S}_1$ is said to be a {\it quasi-free bimodule} finitely generated by
$U_1,\cdots,U_m$, provided  that the morphism
$$(R_{X_1}\otimes_kR_{Y_1})\oplus \cdots \oplus (R_{X_m}\otimes_kR_{Y_m})\rightarrow
\mathcal S_1,\quad 1_{X_i}\otimes_k 1_{Y_i}\mapsto U_i$$
is an isomorphism. In this case, $\{U_1,\ldots,U_m\}$ is called an
{\it $R$-$R$-quasi-free basis} of $\mathcal S_1$, or R-R-quasi-basis of $\mathcal S_1$ for short.\vspace{1mm}

Let $\mathcal S_{p}=R^{\otimes (p+1)}\otimes_{R^{\otimes 2}}\mathcal S_1$, which possesses an $R$-$R$-bimodule
structure:
$$\begin{array}{c}(r\otimes_ks)\otimes_{R^{\otimes 2}}(\alpha\otimes_{R^{\otimes 2}}U)=(r\alpha s)\otimes_{R^{\otimes 2}}U
=\alpha\otimes_{R^{\otimes 2}}((r\otimes_ks)\otimes_{R^{\otimes 2}}U)\end{array}\eqno{(1.1\mbox{-}8)}$$
for $r,s\in R,\alpha\in R^{p+1},U\in\mathcal S$.
Thus $\mathcal S=\sum_{p=1}^\infty\mathcal S_p=\bar\Delta\otimes_{R^{\otimes 2}}\mathcal S_1$
is an $R$-$R$-bimodule, and $\mathcal S_p$ is said to have {\it index $p$} in $\mathcal S$.

\medskip

{\bf Definition 1.1.1}\, Let $T=\{1,2,\cdots,t\}$ be a set of integers,
and let $\sim$ be an equivalent relation on
$T$, such that the set $T/\sim$ is one-to-one corresponding
to the vertex set $\T$ of a minimal algebra $R$. It may be written as $\T=T/\sim$.

\medskip

{\bf Definition 1.1.2}\, (i)\, Define an $R$-$R$-bimodule:
$$\begin{array}{c}\K_0=\{$diag$(s_{11},\cdots,s_{tt})\mid\, \mbox{when}\,\,\,i\in X, s_{ii}\in R_X, $ and $s_{ii}=s_{jj},
\forall\, i\sim j\}.\end{array}$$
Let $E_X\in\mathbb D_{t}(R_X)$ with the entry $s_{ii}=1_X$ if $i\in X$ and $s_{ii}=0$ if $i\notin X$, then
$\{E_X\mid X\in\T\}$ is an $R$-quasi-basis of $\K_0$, and $E=\sum_{X\in\T}E_X$ is the identity matrix of size $t$.

(ii)\, Define a quasi-free $R$-$R$-bimodule
$\K_1\subseteq\mathbb N_t(R\otimes_k R)$ with an $R$-$R$-quasi-basis:
$$\begin{array}{c}\V=\cup_{(X,Y)\in\T\times\T}\V_{XY}=\{V_1,V_2,\cdots,V_m\},
\quad \V_{XY}\subset\mathbb N_t(R_{X}\otimes_k R_{Y}),\end{array}$$
where $V\in\V_{XY}$ if and only if $1_XV1_Y=V$.

(iii)\, Suppose $\K=\K_0\oplus(\bar\Delta\otimes_{_{R^{\otimes 2}}}\K_1)$ possesses an algebra structure, where
the multiplication $\tm: \K\times\K\rightarrow \K$ is the usual
matrix product over $\Delta$ consisting of $\tm_{pq}:\K_p\times\K_q\rightarrow \K_{p+q},\forall\, p,q\geqslant 0$;
the unit is given by the canonical inclusion $\te: R\cong \K_0\hookrightarrow\K$. $\K$ is said to be {\it finitely
generated in index $(0,1)$ over $\Delta$}.

\medskip

Clearly $\K_0\simeq R$ as algebras. The multiplication $\tm_{pq},\forall\, p,q>0$ is determined by
$\tm_{11}:\K_1\times\K_1\rightarrow\K_2$, see Formula
(1.1-7). And for any $r\in R,\alpha\in\Delta, (r\otimes_{_R}E_X)\otimes_{_R}(\alpha\otimes_{_{R^{\otimes 2}}}V_j)
=(r\otimes_{_R}\alpha)\otimes_{_{R^{\otimes 2}}}V_j$ if $1_XV_j=V_j$, or $0$ otherwise,
$(\alpha\otimes_{_{R^{\otimes 2}}}V_j)\otimes_{_R}(r\otimes_{_R}E_X)$ is similar. $\{E_X\mid X\in\T\}$ is a set of orthogonal
primitive idempotents of $\K$, and $E=\te(1_R)$  is the {\it identity}.

\medskip

Let $ T=\{1,2,\ldots,t\}$ and $ T'=\{1,2,\ldots,t'\}$ be two sets
of integers. An {\it order} on $ T\times T'$ is defined as follows:
$(i,j)\preccurlyeq (i',j')$ provided that $i>i'$, or $i=i'$ but
$j\leqslant j'$. Thus an order on the index set of the entries of a matrix in
$\IM_{t\times t'}(\Delta)$ is obtained.
Let $M=(\lambda_{ij})\in \IM_{t\times t'}(\Delta)$. The entry $\lambda_{pq}$ is
said to be the {\it leading entry} of
$M$ if $\lambda_{pq}\ne 0$, and any $\lambda_{ij}\ne 0$ implies that
$(p,q)\preccurlyeq (i,j)$. Let $\bar{M}=(C_{ij})$ be a partitioned
matrix over $\Delta$, one defines similarly the {\it leading
block} of $\bar{M}$.  In both cases, the index $(p,q)$ is
called the {\it leading position} of $M$ resp. $\bar{M}$.

Let ${\cal S}$ be a subspace of
$\IM_t(k)$. An ordered basis $\U=\{U_1,\cdots,U_r\}$ with the
leading positions $(p_1,q_1), \ldots,(p_r,q_r)$ respectively is
called a {\it normalized basis} of ${\cal S}$ provided that
\subitem(i)\, the leading entry of $U_i$ is $1$;

\subitem(ii)\, the $(p_i,q_i)$-th entry of $U_j$ is $0$ for $j\ne i$;

\subitem(iii)\, $U_i\preccurlyeq U_j$ if and only if $(p_i,q_i)\preccurlyeq (p_j,q_j)$

\noindent The basis  $\U$ is a linearly  ordered set.
It is easy to see that ${\cal S}$ has a normalized basis by Linear algebra.
In fact, if $t^2$ variables $x_{ij}$ under the order of matrix indices defined as above are taken, then ${\cal
S}$ will be the solution space of some  system of linear equations
$\sum_{(i,j)\in\T\times\T} a_{ij}^lx_{ij}=0, a^l_{ij}\in k,1\leqslant l\leqslant s$
for some positive integer $s$. Reducing the coefficient matrix to the simplest echelon
form, we assume that $x_{p_1q_1},\cdots,x_{p_rq_r}$
are all the free variables, and $\{U_1,\cdots,U_r\}$ is a basic system of solutions,
whose $(p_i,q_i)$-entry is $1$ and $(p,q)\prec (p_i,q_i)$-entry is $0$ for $i=1,\cdots,r$,
a normalized basis $\U$ of ${\cal S}$ is obtained.

\medskip

{\bf Definition 1.1.3}\, (i)\, Define a quasi-free $R$-$R$-bimodule
$\M_1\subseteq\IM_t(R\otimes_k R)$, such that $E_X\M_1E_Y$ has a
normalized quasi-basis $\A_{XY}\subseteq\IM_t(k1_X\otimes_k1_Yk)\simeq\IM_t(k)$ as
$k$-vector spaces, where $A\in\A_{XY}$
if and only if $1_XA1_Y=A$. Thus there is a {\it normalized quasi-basis}:\vspace{-2mm}
$$\begin{array}{c}\A=\cup_{(X,Y)\in\T\times\T}\A_{XY}=\{A_1,A_2,\cdots,A_n\}.\end{array}$$

(ii)\, Let $\M=\bar\Delta\otimes_{_{R^{\otimes 2}}}\M_1$, and the
algebra $\K$ be given by Definition 1.1.2.
Define a $\K$-$\K$-bimodule structure on $\M$, such that the left module action
$\tl: \K\times\M\rightarrow \M$ consists of $\tl_{pq}:\K_p\times\M_q\rightarrow\M_{p+q},\forall\,p\geqslant 0,q>0$, and the right one
$\tr: \M\times\K\rightarrow \M$ consists of $\tr_{pq}: \M_p\times\K_q\rightarrow\M_{p+q},\forall\,p>0,q\geqslant 0$ given by usual
matrix product respectively. The $\K$-$\K$-bimodule $\M$ is said to be {\it finitely generated in index
$(0,1)$} with $\M_0=\{0\}$.

\medskip

{\bf Definition 1.1.4}\, Let $H=\sum_{X\in\T}H_X\in\IM_t(R)$ be a matrix, where
$H_X=(h_{ij})_{t\times t}\in E_X\IM(R)E_X$ with $h_{ij}\in
R_X$ for $i,j\in X$, and $h_{ij}=0$ otherwise. Define a
derivation $\td: \K\rightarrow \M, U\mapsto UH-HU$, yielded by $H$ consists of
$\td_p:\K_p\rightarrow \M_p,\forall\,p\geqslant 0$.

\smallskip

It is not difficult to see that $\td_0=0$, and $\td_p$ is determined by $\td_1$ for $p>0$ according to the note
stated under Formula (1.1-6).

\medskip

{\bf Definition 1.1.5}\, A quadruple $\mf{A}=(R,\K,\M,H)$ is called a
{\it matrix bimodule problem} provided

(i)\, $R$ is a minimal algebra with a vertex set $\T$ given by Definition 1.1.1;

(ii)\, $\K$ is an algebra given by Definition 1.1.2;

(iii)\, $\M$ is a $\K$-$\K$-bimodule given by Definition 1.1.3;

(iv)\, There is a derivation $\td: \K\rightarrow \M$
given by Definition 1.1.4.

\smallskip

In particular, if $\M=0$, $\mf{A}$ is said to be a {\it minimal matrix
bimodule problem}.

\bigskip
\bigskip

\noindent{\large\bf 1.2 Bi-comodule problems and Bocses}

\bigskip

We define a notion of bi-comodule problems associated to matrix bimodule problems,
which is the transition into bocses. The concepts and the proofs are proposed by Y. Han.

\medskip

Since $\K_1$ and $\M_1$ are both quasi-free $R$-$R$-bimodules, we have their
$R^{\otimes 2}$-dual structures $\C_1$ and $\N_1$ with $R$-$R$-quasi-basis $\V^*$ and
$\A^*$ respectively:
$$\begin{array}{c}\C_1=\Hom_{_{R^{\otimes 2}}}(\K_1,R^{\otimes 2}),\quad \V^*=\{v_1,v_2,\cdots,v_m\};\\
\N_1=\Hom_{_{R^{\otimes 2}}}(\M_1,R^{\otimes 2}),\quad \A^*=\{a_1,a_2,\cdots,a_n\}.\end{array}\eqno{(1.2\mbox{-}1)}$$
Write $v:X\mapsto Y$ (resp. $a:X\mapsto Y$) provided that $V\in\V_{XY}$ (resp. $A\in\A_{XY}$).

The quasi-basis $\V$ of $\K_1$ has a natural partial order, namely,
$V_i\prec V_j$, if their leading positions $(p_i,q_i)\prec (p_j,q_j)$. Thus
$V_iV_j=\sum_{l<i,j}\gamma_{ijl}\otimes_{_{R^{\otimes 2}}}V_l$,
since $\V\subseteq\mathbb N_t(R\otimes_kR)$.
For a fixed pair $(p,q)$, any fixed order on the set $\{V_i\mid (p_i,q_i)=(p,q)\}$ may be defined,
which gives a linear order on $\V$.

\medskip

{\bf Definition 1.2.1}\, Let $\K$ be a $k$-algebra as in Definition 1.1.2.
We define a quasi-free $R$-module $\C_0=\Hom_R(\K_0,R)\simeq\sum_{X\in\T} R_Xe_{_X}\simeq R$ with an $R$-quasi-basis
$\{e_{_X}\}_{X\in \T}$ dual to $\{E_{X}\}_{X\in \T}$ of $\K_0$; and a
quasi-free $R$-$R$-bimodule $\C_1$ with an $R$-$R$-quasi-basis
$\V^\ast$ defined by the first formula of (1.2-1), which has a linear order yielded from that of $\V$.
Write $\C=\C_0\oplus\C_1$, and
define a coalgebra structure with a counit $\varepsilon: \C\rightarrow R,\, e_{_X}\mapsto 1_{_X},v_j\mapsto 0$
and a comultiplication $\mu:\C\mapsto\C\otimes_{R}\C$ dual to $(\tm_{00},\tm_{01},\tm_{10},\tm_{11})$:
$$\begin{array}{c}\mu=\left(\begin{array}{c}\mu_{00}\\
\mu_{10}+\mu_{01}+\mu_{11}\end{array}\right):\left(\begin{array}{c}\C_0\\ C_1\end{array}\right)
\rightarrow \left(\begin{array}{c}\C_0\otimes_{_R}\C_0,\\
\C_1\otimes_{_R}\C_0\oplus\C_0\otimes_{_R}\C_1\oplus\C_1\otimes_{_R}\C_1\end{array}\right)\\
\mu_{00}(e_{_X})=e_{_X}\otimes_{_R}e_{_X},\,\,\mu_{10}(v_l)=v_l\otimes_{_R}e_{t(v_l)},\\
\mu_{01}(v_l)=e_{s(v_l)}\otimes_{_R}v_l,\,\,
\mu_{11}(v_l)=\sum_{i,j>l}\gamma_{ijl}\otimes_{_{R^{\otimes 3}}}(v_i\otimes_{_R} v_j).
\quad\end{array}$$

Since $\A$ is linearly ordered and $\V\subset\mathbb N_t(R\otimes_{k}R)$,
$V_iA_j=\sum_{l>j}\eta_{ijl}\otimes_{_{R^{\otimes 2}}}A_l$ and
$A_iV_j=\sum_{l>i}\sigma_{ijl}\otimes_{_{R^{\otimes 2}}}A_l$.

\smallskip

{\bf Definition 1.2.2}\, Let $\M$ be a $\K$-$\K$-bimodule as in Definition 1.1.3. A quasi free
$R$-$R$-bimodule $\N_1$  with an $R$-$R$-quasi-basis $\A^\ast$ given
by the second formula of (1.2-1) is define. Write $\N=\N_1$,
then $\N$ has a $\C$-$\C$-bi-comodule structure with the left and right co-module actions
dual to $(\tl_{01},\tl_{11})$ and $(\tr_{10},\tr_{11})$ respectively:
$$\begin{array}{c}\iota=(\iota_{01}+\iota_{11}): \N\rightarrow \C\otimes_{_R}\N=\C_0\otimes_{_R}\N\oplus\C_1\otimes_{_R}\N,\\
\iota_{01}(a_l)=e_{s(a_l)}\otimes_{_R}a_l,\,\iota_{11}(a_l)=\sum_{j<l,i}\eta_{ijl}
\otimes_{_{R^{\otimes 3}}}(v_i\otimes_{_R} a_j);\\
\tau=(\tau_{10}+\tau_{11}):\N\rightarrow\N\otimes_{_R}\C=\N\otimes_{_R}\C_0\oplus\N\otimes_{_R}\C_1,\\
\tau_{10}(a_l)=a_l\otimes_{_R}e_{t(a_l)},\,
\tau_{11}(a_l)=\sum_{i<l,j}\sigma_{ijl}\otimes_{_{R^{\otimes 3}}}(a_i\otimes_{_R} v_j).\end{array}$$

{\bf Definition 1.2.3}\, Assume $\td_1(V_i)=\sum_l\zeta_{il}\otimes_{_{R^{\otimes 2}}}A_l$
defined in 1.1.4. There is a co-derivation
$\partial=(\partial_0,\partial_1): \N\rightarrow\C=\C_0\oplus\C_1$ with
$\partial_0=0$ and $\partial_1(a_l)=\sum_i\zeta_{il}\otimes_{_{R^{\otimes 2}}}v_i$
dual to $(\td_0,\td_1)$, such that
$\mu\partial=(\id\otimes\partial)\iota+(\partial\otimes\id)\tau$.

\medskip

{\bf Definition 1.2.4}\, Let $\mf A=(R,\K,\M,H)$ be a matrix bimodule problem.
A quadruple $\mf{C}=(R,\C,\N,\partial)$ is
said to be a {\it bi-comodule problem associated to $\mf A$} provided

(i) $R$ is a minimal algebra with a vertex set $\T$;

(ii)\, $\C$ is a co-algebra given by Definition 1.2.1;

(iii)\, $\N$ is a $\C$-$\C$-bi-comodule
given by Definition 1.2.2;

(iv)\, $\partial: \N\rightarrow\C$ is a co-derivation given by Definition 1.2.3.

\medskip

Now we construct a bocs via the bi-comodule problem $\mf C$ associated to $\mf A$.
Write $\N^{\otimes p}=\N\otimes_{_R} \cdots \otimes_{_R}\N$ with $p$ copies
of $\N$ and $\N^{\otimes 0}=R$. Define a tensor algebra $\Gamma$ of $\N$ over $R$,
whose multiplication is given by the natural isomorphisms:
$$\Gamma=\oplus_{p=0}^{\infty}
\N^{\otimes p}; \quad \N^{\otimes p}\otimes_{_R}\N^{\otimes
q}\simeq\N^{\otimes (p+q)}.
$$
Let $\Xi= \Gamma \otimes_{R} \C
\otimes_{R} \Gamma$ be a $\Gamma$-$\Gamma$-bimodule of co-algebra
structure induced by $R\hookrightarrow\Gamma$, and denoted by $(\Xi,
{\mu}_{\Xi}, {\varepsilon}_{\Xi})$.
Define the following three $R$-$R$-bimodule maps:
$$\begin{array}{c}\kappa_1 : \N \stackrel{{\iota}}{\rightarrow}
\C \otimes_{_R} \N\stackrel{\cong}{\rightarrow} R
\otimes_{_R} \C\otimes_{_R} \N \hookrightarrow \Gamma
\otimes_{_R} \C \otimes_{_R} \Gamma,\\
\kappa_2 : \N \stackrel{{\tau}}{\rightarrow}
\N \otimes_{_R} \C \stackrel{\cong}{\rightarrow} \N
\otimes_{_R} \C \otimes_{_R} R \hookrightarrow \Gamma
\otimes_{_R} \C \otimes_{_R} \Gamma,\\
\kappa_3 : \N
\stackrel{{\partial}}{\rightarrow} \C
\stackrel{\cong}{\rightarrow} R \otimes_{_R} \C \otimes_{_R}
R \hookrightarrow \Gamma \otimes_{_R} \C \otimes_R
\Gamma.\end{array}$$

\medskip

{\bf Lemma 1.2.5}\,  $\Im (\kappa_1 - \kappa_2 +
\kappa_3)$ is a $\Gamma$-coideal in $\Xi$. Thus $\Omega := \Xi / \Im
(\kappa_1 - \kappa_2 + \kappa_3)$ is a $\Gamma$-$\Gamma$-bimodule
of coalgebra structure.

\smallskip

{\bf Proof}\, Recall the law of bi-comodule: $({\mu} \otimes\id)
{\iota} = (\id\otimes {\iota}) {\iota}, (\id\otimes
{\mu}) {\tau} = ({\tau} \otimes\id){\tau}$, $(\id
\otimes {\tau}){\iota} = ({\iota} \otimes\id){\tau}
$ and $(\id\otimes {\partial}){\iota} -{\mu}{\partial} +
({\partial} \otimes\id){\tau} = 0$. Thus, for any $b \in
\N$, we have
$$\begin{array}{ll} &\mu_{\Xi}(\kappa_1- \kappa_2+ \kappa_3)(b)\\
 =& \mu_{\Xi}(1_\Gamma\otimes {\iota}(b)) - \mu_{\Xi}({\tau}(b)
\otimes 1_\Gamma) + \mu_{\Xi}(\id\otimes {\partial}(b) \otimes1_\Gamma)
\\ =&  ({\mu} \otimes\id){\iota}(b) - (\id\otimes
{\mu}){\tau}(b) + {\mu}({\partial}(b)) \\
=&(\id\otimes {\iota}){\iota}(b) - ({\tau} \otimes\id
){\tau}(b) + (\id\otimes{\partial}){\iota}(b)+
({\partial} \otimes\id)
{\tau}(b) \\
=&  (\id\otimes {\iota} - \id\otimes{\tau} + \id\otimes
{\partial}){\iota}(b) + ({\iota} \otimes \id -{\tau}
\otimes\id + {\partial} \otimes\id) {\tau}(b)\\
=&  u_{(1)} \otimes ({\iota} - {\tau} +
{\partial})(b_{(1)}) + ({\iota}- {\tau} +
{\partial})(b_{(2)}) \otimes u_{(2)} \\
=& u_{(1)} \otimes (\kappa_1 - \kappa_2 + \kappa_3)(b_{(1)}) +
(\kappa_1 - \kappa_2 + \kappa_3)(b_{(2)}) \otimes u_{(2)}\\
& \in  \Xi \otimes \Im(\kappa_1 - \kappa_2 + \kappa_3) + \Im(\kappa_1 -
\kappa_2 + \kappa_3) \otimes \Xi \end{array}$$
where ${\iota}(b) := u_{(1)} \otimes b_{(1)} , {\tau}(b) :=
b_{(2)} \otimes u_{(2)} $, and each term in each step is viewed as
an element in $\Xi \otimes_{\Gamma} \Xi$ naturally. \hfill$\Box$

\medskip

Recall from \cite{Ro} and \cite{CB1}, $\mf{B}=(\Gamma,\Omega)$ defined as above is a bocs with a layer
$$
L=(R; \omega; a_1,a_2,\ldots,a_n; v_1,v_2,\ldots,v_m).
$$
Denote by $\varepsilon_{_\Omega}$ and $\mu_{_\Omega}$ the induced co-unit and co-multiplication,
then $\bar\Omega=$ker$\varepsilon_{_\Omega}$ is a $\Gamma$-$\Gamma$-bimodule
freely generated by $v_1, v_2,\ldots,v_m$, and $\Omega=\Gamma\oplus\bar{\Omega}$ as bimodules.
From this, we use the embedding: $\C_0\oplus \C_1\oplus\N\otimes\C_1\oplus\C_1\otimes\N\hookrightarrow
\Gamma\otimes_R(\C_0\oplus \C_1\oplus\N\otimes_R\C_1\oplus\C_1\otimes_R\N)\otimes_R\Gamma\subset\Omega$;
and the isomorphism: $\bar\Omega\otimes_R\bar\Omega\simeq\bar\Omega\otimes_\Gamma\bar\Omega$.

The group-like $\omega: R\rightarrow\Omega, 1_X\mapsto e_{_X}$ is an $R$-$R$-bimodule map.
Recall from \cite[3.3 Definition]{CB1}, and note that $(\iota_0(a_i)+{\iota_1}(a_i))-(\tau_0(a_i)+{\tau_1}(a_i))+
{\partial_1}(a_i)=(\kappa_1-\kappa_2+\kappa_3)(a_i)=0$ in $\Omega$,
the pair of the differentials determined by $\omega$ is given by $\dz_1: \Gamma\rightarrow \bar{\Omega}$:
$$\begin{array}{ccl}&\dz_1(1_X)&=1_Xe_{_X}-e_{_X}1_X=0,\,\,\, X\in\T,\\
&\dz_1(a_i)&=a_i\otimes_Re_{t(a_i)}-e_{s(a_i)}\otimes_Ra_i=\tau_0(a_i)-\iota_0(a_i)\\
&&={\iota_1}(a_i)-{\tau_1}(a_i)+{\partial_1}(a_i),\,\,\,1\leqslant j\leqslant n;\end{array} $$
and $\dz_2:\bar\Omega\mapsto\bar\Omega\otimes_\Gamma\bar\Omega, \dz_2(v_j)=
\mu(v_j)-e_{s(v_j)}\otimes_Rv_j-v_j\otimes_Re_{t(v_j)}=\mu_{11}(v_j),1\leqslant j\leqslant m$.
Then the bocs $\mf B$ is said to be the {\it bocs associated to the matrix bimodule problem $\mf A$}.
Denote by {\it $(\mf{A},\mf{B})$ the pair of a matrix bimodule problem and its
associated bocs}, or just the {\it pair $(\mf A,\mf B)$}.

Let $(\mf{A},\mf B)$ be a pair with the associated bi-comodule problem $\mf{C}$. Then the module actions
$$\tl(\K_1\times A_i),\, \tr(A_i\times\K_1)\subseteq
\oplus_{l=i+1}^n R^{\otimes 3}\otimes_{_{R^{\otimes 2}}}A_l\eqno{(1.2\mbox{-}2)}$$
by the fact stated before Definition 1.2.2, which is called the {\it triangular property}.
The left and the right co-module actions also possess the
{\it triangular property} by 1.2.2:
$$\begin{array}{c}\iota_1(a_l)\in\, \oplus_{i=1}^{l-1}
\C_1\otimes_{_R}a_i,\quad \tau_1(a_l)\in\, \oplus_{i=1}^{l-1}a_i
\otimes_{_R} \C_1.\end{array}\eqno{(1.2\mbox{-}3)}$$
Define a $\K$-$\K$ sub-bimodule of $\M$,
then a $\K$-$\K$-quotient-bimodule of $\M$:
$$\begin{array}{c}\M^{(h)}=\oplus_{i=h+1}^n\bar\Delta\otimes_{_{R^{\otimes 2}}} A_i\subseteq\M,
\quad \M^{[h]}=\M/\M^{(h)}.\end{array}$$
$\mf{A}^{[h]}=(R,\K,\M^{[h]},\bar\td)$ with $\bar\td$ induced from $\td$
is said to be a {\it quotient problem} of $\mf A$, but $\mf{A}^{[h]}$ itself might be no longer a
matrix bimodule problem. If $\N^{(h)}=\oplus_{i=1}^{h}R\otimes_{_R}a_i\otimes_{_R}R$,
then $\mf C^{(h)}=(R,\C,\N^{(h)},\partial\mid_{\N^{(h)}})$
is a sub-bi-comodule problem of $\mf C$. If $\Gamma^{(h)}$ is a tensor algebra freely generated by $a_1,\cdots,a_{h}$,
then the bocs $\mf{B}=(\Gamma,\Omega)$ has a sub-bocs
$\mf B^{(h)}=(\Gamma^{(h)},\Gamma^{(h)}\otimes_R\Omega\otimes_R
\Gamma^{(h)})$.

Note a simple fact: let $(\mf{A},\mf{C},\mf B)$ be a triple defined as above, then
$$\begin{array}{ll}&\tl(\K_1\times\M_1),\tr(\M_1\times\K_1),\td(\K_1)\subseteq \M^{(h)}_1\,\, \mbox{in}\,\, \mf A\\[1mm]
\Longleftrightarrow & \C_1\otimes_{_R}\N^{(h)}_1=0,\N^{(h)}_1\otimes_{_R}\C_1=0,\partial(\N^{(h)}_1)=0\,\, \mbox{in}\,\, \mf C\\[1mm]
\Longleftrightarrow & \delta(\Gamma^{(h)})=0\,\,\mbox{in}\,\,\mf B.\end{array}\eqno{(1.2\mbox{-}4)}$$
In fact, the condition in $\mf A$ is equivalent to $\eta_{jil}=0,\sigma_{ijl}=0, \zeta_{jl}=0$ for
$l=1,\cdots,h$ and any $i,j$, which is equivalent to the conditions in $\mf C$ and $\mf B$.

\medskip

Recall from \cite{CB1} that a representation of a layered bocs $\mf{B}$ is a
left $\Gamma$-module $P$ of dimension vector $\underline d$ consisting of three sets:
$$\begin{array}{c}\{P_X=k^{d_X}\mid X\in\T\},\quad \{P(x):P_X\rightarrow P_X\mid X\in\T\};\\
\{P(a_i):k^{d_{_{X_i}}}\rightarrow k^{d_{_{Y_i}}}\mid\, a_i:X_i\rightarrow Y_i,\,\,
i=1,\cdots,n\}.\end{array}\eqno{(1.2\mbox{-}5)} $$
A morphism from $P$ to $Q$ is given by a
$\Gamma$-map $f:\Omega\otimes_{\Gamma}P\rightarrow Q$. Clearly,
$\Hom_{\Gamma}(\bar{\Omega}\otimes_{\Gamma}P, Q)\simeq
\oplus_{j=1}^m\Hom_{\Gamma}(\Gamma1_{s(v_j)}\otimes_k1_{t(v_j)}P, Q)\simeq
\oplus_{j=1}^m \Hom_k(1_{t(v_j)}P, 1_{s(v_j)}Q)$. Write
$$\begin{array}{c}f=\{f_X;f(v_j)\mid X\in\T,1\leqslant j\leqslant m\},\end{array}\eqno{(1.2\mbox{-}6)}$$
then \cite{BK} shows that $f$ is a morphism
if and only if for all $a_l\in{\cal A}^\ast, 1\leqslant l\leqslant n$:
$$\begin{array}{ll}P(a_l)f_{_{Y_l}}-f_{_{X_l}}Q(a_l)&=
\sum_{j<l,i}\eta_{ijl}
\otimes_{_{R^{\otimes 3}}}(f(v_i)\otimes_{_R}Q(a_j))\\
&-\sum_{i<l,j}\sigma_{ijl}\otimes_{_{R^{\otimes 3}}}(P(a_i)\otimes_{_R} f(v_j))
+\sum_i\zeta_{il}\otimes_{_{R^{\otimes 2}}}f(v_i).\end{array}\eqno{(1.2\mbox{-}7)}$$

\medskip
\bigskip

\noindent{\large\bf 1.3  Representation categories of matrix bimodule problems}

\bigskip

In this subsection, a notion of ``$\ast$-product" and the
operations between $\ast$-products are defined, which will be used frequently throughout the paper. Based on this nation,
the representation category of a matrix bimodule problem is defined.
It is relatively complicated, but
seems to be useful for the proof of the main theorem.

\medskip

{\bf  Definition 1.3.1} Let $ J(\lambda)=J_d(\lambda)^{e_d}\oplus
J_{d-1}(\lambda)^{e_{d-1}}\oplus \cdots\oplus J_1(\lambda)^{e_1}$ be a Jordan form,
where $e_i$ are non-negative integers. Denote
$e_d+e_{d-1}+\cdots+e_j$ by $m_j$ for $j=1,\cdots,d$. The
following partitioned matrix $W(\lambda)$ similar to $J(\lambda)$ is called
a {\it Weyr matrix of eigenvalue $\lambda$}:
$$
W(\lz)=\left(
\begin{array}{cccccc}
\lz I_{m_1} & W_{12} &0&  \cdots & 0&0 \\
  & \lz I_{m_2} & W_{23}&\cdots &0&0 \\
  && \lz I_{m_3} &\cdots & 0 &0\\
&&&\ddots&\vdots&\vdots\\
  & & & & \lz I_{m_{d-1}}& W_{d-1,d} \\
 & & & && \lz I_{m_d}
\end{array}
\right)_{d\times d}, $$
where $W_{j,j+1}= (I_{m_{j+1}}\, 0)^T$ of size $m_j\times m_{j+1}$ with superscript T denoting transpose.
A direct sum $W=W(\lambda_1)\oplus W(\lambda_2)\oplus
\cdots \oplus W(\lambda_s)$ with
distinct eigenvalues $\lambda_i$ is said to be a {\it Weyr matrix}.
An {\it order ``$\prec$" on the base field} $k$ may be defined, so that
each Weyr matrix has a unique form.
Similarly, let $\{Z_{ij}\mid i,j\in\mathbb Z^+\}$ be a finite set of vertices, and
$S=\oplus_{i,j} k1_{Z_{ij}}\oplus k[z,\phi(z)^{-1}]1_Z$ be a
minimal algebra.
$\bar W\simeq \oplus J_{ij}(\lambda_{i})^{e_{ij}}1_{Z_{ij}}$ or
$\oplus J_{ij}(\lambda_{i})^{e_{ij}}1_{Z_{ij}}\oplus (z1_Z)$
with $\{e_{ij}\}\subset \mathbb Z^+$ is said to be a{\it Weyr matrix over $S$}.
It is possible that some summands of $\bar W$ are diagonal blocks
with the diagonal entries being the primitive idempotents of $S$

Let $\mf A=(R,\K,\M,H)$ be a matrix bimodule problem having a set of integers
$T=\{1,2,\cdots,t\}$ with partition $\T$. A Weyr matrix $W$ over $k$ is called {\it
$R_X$-regular} for $X\in \T_1$, if all the eigenvalues $\lambda$ of $W$ have the property that
$\phi_X(\lambda)\ne 0$. An identity matrix $I$ is also called an {\it $R_X$-regular Weyr matrix} for $X\in\T_0$.
A vector of non-negative integers is said to be a {\it size vector} $\m=(m_1,m_2,\ldots,m_t)$ over $\T$,
if $m_i=m_j,\forall\,i\sim j$. And $\sum_{i=1}^tm_i$ is called the {\it size} of $\m$.

\medskip

{\bf Definition 1.3.2}\, Let $\mf A=(R,\K,\M,H)$ be a matrix bimodule problem,
$S$ a minimal algebra, and $\Sigma=\oplus_{p=1}^\infty S^{\otimes p}$, see Formula (1.1-3).

(i)\, Write $H_X=(h_{ij}(x)1_X)_\tt$ with $h_{ij}(x)\in k[x]$ for $X\in\T_1$, and $x=1,h_{ij}(x)\in k$ for $X\in\T_0$.
Let $\bar W_X$ be a Weyr matrix of size $m_{_X}$ over $S$. There exists an
$\m\times \m$-partitioned matrix over $S$:
$$H_X(\bar W_X)=(B_{ij})_{\tt},\quad B_{ij}=\left\{\begin{array}{cc}
h_{ij}(\bar W_X)_{m_i\times m_j},&i,j\in X,
\\(0)_{m_i\times m_j},&i\notin X\,\, \mbox{or}\,\, j\notin X,\end{array}\right.$$

(ii)\, Let $\m=(m_1,\cdots,m_t)$ and
$\n=(n_1,\cdots, n_t)$ be two size vectors over $\T$, and let
$F\in\IM_{m_{_X}\times n_{_X}}(S^{\otimes p}),p=1,2$,
with an $R_X$-$R_X$-module structure. The star product
$\ast$ of $F_X$ and $E_X$ is defined to be a diagonal $\m\times \n$-partitioned matrix:
$$\begin{array}{c}F_{_X}\ast E_X=\mbox{diag}(B_{11},\cdots,B_{tt}),\quad B_{ii}=
\left\{\begin{array}{cc}F_{_X},& i\in X,\\
0& i\notin X,\end{array}\right.\end{array}$$

(iii)\, Let $U=(u_{ij})\in \V_{XY}\cup\A_{XY}$,
and $\{\bar W_X\in\IM_{m_{_X}}(S)\mid X\in\T\}$,
$\{\bar W'_Y\in\IM_{n_{_Y}}(S)\mid Y\in\T\}$ be two sets of regular Weyr matrices.
Suppose there is an $R_X$-$R_Y$-bimodule structure on
$\IM_{m_{_X}\times n_{_Y}}(S^{\otimes p})$ for $p=1,2$,
and $C\in\IM_{m_{_X}\times n_{_Y}}(S^{\otimes p})$:
$$\begin{array}{c}
C\otimes_{_{R^{\otimes 2}}}(x\otimes_ky)=W_XC W'_Y.\end{array}\eqno(1.3\mbox{-}1)$$
The star product
$\ast$ of $C$ and $U$ is defined to be an $(\m\times \n)$-partitioned matrix:
$$\begin{array}{c}C\ast U=(B_{ij})_\tt,\quad\quad B_{ij}=\left\{\begin{array}{cl}
C\otimes_{_{R^{\otimes 2}}}u_{ij}, &i\in X,j\in Y;\\ (0)_{m_{i}\times n_j},
&i\notin X,\,\mbox{or}\, j\notin Y.\end{array}\right.\end{array}$$

\smallskip

{\bf Lemma 1.3.3}\, Let $\mf A=(R,\K,\M,H)$ be a matrix bimodule problem.

(i)\, If $C\in\IM_{m_{_X}\times n_{_Y}}(S^{\otimes 2})$,
$\td(V_i)=\sum_l\zeta_{il}A_l$ with $\zeta_{il}\in R_X\otimes_k R_Y$ are given by Definition 1.2.3, then by the
usual product of $\Sigma$-matrices:
$$\begin{array}{cc}(C\ast V_i)H_Y(\bar W_Y)-H_X(\bar W_X)(C\ast V_i)=
\sum_{l=1}^n(C\otimes_{R^{\otimes 2}}\zeta_{il})\ast A_l.\end{array}$$

(ii)\, If $F_X\in\IM_{m_{_X}\times n_{_X}}(S^{\otimes p}),p=1,2$ and
$C\in\IM_{m_{_X}\times n_{_Y}}(S^{\otimes q}),q=1,2$,
then by the usual product of $\Sigma$-matrices:
$$\begin{array}{c}(F_{_X}\ast E_X)(C\ast U)=\left\{\begin{array}{cl}(F_{_X}C)\ast U,&1_XU=U;\\
0,&1_XU=0.\end{array}\right.\end{array}$$
Similarly, $(C\ast U)(F_{_X}\ast E_X)=(CF_{_X})\ast U$ for $U1_X=U$ and $0$ otherwise. Moreover,
$$\begin{array}{c}(F_X\ast E_X)(F_X'\ast E_X)=(F_XF_X')\ast E_X,\quad F_X,F'_X\in\IM_{m_{_X}}(S^{\otimes p}),p=1,2.\end{array}$$

\smallskip
(iii)\, Let $U\in E_X\IM_t(R\otimes_kR)E_Y, V\in
E_Y\IM_t(R\otimes_kR)E_Z,G_l\in E_X\IM_t(R\otimes_kR)E_Z$, where \vadjust{\vspace{2.5pt}}$UV=\sum_{l=1}^n\epsilon_l\otimes_{{R^{\otimes 2}}}G_l,\quad
\epsilon_l\in R_X\otimes_kR_Y\otimes_kR_Z.$
Let $\m, \n, \l$ be size vectors over $\T$, $C\in \IM_{m_{_X}\times n_{_Y}}(S^{\otimes p}),
D\in\IM_{n_{_Y}\times l_{_Z}}(S^{\otimes q})$ for $(p,q)\in\{(2,2)$, $(1,2),(2,1)\}$.
Then by the usual $\Sigma$-matrix product:
$$(C\ast U)(D\ast V)={\sum}_{l=1}^n
\big((C\otimes_RD)\otimes_{_{R^{\otimes 3}}}\epsilon_l\big)\ast G_l,$$
where the tensor product $\oplus_{(X,Y,Z)\, \in \T\times \T\times\T}\,\IM_{m_{_X}n_{_Y}}(S^{\otimes p})
\otimes_R\IM_{n_{_Y}\times l_{_Z}}(S^{\otimes q})$ has an $R^{\otimes 3}$-module structure yielded from the $R$-$R$-bimodule
structures given by Formula (1.3-1).

\smallskip

{\bf Proof}\, (i)\, Write $H_X=(\alpha_{pq}),\alpha_{pq}\in R_X; H_Y=(\beta_{pq}),
\beta_{pq}\in R_Y; V_i=(v_{pq}),v_{pq}\in R_X\otimes_k R_Y$,
$$\begin{array}{ll}\mbox{The left side}&=(\sum_l(C\otimes_{R^{\otimes 2}}v_{pl})\otimes_{R}\beta_{lq})
-(\sum_l\alpha_{pl}\otimes_R(C\otimes_{R^{\otimes 2}}v_{lq}))\\
&=(C\otimes_{R^{\otimes 2}}\sum_l(v_{pl}\beta_{lq}-\alpha_{pl}v_{lq}))=
C\ast(V_iH_Y-H_XV_i)\\
&=C\ast\td(V_i)=C\ast(\sum_{l}\zeta_{il}A_l)=\mbox{the right side}.\end{array}$$

(ii)\, Write $U=(u_{pq}),u_{pq}\in R_X\otimes_kR_Y$, the left
side$=(F_X1_X(C\otimes_{R^{\otimes 2}}u_{pq}))
=((F_XC)\otimes_{R^{\otimes 2}}u_{pq})=$the right side.

(iii)\, Write $U=(u_{pq}),u_{pq}\in R_X\otimes_kR_Y,V=(v_{pq}),
v_{pq}\in R_Y\otimes_kR_Z$. The left
side$=\big(\sum_l(C\otimes_{R^{\otimes 2}}u_{pl})(D\otimes_{R^{\otimes 2}}v_{lq})\big)
=\big((CD)\otimes_{R^{\otimes 3}}(\sum_lu_{pl}\otimes_Rv_{lq})\big)
=$the right side. \hfill$\Box$

\medskip

{\bf Definition 1.3.4}\, Let $\mf{A}=(R,\K,\M,H)$ be a matrix
bimodule problem, and $\m$ a size
vector over $\T$. Thus a representation $\bar{P}$ of $\mf{A}$ can be written as an
$\m\times\m$-partitioned matrix over $k$:
$$\begin{array}{c}\bar P=\sum_{X\in\T} H_X(W_X)+\sum_{i=1}^n
P(a_i)\ast A_i,\end{array}$$
where $W_X\in\IM_{m_{_X}}(k)$ is an $R_X$-regular Weyr matrix for any $X\in\T$, and
$P(a_i)\in\IM_{m_{_{X_i}}\times m_{_{Y_i}}}(k)$.
Taken $S=k$, the first summand is defined in 1.3.2 (i), and the second one in (iii) .

\medskip

{\bf Definition 1.3.5}\, Let $\bar P,\bar Q$ be two representations of size vectors $\m,\n$
respectively. A {\it morphism} $\bar f:\bar P\rightarrow \bar Q$
is an $\m\times \n$-partitioned matrix obtained from Definition 1.3.2 (ii) and (iii) for $S=k$:
$$\begin{array}{c}\bar f=\sum_{X\in\T}f_{_X}\ast E_X+\sum_{j=1}^m f(v_j)\ast V_j,\quad \end{array}$$
where $f_{_X}\in\IM_{m_{_X}\times n_{_X}}(k),f(v_j)\in\IM_{m_{s(v_j)}\times n_{t(v_j)}}(k)$,
such that $\bar P\bar f=\bar f\bar Q$ under the matrix product given according to Lemma 1.3.3 (i)--(iii).

\medskip

If $\bar U$ is an object and $\bar f':\bar Q\rightarrow\bar U$ a morphism, then
$\bar f\bar f':\bar P\rightarrow\bar U$ calculated according to Lemma 1.3.3 (ii)--(iii) is still a morphism.
In fact, $(\bar f\bar f')\bar P=\bar f(\bar Q\bar f')=(\bar U\bar f)\bar f'=\bar U(\bar f\bar f')$.
We denote by {\it $R(\mf{A})$ the representation category of the
matrix bimodule problem $\mf{A}$}.

\bigskip
\bigskip

\noindent{\large\bf 1.4. Formal Products and Formal Equations}

\bigskip

In this subsection, 1) a concept of ``formal equation" is introduced to build up a
nice connection between a matrix bimodule problem and its associated bocs;
and 2) a special class of bipartite matrix bimodule problems
is noticed, because of the close relation between such problems and finite dimensional algebras.

\medskip

Let $\mf{A}=(R,\K,\M,H)$ be a matrix bimodule
problem, with the associated bi-comodule problem $\mf{C}=(R,\C,\N,\partial)$ and the bocs $\mf B$.
Recall that $\{E_{_X}\}$ and $\{e_{_X}\}$
are dual bases of $(\K_0,\C_0)$; $\{V_1,\cdots,V_m\}$ and $\{v_1,\cdots,v_m\}$
those of $(\K_1,\C_1)$; $\{A_1,\dots,A_n\}$ and $\{a_1,\cdots,a_n\}$ those
of $(\M_1,\N_1)$. Set $S=R,\Sigma=\Delta,\m=(1,\cdots,1)=\n$ in Definition 1.3.2 (ii)--(iii), then
$$\begin{array}{lll} \Upsilon &=&\sum_{X\in \T}e_{_X}\ast E_{_X}\\
\Pi &=&\sum_{j=1}^mv_j\ast V_j\\
\Theta &=&\sum_{i=1}^na_i\ast A_i
\end{array}\eqno{(1.4\mbox{-}1)}$$
are called the {\it formal products} of $(\K_0,\C_0), (\K_1,\C_1)$
and $(\M_1,\N_1)$ respectively.

\medskip

{\bf Lemma 1.4.1}\, Let $\delta$ be the differential in the bocs $\mf B$. We have
$$\begin{array}{c}\big(\sum_{i=1}^mv_i\ast V_i\big)
\big(\sum_{j=1}^mv_j\ast V_j\big)=\sum_{l=1}^m
\mu_{11}(v_l)\ast V_l;\\[1.5mm]
\big(\sum_{i=1}^na_i\ast A_i\big) \big(\sum_{j=1}^mv_j\ast
V_j\big)=\sum_{l=1}^n
\tau_{11}(a_l)\ast A_l;\\[1.5mm]
\big(\sum_{j=1}^mv_j\ast V_j\big) \big(\sum_{i=1}^na_i\ast
A_i\big)=\sum_{l=1}^n
\iota_{11}(a_l) \ast A_l;\\[1.5mm]
\big(\sum_{j=1}^mv_j\ast V_j\big) H-H\big(\sum_{j=1}^mv_j\ast
V_j\big)=\sum_{l=1}^n
\partial_1(a_l)\ast A_l;\\[1.5mm]
\big(\sum_{l=1}^na_l\ast A_l\big) \big(\sum_{_{X\in \T}} e_{_X}\ast
E_{_X}\big)- \big(\sum_{_{X\in  \T}} e_{_X}\ast E_{_X}\big)
\big(\sum_{l=1}^na_l\ast A_l\big)=\sum_{l=1}^n\delta(a_l)\ast A_l.
\end{array}$$

\smallskip

{\bf Proof}\, 1) The second equality is proved first, and the proofs of the first one and the third one
are similar. By Lemma 1.3.3 (iii) for $S=R,p=q=2$, the left side
$=\sum_{l=1}^n \big(\sum_{i,j}\sigma_{ijl}\otimes_{_{R^{\otimes 3}}}
(a_i\otimes_Rv_i)\big)\ast A_l =$ the right side.

2) For the fourth equality, by Lemma 1.3.3 (i) the left side
$=\sum_{l=1}^n\big(\sum_{j=1}^m\zeta_{lj}\otimes_{_{R^{\otimes 2}}}v_j\big)\ast A_l=$
the right side.

3) For the last one, by Lemma 1.3.3 (ii) for $p=1,q=2$, the left side
$=\sum_{l=1}^n (a_l\otimes_{_R}e_{_{Y_l}}-e_{_{X_l}}\otimes_{_R} a_l)\ast A_l=$ the right side.
\hfill$\Box$

\medskip

The matrix equation
$(\Theta+H)(\Upsilon+\Pi)=(\Upsilon+\Pi)(\Theta+H)$, more precisely,
$$\begin{array}{cc}
&\big(\sum_{i=1}^n a_i\ast A_i+H\big) \big(\sum_{X\in \T}
e_X\ast E_X+\sum_{j=1}^mv_j\ast V_j\big)\\[1mm]
=&\big(\sum_{X\in \T} e_X\ast E_X+\sum_{j=1}^mv_j\ast V_j\big)
\big(\sum_{i=1}^n a_i\ast
A_i+H\big)\end{array}\eqno{(1.4\mbox{-}2)}$$ is called the {\it
formal equation of the pair $(\mf A, \mf B)$} due to the following theorem.

\medskip

{\bf Theorem 1.4.2}\, Let $(p_l,q_l)$ be the leading position of $A_l$ for $l=1,\cdots,n$.
Then the $(p_l,q_l)$-entry of the formal equation is
$$\delta(a_l)=\iota_{11}(a_l)-\tau_{11}(a_l)+\partial_1(a_l).$$

\smallskip

{\bf Proof.} According to Formula (1.4-2) and Lemma 1.4.1:
$$\begin{array}{ll}&\sum_{l=1}^n\delta(a_l)\ast A_l=
\sum_{l=1}^n\big(a_le_{t(a_l)}-e_{s(a_l)}a_l\big)\ast A_l\\[1mm]
=&\sum_{j,i}(v_j\ast V_j)(a_i\ast A_i)- \sum_{i,j}(a_i\ast
A_i)(v_j\ast V_j)+ \sum_{j}\big((v_j\ast V_j) H-
H(v_j\ast V_j)\big)\\[1mm]
=&\sum_{l=1}^n\iota_{11}(a_l)\ast A_l -\sum_{l=1}^n\tau_{11}(a_l)
\ast A_l + \sum_{l=1}^n\partial_1(a_l)\ast A_l\\[1mm]
=&\sum_{l=1}^n\big(\iota_{11}(a_l)-\tau_{11}(a_l)
+\partial_1(a_l)\big)\ast A_l.\end{array}$$
The expression at the leading position $(p_l,q_l)$ of the formal equation
is obtained . \hfill$\Box$

\medskip

Moreover, the first formula of Lemma 1.4.1 gives:
$$\begin{array}{c}\big(\sum_{_{X\in \T}} e_{_X}\ast E_{_X}
+\sum_{i=1}^mv_i\ast V_i\big)\big(\sum_{_{X\in \T}}
e_{_X}\ast E_{_X}+\sum_{j=1}^mv_j\ast V_j\big)\\[1mm]
=\sum_{_{X\in\T}} (e_{_X}\otimes_{_R}e_{_X})\ast E_{_X}+\sum_{l=1}^m
\mu(v_l)\ast V_l.\end{array}\eqno{(1.4\mbox{-}3)}$$

Let $(\mf A,\mf B)$ be a pair with an index set $T$ and a vertex set $\T$. A size vector
$\m=(m_1,\cdots,m_t)$ over $T$, and a dimension vector $\underline d=(d_{_X}\mid X\in\T)$ over $\T$
are said to be {\it associated}, if $m_i=d_X$ for $i\in X$.

\medskip

{\bf Corollary 1.4.3}\, Let $(\mf A,\mf B)$ be a pair. Then the representation categories
$R(\mf A)$ and $R({\mf{B}})$ are equivalent.

\smallskip

{\bf Proof}\, Let $P\in R(\mf B)$ with dimension vector $\underline d$. Without loss of generality, the set $\{P(x)=W_X\mid X\in\T\}$
may be assumed to be a set of regular Weyr matrices. Then $\bar P$ of size vector $\m$
associated with $\underline d$ in Definition 1.3.4 and $P$
in Formula (1.2-5) are one-to-one corresponding; $\bar f$ in Definition 1.3.5 and
$f$ in Formula (1.2-6) are one-to-one corresponding.
Moreover, $\bar P\bar f=\bar f\bar Q$ if and only if $f$ satisfies Formula (1.2-7)
by Theorem 1.4.2.\hfill$\Box$

\medskip

Thanks to Corollary 1.4.3, the representations and morphisms in both  categories $R(\mf A)$ and $R(\mf B)$
can be denoted by $P,f$ in a unified manner.
Finally we define a special class of matrix bimodule problems to end the subsection.
Let $\mf{A}=(R,\K,\M,H=0)$ be a matrix
bimodule problem with $R$ trivial. $\mf{A}$ is said to be {\it
bipartite} provided that $\T=\T'\dot{\cup}\T''$; $R=R'\times R''$ and $\K=\K'\times \K''$ are direct
products of algebras; and $\M$ is a $\K'$-$\K''$-bimodule.

\medskip

{\bf Remark 1.4.4}\, Let $\Lambda$ be a finite-dimensional basic $k$-algebra,
$J=\rad(\Lambda)$ be the Jacobson radical of $\Lambda$ with the
nilpotent index $m$, and the top $S=\Lambda/J$. Suppose $\{e_1,\cdots,e_h\}$ is a complete set
of orthogonal primitive idempotents of $\Lambda$.  Taking the pre-images of
$k$-bases of $e_i(J^i/J^{i+1})e_j$ under the canonical
projections $J^i\rightarrow J^i/J^{i+1}$ in turn for $i=m,
\cdots,1$, an ordered basis of $J$ under the length order is obtained,
see \cite[6.1]{CB1} for details. Then we construct the left regular
representation $\bar\Lambda$ of $\Lambda$ under the $k$-basis
$(a_{n},\cdots,a_2,a_1,e_1,\cdots,e_h)$ of $\Lambda$,
which leads to a bipartite matrix bimodule problem $\mf A=(R,\K,\M,H=0)$ with
$$R=S\times S;\quad\K_0\oplus\K_1=\bar\Lambda\times\bar\Lambda;
\quad\M_1=\mbox{rad}(\bar\Lambda);\quad H=0.$$

A simple calculation shows that the {\it row indices} of the leading positions of
the base matrices in $\A$ are {\it pairwise different}, and the {\it column index} of the
leading position of $A\in\A_{XZ}$ equals $j_{_Z}=$max$\{j\in Z\}$ for any $X\in\T$, they are {\it concentrated},
and the $j_{_Z}$-th column is said to be the {\it main column} over $Z$. Such a condition
is denoted by RDCC for short. The condition may not be essential in the proof of the main theorem,
but makes it easier and more intuitive.

\medskip

{\bf Example 1.4.5}\,\cite{D1,R1} Let $Q= $ {\unitlength=1mm
\begin{picture}(20,4) \put(10,1){\circle*{0.6}}
\put(7,1){\circle{4}} \put(13,1){\circle{4}}
\put(9,2){\vector(0,1){0}} \put(11,2){\vector(0,1){0}}
\put(2,0){$a$} \put(16,0){$b$}
\end{picture}}
be a quiver, $I=\langle a^2,ba- ab, ab^2 , b^3\rangle$ be an ideal
of $k Q$, and $\Lambda=kQ/I$. Denote  the residue
classes of $e, a, b$ in $\Lambda$ still by $e, a,b $ respectively.
Moreover,  set $ c=b^2, d=ab$. Then an ordered $k$-basis $\{d, c, b,
a, e\}$ of $\Lambda$ gives a regular representation $\bar \Lambda$.
A matrix bimodule problem $\mf A$ follows from Remark 1.4.4, we may denote by $A,B,C,D$
the $R$-$R$-quasi-basis of $\M_1$, and by $a,b,c,d$ the $R$-$R$-dual basis of $\N_1$.
Then the associated bocs $\mf B$ of $\mf A$ has a layer $L=(R;\omega; a,b,c,d; u_1,u_2,u_3,u_4,$
$v_1,v_2,v_3,v_4)$.
The formal equation of the pair $(\mf A,\mf B)$ can be written as:
$$\small{
\left(
\begin{array}{ccccc}
e& 0 &  u_1    & u_2 & u_4 \\
  & e &  u_2    & 0 & u_3 \\
  &   &   e   & 0 & u_2 \\
  &   &       & e & u_1 \\
  &   &       &   & e
\end{array}
\right) \left(
\begin{array}{ccccc}
0 & 0 &  a    &  b & d \\
  & 0 &  b    & 0 & c \\
  &   &   0   & 0 & b \\
  &   &       & 0 & a \\
  &   &       &   & 0
\end{array}
\right)= \left(
\begin{array}{ccccc}
0 & 0 &  a    &  b & d \\
  & 0 &  b    & 0 & c \\
  &   &   0   & 0 & b \\
  &   &       & 0 & a \\
  &   &       &   & 0
\end{array}
\right) \left(
\begin{array}{ccccc}
f & 0 &  v_1    & v_2 & v_4 \\
  & f &  v_2    & 0 & v_3 \\
  &   &   f   & 0 & v_2 \\
  &   &       & f & v_1 \\
  &   &       &   & f
\end{array}
\right)}
$$
with $e=e_{_{X}}, f=e_{_{Y}}$ for simplicity.
The differentials of the solid arrows of $\mf{B}$ can be read off according to Theorem 1.4.2:
$$
\unitlength=1mm
\begin{picture}(80, 26)
\put(10,0){\circle*{1.00}} \put(10,20){\circle*{1.00}}
\qbezier(7,19)(5,11)(7,3) \put(7,3){\vector(1,-1){1}}
\qbezier(13,19)(15,11)(13,3) \put(13,3){\vector(-1,-1){1}}

\qbezier[10](9,20)(4.5,19)(5,21) \qbezier[10](9,20)(5.5,22)(5,21)
\qbezier[10](9,21)(6,23)(7,24) \qbezier[10](9,21)(9,25)(7,24)
\qbezier[10](11,20)(15.5,19)(15,21)
\qbezier[10](11,20)(14.5,22)(15,21)
\qbezier[10](11,21)(14,23)(13,24) \qbezier[10](11,21)(11,25)(13,24)

\qbezier[10](9,0)(4,1)(5,-1) \qbezier[10](9,0)(5,-2)(5,-1)
\qbezier[10](9,-1)(6,-3)(7,-4) \qbezier[10](9,-1)(9,-5)(7,-4)
\qbezier[10](11,0)(16,1)(15,-1) \qbezier[10](11,0)(15,-2)(15,-1)
\qbezier[10](11,-1)(14,-3)(13,-4) \qbezier[10](11,-1)(11,-5)(13,-4)

\put(17,20){$X$} \put(17,0){$Y$}

\put(9.00,19.00){\vector(0,-1){18}}
\put(11.00,19.00){\vector(0,-1){18}}

\put(4,10){\mbox{$a$}} \put(8,10){\mbox{$b$}}
\put(10.50,10){\mbox{$c$}} \put(14,10){\mbox{$d$}}

\put(35,8){\mbox{$ \left\{\begin{array}{l} \dz(a)=0,\\ \dz(b)=0, \\
\dz(c)=u_2b-bv_2, \\ \dz(d)=u_1b+u_2a-bv_1-av_2.
\end{array}\right.$}}
\end{picture}
$$

\bigskip
\bigskip
\medskip

\centerline{\Large\bf  2 Reductions on matrix bimodule problems}

\bigskip

In this section, the reduction theorem and eight reductions on matrix bimodule problems
corresponding to those on bocses are stated and proved, thus the induced matrix bimodule problems are
constructed. Finally, a concept of defining systems of pairs is defined
in order to help to construct the induced pairs in a sequence of reductions.

\bigskip
\bigskip

\noindent{\large\bf 2.1  Admissible bimodules and induced matrix bimodule problems}

\bigskip

In the subsection we prove the reduction theorem
on matrix bimodule problems via admissible bimodules;
then give the connection to the corresponding admissible functors and the reduction
theorem on bocses.  Before doing so, the following lemma is mentioned first.

\medskip

{\bf Lemma 2.1.1\, }Let $D$ be a commutative algebra, and $\Lambda, \Sigma$ be commutative
$D$-algebras. Suppose $_\Lambda \mathcal{G}$ and $\mathcal{S}_\Sigma$ are finitely generated
projective left $\Lambda$-module and right $\Sigma$-module respectively.

(i) $\mathcal{G}\otimes_D\mathcal{S}$ is a projective $\Lambda\otimes_D\Sigma$-module.

(ii) There exists a $\Lambda\otimes_D\Sigma$-module isomorphism:
$${\rm Hom}_{\Lambda}(\mathcal{G}, \Lambda)\otimes_D {\rm
Hom}_{\Sigma}(\mathcal{S}, \Sigma)\cong {\rm
Hom}_{\Lambda\otimes_D\Sigma}(\mathcal{G}\otimes_D\mathcal{S},
\Lambda\otimes_D\Sigma).$$

{\bf Proof}\, (i)\, Suppose $_\Lambda\mathcal G,\mathcal S_\Sigma$ are both free with the basis
$\{u_1,\cdots,u_m\},\{v_1,\cdots,v_n\}$ respectively. Choose a free $\Lambda\otimes_D\Sigma$-module
$\mathcal{F}$ with basis $\{w_{ij}\mid 1\leqslant i\leqslant m, 1\leqslant j\leqslant n\}$. For $x=\sum_{i=1}^m\lambda_iu_i\in \mathcal{G}$
and $y=\sum_{j=1}^nv_j\sigma_j\in \mathcal{S}$, we define $f:\mathcal G\times\mathcal S\rightarrow\mathcal F,
(x, y)=\sum_{i=1}^m\sum_{j=1}^n\,(\lambda_i\otimes\sigma_j)w_{ij}$.
Then $xr={\sum}_{i=1}^m(\lambda_ir)u_i,$ and $ry=\sum_{j=1}^nv_j(r\sigma_j)$ for any $r\in D$,
and hence $f(xr, y)=f(x,ry)$, $f(x_1+x_2,y)=f(x_1,y)+f(x_2,y)$ and
$f(x,y_1+y_2)=f(x,y_1)+f(x,y_2)$. Thus there exists a unique
$\Lambda\otimes_D\Sigma$-linear map $\tilde f:
\mathcal{G}\otimes_D\mathcal{S}\to \mathcal{F}$ such that $\tilde
f(x\otimes y)=f(x,y),\forall x\in \mathcal{G},y\in \mathcal{S}$.
In particular, $\tilde f(u_i\otimes v_j)=w_{ij}$. Thus $\{u_i\otimes v_j\mid 1\leqslant i\leqslant m, 1\leqslant j\leqslant n\}$ is a
$\Lambda\otimes_D\Sigma$-basis of $\mathcal{G}\otimes_D\mathcal{S}$, and $\mathcal{G}\otimes_D\mathcal{S}$ is free.
If both $_{\Lambda}\mathcal G,\mathcal S_\Sigma$ are projective, then there are some
$_{\Lambda}\mathcal G',\mathcal S_\Sigma'$, such that both
$_{\Lambda}\mathcal G\oplus\mathcal S_\Sigma,_{\Lambda}\mathcal G'\oplus\mathcal S_\Sigma'$
being free. The assertion follows.

(ii) It is stressed, that $\mathcal G\otimes_D\mathcal S$ and
${\rm Hom}_{\Lambda\otimes_D\Sigma}(\mathcal{G}\otimes_D\mathcal{S},
\Lambda\otimes_D\Sigma)$ are projective
$\Lambda\otimes_D\Sigma$-modules by (i).
Consider the following commutative diagram
$$\begin{array}{c}\xymatrix{
{\rm Hom}_{\Lambda}(\mathcal{G}, \Lambda)\times {\rm
Hom}_{\Sigma}(\mathcal{S}, \Sigma)\ar[r] \ar[d]_{\psi}& {\rm
Hom}_{\Lambda}(\mathcal{G}, \Lambda)\otimes_D {\rm
Hom}_{\Sigma}(\mathcal{S}, \Sigma)\ar[ld]^{\wt{\psi}}\\
{\rm Hom}_{\Lambda\otimes_D\Sigma}(\mathcal{G}\otimes_D\mathcal{S},
\Lambda\otimes_D\Sigma) &} \end{array}$$
Let $f\in {\rm Hom}_\Lambda(\mathcal{G}, \Lambda)$ and $g\in {\rm
Hom}_\Sigma(\mathcal{S}, \Sigma)$. Since $f$ and $g$ are $D$-linear,
there exists a $\Lambda\otimes_D\Sigma$-linear map
$\psi(f,g): \mathcal{G}\otimes_D\mathcal{S}\to \Lambda\otimes_D\Sigma$, such that
$(\psi(f, g))(x\otimes y)=f(x)\otimes g(y),$
for $(x, y)\in \mathcal{G}\times \mathcal{S}$. Now
$\psi(fr,g)=\psi(f, r g)$ for $r\in D$, thus there
exists a unique $(\Lambda\otimes_D\Sigma)$-linear map $\tilde\psi$ given by $f\otimes g\mapsto
\psi(f, g)$, which is clearly natural in
both $\Hom_{\Lambda}(\mathcal{G},\Lambda)$ and $\Hom_{\Sigma}(\mathcal{S},\Sigma)$. $\tilde\psi$ is an isomorphism
if $_\Lambda\mathcal G,\mathcal S_\Sigma$ are free, consequently,  $\tilde\psi$ is an isomorphism
if $_\Lambda\mathcal G,\mathcal S_\Sigma$ are projective. \hfill$\Box$

\medskip

Next we introduce a notion of admissible bimodules which is a module-theory version of admissible
functor \cite{CB1}. Some preliminaries are needed.

Let $(\mf A, \mf B)$ be a pair with a minimal algebra $R$. Recall from Formula (1.2-4), that $\delta(a_i)=0$
for the first $h$ arrows $a_i,i=1,\cdots,h,$ of $\mf B$, if and only if
$\tl(\K_1\times\M_1),\tr(\M_1\times\K_1),\td(\M_1)\subseteq\M^{(h)}_1$.
The algebra $\bar R=R[a_1,\cdots,a_h]$ is said to be {\it pre-minimal}.

Let $\underline d=(n_{_X}\mid X\in\T)$ be a dimension vector over $\T$,
$R'$ a minimal algebra with the vertex set $\T'$ and the algebra $\Delta'=\sum_{p=1}^\infty {R'}^{\otimes p}$, see
Formula (1.1-3). Define an $R'$-$\bar R$-bimodule $L$ (or an $\bar R$-module over $R'$)
of dimension vector $\underline d$ as follows:
$L=\oplus_{X\in\T}L_X$, where $L_X=\oplus_{p=1}^{n_{_X}}R'_{Z_{(X,p)}}$ with $Z_{(X,p)}\in\T'$,
be an $R'$-$\bar R$-bimodule. Let $L^\ast=\oplus_{X\in\T}L^\ast_X$ be
the $R'$-dual module of $L$, where $L^\ast_X=\Hom_{R'}(L_X,R')=\oplus_{p=1}^{n_{_X}}{R'}^\ast_{Z_{(X,p)}}$.
Clearly, $L^\ast$ is an $\bar R$-$R'$-bimodule.

Denote by $\textsf e_{_{Z_{(X,p)}}}$ the
$(1\times n_{_X})$-matrix with the $p$-th
entry $1_{{Z_{(X,p)}}}$ and others zero. Then the set $\{\textsf e_{_{Z_{(X,p)}}}\mid 1\leqslant p\leqslant n_{_X}\}$
 forms an $R'$-quasi-basis of $L_X$. Similarly, the set $\{\f_{_{Z_{(X,p)}}}=\textsf e_{_{Z_{(X,p)}}}^\ast=\textsf e_{_{Z_{(X,p)}}}^T
 \mid 1\leqslant p\leqslant n_{_X}\}$ of $(n_{_X}\times 1)$-matrices
forms an $R'$-quasi-basis of $L^\ast$, where the superscript ``$T$" stands for the transpose of matrix.
Note that $\textsf f_{_{Z_{(Y,q)}}}(\textsf e_{_{Z_{(X,p)}}})=\textsf e_{_{Z_{(X,p)}}}\textsf f_{_{Z_{(Y,q)}}}
=1_{Z_{(X,p)}}$ for $(X,p)=(Y,q)$, or zero otherwise. And $\textsf f_{_{Z_{(X,p)}}}\textsf e_{_{Z_{(X,q)}}}$
is a matrix unit with the $(p,q)$-th entry $1_{Z_{(X,p)}}$ for $Z_{(X,p)}=Z_{(X,q)}$, or a zero matrix otherwise. Define
$$
\begin{array}{ll}
\hat E_d=\oplus_{X\in\T}(\hat E_d)_X, & (\hat E_d)_X=\oplus_{p=1}^{n_{_X}}R'(\textsf f_{_{Z_{(X,p)}}}
\textsf e_{_{Z_{(X,p)}}})R'\subseteq\mathbb D_{n_{_X}}(R'),\\
\hat E_u=\oplus_{X\in\T}(\hat E_u)_X, & (\hat E_u)_X=
\oplus_{1\leqslant p<q\leqslant n_{_X}}R'(\textsf f_{_{Z_{(X,p)}}}\otimes_k
\textsf e_{_{Z_{(X,q)}}})R'\subseteq\mathbb N_{n_{_X}}(R'\otimes_kR'),\\
\hat E_l=\oplus_{X\in\T}(\hat E_l)_X, & (\hat E_l)_X
=\oplus_{1\leqslant q<p\leqslant n_{_X}}R'(\textsf f_{_{Z_{(X,p)}}}\otimes_k\textsf e_{_{Z_{(X,q)}}})R'\subseteq\mathbb N_{n_{_X}}(R'\otimes_kR')^T.
\end{array}
$$
In the following, $c:X\rightarrow Y$ is an arrow of $\bar R$, and the set $\mathscr S=\{x\mid X\in\T_1\}\cup\{a_1,\cdots,a_h\}$.
$$
\begin{array}{l}
\bar E_0=\{(B_X)_{X\in\T}\in\hat E_d\mid B_XL(c)=L(c)B_Y,
\forall\, c\in\mathscr S\},\\
\bar E_1=\{(B_X)_{X\in\T}\in\hat E_u\mid B_XL(c)=L(c)B_Y,
\forall\, c\in\mathscr S\},\\
\bar E_l=\{(B_X)_{X\in\T}\in\hat E_l\mid B_XL(c)=L(c)B_Y,
\forall\, c\in\mathscr S\}.
\end{array}
$$

{\bf  Definition 2.1.2}\,  With the notations as above. The
$R'$-$\bar R$-bimodule $L$ of dimension vector $\underline d$ is said to be {\it admissible} provided that

\subitem\hspace{-7mm}{\bf (a1)} $L$ is sincere over $R'$;
\subitem\hspace{-7mm}{\bf (a2)} $\bar E_0\simeq R'$ with an $R'$-quasi-basis:
$$\begin{array}{c}\F_0=\{F_Z=(F_{ZX})_{X\in\T}\mid Z\in\T'\},\quad F_{ZX}=\sum_{Z_{(X,p)}=Z}\textsf f_{_{Z_{(X,p)}}}
\textsf e_{_{Z_{(X,p)}}};\end{array}$$
\subitem\hspace{-7mm}{\bf (a3)} $\bar E_1$ is a quasi-free $R'$-$R'$-bimodule
with a quasi-basis $\F_1$ for $i=1,\cdots,l$:
$$\begin{array}{c}\F_1=\{F_{i}=\sum_{p_{i\iota}<q_{i\iota}}\varepsilon_{_{X_i},p_{i\iota}q_{i\iota}}(\textsf f_{_{Z_{(X_i,p_{i\iota})}}}\otimes_k
\textsf e_{_{Z_{(X_i,q_{i\iota})}}})\mid \varepsilon_{_{X_i},p_{i\iota}q_{i\iota}}=0\mbox{ or }1\},
\end{array}$$
where $\{(p_{i\iota},q_{i\iota})\mid \varepsilon_{_{X_i},p_{i\iota}q_{i\iota}}=1\}
\cap\{(p_{j\kappa},q_{j\kappa})\mid \varepsilon_{_{X_j},p_{j\kappa}q_{j\kappa}}=1\}=\emptyset, \forall\, X_i=X_j$;
\subitem\hspace{-7mm}{\bf (a4)} $\bar E_l=\{0\}$.

\vspace{1mm}
\noindent The $k$-algebra $\bar E=\bar E_0\oplus(\bar\Delta'\otimes_{{R'}^{\otimes 2}}\bar E_1)
\subseteq\Pi_{X\in\T}\mathbb T_{n_{_X}}(\Delta')$ may be
called a {\it pseudo endomorphism algebra} of $L$, which is finitely generated in index $(0,1)$.

\medskip

Let $M=\oplus_{X\in\T}M_X\subseteq\oplus_{X\in\T}\mathbb N_{n_{_X}}(R'\otimes_kR')$ be an $R'$-$R'$-bimodule,
where $M_X$ has an $R'$-$R'$-quasi-basis $\{E_{(Xpq)}\mid p<q\}$, the matrix units of size $n_{_X}$
with the $(p,q)$-entry $1_{Z_{(X,p)}}\otimes_k 1_{Z_{(X,q)}}$ and others zero. There is an $R'$-$R'$-isomorphism
$\kappa:\hat E_u\rightarrow M, \textsf f_{_{Z_{(X,p)}}}\otimes_k\textsf e_{_{Z_{(X,q)}}}\mapsto E_{(Xpq)}$.
Furthermore, $M_{X},\forall\,X\in\T$,
possesses an $\bar R$-$\bar R$-bimodule structure as follows: if
$b, c\in\Lambda$ with $t(b)=X=s(c)$, then
$b\otimes_{\bar R}E_{Xpq}\otimes_{\bar R}c
=L(b)E_{Xpq}L(c)$. It is clear that $\kappa$ is also an $\bar R$-$\bar R$-bimodule
isomorphism, and $\hat E_u$ may be identified with $M$.
Thus $\bar E_1$ can be viewed as a submodule of $M$.

\medskip

Write the quasi-free $R'$-$R'$-bimodule
$L^\ast\otimes_kL=\oplus_{(X,Y)\in\T\times\T}L_X^\ast\otimes_kL_Y$.
An induced matrix bimodule $\mf A'$
of $\mf A$ based on an admissible bimodule is described below.

\medskip

{\bf Construction 2.1.3}\, Let $\mf{A}=(R,\K,\M,H)$ be a matrix
bimodule problem. Suppose $\bar R,R',L, \underline d,\bar E$ are given as in Definition 2.1.2.
Then there exists an induced matrix bimodule problem $\mf A'=(R',\K',\M',H')$
in the following sense.

(i)\, Let $\n=(n_1,\cdots,n_t)$ be the size vector
associated with $\underline d$ stated before Corollary 1.4.3. Define a set of integers
$T'=\{1,\cdots,t'\}$ with $t'=\sum_{i\in T}n_{i}$.
Then the set of vertices of $R'$ is the partition $\T'$ of $T'$,
and the matrices in $\K',\M'$ and $H'$ are of size $\n\times \n$ partitioned under $\T$.

(ii)\, The $k$-algebra $\K'$ is given as follows.
First, let $\F'_0=\{E'_Z=\sum_{X\in\T}F_{Z,X}\ast E_X\in\mathbb D_{t'}(R')\mid Z\in\T'\}$
by Definition 1.3.2 (ii) for $S=R',p=1$, and $\K_0'$ be an algebra
generated by $\F'_0$ over $R'$. $\F'_0$ is a quasi-basis of $\K_0'$
via the isomorphism $\bar E_0\stackrel{\nu_0}\rightarrow\K_0',F_Z\mapsto E'_Z$, and $\K_0'\simeq R'$.

Second, let $\mathcal F'_1=\{F_i'=F_{i}\ast E_{X_i} \mid i=1,\cdots,l\}$
by 1.3.2 (ii) for $S=R',p=2$, and $\K_{10}'$ be an $R'$-$R'$-bimodule generated by $\F'_1$, where $\F_1'$
is a quasi-basis via the isomorphism $\bar E_1\stackrel{\nu_1}\rightarrow\K_{10}',
F_i\mapsto F_i'$. There exists a natural order on $\mathcal F_1'$
according to the leading positions of matrices.
Let $\K_{11}'=(L^\ast\otimes_kL)\otimes_{R^{\otimes 2}}\K_1\subseteq
\mathbb N_{t'}(R'\otimes_k R')$ be an $R'$-$R'$-bimodule with a quasi-basis:
$$\begin{array}{c}\U'=\{(\f_{_{Z_{(X'_{j},p)}}}\otimes_k\textsf e_{_{Z_{(Y'_{j},q)}}})\ast
V_j\mid V_j\in\V_{X'_{j}Y'_j},\,1\leqslant p\leqslant n_{_{X_j'}},\,
1\leqslant q\leqslant n_{_{Y_j'}}, 1\leqslant j\leqslant m\}\end{array}$$
given by Definition 1.3.2 (iii) for $S=R',p=2$. Finally, set
$\K_1'=\K_{10}'\oplus\K_{11}'$ with a quasi-basis $\V'=\mathcal F_1'\cup\U'$.

(iii)\, Let $\M'_1=(L^\ast\otimes_kL)\otimes_{R^{\otimes 2}}\M_1$ be an $R'$-$R'$-bimodule
with a normalized quasi-basis
$$\begin{array}{c}\A'=\{(\f_{_{Z_{(X_{i},p)}}}\otimes_k\textsf e_{_{Z_{(Y_{i},q)}}})\ast
A_i\mid A_i\in\A_{X_{i}Y_{i}},\,1\leqslant p\leqslant n_{_{X_j}},1\leqslant q\leqslant n_{_{Y_j}}, h<i\leqslant n\}\end{array}$$
given by Definition 1.3.2 (iii) for $S=R',p=2$.

(iv)\, Let $H'=\sum_{X\in\T}H_X(L_X(x))+\sum_{i=1}^hL(a_i)\ast A_i$ be a matrix over $R'$, where
$L_X(x)=\bar W_X$, $H_X(\bar W_X)$ is defined in 1.3.2 (i);
and $\ast$ is given by 1.3.2 (iii) for $S=R',p=1$.

(v)\, The product $\tm_{11}':(\K_{10}'\oplus\K_{11}')\times(\K'_{10}\oplus\K_{11}')\rightarrow \K_2'$ is
given by
$$\begin{array}{c}\big((\f_{Z_{(X'_{i},p_1)}}\otimes_k\textsf e_{Z_{(Y'_{i},q_1)}})\ast V_i\big)\big((\f_{Z_{(X'_{j},p_2)}}
\otimes_k\textsf e_{Z_{(Y_{j}',q_2)}})\ast V_j\big)\\
=\sum_l\Big(\big((\f_{Z_{(X'_{i},p_1)}}\otimes_k\textsf e_{Z_{(Y'_{i},q_1)}})
(\f_{Z_{(X'_j,p_2)}}\otimes_k\textsf e_{Z_{(Y'_j,q_2)}})\big)\otimes_{R^{\otimes 3}}\gamma_{ijl}\Big)\ast V_l,\\[1.5mm]
(F_{i}\ast E_{X_i})\big((\f_{X'_l,p}\otimes_k\textsf e_{Z_{(Y'_l,q)}})\ast V_l\big)
=\big((F_{i}(\f_{Z_{(X'_l,p)}}\otimes_k\textsf e_{Z_{(Y'_l,q)}})\big)\ast(E_{X_i}V_l),\\[1mm]
\big((\f_{_{Z_{(X'_l,p)}}}\otimes_k\textsf e_{_{Z_{(Y'_l,q)}}})\ast V_l\big)(F_{i}\ast E_{X_i})
=\big((\f_{_{Z_{(X'_l,p)}}}\otimes_k\textsf e_{_{Z_{(Y'_l,q)}}}) F_{i}\big)\ast (V_lE_{X_i}),\\[1mm]
(F_{i}\ast E_{X_i})(F_{j}\ast E_{X_j})=(F_iF_j)\ast(E_{X_i}E_{X_j})
\end{array}$$
according to 1.3.3 (iii) and (ii) for $p=2=q$. The left module action $\tl_{11}':(\K_{10}'\oplus\K_{11}')\times\M'_1\rightarrow \M_2'$ is
similar to the first and second formulae above, the right one
$\tr_{11}':\M'_1\times(\K_{10}'\oplus\K_{11}')\rightarrow \M_2'$ to the first and third ones.
Finally, the derivation $\td_1=(\td_{10},\td_{11})$ with $\td_{10}':\K'_{10}\rightarrow\{0\}$ and
$$\begin{array}{c}\td'_{11}:\K_{11}'\rightarrow \M_1',(\f_{_{Z_{(X_j',p)}}}\otimes_k\textsf e_{_{Z_{(Y'_j,q)}}})\ast V_j
\mapsto\sum_l\big(\zeta_{jl}\otimes_{R^{\otimes 2}}(\f_{_{Z_{(X_j',p)}}}\otimes_k\textsf e_{_{Z_{(Y'_j,q)}}})\big)\ast A_l.\end{array}\eqno{\Box}$$

Admissible bimodules can be transferred to admissible functors as follows.
A minimal algebra $R$ can be viewed as a minimal category $A'=\prod_{X\in\T_0}$mod$R_X
\times\prod_{X\in\T_1}P(R_X)$, \cite[2.1]{CB1} by the one-to-one correspondence
between the vertex set of $R$ and the set of indecomposable objects of $A'$, the two sets may
be identified for the sake of convenience.
Then $\bar R$ determines a pre-minimal category $A''$ by adding some
morphisms $a_i:X_i\mapsto Y_i,i=1,\cdots,h,$ into $A'$. A similar transfer holds from
$R'$ to $B'$. Thus $L$ can be viewed as a functor $\theta':A'\rightarrow B'$, where
$$\begin{array}{c}\theta'(X)=\oplus_{p=1}^{n_{_X}}Z_{(X,p)},
\quad \theta'(c)=L(c),\,\,\,\forall\,c\in\Lambda.
\end{array}\eqno{(2.1\mbox{-}1)}$$
We stress, that the opposite construction is usually impossible.
Throughout the paper, the right module structure and upper triangular matrix are mainly used, which is
opposite to the left module and lower triangular matrix used in \cite{CB1}.

\medskip

{\bf Proposition 2.1.4}\, The functor $\theta'$ is admissible in the sense of \cite[4.3 Definition]{CB1}.

\smallskip

{\bf Proof}\, 1)\, (A1) is clear; (a1) implies (A2); the finite set
$\Lambda=\{Z_{(X,p)}\mid 1\leqslant p\leqslant n_{_X},X\in\T\}$ has a partial order: $Z_{(X,p)}\prec Z_{(X,q)}$ if
$p<q$ in (A3); (A4) follows from $\textsf e_{_{Z_{(X,p)}}}\textsf f_{_{Z_{(Y,q)}}}=1_{Z_{X,p}}$
for $(X,p)=(Y,q)$, or $0$ otherwise; (A5) from $\sum\textsf f_{_{Z_{(X,p)}}}\textsf e_{_{Z_{(X,p)}}}=1_{L_X}$.

2)\, Let $\bar E^\ast_0=$Hom$_{R'}(\bar E_0,R')$, then $\bar E^\ast_{0}\simeq$Hom$_{R'}(R',R')\simeq R'$ by (a2).
Since $(\hat E_d)_X\simeq\oplus_{p=1}^{n_{_X}}R'(\textsf f_{_{Z_{(X,p)}}}
\textsf e_{_{Z_{(X,p)}}}) \simeq\oplus_{p=1}^{n_{_X}}R'(\textsf f_{_{Z_{(X,p)}}}\otimes_{_{R'}}R'
\textsf e_{_{Z_{(X,p)}}})\simeq\oplus_{p=1}^{n_{_X}}R'(\textsf f_{_{Z_{(X,p)}}}\otimes_{_{R'}}
\textsf e_{_{Z_{(X,p)}}})R'$, Lemma 2.1.1 (ii) shows
$$\begin{array}{ll}&\Hom_{R'}(L^\ast\otimes_{R'}L,R')\simeq
\Hom_{R'\otimes_{R'}R'}(L^\ast\otimes_{R'}L,R'\otimes_{R'}R')\\
\simeq&\Hom_{R'}(L^\ast,R')\otimes_{R'}\Hom_{R'}(L,R')
\simeq L\otimes_{R'}L^\ast.\end{array}$$
We claim, that $\bar E^\ast_0=\sum_{X\in\T,p} R'(\textsf e_{_{Z_{(X,p)}}}
\otimes_{_{\bar R}}\f_{_{Z_{(X,p)}}})R'$. In fact, ${\hat E}^\ast_d\subseteq L\otimes_{R'} L^\ast$, and ${\hat E}^\ast_d$ has
an $\bar R$-$\bar R$-bimodule structure given
by $b(\textsf e_{_{Z_{(X,p)}}}\otimes_{_{R'}}
\textsf f_{_{Z_{(X,p)}}})c=(\textsf e_{_{Z_{(X,p)}}}L(b))\otimes_{_{R'}}
(L(c)\textsf f_{_{Z_{(X,p)}}}),\forall\, b, c\in\Lambda$ with $s(b)=X=t(c)$.
The $\bar R$-$\bar R$-structures on $\hat E_d$ and $\hat E_d^\ast$
ensure that ${\bar E}^\ast_0=\{(B^\ast_X)_{X\in\T}\in{\hat E}^\ast_d\mid cB^\ast_X=B^\ast_Yc,
\forall\, c\in S\}$ by \cite[Proposition 3.4, 3.5]{J} as claimed.

3) Let $\bar E^\ast_1=$Hom$_{R'^{\otimes 2}}(\bar E_1,R'^{\otimes 2})$.
Lemma 2.1.1 (ii) shows:
$$\Hom_{R'\otimes_{k}R'}(L^\ast\otimes_{k}L,R'\otimes_{k}R')
\simeq\Hom_{R'}(L^\ast,R')\otimes_{k}\Hom_{R'}(L,R')\simeq L\otimes_{k}L^\ast.\eqno(2.1\mbox{-}2)$$
By (a3) and a similar argument as in 2),
$\bar E_1^\ast\simeq\sum_{X\in\T,p<q}R'(\textsf e_{_{Z_{(X,p)}}}\otimes_{_{\bar R}}\f_{_{Z_{(X,q)}}})R'$.

4) Combining 2), 3) and noting $\bar E_l^\ast=\sum_{X\in\T,p>q}
R'(\textsf e_{_{Z_{(X,p)}}}\otimes_{_{\bar R}}\f_{_{Z_{(X,q)}}})R'=\{0\}$ by (a4),
$\bar E_0^\ast\oplus\bar E_1^\ast=L\otimes_{\bar R}L^\ast$.
Recall from Formula (2.1-1) that $L\otimes_{\bar R}L^\ast$
corresponds to $B'\otimes_{A'}B'$.
There is an exact sequence $0\rightarrow \bar E_1^\ast\rightarrow\bar E^\ast_1\oplus\bar E^\ast_0\stackrel{p}
\rightarrow\bar E_0^\ast\rightarrow 0$,
$p(\textsf e_{_{Z_{(X,p)}}}\otimes_{_{\bar R}}\f_{_{Z_{(Y,q)}}})=\textsf e_{_{Z_{(X,p)}}}\f_{_{Z_{(Y,q)}}}$,
which corresponds to the map $B'\otimes_{A'}B'\rightarrow B'$ in \cite[4.3 (A3)]{CB1}.
Thus the kernel $J'$ of the map corresponds to
$\bar E^\ast_1$, and $J'$ is projective from $\bar E^\ast$ being so. (A3) follows.

5)\, (A6) concerns the bimodule action on $\bar E_0^\ast\simeq R'$.
Since $\textsf e_{_{Z_{(X,p)}}}\otimes_{_{\bar R}}\f_{_{Z_{(X,q)}}}=0$ for $p>q$ in $\bar E_l^\ast$, (A7) follows.  \hfill$\Box$

\medskip

{\bf Proposition 2.1.5}\, Let $(\mf A,\mf B)$ be a pair, and let $\mf A'$ be given by
Construction 2.1.3. Then the associated bocs $\mf B'$ of $\mf A'$
is the induced bocs of $\mf B$ given by \cite[4.5 Proposition]{CB1}.

\smallskip

{\bf Proof}\, Denote by $\mf C'=(R',\C',\N',\partial')$ the associated
bi-comodule problem of $\mf A'$.

1) $\C'_0=\Hom_{R'}(\K_0',R')$. The isomorphism
$\bar E_0^\ast=\Hom_{R'}(\bar E_0,R')\stackrel{\nu_0^\ast}\rightarrow\Hom_{R'}(\K_0',R')=\C_0'$ with
$\nu_0^\ast$ being the $R'$-dual map of $\nu_0$ in 2.1.3 (ii)
gives the $R'$-quasi-basis ${\F_0'}^\ast=\{e'_{_Z}=\nu_0^\ast(F_Z^\ast)\mid Z\in\T'\}$ of $\C_0'$, see
the proof 2) of Proposition 2.1.4. And ${\F_0'}^\ast$ is $R'$-dual to $\F_0'$ of $\K'_0$.

2)\, $\C_1'=\Hom_{R'^{\otimes 2}}(\K'_{10}\oplus\K'_{11},R'^{\otimes 2})=\C_{10}'\oplus\C_{11}'$. Since
$$\begin{array}{c}\bar E_1^\ast=\Hom_{R'^{\otimes 2}}(\bar E_1,R'^{\otimes 2})\stackrel{\nu_1^\ast}\rightarrow
\Hom_{R'^{\otimes 2}}(\K_{10}',R'^{\otimes 2})=\C_{10}'\end{array}$$
is an isomorphism with $\nu_1^\ast$ being the ${R'}^{\otimes 2}$-dual map of $\nu_1$ in 2.1.3 (ii),
${\F'_1}^\ast=\{{F_i'}^\ast=\nu_1^\ast(F_i^\ast)\mid i=1,\cdots,l\}$
forms an $R'$-$R'$-basis of $\C_{10}$. And ${\F'_1}^\ast$ inherits a linear order from $\mathcal F'_1$.
Since $R'$ is an $R$-$R$-bimodule via the isomorphism
$R'\simeq \K_0'$, by Lemma 2.1.1 (ii) and Formula (2.1-2):
$$\begin{array}{ll}&\C_{11}'=\Hom_{R'^{\otimes 2}}(\K_{11}',R'^{\otimes 2})
=\Hom_{R'^{\otimes 2}}\big((L^\ast\otimes_kL)\otimes_{R^{\otimes 2}}\K_1,R'^{\otimes 2}\big)\\[0.5mm]
\simeq&\Hom_{R'^{\otimes 2}\otimes_{R^{\otimes 2}}R^{\otimes 2}}\big((L^\ast\otimes_kL)\otimes_{R^{\otimes 2}}\K_1,
R'^{\otimes 2}\otimes_{R^{\otimes 2}}R^{\otimes 2}\big)\\[0.5mm]
\simeq&\Hom_{R'^{\otimes 2}}(L^\ast\otimes_kL,R'^{\otimes 2})\otimes_{R^{\otimes 2}}
\Hom_{R^{\otimes 2}}(\K_1,R^{\otimes 2})\simeq (L\otimes_kL^\ast)\otimes_{R^{\otimes 2}}\C_1.\end{array}$$
Write $(\textsf e_{_{Z_{(X'_{j},p)}}}\otimes_{k}\f_{_{Z_{(Y'_{j},q)}}})\otimes_{_{R^{\otimes 2}}}v_j
=\textsf e_{_{Z_{(X'_{j},p)}}}\otimes_{_R}v_j\otimes_{_R}\f_{_{Z_{(Y'_{j},q)}}}=v_{jpq}$.
${\U'}^\ast=\{v_{jpq}\mid 1\leqslant p\leqslant n_{s(v_j)},1\leqslant q\leqslant n_{t(v_j)}; 1\leqslant j\leqslant m\}$
is an ${R'}$-$R'$-quasi-basis of $\C'_{11}$ dual to $\U'$ of $\K_{11}$ given in Construction 2.1.3 (ii).
The $R'$-$R'$-quasi basis of $\C_1'$ is ${\V'}^\ast={\F'_1}^\ast\cup{\U'}^\ast$
dual to  $\V'$ of $\K_1'$.

3)\, $\N_1'=\Hom_{R'^{\otimes 2}}(\M_1',R'^{\otimes 2})\simeq (L\otimes_kL^\ast)\otimes_{R^{\otimes 2}}\M_1$
can be proved in a similar manner as that for $\C_{11}'$.
Write $(\textsf e_{_{Z_{(X_{i},p)}}}\otimes_k\f_{_{Z_{(Y_{i},q)}}})
\otimes_{_{R^{\otimes 2}}}a_i=\textsf e_{_{Z_{(X_{i},p)}}}\otimes_{_R}a_i\otimes_{_R}\f_{_{Z_{(Y_{i},q)}}}=a_{ipq}$.
$\mathcal A'^\ast=\{a_{ipq}\mid 1\leqslant p\leqslant n_{s(a_i)},1\leqslant q\leqslant n_{t(a_i)};h<i\leqslant n\}$
is an ${R'}$-$R'$-quasi basis of $\N'_1$ dual to $\A'$ of $\M_1'$.

4) Formula (1.4-1) shows the formal products of $(\K_0',\C_0')$,  $(\K_1',\C_1')$,
$(\M_1',\N_1')$ respectively:
$$\begin{array}{c}\Upsilon'=\sum_{Z\in\T'}e_{_{Z'}}\ast E'_{Z};\\[0.5mm]
\Pi'=\sum_{j,p,q}v_{jpq}\ast\big((\textsf f_{_{Z_{(X_j'p)}}}\otimes_k\textsf e_{_{Z_{(Y_j',q)}}})
\ast V_j\big)+\sum_{i=1}^l{F_i'}^\ast\ast F'_i\\
=\sum_{j=1}^m(v_{jpq})_{_{n_{_{X'_j}}\times n_{_{Y'_j}}}}
\ast V_j+\sum_{i=1}^l{F_i'}^\ast F'_i=\Pi_1'+\Pi_0';\\
\Theta'=\sum_{i,p,q}a_{ipq}\ast\big((\textsf f_{_{Z_{(X_i,p)}}}\otimes_k\textsf e_{_{Z_{(Y_i,q)}}})\ast A_i\big)=
\sum_{i=1}^n(a_{ipq})_{_{n_{_{X_i}}\times n_{_{Y_i}}}}\ast A_i\end{array}$$
Exhibit the formal
equation $\Theta'\Upsilon'-\Upsilon'\Theta'=(\Pi'_1\Theta'-\Theta'\Pi'_1+\Pi_1'H'-H'\Pi_1')+\Pi'_0\Theta'-\Theta'\Pi'_0$:
$$\begin{array}{c}
\sum_l(a_{lpq}\otimes_{_{R'}} e_{_{t(a_{lpq})}})\ast A_l
-\sum_l(e_{_{s(a_{lpq})}}\otimes_{_{R'}}a_{lpq})\ast A_l\\
=\sum_{l,(i,j)}\Big(\big((v_{ipq})(a_{jpq})\big)\otimes_{R^{\otimes 3}}\eta_{ijl}\Big)\ast A_l-
\sum_{l,(i,j)}\Big(\big((a_{ipq})(v_{jpq})\big)\otimes_{R^{\otimes 3}}\sigma_{ijl}\Big)\ast A_l\\
+\sum_{l,j}\big(\zeta_{jl}\otimes_{R^{\otimes 2}}(v_{jpq})\big)\ast A_l+
\sum_{l,i}\big({F_i'}^\ast(a_{lpq})\big)\ast (F_i'A_l)
-\sum_{l,i}\big((a_{lpq}){F_i'}^\ast\big)\ast (A_lF_i'),\\\end{array}$$
where $(v_{ipq})(a_{jpq})=(\sum_h v_{iph}\otimes_{R'}a_{jhq})$, and other matrix products are similar. Thus
Theorem 1.4.2 shows the differential $\dz'$ in $\mf B'$:
$$\begin{array}{c}\dz'(a_{lpq})=\textsf e_{_{Z_{(X_l,p)}}}\otimes_{_R}\dz(a_l)\otimes_{_R}\f_{_{Z_{Y_l,q}}}\\
+\sum_{p'}\nu^\ast_1(\textsf e_{_{Z_{(X_l,p)}}}
\otimes_{_{\bar R}}\f_{_{Z_{(Y_l,p')}}})\otimes_{_{R'}} a_{lp'q}
-\sum_{q'}a_{lpq'}\otimes_{_{R'}}\nu^\ast_1(\textsf e_{_{Z_{(X_l,q')}}}
\otimes_{_{\bar R}}\f_{_{Z_{(Y_l,q)}}}).\end{array}$$
This coincides with the formula in \cite[4.5 Proposition]{CB1}, i.e. $\mf B'$ is the induced bocs of $\mf B$. \hfill$\Box$

\bigskip
\bigskip
\noindent{\large\bf 2.2  Eight reductions}
\bigskip

In this subsection seven reductions of matrix bimodule problems based on
Definition 2.1.2 and Construction 2.1.3 are introduced, where the last two
do not occur in any references on bocses. And finally, a regularization is presented
as the eighth reduction.

\medskip

{\bf Proposition 2.2.1}\, (Localization) Let $(\mf{A},\mf B)$ be a
pair with $R_X=k[x,\phi(x)^{-1}]$ and
$R'_X=k[x,\phi(x)^{-1}c(x)^{-1}]$ a finitely generated localization
of $R_X$. Define two algebras $\bar R=R$, $R'=R'_X\times\prod_{Y\in\T\setminus\{X\}} R_Y$,
and an $R'$-$\bar R$-bimodule $L=R'$. Then $L$ is admissible.

(i)\, There exists an induced matrix bimodule problem
$\mf{A}'=(R',\K',\M',H')$ of $\mf{A}$  and a fully faithful functor
$\vartheta: R(\mf{A}')\rightarrow R(\mf{A})$.

(ii) The induced bocs $\mf{B}'$ of $\mf B$ given by localization
\cite[4.8]{CB1} is the associated bocs of $\mf A'$.

\medskip

{\bf Proposition 2.2.2}\, (Loop mutation) Let $(\mf{A},\mf B)$
be a pair with the first arrow $a_1:X\mapsto X$, such that $\dz(a_1)=0,X\in\T_0$.
Define a pre-minimal algebra $\bar R=R[a_1]$, a minimal algebra
$R'=R'_X\times\prod_{Y\in\T\setminus\{X\}}R_Y$ with $R'_X=k[x]$, and an $R'$-$\bar R$-bimodule $L=R'$.
Then $L$ is admissible.

(i)\, There exists an induced
matrix bimodule problem ${\mf A}'=(R',\K',\M',\td')$ of $\mf A$, and an
equivalent functor $\vartheta: R({\mf A}')\rightarrow
R(\mf{A})$.

(ii) The induced bocs $\mf{B}'$ of $\mf B$ given by the functor $\theta': A'\rightarrow B'$, with
$\theta'(Y)=Y,\forall\, Y\,\in\T$, $\theta'(a_1)=x$, is the  associated bocs of $\mf A'$ by Proposition 2.1.5.

\medskip

{\bf Proposition 2.2.3} (Deletion) Let $(\mf{A},\mf B)$ be a
pair, $\T'\subset \T$. Define two algebras $\bar R=R$, $R'=\prod_{X\in \T'}R_X$, and an
$R'$-$\bar R$-bimodule $L=R'$. Then $L$ is admissible.

(i)\,  There exists an induced
matrix bimodule problem $\mf{A}'=(R',\K',\M',H')$ of
$\mf{A}$, and a fully faithful functor $\vartheta:
R(\mf{A}')\rightarrow R(\mf{A})$.

(ii) The induced bocs $\mf{B}'$ of $\mf B$ obtained by deletion of $\T\setminus\T'$
\cite[4.6]{CB1} is the associated bocs of $\mf A'$.
\medskip

Let the algebra $R_X=k[x,\phi(x)^{-1}]$, $r\in\mathbb Z^+$, and
$\lambda_1, \ldots, \lambda_s\in k$ with
$\phi(\lambda_i)\ne 0$. Write $g(x)=(x-\lambda_1)\cdots (x-\lambda_s)$.
Define a minimal algebra $S$ and an $S$-$R_X$-bimodule $K$:
$$\begin{array}{c}
S=\big(\prod_{i=1}^s\prod_{j=1}^r
k1_{Z_{ij}}\big)\times k[z,\phi(z)^{-1}g(z)^{-1}];\\[1.5mm]
K=\big(\oplus_{i=1}^s\oplus_{l=1}^r\oplus_{j=1}^l
k1_{Z_{ijl}}\big)\oplus k[z,\phi(z)^{-1}g(z)^{-1}],\,\,Z_{ijl}=Z_{ij};\\[1.5mm]
K(x)=\bar W: K\rightarrow K,\,\, \bar W\simeq \oplus_{i=1}^s
\oplus_{j=1}^r J_j(\lambda_j)1_{Z_{ij}}\oplus (z1_{Z_0}),\end{array}\eqno{(2.2{\mbox{-}}1)}$$
where $\bar W$ is a Weyr matrix over $S$. Let $n=\frac{1}{2}sr(r+1)+1$.
Denote by $\{(i,j,l)\mid 1\leqslant j\leqslant r,1\leqslant l\leqslant j,1\leqslant i\leqslant s\}\cup\{n\}$
the index set of the direct summands of $K$. There is a partition given by
$Z_{ij}=\{(i,j,l)\mid l=1,\cdots,j\}$,
$Z=\{n\}$. An order on the set is defined as
$$\begin{array}{c}(i,j, l)\prec (i',j', l')\Longleftrightarrow i<i'; \,\mbox
{or}\,\, i=i', l<l';\,\mbox {or}\,\, i=i',l=l',
j>j',\end{array}$$
and $n\prec(i,j,l)$. Let $\textsf e_{_{(ijl)}}$ be a $1\times n$, (resp. $\f_{_{(i,j,l)}}$
an $n\times 1$) matrix with $1_{_{Z_{ij}}}$ at the $(i,j,l)$-th component and $0$ at others.
Then $K$ has an $S$-$S$-quasi-basis $\{\textsf e_{_{(ijl)}},\textsf e_{_{n}}
\mid 1\leqslant j\leqslant r,1\leqslant l\leqslant j,1\leqslant i\leqslant s\}$, and $K^\ast=\Hom_S(K,S)$ has
$\{\textsf f_{_{(ijl)}},\textsf f_{_{n}}\mid 1\leqslant j\leqslant r,1\leqslant l\leqslant j,1\leqslant i\leqslant s\}$.
The $S$-quasi-free-module $\bar E_0$, and the $S$-$S$-quasi-free bimodule $\bar E_1$ have the quasi-basis respectively:
$$\begin{array}{c}\{F_{ij}=\sum_{l=1}^{j}\f_{_{(ijl)}}\otimes_{_S}\textsf e_{_{(ijl)}};
F_{n}=\f_{_{n}}\otimes_{_S}\textsf e_{_{n}}\mid
1\leqslant j\leqslant r, 1\leqslant i\leqslant s\},\\[1.5mm]
F_{ijj'l}=\sum_{h=1}^{j'-l+1}\f_{_{(ijh)}}\otimes_k\textsf e_{_{(ij',l+h-1)}},\quad
l=\left\{\begin{array}{c}1,\cdots,j',\,\,\mbox{if}\,\,j>j';\\
2,\cdots,j',\,\,\mbox{if}\,\,j=j';\\
j,\cdots,j',\,\,\mbox{if}\,\,j<j'.\end{array}\right.\end{array}\eqno{(2.2{\mbox{-}}2)}
$$

{\bf Proposition 2.2.4} (Unraveling) Let $(\mf{A},\mf B)$ be a
pair with $R_X=k[x,\phi(x)^{-1}]$. Define two algebras $\bar R=R$, $R'=S\times\prod_{Z\in\T\setminus\{X\}} R_Z$,
and an $R'$-$\bar R$-bimodule $L=K\oplus(\oplus_{Z\in\T\setminus\{X\}} R_Z)$ with $S$ and $K$ given by Formula (2.2-1).
Then $L$ is admissible.

(i)\, There exists an induced matrix bimodule
problem $\mf{A}'=(R',\K',\M',H')$ and a fully faithful
functor $\vartheta: R(\mf{A}')\rightarrow R(\mf{A})$.

(ii) The induced bocs $\mf{B}'$ of $\mf B$ given by unraveling \cite[4.7]{CB1}
is the associated bocs of $\mf A'$.

\smallskip

The picture below shows $\textsf e_{_{(ij1)}}\otimes_{_{\bar R}}\f_{_{(ij'l)}}$ in
$\bar E_1^\ast$ as dotted arrows for $s=1, r=3$:

\vskip 3mm

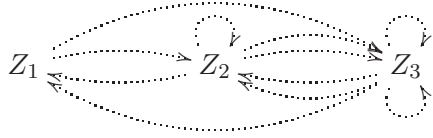
$$\xymatrix{
Z_1\ar@{.>}@/^2pc/[rrrr] \ar@{.>}@/^/[rr]   &&
Z_2\ar@{.>}@/^/[ll]\ar@{.>}@/^1pc/[rr]\ar@{.>}@/^/[rr]\ar@{.>}@(ul,ur) &&
Z_3\ar@{.>}@/^/[ll]\ar@{.>}@/^1pc/[ll] \ar@{.>}@/^2pc/[llll]\ar@{.>}@(ul,ur)\ar@{.>}@(dl,dr) }$$

\vskip 9mm

\noindent

Let an algebra $\bar R_{XY}$, a minimal algebra $S$ and an $S$-$\bar R_{XY}$-module $K$ be defined as follows:
$$\begin{array}{c}\bar R_{XY}:X\stackrel{a_1}\rightarrow Y;\quad S=\prod_{i=1}^3
S_{Z_i},\quad S_{Z_i}=k1_{Z_i}, i=1,2,3;\\[1.5mm]
K_{X}=k1_{Z_2}\oplus k1_{Z_1},\quad K_{Y}=k1_{Z_3}\oplus k1_{Z_2},
\quad K(a_1)=\left(\begin{array}{cc} 0& 1_{Z_2}\\[1.5mm] 0&
0\end{array}\right): K_{X}\rightarrow K_{Y}.\end{array}\eqno{(2.2\mbox{-}3)}$$
Let $Z_{(X,1)}=Z_2=Z_{(Y,2)}$, and $Z_{(X,2)}=Z_1,Z_{(Y,1)}=Z_3$, then
$\{\textsf e_{_{Z_{(X,1)}}},\textsf e_{_{Z_{(X,2)}}}\}$ is an $S$-quasi-basis of $K_X$,
and $\{\textsf f_{_{Z_{(X,1)}}},\textsf f_{_{Z_{(X,2)}}}\}$ is that of $K^\ast_X=\Hom_S(K_X,S)$.
There is a similar observation on $K_Y$.
The $S$-quasi-free module $\bar E_0$, and the $S$-$S$-quasi-free bimodule $\bar E_1$ have
the quasi-basis respectively:
$$\begin{array}{c}F_{Z_1}=\f_{_{Z_{(X,2)}}}\otimes_{_S}\textsf e_{_{Z_{(X,2)}}},\quad
F_{Z_3}=\f_{_{Z_{(Y,1)}}}\otimes_{_S}\textsf e_{_{Z_{(Y,1)}}},\\
F_{Z_2}=(\f_{_{Z_{(X,1)}}}\otimes_{_S}\textsf e_{_{Z_{(X,1)}}},
\f_{_{Z_{(Y,2)}}}\otimes_{_S}\textsf e_{_{Z_{(Y,2)}}});\\
F_{Z_2Z_1}=\f_{_{Z_{(X,1)}}}\otimes_{k}\textsf e_{_{Z_{(X,2)}}},
\quad F_{Z_3Z_2}=\f_{_{Z_{(Y,1)}}}\otimes_{k}\textsf e_{_{Z_{(Y,2)}}}.
\end{array}\eqno{(2.2{\mbox{-}}4)}$$

{\bf Proposition 2.2.5} (Edge reduction)\, Let $(\mf{A},\mf B)$
be a pair with the first arrow $a_1:X\mapsto Y$, such that $X,Y\in\T_0,\dz(a_1)=0$.
Define a pre-minimal algebra $\bar R=R[a_1]$, a minimal algebra $R'=S\times\prod_{Z\in\T\setminus\{X,Y\}} R_Z$,
and an $R'$-$\bar R$-bimodule $L=K\oplus(\oplus_{Z\in\T\setminus\{X,Y\}} R_Z)$ with $S$ and $K$ defined in Formula
(2.2-3). Then $L$ is admissible.

(i)\, There exists an induced matrix bimodule problem $\mf{A'}=(R',\K',\M',H')$, and an equivalent functor $\vartheta:
R(\mf{A}')\rightarrow R(\mf{A})$.

(ii)\, The induced bocs $\mf{B}'$ of $\mf B$ given by edge reduction \cite[4.9]{CB1}
is the associated bocs of $\mf A'$.

\medskip

{\bf Proposition 2.2.6}\, Let $(\mf A,\mf B)$ be a pair with the first arrow $a_1:X\mapsto Y$, such that
$X,Y\in\T_0$, $\dz(a_1)=0$. Set two algebras $\bar R=R[a_1]$, $R'=R$, and an $R'$-$\bar R$-bimodule
$L=K\oplus(\oplus_{U\in\T\setminus\{X,Y\}}R_U)$ with $K:R_X\stackrel{(0)}\rightarrow R_Y$.
Then $L$ is admissible

(i)\, There are an induced matrix bimodule problem $\mf A'=(R',\K', \M', \td')$
with $\K'=\K$, $\M'=\M^{(1)}$, $H'=H$, and an induced fully faithful functors $\vartheta:
R(\mf{A}')\rightarrow R(\mf{A})$.
The subcategory of $R(\mf A)$ consisting of
representations $P$ with $P(a_1)=0$ is equivalent to $R(\mf A')$.

(ii)\,  The induced bocs $\mf{B}'$ of $\mf B$ given by the admissible functor $\theta':A'\rightarrow B'$ with
$\theta'(U)=U,\forall\, U\in\T$ and $\theta'(a_1)=0$ is the associated bocs of $\mf A'$.

\medskip

Let $R_{XY}=R_X\times R_Y$ be a minimal algebra with $R_X=k[x,\phi(x)^{-1}]1_X,R_Y=k1_Y$.
Define an algebra $\bar R_{XY}=R_{XY}[a_1]$ with $a_1:X\rightarrow Y$,
an algebra $S=k[z,\phi(z)^{-1}]$, and an $S$-$\bar R_{XY}$-bimodule $K$ with
$K_{X}=S,\, K_{Y}=S,K(x)=(z),K(a_1)=(1_Z)$. Then there are a $S$-module
$\bar E_0=SF_Z$ with $F_Z$ defined below, and a $S$-$S$-bimodule $\bar E_1=0$.
$$\begin{array}{c}\bar R_{XY}:\xymatrix{X\ar@(ul,dl)[]_x\ar[r]^{a_1}&Y};
\quad S:\xymatrix{Z\ar@(ul,dl)[]_z};\quad K:\xymatrix{S\ar@(ul,dl)[]_z\ar[r]^{(1)}&S};\\[2mm]
F_{Z}=(\f_{_{Z_{X}}}\otimes_{_S}\textsf e_{_{Z_{X}}}, \f_{_{Z_{Y}}}\otimes_{_S}\textsf e_{_{Z_{Y}}})=(1_Z,1_Z).
\end{array}\eqno{(2.2\mbox{-}5)}$$

{\bf Proposition 2.2.7}\,  Let $(\mf A,\mf B)$ be a pair with the first arrow $a_1:X\mapsto Y$ and $\dz(a_1)=0$.
Define a pre-minimal algebra $\bar R=R[a_1]$, a minimal algebra $R'=S\times\prod_{U\in\T\setminus\{X,Y\}} R_U$,
and an $R'$-$\bar R$-bimodule $L=K\oplus(\oplus_{U\in\T\setminus\{X,Y\}} R_U)$,
where $S$ and $K$ are given by Formula (2.2-5). Then $L$ is admissible.

(i)\, There are an induced matrix bimodule problem $\mf A'$, and an induced fully faithful functors $\vartheta:
R(\mf{A}')\rightarrow R(\mf{A})$. The subcategory of $R(\mf A)$ consisting of representations $P$
with $P(a_1)$ invertible is equivalent to $R(\mf A')$.

(ii)\,  The induced bocs
$\mf{B}'$ given by the admissible functor $\theta':A'\rightarrow B'$ with $\theta'(X)=Z,\theta'(Y)=Z$;
$\theta'(x)=z,\theta'(a_1)=(1)$, is the associated bocs of $\mf A'$.

\medskip

Let $\mf{A}=(R,\K,\M,\td)$ be a matrix
bimodule problem, $\mf C=(R,\C,\N,\partial)$ be the associated
bi-comodule problem, and $\mf B$ the bocs of $\mf C$. Then
$$\begin{array}{c}\td(V_1)=A_1+\sum_{l>1}\zeta_{1l} A_l\,\,\mbox{and}\,\,\td(V_j)\in
\M^{(1)}_1\,\,\mbox{for}\,\,j\geqslant 2\,\, \mbox{in}\,\, \mf A\\[1mm]
\Longleftrightarrow  \partial (a_1)=v_1\,\, \mbox{in}\,\, \mf C
\Longleftrightarrow \delta(a_1)=v_1\,\,\mbox{in}\,\,\mf B.\end{array}\eqno{(2.2\mbox{-}6)}$$
In fact, since $\partial(a_1)=\sum_{j=1}^m\zeta_{j1}v_{j}$,
we have $\partial(a_1)=v_1$, if and only if $\zeta_{11}=1_{s(a_1)}\otimes_k1_{t(a_1)}$ and
$\zeta_{j1}=0$ for all $j\geqslant 2$, if and only if
$\td(V_1)=A_1+\sum_{l>1}\zeta_{1l} A_l$ and $\td( V_j)\in
\M^{(2)}$ for all $j\geqslant 2$, since $\td(V_{j})=\zeta_{j1}A_1+
\sum_{i>1}\zeta_{ji}A_i$. On the other hand, noting $\iota_1(a_1)=0$ and $\tau_1(a_1)=0$
by triangularity of $\mf C$, thus $\delta(a_1)=v_1$ in $\mf B$, if and
only if $\partial(a_1)=v_1$ in $\mf C$ by the definition of $\dz_1:\Gamma\rightarrow\bar\Omega$
given below Lemma 1.2.5.

\medskip

{\bf Remark}\, Let $\mf{A}, \mf C$ be given as above with $\partial(a_1)=v_1$. Then

(i)\, $\K^{(1)}=\K_0\oplus(\oplus_{j=2}^m
\bar\Delta\otimes_{_{R^{\otimes 2}}} V_j)$ is a sub-algebra of $\K$,
and $\M^{(1)}=\oplus_{i=2}^n\bar\Delta\otimes_{_{R^{\otimes 2}}}a_i$ is a
$\K^{(1)}$-$\K^{(1)}$-sub-bimodule;

(ii)\, $\C^{(1)}=\bar\Delta\otimes_{_{R^{\otimes 2}}}v_1$ is a
coideal of $\C$, and $\C^{[1]}=\C/\C^{(1)}$ is a
quotient coalgebra;
$\N^{(1)}=\bar\Delta\otimes_{_{R^{\otimes 2}}}a_1$ is a
$\C$-$\C$-sub-bi-comodule, and $\N^{[1]}=\N/\N^{(1)}$ is a
$\C^{[1]}$-$\C^{[1]}$-quotient bi-comodule.

\smallskip

{\bf Proof}\, (i)\, By the triangularity (1.2-2), $\td(V_iV_j)=\td(V_i)V_j+V_i\td(V_j)\in \M^{(1)}$
and $\forall\, V_i,V_j\in\V$:
$$\begin{array}{l}\td(V_iV_j)= \td(\sum_{l=1}^m
\gamma_{ijl}\otimes_{_{R^{\otimes 2}}}V_l)
=\sum_{l=1}^m \gamma_{ijl}\otimes_{_{R^{\otimes 2}}}\td(V_l)
=\sum_{l,p} (\gamma_{ijl}\otimes_{_{R^{\otimes 2}}}\zeta_{lp})\otimes_{_{R^{\otimes 2}}}
A_p.\end{array}$$
The coefficient of $A_1$ is $\sum_{l}
(\gamma_{ijl}\otimes_{_{R^{\otimes 2}}}\zeta_{l1})=0$, where $\zeta_{11}=1$, $\zeta_{l1}=0$ for
$l>1$ by hypothesis, so that $\gamma_{ij1}=0$ for all
$1\leqslant i,j\leqslant m$. Therefore
$V_iV_j=\sum_{l>1}\gamma_{ijl}\otimes_{_{R^{\otimes 2}}}V_l\in\K^{(1)}$
and hence $\K^{(1)}$ is a subalgebra of $\K$. $\M^{(1)}$ is a
$\K^{(1)}$-$\K^{(1)}$-bimodule deduced from the triangularity (1.2-2) easily.

(ii) Since $\mu(v_1)=\mu(\partial(a_1))
=(\partial\otimes\id)(\iota(a_1))+(\id\otimes\partial)(\tau( a_1))
\stackrel{(1.2\mbox{-}3)}=(\partial\otimes \id)(e_{s(a_1)}\otimes_{_R}a_1)
+(\id\otimes\partial)(a_1\otimes_{_R}e_{t(a_1)})=0$,
$\C^{(1)}$ is a coideal of $\C$. $\N^{(1)}$ is a
$\C$-$\C$-sub-bi-comodule deduced from (1.2-3) easily. \hfill$\Box$

\medskip

{\bf Proposition 2.2.8} (Regularization)\, Let
$(\mf{A},\mf B)$ be a pair with $\delta(a_1)=v_1$.

(i) There is an induced matrix bimodule problem $\mf{A}'=(R,
\K^{(1)}, \M^{(1)}, H)$ of $\mf A$, and an
equivalent functor $\vartheta: R(\mf{A}') \rightarrow
R(\mf{A})$.

(ii) The induced bocs $\mf{B}'$ of $\mf B$ given by regularization
\cite[4.2]{CB1} is the associated bocs of $\mf A'$.\hfill$\Box$

\smallskip

{\bf Proof}\, (i) $\mf A'$ is a matrix bimodule problem by Remark (i) above.
Note that $R'=R,\T'=\T$. Taken any $P\in R(\mf A)$
of size vector $\m$, let $f=\sum_{X\in\T}I_{m_X}\ast E_X+P(a_1)\ast V_1$, then
$P'=f^{-1}Pf\in R(\mf A')$. Therefore, $\vartheta$ is an equivalent functor.

(ii) $\mf C'=(R,\C^{[1]}, \N^{[1]}, \bar\partial)$ with $\bar\partial$ induced
from $\partial$ is the associated bi-comodule problem of $\mf A'$ by
Remark (ii) above. Thus the associated
bocs $\mf B'$ of $\mf A'$ is given by regularization from $\mf B$.\hfill$\Box$

\medskip

Let $(\mf A,\mf B)$ be a pair with a layer
$L=(R;\omega;a_1,\cdots,a_n;v_1,\cdots,v_m)$ in $\mf B$. Suppose the first arrow $a_1: X\mapsto Y$ with
$\dz(a_1)=\sum_{j=1}^m f_{j}(x,y)v_{j}\ne 0$. In order to obtain
$\dz(a_1)=h(x,y)v'_1$, we make the following base change:
$$\begin{array}{c}(v_{1}',\cdots,v_{m}')=(v_{1},\cdots,v_{m})F(x,y),
\end{array}\eqno{(2.2{\mbox{-}}7)}$$
where $F(x,y)\in\IM(R\otimes_kR)$ is invertible. When $X\in \T_0$ or $Y\in \T_0$, $R$ is preserved;
but when $X, Y\in \T_1$, some localization $R_{X}'=R_{X}[c(x)^{-1}]$
(resp. $R_{Y}'=R_{Y}[c(y)^{-1}]$) is needed \cite[$\S
5$]{CB1}. Consequently, we have a base change of $\K_1$ dually given by
$$\begin{array}{c}(V'_{1},\cdots,V'_{m})=(V_{1},\cdots,V_{m})F(x,y)^{-T}.\end{array}\eqno{(2.2{\mbox{-}}8)}$$

Finally, a simple fact according to all the reductions defined above is mentioned
to end the subsection. Let us start from a matrix bimodule problem
$\mf{A}^0=(R^0,\K^0,\M^0,H=0)$ with $\T^0$ trivial, and after a series of reductions,
an induce matrix bimodule problem $\mf{A}'=(R',\K',\M',H')$ is obtained. Then for any $X\in \T'$,
$H_X'=(h_{ij}(x))$, any entry $h_{ij}(x)=a_{ij}+b_{ij}x\in k[x]$.

\bigskip
\bigskip
\noindent{\large\bf 2.3 Canonical forms}
\bigskip

In this subsection, a canonical form (cf.\cite{S}) for each representation of a matrix
bimodule problem is calculated; and  a notion of reduction blocks is defined.

\medskip

{\bf Convention 2.3.1}\, Suppose $\mf{A}$ is a matrix bimodule problem,
$\mf{A}'$ an induced matrix bimodule problem and
$\vartheta: R(\mf A')\rightarrow R(\mf A)$ an induced functor.
Let $\m'$ be a size vector over $\T'$ of $\mf A'$.
A size vector $\m=(m_1,m_2,\ldots,m_t)$ over $\T$ of $\mf A$ based on $\m'$ is defined:

(i) for regularization, loop mutation, localization, and Proposition 2.2.6, set
$\m=\m'$;

(ii) for deletion, set $m_i=m'_i$ if $i\in X, X \in \T'$, and $m_i=0$ if
$i\in X, X\in \T\setminus \T'$;

(iii) for edge reduction, set $m_i=m'_i$ if $i\in Z, Z\ne X,Y$,
$m_i=m'_{Z_1}+m'_{Z_2}$ if $i\in X$, and $m_i=m'_{Z_2}+m'_{Z_3}$ if $i\in
Y$; For proposition 2.2.7, set $m_{_X}=m_{_{Z}}',m_{_Y}=m'_{_Z}$;

(iv) for unraveling, set  $m_i=m'_i$ if $i\not\in X$, and
$m_i=\sum_{i=1}^s\sum_{j=1}^r jm'_{Z_{ij}}+m'_{Z_0}$ if $i\in X$.

Then $\m$ is said to be the {\it size vector determined by $\m'$},
and is denoted by $\vartheta(\m')$.

\medskip

Let $\mf A=(R,\K,\M,H)$ be a matrix bimodule problem with $\T$ being trivial,
and $\m$ be a size vector. For the sake of simplicity,
we write
$$\begin{array}{c}H_\m(k)=\sum_{X\in \T} H_X(I_{m_{_X}}),\quad
H(k)=\sum_{X\in \T}H_X(1).\end{array}\eqno{(2.3\mbox{-}1)}$$

Let $P$ be a representation of size vector $\m$ in $R(\mf A)$. Then Definition 1.3.4 shows:
$$\begin{array}{c}P=H_{\m}(k)+\sum_{i=1}^nP(a_i)\ast A_i.\end{array}\eqno{(2.3\mbox{-}2)}$$
Let $\T'=\{i\mid m_i\ne 0\}$, and the induced bimodule problem $\mf A'$ be
given by a deletion of $\T\setminus \T'$. Then the size vector
$\m'=(m_i\mid m_i\ne 0)$ is sincere over $\T'$. It may be assumed that $\m$ is sincere over $\T$ in the sequel.
We will find an induced matrix bimodule problem $\mf A'$ given by minimal steps of reducions, and an object
$P'\in R(\mf A')$ of sincere size vector $\m'$ over $\T'$, such that $\vartheta(P')\simeq P$ under the induced functor
$\vartheta:R(\mf A')\rightarrow R(\mf A)$.
Let $\mf B$ be the associated bocs of $\mf A$ with the first arrow $a_1: X\rightarrow Y$.
There are three possibilities.

(i)\, If $\dz(a_1)=v_1$, we proceed with a regularization, and obtain an induced matrix bimodule problem $\mf{A}'$. Set
$$B=(\emptyset)_{m_{_X}\times m_{_Y}},\quad G=(\emptyset)_{1\times 1},\eqno{(2.3\mbox{-}3)}$$
where $\emptyset$ indicates a distinguished zero entry or block. Suppose
$P'$ is given in the proof (i) of Proposition 2.2.8 with $\T'=\T,\m'=\m$. Then $P'(a_1)=B$ and $\vartheta(P')\simeq P$ in $R(\mf A)$.

(ii)\, If $\dz(a_1)=0$ and $X=Y$, suppose $P(a_1)\simeq J=
\oplus_{i=1}^s(\oplus_{j=1}^r J_j(\lambda_i)^{e_{ij}}),e_{ij}\geqslant 0$,
a Jordan form over $k$ with the maximal
size $r$ of the Jordan blocks. We first proceed with a loop
mutation $a_1\mapsto (x)$, then with an unraveling for the polynomial
$g(x)=(x-\lambda_1)\cdots(x-\lambda_s)$ and the integer
$r$, thus an induced matrix bimodule problem $\mf A_1$ of $\mf A$ is obtained.
Let $f_{X}\in\IM_{m_{_X}}(k)$ be invertible, such that
$$B=f_{X}^{-1}P(a_1)f_{X}=W;\quad G=\bar W,\eqno{(2.3\mbox{-}4)}$$
where $W$ is a Weyr matrix over $k$, $\bar W$ is that over $R'$ similar
to $\oplus_{i=1}^s\oplus_{e_{ij}>0} J_j(\lambda_i)1_{Z_{ij}}$.
Deleting a set of vertices $\{Z_0\}\cup\{Z_{ij}\mid e_{ij}=0\}$ from $\mf A_1$,
an induced problem $\mf{A}'$ of $\mf A$ is obtained. Let $\m'=(m'_i)_{i\in T'}$ be a size vector over $\T'$ with $m'_{_Z}=m_{_Z}$
for $Z\in\T\setminus\{X\}$, and $m'_{_{Z_{ij}}}=je_{ij}$, then $\m'$ is sincere. Let $f=f_{X}\ast E_{X}+
\sum_{Z\in\T\setminus\{X\}}I_{m_{_Z}}\ast E_Z$, and $P'=f^{-1}Pf$ with the size vector $\m'$ in $R(\mf A')$. Then
$P'(a_1)=B$ and $P\simeq\vartheta(P')$ in $R(\mf A)$.

(iii) If $\dz(a_1)=0$ and $X\ne Y$, we
proceed with an edge reduction for $\mf{A}$ and obtain an induced
problem $\mf A_1$ with the vertex set $\T_1$. If $\rank(P(a_1))=r$, let
$f_X\in\IM_{m_{_X}}(k),f_Y\in \IM_{m_{_Y}}(k)$ be invertible, such that
$$\begin{array}{c}
B=f_{X}^{-1}P(a_1) f_{Y}={{\,0\,\, I_r}\choose{0\,\,\
0}}_{m_{_X}\times m_{_Y}}.\\[2mm]
\mbox{Set}\,\,\,G=\mbox{\ding{172}}\, (0),\,\,\mbox{\ding{173}}\, (1_{Z_2}),\,\,\mbox{\ding{174}}\, (0\, 1_{Z_2}),\,\,
\mbox{\ding{175}}\, {{1_{Z_2}}\choose{0}},\,\,
\mbox{or \ding{176}}\, {{\,0\,\, 1_{Z_2}}\choose{0\,\,\,\,
0\,\,}},\end{array}\eqno{(2.3\mbox{-}5)}
$$
where the five cases of $G$ are obtained by deleting a subset $\hat\T\subset\T_1$ from $\mf A_1$:
\ding{172} $\hat\T=\{Z_2\}$ for $r=0$; now suppose $r>0$, \ding{173} $\hat\T=\{Z_1,Z_3\}$ for $m_{_X}=r=m_{_Y}$;
\ding{174} $\hat\T=\{Z_1\}$ for $m_{_X}=r, m_{_Y}>r$; \ding{175} $\hat\T=\{Z_3\}$
for $m_{_X}>r, m_{_Y}=r$; \ding{176} $\hat\T=\emptyset$ for $m_{_X},m_{_Y}>r$.
An induced matrix bimodule  problem $\mf{A}'$ of $\mf A$ given by $a_1\mapsto G$ is obtained.
Let $\m'=(m'_i)_{i\in T'}$ be a size vector over $\T'$, with $m'_{_Z}=m_{_Z}$ for $Z\in
\T\setminus \{X,Y\}$, $m_{_{Z_1}}'=m_{_X}-r$, $m_{_{Z_2}}'=r$ and
$m_{_{Z_3}}'=m_{_Y}-r$, thus $\m'$ is sincere over $\T'$. Let $f=f_{X}\ast E_{X}+f_{Y}\ast E_{Y}+
\sum_{Z\in\T\setminus\{X,Y\}}I_{m_{_Z}}\ast E_Z$,
and $P'=f^{-1}Pf$ with the size vector $\m'$ in $R(\mf A')$. Then $P'(a_1)=B$ and
$P\simeq\vartheta(P')$ in $R(\mf A)$.

\medskip

{\bf Lemma 2.3.2} (cf.\cite{S})\, Let $\mf A=(R,\K,\M,H)$ be a matrix bimodule problem
with $\T$ trivial, and let $P$ be given in Formula (2.3-2).
Then there exists an induced matrix bimodule problem
$\mf{A}'=(R',\K',\M',H')$ given by one of the following three reductions:

(i) Regularization,

(ii) Loop reduction: a loop mutation, then a unraveling, followed by a deletion.

(iii) Edge reduction: first an edge reduction, followed by a deletion,

\noindent There is a representation $P'$ of sincere size vector $\m'$ over $\T'$ in
$R(\mf A')$, such that $P\simeq\vartheta(P')$ in $R(\mf A)$ under the fully faithful functor $\vartheta:
R(\mf{A}')\rightarrow R(\mf{A})$. According to Formulae (2.3-3)--(2.3-5):
$$\begin{array}{c}P'=H_{\m}(k)+
B\ast A_1+\sum_{i=1}^{n'}P(a'_i)\ast A'_i.\end{array}$$

The procedure may be
said that {\it the reduction is given by $a_1\mapsto G$} in one of Formulae (2.3-3)--(2.3-5),
and $G$ is called the {\it reduction block from $\mf A$ to $\mf A'$}.

Applying Lemma 2.3.2 repeatedly, the following theorem is obtained by induction.

\medskip

{\bf Theorem 2.3.3} (cf.\cite{S})\, Let $\mf{A}=(R,\K,\M,H=0)$ be a matrix bimodule
problem with $\T$ trivial. Let $P\in R(\mf A)$ be a representation of sincere size vector $\m$.
Then there exists a unique sequence of matrix bimodule problems:
$$
\begin{array}{llllllll}&\mf{A}=\mf{A}^0,&\,\, \mf{A}^1, &\cdots,
&\mf{A}^{i}, &\mf{A}^{i+1},& \ldots,&\mf{A}^s \\
&&(G^1,&\ldots,&G^{i},&G^{i+1},&\cdots,&G^s\,\,)\end{array}\eqno{(\ast)}$$
where $\mf A^{i+1}$ is obtained from $\mf A^{i}$ by $a^i_1\mapsto G^{i+1}$
defined by one of Formula (2.3-3)--(2.3-5).
There is also a unique sequence of representations:
$$\begin{array}{llllllll}&P^0=P,&\,\, P^1,&\ldots,&P^{i},&P^{i+1},&\ldots,&P^s\\
&&(B^1,&\ldots,& B^{i},&B^{i+1},&\ldots,&B^s\,\,)\end{array}\eqno{(\ast\ast)}$$
where $P^i(a^i_1)=B^{i}$ is defined by one of Formulae (2.3-3)--(2.3-5).
Let $\vartheta^{i,i+1}: R(\mf{A}^{i+1})\rightarrow R(\mf{A}^{i})$ be
the induced functor. There is a representation $P^{i+1}\in R(\mf A^{i+1})$ of sincere size vector $\m^{i+1}$
with $\vartheta^{i,i+1}(P^{i+1})\simeq P^{i}$ for $0\leqslant i<s$.

\medskip

Write for $i<j$ the composition of induced functors $\vartheta^{ij}
=\vartheta^{i,i+1}\cdots\vartheta^{j-1,j}:R(\mf A^j)\rightarrow R(\mf
A^{i})$. Denote by $A_1^i$ the first quasi-basis matrix of $\M^i_1$ in $\mf A^i$, then
$P^{i+1}=H^i_{\m^{i}}(k)+B^{i+1}\ast A^{i}_1+\sum_{j=2}^{n^{i+1}}M^{i+1}(a^{i+1}_j)\ast A^{i+1}_j$.
Using the formula inductively:
$$\begin{array}{c}\vartheta^{0s}(H^s_{\m^{s}}(k))=
\sum_{i=0}^{s-1} B^{i+1}\ast A^i_1\in R(\mf A).\end{array}\eqno{(2.3\mbox{-}6)}$$
In particular, if $\mf A^s$ is minimal, then
$P^s=H^s_{\m^s}(k)$.
In this case, the matrix $\vartheta^{0s}(P^s)$ is called the {\it canonical form} of
$P$, and denoted by $P^{\infty}$.
The entry $``1"$ appearing in $B^{i+1}$ of $P^\infty$, which is not an
eigenvalue when $B^{i+1}$ being a Weyr matrix, is called a
{\it link} of $P^\infty$. And denote by $l(P^\infty)$ the number of the links in $P^\infty$.

\medskip

{\bf Corollary 2.3.4} \cite{S,XZ}\, The canonical form of any representation $P$ over
a matrix bimodule problem $\mf{A}=(R, \K,\M,H=0)$ with $R$ trivial is uniquely determined.
Moreover,

(i)\, for any $P, Q\in R(\mf{A})$, $P\simeq Q$ if and only if $P$
and $Q$ have the same canonical form;

(ii)\, $P$ is indecomposable if and only if $l(P^\infty)=\mbox{dim}(P)-1$.

\medskip

{\bf Corollary 2.3.5}\,  Let $\mf{A}=(R,\K,\M,H)$ be a matrix
bimodule problem with $\T$ trivial, let $\mf{A}'=(R',
\K',\M',H')$ be an induced matrix bimodule problem
obtained by a series of reductions with an induced functor $\vartheta: R(\mf{A}')\rightarrow
R(\mf{A})$. If $\T'$ is trivial and $P=\vartheta(H'(k))$ with a sincere size vector $\m$ over $\T$ in $R(\mf A)$,
then there is a unique reduction sequence $(\ast)$ performed for $P$ by Theorem 2.3.3.

\medskip

We conclude this subsection with a definition of some reduction block $G^{i+1}_s(R^s)$ (or $G^{i+1}_s$ for short) over $R^s$.
Under the hypothesis of Corollary 2.3.5, set $\mf A'=\mf A^s$ in the sequence $(\ast)$.
Let $\m^s=(1,\cdots,1)$ and $\m^{i}=\vartheta^{is}(\m^s)$.
Write $s(a^i_1)=X^i$ and $t(a^i_1)=Y^i$.
A matrix $G_s^{i+1}\in\IM_{m^i_{X^i}\times m^i_{Y^i}}(R^s)$ is determined by

\subitem\hspace{-6mm}(i) $G_s^{i+1}(k)=G_s^{i+1}\otimes 1\in \IM_{m^i_{X^i}\times m^i_{Y^i}}(R^s)\otimes_{R^s}k\simeq
\IM_{m^i_{X^i}\times m^i_{Yi}}(k)$ is equal to $B^{i+1}$ given in Theorem 2.3.3;

\subitem\hspace{-6mm}(ii) write the matrix $G^{i+1}_s\ast A^i_1=(g_{pq})\in\IM_{t^s}(R^s)$,
if $g_{pq}\ne 0$, and $p\in X^s$ (or equivalently, ${q}\in X^s$), then $g_{pq}\in R^s_{X^s}$.

\noindent And $G^{i+1}_s$ for $i=0,\cdots,s-1$ are said to be the
{\it reduction blocks of $H^s$}. Furthermore,
$$\begin{array}{c}H^s=\sum_{i=0}^{s-1}G^{i+1}_s\ast A^i_1\quad \mbox{and}\quad
\vartheta^{0s}(H^s(k))=\sum_{i=0}^{s-1}G^{s+1}_i(k)\ast A^i_1.\end{array} \eqno{(2.3\mbox{-}7)}$$
For the sake of convenience, a links of $\vartheta^{0s}(H^s(k))$ is also said to be a link of $H^s$.
Thus $\mf A^s$ is local if and only if
$l(H^s)=$dim$\big(\vartheta^{0s}(H^s(k))\big)-1$.

\bigskip
\bigskip
\noindent{\large\bf 2.4 Defining systems}
\bigskip

We introduce a concept of defining systems in this subsection. There exist
two sorts of systems used in different situations in order to construct induced
matrix bimodule problems in a reduction sequence.

\medskip

Let $B=(b_{ij})_{t\times t}$ and $C=(c_{ij})_{t\times t}$ be two
$t\times t$ matrices over $k$. Given $1\leqslant p,q\leqslant t$, the notation
$B\equiv_{\prec(p,q)}C$ (resp. $B\equiv_{(p,q)}C$, $B\equiv_{\preccurlyeq(p,q)}C$) means that
$b_{ij}=c_{ij}$ for any $(i,j)\prec(p,q)$ (resp.
$(i, j)=(p,q),\, (i,j)\preccurlyeq(p,q)$). One can define the similar notations
for partitioned matrices.

Let $(\mf A,\mf B)$ be a pair with $\T$ trivial, $H=0$, the $R$-$R$-quasi-basis $\V=\{V_1,\cdots,V_m\}$ of $\K_1$
and $\A=\{A_1,\cdots,A_n\}$ of $\M_1$. Denote by $(\textsf p_j,\textsf q_j)$ the leading position
of $A_j$ for $j=1,\cdots,n$. Suppose there exists a sequence of reductions in the sense of Lemma 2.3.2:
$$
(\mf{A},\mf{B})=(\mf{A}^0,\mf{B}^0), (\mf{A}^1,\mf{B}^1), \cdots, (\mf{A}^{i},\mf{B}^i), (\mf{A}^{i+1},\mf B^{i+1}), \cdots,
(\mf{A}^r,\mf B^r),\cdots,({\mf{A}^s},\mf B^s).\eqno{(2.4\mbox{-}1)}
$$
For each $0\leqslant i\leqslant s$ in the sequence, a matrix equation is defined by
$$\begin{array}{c}\mbox{\IE}^{i}: \,\,\ \Phi_{\m^{i}}
H^{i}(k)\equiv_{\prec(p^{i},q^{i})}H^{i}(k)\Phi_{\m^{i}}\\[1mm]
\mbox{with}\quad
\Phi_{\m^i}=\sum_{X\in\T}Z_X\ast E_X +\sum_{j=1}^m Z_j\ast V_j,\end{array}\eqno
{(2.4\mbox{-}2)}$$
where $(p^i,q^i)$ is the leading position of $A_1^i$ of $\M^i_1$ in the $i$-th pair, $\m^i=\vartheta^{0i}(1,\cdots,1)$ is the
size vector of $\vartheta^{0i}(H^i(k))\in R(\mf A)$ over $\T$, $Z_X=(z_{pq}^X)_{m_{_X}\times m_{_X}}$, $\forall\, X\in \T$, and
$Z_j=(z_{pq}^j)_{m_{s(v_j)}\times m_{t(v_j)}}$ for all quasi-base matrices $V_j$ of $\K_1$ in $\mf A$,
$z^X_{pq},\; z^j_{pq}$ are pairwise different variables over
$k$. $\Phi_{\m^i}$ is called a {\it variable matrix}.
The system of linear equations in $\IE^i$, which consists of  equations locating in the
$(\textsf p_j,\textsf q_j)$-th block for $j=1,\cdots,n$, is said to be a {\it  defining system of $\K^i$},
and is denoted still by $\IE^i$.

\medskip

{\bf Theorem 2.4.1}\, The solution space of the defining system $\IE^i$ in Formula (2.4-2)
is the $k$-vector space spanned by the quasi-basis of $\K^{i}_0\oplus\K^i_1$ for $0\leqslant i\leqslant s$
in the sequence (2.4-1).

\smallskip

{\bf Proof}\, Since $H^0=H=0$, our theorem holds true for $i=0$. Suppose the theorem is true
for the defining system $\IE^i$, now consider the defining system $\IE^{i+1}$.

In the case of Regularization, $\m^{i+1}=\m^{i}$, the $k$-vector space spanned by the quasi-basis of
$\K_0^{i+1}\oplus \K_1^{i+1}$ is just the solution space of the equations of
$\IE^{i}$ and the equation $\Phi_{\m^{i}}H^{i}(k)\equiv_{(p^{i},q^{i})}H^{i}(k)\Phi_{\m^{i}}$,
which form the equation system $\IE^{i+1}:\Phi_{\m^{i+1}}H^{i+1}(k)\equiv_{\prec(p^{i+1},q^{i+1})}H^{i+1}(k)\Phi_{\m^{i+1}}$.

In the case of Loop or Edge reduction, set $a_1:X\rightarrow Y$ with $X=Y$ for loop reduction;
denote by $\n=\vartheta^{i,i+1}(1,\cdots,1)$, the size vector of $\vartheta^{i,i+1}(H^{i+1}(k))$
over $\T^i$. Then the size vector $\m^{i+1}=\vartheta^{0i}(\n)$ over $\T$. Since $G^{i+1}=L(a_1)$,
write $G^{i+1}(k)=L(a_1)(k)$. The $k$-vector space spanned by the quasi-basis of
$\K_0^{i+1}\oplus \K_1^{i+1}$ is the solution space of the matrix equations partitioned under $\T^i$:
$$
\left\{\begin{array}{l}(\Phi_{\m^{i}})_{\n}H^{i}_{\n}(k)
\equiv_{\preccurlyeq(p^{i},q^{i})}H^{i}_{\n}(k)(\Phi_{\m^{i}})_{\n},\\[1.5mm]
(\Phi_{\m^{i}})_{\n, X} L(a_1)(k)=
L(a_1)(k)(\Phi_{\m^{i}})_{\n, Y},\end{array}\right.$$
where if $\Phi_{\m^i}=(x_{pq})$, then $(\Phi_{\m^i})_{\n}=(X_{pq})$ with
$X_{pq}=(x_{pq,\alpha\beta})_{n_p\times n_q}$; since $\dz(a_1)=0$, the $(p^i,q^i)$-th equation
of $\IE^i$ is a linear combination of previous equations, $``\preccurlyeq(p^{i},q^{i})"$
can be used in the first formula; in the second one $(\Phi_{\m^{i}})_{\n, X}$ stands for the $(j,j)$-th
block of $(\Phi_{\m^{i}})_{\n}$ with $j\in X$. Since
$(\Phi_{\m^i})_{\n}=\Phi_{\m^{i+1}}$ and $H^{i+1}(k)=H^i_{\n}(k)+L(a_1)(k)\ast A^i_1$,
the equation system above is just
$\IE^{i+1}:\Phi_{\m^{i+1}}H^{i+1}(k)\equiv_{\prec(p^{i+1},q^{i+1})}H^{i+1}(k)\Phi_{\m^{i+1}}$,
where $(p^{i+1},q^{i+1})$ is the index followed by the biggest index over $\T^{i+1}$ of
the $(p^i,q^i)$-block partitioned under $\T^i$. Our theorem is proved by induction.\hfill$\Box$

\medskip

For an example, see 2.4.5 (iv) below. The theorem implies the following fact obviously.

\smallskip

{\bf Corollary 2.4.2}\, $\delta(a_1^i)=0$ in ${\mf B}^i$ if and only if the $(p^i,q^i)$-equation of $\IE^i$,
is a linear combination of the equations of ${\IE}^i$, namely, the equations locating
before $(p^i,q^i)$.

\medskip

Next, we give a deformed system based on the defining system.
Fix some $0<r<s$ in the sequence (2.4-1). Suppose $\mf A^r=(R^r,\K^r,\M^r, H^r)$, and $\{V^r_1,\ldots,V^r_{m^r}\}$
is a normalized quasi-basis of $\K^r_1$. If $i\geqslant r$, set a size vector over $\m^{ri}=\vartheta^{ri}(1,\cdots,1)$
over $\T^r$, where $(1,\cdots,1)$ is a size vector over $\T^i$. Let
$Z^{ri}_{Y^r}$ be variable matrices of size $m^{ri}_{_{Y^r}},\forall\,Y^r\in \T^r$, and $Z^{ri}_j$ be those
of size $m^{ri}_{s(v^r_j)}\times m^{ri}_{t(v^r_j)}, 1\leqslant j\leqslant m^r$, then set a variable matrix
$\Psi_{\m^{ri}}=\sum_{Y^r\in\T^r}Z^{ri}_{Y^r}\ast E_{Y^r}^r +\sum_{j=1}^{m^r}
Z^{ri}_j\ast V^r_j$. Let
$H^i=H^i_1+H^i_2$ with $H_1^i=\sum_{j=1}^{r-1}G_i^{j+1}\ast A^j_1,
H^i_2=\sum_{j=r}^{i-1}G_i^{j+1}\ast A^j_1$. Then the matrix equation $\Psi_{\m^{ri}}
H^i(k)\equiv_{\prec(p^i,q^i)} H^i(k) \Psi_{\m^{ri}}$ can be rewritten as
$$\begin{array}{c}{\IF}^{ri}: \Psi_{\m^{ri}}
H^i_2(k)\equiv_{\prec(p^i,q^i)} \Psi^0_{\m^{ri}}+H^i_2(k) \Psi_{\m^{ri}},\\[1mm]
\Psi_{\m^{ri}}^0=H^i_1(k)\Psi_{\m^{ri}}-\Psi_{\m^{ri}} H^i_1(k).\end{array} \eqno{(2.4\mbox{-}3)}$$

{\bf Corollary 2.4.3}\, The equation system $\IF^{ri}$ is equivalent to $\IE^i$.
And $\delta(a_1^i)=0$ in ${\mf B}^i$ if and only if the $(p^i,q^i)$-th equation
${\IF}^i$ is a linear combination of equations of
${\IF}^i$, namely, the equations locating before $(p^i,q^i)$.

\medskip

The above theorem and corollaries will be used in Subsections 5.2--5.4 to calculate the differentials
of bocses given by some bordered matrices.
Sometimes, it is difficult to determine the dotted arrows in the induced bocs
after some reductions. Instead, we may consider a system of equations on ``dotted elements"
(see the definition below), and give explicitly the linear relations on those elements,
which will be used in Subsections 4.1 and 4.3--4.5.

\medskip

{\bf Theorem 2.4.4}\, For each $0\leqslant i\leqslant s$ in the sequence (2.4-1), there exists a system of
equations $\bar{\IE}^i$ over $R^i\otimes_kR^i$,
whose general solution can be expressed as the formal product
$\Pi^i=\sum_{j}v_j^i\ast V_j^i$, namely

(i)\, the $R^i\otimes_kR^i$-quasi-basis $\{V_j^i\}_j$ of $\K_1^i$ forms a basic system of solutions of
$\bar\IE^i$;

(ii)\, the $R^i\otimes_kR^i$-quasi-basis $\{v_j^i\}_j$ of $\C_1^i$ forms a set of free variables.

\smallskip

{\bf Proof}\, For $i=0$, let $\bar\Phi_{\m^0}=\sum_{j=1}^mv_j\ast V_j=\Pi$
and $\bar\IE^0: \bar\Phi_{\m^0} H^0\equiv_{\prec (p,q)} H^0\bar\Phi_{\m^0}$
be a matrix equation with $(p,q)$ being the leading position of $A_1$ in $\M_1$,
$H^0=0$. Then $\Pi$ is a general solution of $\bar\IE^0$.

Suppose a system $\bar\IE^i$ of the pair $(\mf A^i,\mf B^i)$ satisfying condition (i)--(ii) has been obtained:
$$\begin{array}{c}\bar\IE^i: \bar\Phi_{\m^i}H^i\equiv_{\prec(p^{i},q^{i})} H^i\bar\Phi_{\m^i},
\end{array}\eqno{(2.4\mbox{-}4)}$$
where $\bar \Phi^i=(u_{pq})$ is strictly upper triangular with $u_{pq}$ being a $k$-linear combination
of some variables over $R^i_{X}\times R^i_{Y}$ for $p\in X,q\in Y,X,Y\in\T^i$.
We now construct a system $\bar\IE^{i+1}$. For the sake of convenience,
$\bar\IE^i_{\prec(p^i,q^i)}$ is used for the equation system $\bar\IE^i$, and $\bar\IE^i_{(p^i,q^i)}$ stands for the
$(p^i,q^i)$-th equation of $\bar\IE^i$.

1)\, If $\bar\IE^i_{(p^i,q^i)}$ is not a linear combination of the equations of $\bar\IE^i_{\prec(p^{i},q^{i})}$,
we proceed with a regularization. Thus $\m^{i+1}=\m^{i}$, $T^{i+1}=T^i$ and $\T^{i+1}=\T^i$.
The combination of the equation $\bar\Phi_{\m^{i}} H^{i}\equiv_{(p^{i},q^{i})}H^{i}\bar\Phi_{\m^{i}}$,
and the equations of $\bar\IE^i_{\prec(p^i,q^i)}$ forms an equation system $\bar\IE^{i+1}:\bar\Phi_{\m^{i+1}}H^{i+1}
\equiv_{\prec(p^{i+1},q^{i+1})} H^{i+1}\bar\Phi_{\m^{i+1}}$, where $H^{i+1}=H^i+\emptyset\ast A^i_1$
by Proposition 2.2.8. And $\bar\IE^{i+1}$ satisfies assertions (i)--(ii).

\smallskip

2)  If $\bar\IE^i_{(p^i,q^i)}$ is a linear combination of the equations of $\bar\IE^i_{\prec(p^{i},q^{i})}$,
we proceed with a loop or an edge reduction. There are a pre-minimal algebra
$\bar R^i=R^i$ in a loop reduction, or $\bar R^i=R^i[a^i_1]$ in an edge reduction; a minimal algebra
$R^{i+1}$; an admissible $R^{i+1}$-$\bar R^i$-bimodule $L^i$.
Set a size vector $\n=\vartheta^{i,i+1}(1,\cdots,1)$ over $\T^i$
with $(1,\cdots,1)$ being a size vector over $\T^{i+1}$, then $\m^{i+1}=\vartheta^{0i}(\n)$ is a size vector over $\T$.
Denote by $\bar\Phi_{XY}=(u_{pq}')$ for any $(X,Y)\in\T^i\times\T^i$ a submatrix of $\bar\Phi_{\m^i}$, such that $u'_{pq}=u_{pq}$
for $p\in X,q\in Y$, or $0$ otherwise. Define a variable matrix over $R^{i+1}\otimes_k R^{i+1}$,
and a matrix in $\IM_{\m^{i+1}}(R^{i+1})$:
$$\begin{array}{c}(\bar\Phi_{\m^{i+1}})_{\n}=\sum_{(X,Y)\in \T^i\times\T^i}\sum_{1\leqslant p\leqslant n_{_X},1\leqslant q\leqslant n_{_Y}}
(\textsf f_{_{Z_{(X,p)}}}\otimes_k\textsf q_{_{Z_{(Y,q)}}})
\ast\bar\Phi_{XY};\\
H^i_{\n}=\sum_{X\in\T^i}\sum_{1\leqslant p\leqslant n_{_X}}(\textsf f_{_{Z_{(X,p)}}}
\otimes_k\textsf q_{_{Z_{(X,p)}}})\ast H_X^i.\end{array}$$
where the vertices $Z$ and the matrices $\textsf f\otimes_k\textsf e$ are  given before
Definition 2.1.2. If the $R^i$-$R^i$-quasi-basis $\V^i=\{V^i_j\mid j=1,\cdots,m^i\}$ is a basic solutuion
of $\bar\IE^i_{\prec(p^i,q^i)}$, then the $R^{i+1}$-$R^{i+1}$-quasi-basis
$\{(\textsf f_{_{Z_{(s(v_j^i),p)}}}\otimes_k\textsf q_{_{Z_{(t(v^i_j),q)}}})
\ast V^i_j\mid 1\leqslant p\leqslant n_{s(v^i_j)},1\leqslant q\leqslant n_{t(v^i_j)};j=1,\cdots,m^i\}$
of $\K^{i+1}_{11}$ is a basic system of solutions of the matrix equation
$(\bar\Phi_{\m^{i+1}})_{\n}H^i_{\n}\equiv_{\preccurlyeq (p^{i},q^{i})} H^i_{\n}(\bar\Phi_{\m^{i+1}})_{\n}$
partitioned under $\T^i$, since the $(p^i,q^i)$-th block is a $k$-linear combination of the others.
In other words, the formal product $\Pi^{i+1}_1$ of $(\K^{i+1}_{11},\C^{i+1}_{11})$ is a general solution of
the matrix equation.

We may assume that $\bar E_1$ has a $R^{i+1}$-$R^{i+1}$-quasi-basis
$\{F_1,\cdots,F_l\}$ by Definition 2.1.2 (a3), where $0\leqslant l\leqslant \frac{1}{2} sr(r+1)$
after some deletion in Formula (2.2-2) for $a_1^i: X\rightarrow X$, and $F_jL(a^i_1)=L(a_1^i)F_j$;
or $l=0,1,2$ after some deletion in Formula (2.2-4) for $a^i_1: X\rightarrow Y$, and $F_1L(a^i_1)=0=L(a^i_1)F_2$.
Then either $\K_{10}^{i+1}=\{0\}$, or $R^{i+1}$-$R^{i+1}$- quasi-basis of $\K_{10}^{i+1}$ is
$\{F_j'=F_j\ast E^i_{X}\mid j=1,\cdots,l\}$; or $F_1'=F_1\ast E^i_X$
or $F_2'=F_2\ast E^i_Y$, or both of them given by Construction 2.1.3 (ii), and that
of $\C^{i+1}_{10}$ is $\{{F_1'}^\ast,\cdots,{F_l'}^\ast\}$ given by Proof 2) of Proposition 2.1.5.
Thus the formal product $\Pi^{i+1}_0=\sum_{j=0}^l{F_j'}^\ast\ast F'_j$ of $(\K_{10}^{i+1},\C_{10}^{i+1})$ is
a general solution of the matrix equation $\Pi_0^{i+1}\big(L(a_1)\ast A^i_1\big)
\equiv_{\preccurlyeq(p^i,q^i)}\big(L(a_1)\ast A^i_1\big)\Pi_0^{i+1}$ partitioned under $\T^i$, since the $(p^i,q^i)$-th block
is $\sum_1^l(F_j^\ast\ast F_j)L(a_1^i)=F_j^\ast\ast(L(a_1^i)F_j)=\sum_1^lL(a_1^i)(F_j^\ast\ast F_j)$. Define
$$\begin{array}{c}\bar\Phi_{\m^{i+1}}=\Pi^{i+1}_0+(\bar\Phi_{\m^i})_{\n},\quad H^{i+1}=H^i_{\n}+L(a)\ast A_1^i;\\[1mm]
\bar\IE^{i+1}: \bar\Phi_{\m^{i+1}}H^{i+1}\equiv_{\prec(p^{i+1},q^{i+1})} H^{i+1}\bar\Phi_{\m^{i+1}},\end{array}$$
where $(p^{i+1},q^{i+1})$ is the index followed by the biggest index over $\T^{i+1}$ of
the $(p^i,q^i)$-block partitioned under $\T^i$.
We claim that the formal product $\Pi^{i+1}=\Pi^{i+1}_0+\Pi^{i+1}_1$ of $(\K^{i+1}_1,\C^{i+1}_1)$
is a general solution of $\bar\IE^{i+1}$.
First, both left and right sides of each block equation of $(\bar\Phi_{\m^{i}})_{\n}\big(L(a)\ast A_1^i\big)
\equiv_{\preccurlyeq(p^{i},q^{i})}\big(L(a)\ast A_1^i\big)(\bar\Phi_{\m^{i}})_{\n}$ partitioned under $\T^i$
are zero blocks, since $(\bar\Phi_{\m^{i}})_{\n}$ is a strict upper triangular partitioned matrix
and the index of the leading block of $L(a^i_1)\ast A^i_1$ is $(p^i,q^i)$.
Second, at the left and right sides of each block equation of $\Pi^{i+1}_0H^i_{\n}\equiv_{\preccurlyeq(p^{i},q^{i})}H^i_{\n}\Pi^{i+1}_0$
under $\T^i$ are equal blocks. In fact, suppose $H_X^i=(h_{\alpha\beta}1_X)_{t^i}$ with
$h_{\alpha\beta}=0$ if $\alpha\notin X$ or $\beta\notin X$, then $(H_{\n}^i)_X=(H_{X,\alpha\beta})$ with
$H_{X,\alpha\beta}=h_{\alpha\beta}$diag$(1_{Z_{(X,1)}},\cdots,1_{Z_{(X,n_{_X})}})$,
therefore $F_jH_{X,\alpha\beta}=H_{X,\alpha\beta}F_j$. And the same assertion
is valid for $Y\in\T^i$. Our Theorem follows by induction. \hfill$\Box$

\medskip

For an example, see 2.4.5 (iv) below. $\bar\IE^i$ in Formula (2.4-4) is also called a {\it defining system of the pair $(\mf A^i,\mf B^i)$}.
The matrix $\bar\Phi_{\m^{i}}$ is called a {\it matrix of dotted elements}.
The concept of the dotted elements possesses two folds of meanings: 1) as variables in the equation system
$\bar{\IE}$; 2) as the elements with a series of linear relations after a sequence of reductions.
Different meanings will be used for different cases frequently in Section 4.

Next, we define a deformed system $\bar\IF^{ri}$ for some fixed $0<r<s$, which is equivalent to $\bar\IE^i$.
Like the discussion stated before Formula (2.4-3), a matrix equation and a variable matrix
of size vector $\m^{ri}$ over $\T^{r}$ are defined:
$$\begin{array}{c}{\bar\IF}^{ri}: \bar\Psi_{\m^{ri}}
H^i_2\equiv_{\prec(p^i,q^i)} \bar\Psi^0_{\m^{ri}}+H^i_2 \bar\Psi_{\m^{ri}},\\[2mm]
\bar\Psi_{\m^{ri}}^0=H^i_1\bar\Psi_{\m^{ri}}-\bar\Psi_{\m^{ri}} H^i_1,\\[1mm]
\bar\Psi_{\m^{ri}}=\sum_{X^r\in\T^r}\bar w_{_{X^r}}^{ri}\ast E_{X^r}+
\sum_j\bar v^{ri}_j\ast V^r_j.
\end{array}\eqno(2.4\mbox{-}5)$$
where the definition of $\bar v_j^{ri}=(v^{ri}_{jpq})$ and $\bar w_{_{X^r}}^{ri}=(w_{{_{X^r}}pq}^{ri})$
is analogous to that of Formula (2.4-4).

At the end of the subsection, we perform reduction procedure for the matrix bimodule
problem given in Example 1.4.5 in order to show some concrete calculations.

\medskip

{\bf Example 2.4.5}\, (i) Making an edge reduction for the first arrow $a:X\rightarrow Y$ by
$a\mapsto G^1=(1_{Z})$, an induced local pair $(\mf{A}^1,\mf B^1$) with $R^1=k1_Z$;
$H^1=(1_Z)\ast A$ is obtained.

(ii) Making a loop reduction for $b: Z\rightarrow Z$ by
$b\mapsto G^2=J_2(0)1_X$, an induced local pair $(\mf{A}^2,\mf{B}^2)$ with
$R^2=k1_X, H^2={{1_X\,\,\,0\,}\choose{\,\,0\,\,1_X}}\ast A
+{{\,0\,1_X}\choose{0\,\,\,0\,}}\ast B$ is obtained. There are two matrix equalities
in the formal equation of the pair $(\mf{A}^2,\mf{B}^2)$:
$$\small{~_{\begin{array}{ll}  \left(\begin{array}{cc} e& v\\
0&e\end{array}\right)\left(\begin{array}{cc} c_{11}& c_{12}\\
c_{21}&c_{22}\end{array}\right)+\left(\begin{array}{cc} u^2_{11}& u^2_{12}\\
u^2_{21}&v^2_{22}\end{array}\right)\left(\begin{array}{cc} 0& 1_{X}\\
0&0\end{array}\right)\\[3mm]
= \left(\begin{array}{cc} 0& 1_{X}\\
0&0\end{array}\right)\left(\begin{array}{cc} v^2_{11}& v^2_{12}\\
v^2_{21}&v^2_{22}\end{array}\right)+\left(\begin{array}{cc} c_{11}& c_{12}\\
c_{21}&c_{22}\end{array}\right)\left(\begin{array}{cc} e& v\\
0&e\end{array}\right),    \end{array}}}$$
$$\small{~_{\begin{array}{ll}  \left(\begin{array}{cc} e& v\\
0&e\end{array}\right)\left(\begin{array}{cc} d_{11}& d_{12}\\
d_{21}&d_{22}\end{array}\right)+\left(\begin{array}{cc} u^1_{11}& u^1_{12}\\
u^1_{21}&u^1_{22}\end{array}\right)\left(\begin{array}{cc} 0& 1_{X}\\
0&0\end{array}\right)+\left(\begin{array}{cc} u^2_{11}& u^2_{12}\\
u^2_{21}&u^2_{22}\end{array}\right)\left(\begin{array}{cc} 1_X& 0\\
0&1_X\end{array}\right) \\[3mm]
= \left(\begin{array}{cc} 1_{X}&0 \\
0&1_X\end{array}\right)\left(\begin{array}{cc} v^2_{11}& v^2_{12}\\
v^2_{21}&v^2_{22}\end{array}\right)+\left(\begin{array}{cc} 0& 1_{X}\\
0&0\end{array}\right)\left(\begin{array}{cc} v^1_{11}& v^1_{12}\\
v^1_{21}&v^1_{22}\end{array}\right)+\left(\begin{array}{cc} d_{11}& d_{12}\\
d_{21}&d_{22}\end{array}\right)\left(\begin{array}{cc} e& v\\
0&e\end{array}\right), \end{array}}}$$
where $(c_{pq})_{2\times 2},(d_{pq})_{2\times 2}$ are splits from $c,d$
respectively, $e\in\C_0^2$ is dual to $E_X=(1_XI_{10},1_XI_{10})\in\K^2_0$,
and $v\in \C_1^2$ is dual to $V={{0\,\,1_X\otimes_k1_X}\choose{0\qquad 0\quad}}
\ast E_X\in \K^2_1$ respectively.

(iii) Making a loop mutation $c_{21}\mapsto (x)$, followed by three regularizations, such
that $c_{22}\mapsto\emptyset, u^2_{21}= xv$;
$c_{11}\mapsto\emptyset, v^2_{21}=vx $; $c_{12}\mapsto\emptyset,
u^2_{11}=v^2_{22}$, an induced pair
$(\mf A^3,\mf B^3)$ is obtained, and the differentials of the solid arrows in $\mf B^3$ are:
$$
\small{\left\{\begin{array}{l} \delta(d_{21})= xv-vx\\
\delta(d_{22})=u^1_{21}+u^2_{22}-v^2_{22}-d_{21}v\\
\delta(d_{11})= u^2_{11}-v^2_{11}-v^1_{21}+vd_{21}\\
\delta(d_{12})= u^1_{11}+u^2_{12}-v^2_{12}-v^1_{22} -d_{11}v +vd_{22}.
\end{array}\right.}
$$

(iv)\, Finally, we describe the defining systems of Theorem 2.4.1 and 2.4.4 for the pair $(\mf A^2,\mf B^2)$.
Since $(\mf A,\mf B)$ is bipartite, $\T=\T'\times \T''$. Thus $\Phi_{\m^2}=\Psi_{\l^2}^1\times\Phi_{\n^2}^2$
and $\bar\Phi_{\m^2}=\bar\Phi_{\l^2}^1\times\bar\Phi_{\n^2}^2$, where
the size vector $\l^2=(2,2,2,2,2)$ is over $\T'$ and $\n^2=(2,2,2,2,2)$ over $T''$.
Suppose the systems are $\Phi_{\l^2}^1H^2(k)=H^2(k)\Phi_{\n^2}^2$
given by 2.4.1, and  $\bar\Phi_{\l^2}^1H^2=H^2\bar\Phi_{\n^2}^2$ by 2.4.4 respectively, where
$$\small{\Phi^1_{\l^2}\,(\mbox{or\,\,}\bar\Phi^1_{\l^2})=
\left(
\begin{array}{ccccc}
\Phi_0 & 0 &  \Phi_1    &  \Phi_2 & \Phi_4 \\
  & \Phi_0 &  \Phi_2    & 0 & \Phi_3 \\
  &   &   \Phi_0   & 0 & \Phi_2 \\
  &   &       & \Phi_0 & \Phi_1 \\
  &   &       &   & \Phi_0
\end{array}
\right),\quad\Phi^2_{\n^2}\,(\mbox{or}\,\,\bar\Phi^2_{\n^2})= \left(
\begin{array}{ccccc}
\Phi_0' & 0 &  \Phi_1'    &  \Phi'_2 & \Phi'_4 \\
  & \Phi'_0 &  \Phi'_2    & 0 & \Phi'_3 \\
  &   &   \Phi'_0   & 0 & \Phi'_2 \\
  &   &       & \Phi'_0 & \Phi'_1 \\
  &   &       &   & \Phi'_0
\end{array}
\right)}.
$$
Then $\Phi_i={{x^i_{11}\,x_{12}^i}\choose{x^i_{21}\,x^i_{22}}}$ in $\Phi^1_{\l^2}$,
$\Phi'_i={{y^i_{11}\,y_{12}^i}\choose{y^i_{21}\,y^i_{22}}}$ in $\Phi^2_{\n^2}$ for $i=0,1,2,3,4$ by Theorem 2.4.1;
and $\Phi_0={{0\, v}\choose{0\, 0}}=\Phi_0'$ in $\bar\Phi_{\m^2}$, which is obtained from a loop reduction $b\mapsto J_2(0)1_X$;
$\Phi_i={{u^i_{11}\,u_{12}^i}\choose{u^i_{21}\,u^i_{22}}}$ in $\bar\Phi^1_{\l^2}$,
$\Phi'_i={{v^i_{11}\,v_{12}^i}\choose{v^i_{21}\,v^i_{22}}}$ in $\bar\Phi^2_{\n^2}$ for $i=1,2,3,4$ by Theorem 2.4.4.

\bigskip
\bigskip
\bigskip

\centerline{\Large\bf  3 Classification of minimal wild bocses}
\bigskip

Based on the well known Drozd's wild configurations, this section is devoted to classifying so-called minimal wild
bocses, which are divided into five classes. Then the non-homogeneity of bocses in the first four classes
is proved. But those in the last class have been proved to be strongly homogeneous.
Some preliminaries are stated in subsections 3.1 and 3.2.

\bigskip
\bigskip
\noindent {\large\bf 3.1 An exact structure on representation categories of bocses}
\bigskip

In this subsection the concept on exact structure
of categories is recalled, especially the exact structure on representation
categories of bocses.

\medskip

Let $\mathscr{A}$ be an additive category with Krull-Schmidt
property. We recall from \cite{GR}  and \cite{DRSS}  the following
notions. A pair $(\iota, \pi)$ of composable morphisms
$$ (e) \qquad M\stackrel{\iota}{\longrightarrow}E\stackrel{\pi}{\longrightarrow}N\eqno{(\ast)}$$ in
$\mathscr{A}$ is called {\it exact} if $\iota$ is a kernel of $\pi$
and $\pi$ is a cokernel of $\iota$.

Let ${\cal E}$ be a class of exact pairs
which is closed under isomorphisms. The morphisms $\iota$ and $\pi$
appearing in a pair $(e)$ are called an {\it
inflation} and a {\it deflation} of ${\cal E}$ respectively, the
pair itself is called a {\it conflation}, and is denoted by $(\iota, \pi)$.

\medskip

{\bf Definition 3.1.1} The class ${\cal E}$ is said to be an {\it
exact structure} on $\mathscr{A}$, and $(\mathscr{A}, {\cal E})$ an
{\it exact category} if the following axioms are satisfied:

{\bf E1}\, The composition of two deflations is a deflation.

{\bf E2}\, For each $\varphi$ in $\mathscr{A}(N', N)$ and each
deflation $\pi$ in $\mathscr{A}(E,N)$, there are some $E'$ in
$\mathscr{A}$, an $\varphi'$ in $\mathscr{A}(E',E)$ and a deflation
$\pi': E'\rightarrow N' $ such that $\pi'\varphi=\varphi' \pi$.

{\bf E3}\, Identities are deflations. If $\varphi\pi$ is a deflation,
then so is $\pi$.

(Or {\bf E3$^{op}$}\, Identities are inflations, if $\iota\varphi$ is an
inflation, then so is $\iota$.)

\medskip

An object $P$ in $\mathscr{A}$ is said to be {\it ${\cal
E}$-projective} (or {\it projective} for short ) if any conflation
ending at $P$ is split. Dually an object $I$ in $\mathscr{A}$ is
said to be {\it ${\cal E}$-injective} (or {\it injective} for short
) if any conflation starting at $I$ is split.

Let $\mathscr{A}$ be a Krull-Schmidt category. A morphism $\pi:
E\rightarrow N $ in $\mathscr{A}$ is called {\it right almost split}
if it is not a retraction and for any non-retraction $\varphi:
L\rightarrow N$, there exists a morphism $\psi: L\rightarrow E$ such
that $\varphi=\psi\pi$. It is said that $\mathscr{A}$ has {\it right
almost split morphisms} if for all indecomposable $N$ there exist
 right almost split morphisms ending at $N$. Dually,
{\it left almost split morphisms} are defined. It is said that {\it $\mathscr{A}$
has almost split morphisms} if $\mathscr{A}$ has right and left
almost split morphisms.

A morphism $\pi: E\rightarrow N$ is called {\it right minimal} if
every endomorphism $\eta: E\rightarrow E$ with the property that
$\pi=\eta\pi$ is an isomorphism. A {\it left minimal morphism}
$\iota: M\rightarrow E$ is defined dually.

\medskip

{\bf Proposition 3.1.2} Suppose that the Krull-Schmidt category
$\mathscr{A}$ carries an exact structure ${\cal E}$. Let
$(e)$ given in Formula $(\ast)$ be a conflation. Then the following
assertions are equivalent.

(i) $\iota$ is minimal left almost split;

(ii) $\pi$ is minimal right almost split;

(iii) $\iota$ is left almost split and $\pi$ is right almost split.

\medskip

The conflation $(e)$ in the above proposition
is said to be an {\it almost split conflation}. The exact category
$(\mathscr{A},{\cal E})$ is said to {\it have almost split
conflations} if (i) $\mathscr{A}$ has almost split morphisms;
(ii) for any indecomposable non-projective $N$, there exists an
almost split conflation $(e)$ ending at $N$; (iii)
for any indecomposable non-injective $M$, there exists an almost
split conflation $(e)$ starting at $M$.

\medskip

Now we turn to the representation category of a bocs.
Let $\mathfrak{B}=(\Gamma,\Omega)$ be a bocs with a layer
$L=(\Gamma';\omega;a_1,\cdots,a_n;v_1,\cdots,v_m)$. From now on
it is always assumed that $\mf B$ is {\it triangular on the dotted arrows},
i.e. $\delta(v_j)$ involves only $v_1, \cdots, v_{j-1}$.
In particular, the bocs $\mf B$ associated to a matrix bimodule problem $\mf A=(R,\K,\M,H)$ is triangular
by Definition 1.2.1.

The bocs $\mathfrak{B}_0=(\Gamma,\Gamma)$
is called a {\it principal bocs} of $\mf{B}$.
The representation category $R(\mf{B}_0)$ is just the
module category of $\Gamma$.

\medskip

{\bf Lemma 3.1.3}  \cite{O}\, Let $\mathfrak{B}=(\Gamma,\Omega)$
be a layered bocs with a principal bocs
$\mf{B}_0$. Suppose $\mf B$ is triangular on the dotted arrows.

(i) If $\iota: M\rightarrow E$ is a morphism of $R(\mathfrak{B})$
with $\iota_0$ injective, then there exist an isomorphism $\eta$
and a commutative diagram in $R(\mathfrak{B})$, such that  the bottom row
is exact in $R(\mathfrak{B}_0)$.
Dually, if $\pi: E\rightarrow N$ is a morphism of $R(\mathfrak{B})$
with $\pi_0$ surjective, then there exist an isomorphism $\eta$ and
a commutative diagram in $R(\mathfrak{B})$, such that  the bottom row
is exact in $R(\mathfrak{B}_0)$.
$$\small{
\begin{CD}
 &  & M @>\iota>> E\\
 && @VidVV @VV\eta V\\
0@>>>M @>>\iota^{\prime}>E^{\prime}
\end{CD}}\qquad\qquad
\small{\begin{CD}
  E @>\pi>> N\\
  @V\eta VV @VVidV\\
E^{\prime}@>>\pi^{\prime}>N @>>>0
\end{CD}}
$$

(ii) If $(e):  M\stackrel{\iota}{\longrightarrow}
E\stackrel{\pi}{\longrightarrow} N$ with $\iota\pi=0$ is a pair of
composable morphisms  in $R(\mathfrak{B})$ and
$(e_0):0\longrightarrow M\stackrel{\iota_0}{\longrightarrow}
E\stackrel{\pi_0}{\longrightarrow} L\longrightarrow 0$ is exact in
the category of vector spaces, then there exists an isomorphism
$\eta$ and a commutative diagram in $R(\mathfrak{B})$:
$$\small{\begin{CD} (e) \qquad &     & M@>\iota>> E @>\pi>> N  \\
& &@VidVV @VV\eta V@VVidV \\
(e') \qquad   0 @>>> M @>>\iota^{\prime}> E^{\prime}@>>\pi^{\prime}>
N @>>>0
\end{CD}}
$$
such that $(e')$ is an exact sequence in $R(\mathfrak{B}_0)$.
Moreover, by choosing a suitable basis of $M,E',N$, we are able to obtain
$\iota'_X=(0, I)$ and $\pi'_X=(I,0)^T$ for all $X\in \T$.\hfill$\Box$

\medskip

{\bf Lemma 3.1.4} Let $\mf{B}=(\Gamma,\Omega)$ be a layered
bocs, which is triangular on the dotted arrows.

(i) $\iota: M\rightarrow E$ is monic in $R(\mf B)$ if $\iota_0: M\rightarrow E$ is
injective. Dually, $\pi: E\rightarrow N$ is epic in $R(\mf B)$ if $\pi_0: E\rightarrow N$ is
surjective.

(ii)  A pair of composable morphisms $(e):
M\stackrel{\iota}{\longrightarrow} E\stackrel{\pi}{\longrightarrow}
N$ with $\iota\pi=0$ is exact in $R(\mf{B})$, if $(e_0):
0\longrightarrow M\stackrel{\iota_0}{\longrightarrow}
E\stackrel{\pi_0}{\longrightarrow} N\longrightarrow 0$ is exact as
a sequence of vector spaces.

\smallskip

{\bf Proof}\, (i) If $\iota_0$ is injective, then Lemma 3.1.3 (i) gives a
commutative diagram with $\iota': M\rightarrow E'$ in
$R(\mf{B}_0)$. Given any morphism $\varphi: L\rightarrow M$
with $\varphi\iota=0$, there is $\varphi\iota\eta=\varphi\iota'=0$.
Then $\varphi_0\iota'_0=0$ yields $\varphi_0=0$. And for any dotted arrow $v_l: X
\rightarrow Y$, suppose $\delta(v_l)=\sum_{i,j}u_i
\otimes_{\Gamma} u_{j}$ with $u_i,u_j\in \oplus_{l'<l}\Gamma v_{l'}\Gamma$.
There is inductively, $0=(\varphi\iota')(v_l)=\varphi(v_l)\iota'_Y+\varphi_X\iota'(v_l)+\sum_{i,j}
\varphi(u_i)\iota'(u_j)=\varphi(v_l)\iota'_Y$, which yields $\varphi(v_l)=0$ by
the injectivity of $\iota'_Y$. Thus $\varphi=0$ and $\iota$ is monic. The second assertion on $\pi$
is proved dually.

(ii) It is proved first that $\iota$ is the kernel of $\pi$.  \ding{172} If
$(e_0)$ is exact, then (i) shows that $\iota$ is monic.
\ding{173} It is known that
 $\iota\pi=0$. \ding{174} Lemma 3.1.3 (ii) gives a commutative diagram.
If $\varphi: L\rightarrow E $ with $\varphi\pi=0$, then
$\varphi(\eta\pi^{\prime})=0$. Let $\xi=\varphi\eta$, then
$\xi\pi^{\prime}=0$ implies that $\xi_X\pi_X^{\prime}=0$ and
$\xi(v)\pi_Y^{\prime}=0$ for any vertex $X$ and any dotted arrow $v:
X\rightarrow Y$. Let $\varphi^{\prime}:L\rightarrow M$ be given by
$\varphi^{\prime}_X\iota^{\prime}_X=\xi_X$,
$\varphi^{\prime}(v)\iota^{\prime}_Y=\xi(v)$,  then
$\varphi^{\prime}\iota^{\prime}=\xi$ is obtained. Thus
$\varphi^{\prime}\iota^{\prime}\eta^{-1}=\varphi$, i.e.
$\varphi^{\prime}\iota=\varphi$. Therefore $\iota$ is a kernel of
$\pi$. Second, it can be proved dually that $\pi$ is a cokernel of $\iota$.
\hfill$\Box$

\medskip

Let a layered bocs $\mathfrak{B}=(\Gamma,\Omega)$ be triangular on the dotted arrows. A class
${\cal E}$ of composable morphisms
in $R(\mf{B})$ is defined, such that $M\stackrel{\iota}{\longrightarrow}
E\stackrel{\pi}{\longrightarrow} L$ in ${\cal E}$, provided that
$\iota\pi=0$ and
$$0\longrightarrow M\stackrel{\iota_0}{\longrightarrow}
E\stackrel{\pi_0}{\longrightarrow} L\longrightarrow 0\eqno{(3.1\mbox{-}1)}$$ is exact as
a sequence of vector spaces. It is clear that ${\cal E}$ is closed
under isomorphisms.

\medskip

{\bf Proposition 3.1.5}\, \cite[Theorem 4.4.1]{O} and \cite{BBP})
Suppose a layered bocs $\mf B$ is triangular on the dotted arrows.
Then the class ${\cal E}$ defined by Formula (3.1-1) is an exact structure on
$R(\mathfrak{B})$, and $(R(\mathfrak{B}),{\cal E})$ is an exact category.

\medskip

{\bf Corollary 3.1.6} (\cite{B1}, \cite [Lemma 7.1.1]{O})\,
Let $\mf{B}=(\Gamma, \Omega)$ be a layered bocs.

(i)\, For any $M\in R(\mf{B})$ with
$\dim M=\m$, if $m_X\ne 0$ for some vertex $X\in \T_1$, then $M$ is
neither projective nor injective.

(ii)\, For any positive integer $n$, there are only finitely many
iso-classes of indecomposable projectives and injectives in $R(\mf B)$
of dimension at most $n$.

\medskip

{\bf Remark 3.1.7} (\cite{BCLZ}, \cite[Definition 4.4.1]{O})\,
Let $\mf{B}=(\Gamma, \Omega)$ be a layered bocs, such that $(R(\mf{B}),
{\cal E})$ is an exact structure. The almost split conflations have been defined
in a general exact category, particularly in $R(\mf B)$.

(i) An indecomposable representation $M\in R(\mf{B})$ is said to be
{\it homogeneous} if there is an almost split conflation
$M\stackrel{\iota}{\longrightarrow}E\stackrel{\pi}{\longrightarrow}M$.
The iso-class of $M$ is also said to be {\it homogeneous}.

(ii) The category $R(\mf{B})$ (or bocs $\mf B$) is said to be {\it homogeneous} if for
each positive integer $n$, almost all (except finitely many)
iso-classes of indecomposable representations in $R(\mf{B})$ with
size at most $n$ are homogeneous. For example,
If $\mf B$ is of representation tame type,
then $R(\mf B)$ is homogeneous \cite{CB1}.

(iii) The category $R(\mf{B})$ (or bocs $\mf B$) is said to be {\it strongly
homogeneous} if  there exists neither projectives nor injectives,
and all indecomposable representations in $R(\mf{B})$ are
homogeneous. For example, if $\mf B$ is a local bocs with a layer
$(R;\omega; a; v), R=k[x,\phi(x)^{-1}]$, and the differential
$\delta(a)= xv-vx$. Then $R(\mf B)$ is strongly homogeneous and
representation wild type \cite{BCLZ}. In particular the induced bocs given in
Example 2.4.5 (iii) is strongly homogeneous.

\medskip

Note that $(R(\mathfrak{B}),{\cal E})$ may not have any almost split
conflation. For example, set quiver $Q=\xymatrix{\cdot\ar@(ul,dl)_{a}\ar@(ur,dr)^{b}}$, the path algebra
$\Gamma=kQ$, and the principal bocs $\mf{B}=(\Gamma,\Gamma)$.  Then $R(\mf{B})$ has
no almost split conflations, see \cite{V,ZL} for details.

Recalling from \cite{CB1}, let
$\mf{B}=(\Gamma,\Omega)$ be a minimal bocs. Then
for any $X\in \T_1$ with $R_X=k[x, \phi_{_X}(x)^{-1}]$, and for
any $\lambda\in k$ with $\phi_{_X}(\lambda)\ne 0$, there is an almost split conflation:
$$
\begin{array}{c}
S(X,1,\lambda)\stackrel{(0\, 1)}\rightarrow S(X,2,\lambda)
\stackrel{{1}\choose{0}}\rightarrow
S(X,1,\lambda)\quad\mbox{in}\,\,\,R(\mf{B}),
\end{array}\eqno{(3.1\mbox{-}2)}
$$
where
$S(X,1,\lambda)$ (resp. $S(X,2,\lambda)$) is given by $\xymatrix{k\ar@(ur,dr)^{J_1(\lambda)}}$
(resp. $=\xymatrix{k^2\ar@(ur,dr)^{J_2(\lambda)}}$) at $X$, and $\{0\}$ at other vertices.

\bigskip
\bigskip
\noindent{\large\bf 3.2 Almost split conflations in the process of reductions}
\bigskip

In order to prove the non-homogeneity of some wild bocses, we must understand the behavior
of almost split conflations during reduction procedures.
We study under what conditions the homogeneous property is preserved after a sequence of reductions
in this subsection.

\medskip

{\bf Lemma 3.2.1} \cite{B1}\, Let $\mf B'=(\Gamma',\Omega')$ be
the induced bocs of a bocs $\mf B$ given by one of eight reductions in the subsection 2.2, and
$N'$ be an indecomposable representation in
$R(\mf B')$. If $N'$ is non-projective (resp. non-injective) in
$R(\mf B')$, then so is $\vartheta(N')$ in $R(\mf B)$.

\medskip

{\bf Lemma 3.2.2} \cite{B1}\, Let $\mf B'=(\Gamma',\Omega')$ be
the induced bocs of $\mf B$ given by one of eight reductions in the subsection 2.2.

(i) If $\iota': M'\rightarrow E'$ is a morphism
in $R(\mf B')$ with $\vartheta(\iota'): \vartheta(M')\rightarrow
\vartheta(E')$ being a left minimal almost split inflation in
$R(\mf B)$, then so is $\iota'$ in $R(\mf B')$.
Dually if $\pi': E'\rightarrow N'$ is a morphism in $R(\mf B')$
with $\vartheta(\pi'): \vartheta(E')\rightarrow \vartheta(N')$
being a right minimal almost split deflation in $R(\mf B)$,  then so is
$\pi'$ in $R(\mf B')$.

(ii) If $(e'): M' \stackrel {\iota'} \rightarrow
E'\stackrel {\pi'} \rightarrow M'$ is a conflation in
$R(\mf B')$ with $\vartheta (e'): \vartheta (M')
\stackrel {\vartheta (\iota')} \rightarrow \vartheta
(E')\stackrel {\vartheta(\pi')} \rightarrow \vartheta
(M')$ being an almost split conflation in $R(\mf B)$, then
so is $(e')$ in $R(\mf B')$.

\medskip

Let $(\mf A,\mf B)$ be a pair with trivial $R$, $M\in R(\mf A)$ be an indecomposable object of
size vector $\m$. Set $\T^1=\{X\in\T\mid m_{_X}\ne 0\}$,
suppose $\mf A^1$ is obtained by deleting $\T\setminus\T^1$ from $\mf A$,
and $M^1\in R(\mf A^1)$ with $\vartheta^{01}(M^1)\simeq M$. Suppose a sequence of reductions
in the sense of Lemma 2.3.2 is given
by Theorem 2.3.3 with respect to $M^1$:
$$(\mf A,\mf B),(\mf A^1,\mf B^1), \cdots,(\mf A^i,\mf B^i),(\mf A^{i+1},\mf B^{i+1}),\cdots,
(\mf A^s,\mf B^s).\eqno {(3.2\mbox{-}1)}$$
Then there is some $M^s\in R(\mf A^s)$ of sincere size vector $\m^s$, such that $\vartheta^{0s}(M^s)\simeq M$.

\medskip

{\bf Theorem 3.2.3}\,  Suppose the first arrow $a^s_1$ is a loop at $X^s$ with $\delta(a_1^s)=0$; and
$M^s_{X^s}=k$, $M^s(a_1^s)=(\lambda)$ in the last term $\mf B^s$ of the sequence (3.2-1).
If $M$ is homogeneous and $(e): M\rightarrow E\rightarrow M$ is an almost
split conflation in $R(\mf A)$,
then for $i=1,\cdots,s$ there exits
an almost split conflation $(e^i): M^i\rightarrow E^i\rightarrow M^i$ in $R(\mf A^i)$,
such that $\vartheta^{0i}(e^i)\simeq (e)$.

\smallskip

{\bf Proof}\, Induction is used for the proof.  The assertion is obviously true for $i=1$,
since the size vector of $E$ is $2\m$ by Formula (3.1-1).
According to Definition 1.3.4 and Formula (2.3-6):
$$\begin{array}{c}M^s=H^s_{\m^s}(k)+\sum_{j}M^s(a^s_j)\ast A^s_j,\quad
H^s_{\m^s}(k)=\sum_{j=1}^{s-1}B^{j+1}\ast A^j_1.\end{array}$$

Suppose the assertion is valid
for some $1\leqslant i<s$. The formula below gives
as $k$-matrices of size vector $\m^i=\vartheta^{i0}(\m^s)$:
$$\begin{array}{c}M^i=M^s=H^i_{\m^i}(k)+B^{i+1}\ast A^i_1+
\sum_{j=2}^{n^i}M^i(a^i_j)\ast A^i_j.\end{array}\eqno {(3.2\mbox{-}2)}$$
There exits an object
$E^i= H^i_{2\m^i}(k)+\sum_{j=1}^{m^i}E^i(a_j^i)\ast A_j^i\in R(\mf A^i)$ and
an almost split conflation $(e^i): M^i\stackrel{\iota^i}\rightarrow E^i
\stackrel{\pi^i}\rightarrow M^i$ in $R(\mf A^i)$,
such that $\vartheta^{0i}(e^i)\simeq(e)$.  We now treat the $(i+1)$-th stage
via proving the existence of an isomorphism $\eta:E^i\rightarrow \hat E^i$ in $R(\mf A^i)$ with
$\hat E^i(a^i_1)=B^{i+1}\oplus B^{i+1}$.

If this is the case, suppose $a^i_1:X\rightarrow Y$, $S_{X}$ and $S_{Y}$ are
invertible matrices determined by changing certain rows and columns of
$B^{i+1}\oplus B^{i+1}$, such that $S_X^{-1}(B^{i+1}\oplus B^{i+1})S_Y
=I_2\otimes B^{i+1}$, the usual Kronecker product of two matrices. Define
a matrix $S=\sum_{Z\in\T^i}S_Z\ast E_Z$ with
$S_Z=I_{m_{_Z}}$ for $Z\in\T^i\setminus \{X,Y\}$. Then there are $k$-matrices:
$$\begin{array}{ll}
R(\mf B^i)\ni \xi(\hat E^i)&:=S^{-1}\hat E^i S\\
&\,\,=H^i_{2\m^i}(k)+(I_2\otimes B^{i+1})\ast A^i_1+
\sum_{j=2}^{n^i}S_{s(a^i_j)}^{-1}\hat E^i(a^i_j)S_{t(a^i_j)}\ast A^i_j\\
&:=H^{i+1}_{2\m^{i+1}}(k)+\sum_{j=1}^{n^{i+1}}E^{i+1}(a^{i+1}_j)\ast A^{i+1}_j
=E^{i+1}\in R(\mf A^{i+1}).\end{array}$$
Thus an almost split conflation $(e^{i+1})$ in $R(\mf A^{i+1})$ is obtained
by Lemma 3.2.2 (ii), such tat $\vartheta^{i,i+1}(e^{i+1})\simeq(\hat e^i)\simeq (e^i)$
via the isomorphisms $E^{i+1}\stackrel{\xi^{-1}}\rightarrow\hat E^i\stackrel{\eta^{-1}}\rightarrow E^i$ in $R(\mf A^i)$.
Consequently $\vartheta^{0,i+1}(e^{i+1})\simeq\vartheta^{0i}(e^i)\simeq(e)$.

The existence of such an isomorphism $\eta$ is established below.

If $\dz(a_1^i)=v^i_1\ne 0$, then $B^{i+1}=M^i(a^i_1)=(0)$. By the proof (i) of Proposition 2.1.8,
there exists an isomorphism $\eta$, such that $\hat E^i=\eta(E^i)\in R(\mf A^i)$ with $\hat E^i(a^i_1)=(0)
$ as desired.

If $\dz(a_1^i)=0$ in the case of loop or edge reduction, the proof is divided into three parts.

1)\, We define an object $L^s=H^s_{2\m^s}(k)+\sum_{i=1}^{n^s}L^s(a^s_i)\in R(\mf A^s)$
with $L^s(a^s_1)={{\lambda\,1}\choose{0\,\lambda}}$.
Let $\varphi^s: L^s\rightarrow M^s$ be a morphism in $R(\mf A^s)$,
such that $\varphi^s_{Y^s}={I_{Y^s}\choose 0},\forall\,Y^s\in \T^s$,
and $\varphi^s({v^s_j})=0$
for any dotted arrow $v^s_j$ in $\mf{B}^s$. Clearly, $\varphi^s$ is
not a retraction. Thus $\vartheta^{is}(\varphi^s):\vartheta^{is}(L^s)\mapsto\vartheta^{is}(M^s)=M^i$
is not a retraction, since the functor $\vartheta^{is}$ is fully faithful.

2)\, Because $\vartheta^{is}(L^s)=H^i_{2\m^i}(k)+(I_2\otimes B^{i+1})\ast A^i_1+
\sum_{j=2}^{n^i}\vartheta^{is}(L^s)(a_j^i)\ast A^i_j$, it is possible to construct
an object $L^i$ with $L^i(a^i_1)=B^{i+1}\oplus B^{i+1}$
by changing certain rows and columns in $I_2\otimes B^{i+1}$, and an isomorphism
$\vartheta^{is}(L^s)\stackrel{\varsigma}\rightarrow L^i$. Let
$\varphi^i=\vartheta^{is}(\varphi^s)\varsigma^{-1}:L^i\rightarrow M^i$, which is not a retraction by 1).
Thus there is a lifting $\wt \varphi^i:L^i\rightarrow E^i$ of $\varphi^i$ with
$\varphi^i=\wt \varphi^i \pi^i$, since $\pi^i: E^i\rightarrow M^i$ is
right almost split in $R(\mf{A}^i)$ by the assumption on $(e^i)$. The
triangle and the square below are both commutative:
$$\xymatrix {&L^i\ar[dl]_{\wt
\varphi^i}\ar[dr]^{\varphi^i}&\\E^i\ar[rr]_{\pi^i}&&
M^i}\quad\quad\quad \xymatrix {L^i_X\ar[r]^{L^i(a^i_1)}\ar[d]_{\wt
\varphi^i_X}&L^i_Y\ar[d]^{\wt \varphi^i_Y}\\
E^i_X\ar[r]_{E^i(a^i_1)}&E^i_Y}$$

3)\, According to Lemma 3.1.3, it may be assumed that the sequence $(e^i)\in R(\mf B^i_0)$ with
$\iota^i_{_Z}=(0\,\, I_Z),\pi^i_{_Z}={I_Z\choose 0},\forall\,Z\in \T^i$, then
$E^i(a^i_j)={{M^i(a^i_j)\,\,\,K^i_j}\,\,\,\choose{\quad 0\quad M^i(a^i_j)}}$.
The commutative triangle
forces $\wt \varphi_Z={{I_Z\,\,C_Z}\choose {0_Z\,\, D_Z}}$ for each $Z\in \T^i$. The commutative
square for $j=1$ gives an equality
$$\left
(\begin{array}{cc}I_X&C_X\\&D_X\end{array}\right )\left
(\begin{array}{cc}B^{i+1}&K^i_1\\& B^{i+1}\end{array}\right
)=\left (\begin{array}{cc}B^{i+1}&\\&B^{i+1}\end{array}\right
)\left (\begin{array}{cc}I_Y&C_Y\\&D_Y\end{array}\right ).$$
Let $\hat E^i=\{\hat E^i_Z\mid$ dim$(E^i_Z)=2m^i_{_Z}, Z\in\T^i\}$ be a set of vector spaces.
Define a set of maps $\eta:
E^i\rightarrow\hat E^i$, such that $\eta_{_{X}}={{I_{X}\,\,C_{X}}\choose
{\,0\,\,\,\,\,I_{X}}}$,
$\eta_{_{Y}}={{I_{Y}\,\,C_{Y}}\choose {\,0\,\,\,\,\,I_{Y}}}$,
and $\eta_{_Z}=I_{2m^i_{_Z}}$ for $Z\in \T^i\setminus\{X,Y\}$; $\eta(v_j)=0$
for any $j=1,\cdots,m^i$.
Let $\hat E^i(a^i_j)=\eta_{_{s(a^i_j)}}E^i(a^i_j)\eta_{_{t(a^i_j)}}^{-1}$ for $j=1,\cdots,n^i$,
an object $\hat E^i=\eta E^i\eta^{-1}\simeq E^i$ in $R(\mf A^i)$ with $\hat
E^i(a^i_1)=B^{i+1}\oplus B^{i+1}$ is obtained as desired.\hfill$\Box$

\medskip

Suppose in the sequence below, the first part from the $0$-th pair up to the $s$-pair is
given by Formula (3.2-1) with respect to the indecomposable object $M=\vartheta^{0i}(M^s)\in R(\mf A)$:
$$(\mf A,\mf B),(\mf A^1,\mf B^1),\cdots,
(\mf A^s,\mf B^s),(\mf A^{s+1},\mf B^{s+1})\cdots,(\mf A^\varepsilon,\mf B^\varepsilon),
\cdots,(\mf A^\tau,\mf B^\tau).\eqno {(3.2\mbox{-}3)}$$

Firstly, it is assumed that in the sequence (3.2-3), $\mf B^s$ is local, $\T^s=\{X\}$; $\mf B^{s+1}$ is induced from $\mf B^s$
by a loop mutation; the reduction from $\mf B^i$ to $\mf B^{i+1}$
is given by a localization followed by a regularization, such that $R^{i+1}=k[x,\phi^{i+1}(x)^{-1}]$
for $s<i<\tau$; and $\mf B^\tau$ is minimal.

\medskip

{\bf Corollary 3.2.4}\, Suppose $M^\tau$ is an object of $R(\mf B^\tau)$ with $M^\tau_{X}=k,M^\tau(x)=(\lambda),\phi^\tau(\lambda)\ne 0$.
If $\vartheta^{0\tau}(M^\tau)=M\in R(\mf B)$ is homogeneous with an almost split conflation $(e)$, then
there exists an almost split conflation $(e^\tau_\lambda)$ given by Formula (3.1-2) in $R(\mf B^\tau)$,
such that $\vartheta^{0\tau}(e^\tau_\lambda)\simeq (e)$.

\smallskip

{\bf Proof}\, Set $M^s=\vartheta^{s\tau}(M^\tau)$, then $M^s(a^s_1)=(\lambda)$. Theorem 3.2.3 gives an
almost split conflation $(e^s)$ in $R(\mf B^s)$ with $\vartheta^{0s}(e^s)\simeq (e)$.
Furthermore, $R(\mf A^{s+1})$ is equivalent to $R(\mf A^{s})$, and for
$i> s$, $R(\mf A^{i+1})$ is equivalent to a subcategory of $R(\mf A^{i})$ consisting of the
objects $M^{i}$ with $M^i(x)=(\lambda),\phi^{i+1}(x)$.
The assertion follows by the fact that $\phi^{i+1}(x)\mid\phi^\tau(x)$, and induction on $i=s+1,\cdots,\tau-1$.\hfill$\Box$

\medskip

Secondly, it is assumed that in the sequence (3.2-3) the bocs $\mf B^s$ has two vertices
$X,Y$, and the first arrow $a^s_1:X\rightarrow X$ with $\dz(a_1^s)=0$. The bocs $\mf B^{s+1}$ is
induced from $\mf B^s$ by a loop mutation $a^s_1\mapsto (x)$, such that $a^{s+1}_j$ is either a
loop at $X$, or an edge from $X$ to $Y$ for $j=1,\cdots,\tau-s$. In particular,
there exists a certain index $s<\varepsilon<\tau$, such that $a^{s+1}_{\varepsilon-s}:X\rightarrow Y$ is an edge.
The reduction from $\mf B^i$ to $\mf B^{i+1}$ is given by one of the following three cases:

\subitem\hspace{-6mm}{(i)}\, when $s< i<\varepsilon$, if $a^i_1:X\rightarrow X$, a localization followed
by a regularization are made with $R^{i+1}=k[x,\phi^{i+1}(x)^{-1}]$; if $a_1^i:X\rightarrow Y$, a
regularization, or a reduction given by proposition 2.2.6 is made;

\subitem \hspace{-6mm}(ii) when $i=\varepsilon$, a reduction given by proposition 2.2.7 is made,
the induced bocs $\mf B^{\epsilon+1}$ is local with a vertex $Z$;

\subitem \hspace{-6mm}(iii) when $\varepsilon<i<\tau$, then $a^i_1:Z\rightarrow Z$, a localization followed by a regularization are made.
Finally, $\mf B^\tau$ is minimal with $R^\tau=k[x,\phi^\tau(x)^{-1}]$.

\medskip

{\bf Corollary 3.2.5}\, Suppose $M^\tau$ is an object of $R(\mf B^\tau)$, $M^\tau_{Z}=k,M^\tau(z)=(\lambda)$ with $\phi^\tau(\lambda)\ne 0$.
If $\vartheta^{0\tau}(M^\tau)=M\in R(\mf B)$ is homogeneous with an almost split conflation $(e)$, then
there exists an almost split conflation $(e^\tau_\lambda)$ given by Formula (3.1-2) in $R(\mf B^\tau)$,
such that $\vartheta^{0\tau}(e^\tau_\lambda)\simeq (e)$.

\smallskip

{\bf Proof}\, Set $M^s=\vartheta^{st}(M^t)$, then $M^s_X=k,M^s_Y=k$, and
$M^s(a^s_1)=(\lambda)$. Theorem 3.2.3 gives an
almost split conflation $(e^s):M^s\rightarrow E^s\rightarrow M^s$ in $\mf B^s$ with $\vartheta^{0s}(e^s)\simeq(e)$.
Since $R(\mf B^{s+1})$ is equivalent to $R(\mf B^s)$, we may suppose for some $i>s$, there exists an almost split conflation
$$\begin{array}{c}(e^i):M^i\stackrel{\iota^i}\rightarrow
E^i\stackrel{\pi^i}\rightarrow M^i\in R(\mf B^i)\quad \mbox{with}\quad
M^i=\vartheta^{it}(M^t),\,\, \vartheta^{si}(e^i)\simeq(e^s). \end{array}$$
An almost split conflation $(e^{i+1})$ in $R(\mf B^{i+1})$ with $\vartheta^{i,i+1}(e^{i+1})\simeq(e^i)$
will be constructed according to cases (i)--(iii) stated before the corollary.

(i) A regularization for an edge gives an equivalence $R(\mf B^{i+1})\simeq R(\mf B^i)$. And the proof of a regularization
for a loop is similar to that of Corollary 3.2.4.
Suppose $\dz(a^i_1)=0$ in $\mf B^i$, and a reduction of Proposition 2.2.6 is made by $a^i_1\mapsto (0)$.
Then $M^i(a^i_1)$ must be $(0)$. By Lemma 3.1.3 (ii), it may be assumed that
$\iota^i=(0\, 1),\pi^i={{1}\choose{0}}$, thus $E^i(a^i_1)={{0\,b}\choose{0\, 0}}$. Define an object
$L\in R(\mf B^i)$ of size $2\m^i$, such that $L(x)=J_2(\lambda), L(a^i_j)=0,\forall j$;
and a morphism $g: L\rightarrow M^i$, such that
$g_{_X}={{1}\choose{0}}=g_{_Y}$, $g(v)=0$ for any dotted arrow $v$.
Then $g$ is not a split epimorphism. Thus there exists a lifting $\tilde g:L\rightarrow E^{i}$
with $\tilde g_X={{1\,  c}\choose{0\,d}}, \tilde g_Y={{1\,  c'}\choose{0\,d'}}$. Since $\tilde g$
is a morphism, $L^i(a^i_1)\tilde g_{X}=\tilde g_{Y}E^i(a^i_1)$, which leads to
$0={{0\, b}\choose{0\,0}}$.
Therefore $b=0,E^i(a^i_1)={{0\,0}\choose{0\,0}}$. Set $E^{i+1}\in R(\mf B^{i+1})$ with
$E^{i+1}(x)=E^i(x), E^{i+1}(a^{i+1}_{j-1})=E^i(a^i_j)$
for $j=2,\cdots,n^i$, then
$\vartheta^{i,i+1}(E^{i+1})=E^i$.

(ii) Proposition 2.2.7 ensures a possibility that $M^\varepsilon(a^\varepsilon_1)=(1)$, so
$E^\varepsilon(a^\varepsilon_1)={{1\, b}\choose{0\,1}}$. Define a set of matrices: $\{\xi_{_Y}={{1\, -b}\choose{0\,\,1\,}},\xi_{_X}=I_2,
\xi(v^\varepsilon_j)=0\}$, and an object $\bar E^\varepsilon\in R(\mf B^\varepsilon)$: $\bar E^{\varepsilon}_X=k^2=\bar E^\varepsilon_Y$,
$\bar E^{\varepsilon}(x)=E^\varepsilon(x)$, $\bar E^{\varepsilon}(a^\varepsilon_j)=\xi_{s(a_j^\varepsilon)}^{-1}E^i(a_j^\varepsilon)\xi_{t(a^\varepsilon_j)}$,
then $\bar E^{\varepsilon}(a_1^{\varepsilon})=I_2$.
Let $E^{\varepsilon+1}\in R(\mf B^{\varepsilon+1})$ with $E^{\varepsilon+1}_Z=k^2$; $E^{\varepsilon+1}(z)=E^\varepsilon(x)$,
$E^{\varepsilon+1}(a_{j-1}^{\varepsilon+1})=\bar E^\varepsilon(a^\varepsilon_{j})$ for $j>1$. Then
$\vartheta^{\varepsilon,\varepsilon+1}(E^{\varepsilon+1})=\bar E^\varepsilon\simeq E^\varepsilon$ in $R(\mf B^\varepsilon)$.

(iii)\, Since $\phi^\tau(\lambda)\ne 0$ and $\phi^{i+1}(x)\mid\phi^t(x)$, $\phi^{i+1}(\lambda)\ne 0$.
There is an object $E^{i+1}\in R(\mf B^{i+1})$ with $\vartheta^{i,i+1}(E^{i+1})\simeq E^i$.

By induction, there is $(e^\tau_\lambda)\in R(\mf B^\tau)$
with $\vartheta^{s\tau}(e^\tau_\lambda)\simeq e^s_\lambda$,
thus $\vartheta^{0\tau}(e^\tau_\lambda)\simeq (e)$. \hfill$\Box$

\medskip

{\bf Lemma 3.2.6}\, (i)\, Suppose that $f(x,y)=\sum_{i,j\geqslant 0}\alpha_{ij}x^iy^j\in k[x,y]$ with $f(\lambda,
\mu)\ne 0$. Let $W_{\lambda}, W_{\mu}$ be Weyr matrices of size
$m,n$ and eigenvalues $\lambda,\mu$ respectively, and
$V=(v_{ij})_{m\times n}$ with $\{v_{ij}\mid 1\leqslant
i\leqslant m, 1\leqslant j\leqslant n\}$ being $k$-linearly
independent. Let $f(W_\lambda,W_\mu)V=\sum_{i,j\geqslant
0}\alpha_{ij}W_\lambda^iVW_\mu^j= (u_{ij})_{m\times n}$. Then
$\{u_{ij}\mid 1\leqslant i\leqslant m, 1\leqslant j\leqslant n\}$
is also $k$-linearly independent.

(ii)\, Let $\mf B$ be a bocs with $R=R_X\times R_Y$, where
$R_X=k[x,\phi_X(x)^{-1}]$, $R_Y=k[y,\phi_Y(y)^{-1}]$,
and $a_i:X\rightarrow Y$. Define $\dz^0(a_i)$ to be a {\it part
of $\delta(a_i)$} without terms involving any solid arrow.
It is possible that $X=Y$, in this case $x$ stands for the
multiplying a dotted arrow from the left and $y$ from the right. Suppose
$$\left\{\begin{array}{l}\dz^0(a_1)=f_{11}(x,y)v_{1}\\
\dz^0(a_2)=f_{21}(x,y)v_1+f_{22}(x,y)v_2,\\
\quad \quad\cdots\quad \cdots\\
\dz^0(a_n)=f_{n1}(x,y)v_1+f_{n2}(x,y)v_2+\cdots
f_{nn}(x,y)v_n,\end{array}\right.$$ where $f_{ii}(x,y)\in R_X\times
R_Y$ are invertible for $i=1,2,\cdots,n$. If $x\mapsto W_X$ of size $m$
with eigenvalue $\lambda$ and $\phi_{_X}(\lambda)\ne 0$,
$y\mapsto W_Y$ of size $n$ with eigenvalue $\mu$ and $\phi_{_Y}(\mu)\ne 0$,
then the solid arrows splitting from $a_1,\cdots, a_n$ are all going
to $\emptyset$ by regularizations in further reductions.

\smallskip

{\bf Proof}\, (i)\, Since $u_{ij}=f(\lambda,\mu)v_{ij}+
\sum_{(i',j')\prec(i,j)}f_{i'j'}(\lambda,\mu)v_{i'j'}$ with
$f_{i'j'}(x,y)\in k[x,y]$, the assertion follows by induction on
the ordered index set $\{(i,j)\mid 1\leqslant i\leqslant m,
1\leqslant j\leqslant n\}$.

(ii) The conclusion follows by (i) inductively on $1,\cdots,n$.\hfill$\Box$

\bigskip
\bigskip
\noindent {\large\bf 3.3\, Minimal wild bocses}
\bigskip

In this subsection five classes of minimal wild bocses  is defined
in order to prove the main theorem.
Our classification relies on the Drozd's wild
configurations with refinements at some last reduction stages.

\medskip

{\bf Proposition 3.3.1} (\cite{D1},\cite{CB1})\, Let
$\mathfrak{B}=(\Gamma, \Omega)$ be a bocs with a layer $L=(R;\omega;a_1,\ldots,a_n;\\v_1,\ldots,v_m)$
and suppose $a_1:X\rightarrow Y$. If $\mf B$ is of representation wild type, then
it is bound to meet one of the following configurations at some stage of reductions:

{\bf Case 1}\, $X\in \T_1,\, Y\in \T_0$ (or dually $X\in \T_0,\, Y\in \T_1$), $\dz(a_1)=0$.

{\bf Case 2}\, $X, Y\in \T_1$ (possibly $X=Y$), $\dz(a_1)=f(x,y)v_1$ for some
non-invertible $f(x,y)$ in $k[x,y,\phi_X(x)^{-1},\phi_Y(y)^{-1}]$.\hfill$\Box$

\medskip

Some notations will be fixed first before the classification. There is a decomposition for any
non-zero polynomial $f(x,y)\in k[x,y]$:
$$\begin{array}{c}f(x,y)=\alpha(x)h(x,y)\beta(y),\quad\mbox{where}\,\,\,\alpha(x)\in k[x],\,\,
\beta(y)\in k[y];\end{array}\eqno {(3.3\mbox{-}1)}$$
and every irreducible factor of $h(x,y)$ contains both $x$ and $y$, or $h(x,y)\in k^\ast$ with $k^\ast=k\setminus\{0\}$.
Sometimes $\bar x$ is used instead of $y$.

Let $k(x,y,z)$ be
the fractional field of the polynomial ring $k[x,y,z]$
of three indeterminates. Consider a vector space $\mathcal S$ generated by the dotted arrows
$\{v_1,\cdots,v_m\}$ of a bocs $\mf B$ over $k(x,y,z)$. Suppose there is a linear combination:
$$\begin{array}{c}G=f_1(x,y,z)v_1
+\cdots+f_m(x,y,z)v_m,\quad f_i(x,y,z)\in k[x,y,z].\end{array}\eqno{(\ast)}$$
Let $h(x,y,z)$ be the greatest common factor of $f_1,\cdots,f_m$, then $f_1/h,$ $\cdots,$ $f_m/h$
are co-prime, and $G=h\sum_{i=1}^m(f_i/h)v_j$. There exists some $s_i
\in k[x,y,z]$ for $i=1,\cdots,m$, such that $\sum_{i=1}^ms_i(f_i/h)=c(x,y)\in k[x,y]$.
Since  $S=k[x,y,z,c(x,y)^{-1}]$ is a Hermite ring, there exists
some invertible $F(x,y,z)\in \IM_{m}(S)$ with the first column $(f_1/h,\cdots,f_m/h)$.
A base change of the form $(w_1,\cdots,w_m)=(v_1,\cdots,v_m)F$ is made, thus
$G=h(x,y,z)w_1$.

\medskip

{\bf Classification 3.3.2}\, Let $\mf B^0$ be a wild bocs
given by  Proposition 3.3.1.
Then we are bound to meet an induced bocs $\mf B$ with a
layer $L=(R;\omega;a_1,\ldots,a_n;v_1,\ldots,v_m)$ in one of the five classes at
some stage of reductions. And a bocs in those classes is said to be {\it minimal wild},
which might be written briefly by MW.

\smallskip

Suppose the bocs $\mf B$ has two vertices $\T=\{X,Y\}$, such that the
induced local bocs $\mf B_X$ is tame infinite with $R_X=k[x,\phi_{_X}(x)^{-1}]$.

\smallskip

{\bf MW1}\, $\mf B_Y$ is finite with $R_Y=k1_Y$, and $\dz(a_1)=0$:
$$\xymatrix {X\ar[rr]^{a_1}\ar@(ul, dl)[]_{x}&&Y}.$$

\medskip

{\bf MW2}\, $\mf B_Y$ is tame infinite with $R_Y=k[y,\phi_{_Y}(y)^{-1}]$, and
$\dz(a_1)=f(x,y)v_1$, such that $f(x,y)$ is non-invertible in $k[x,y,
\phi_{_X}(x)^{-1},\phi_{_Y}(y)^{-1}]$:
$$\xymatrix {X\ar@(ul, dl)[]_{x}\ar[rr]^{a_1}&& Y\ar@(ur,
dr)[]^{y}}.$$

Suppose now we have a local bocs $\mf B$ with $R=k[x,\phi(x)^{-1}]$:
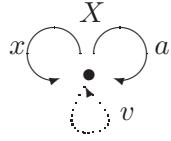
\begin{center}
\unitlength=0.5pt
\begin{picture}(60,60)
\put(-10,30){\oval(40,40)[t]} \put(-10,30){\oval(40,40)[bl]}
\put(-10,10){\vector(3,1){5}} \put(40,30){\oval(40,40)[t]}
\put(40,30){\oval(40,40)[br]} \put(40,10){\vector(-3,1){5}}

\put(11,9){$\bullet$} \qbezier[10](13,5)(-8,-28)(16,-30)
\qbezier[10](13,5)(50,-28)(13,-30) \put(17,-1){\vector(-1,2){3}}

\put(9,54){\makebox{$X$}}\put(40,-20){$v$}
\put(-44,30){\makebox{$x$}} \put(66,30){\makebox{$a$}}
\end{picture}
\end{center}\vskip 3mm

{\bf MW3}\, The differential $\dz^0$ of the solid arrows of $\mf B$ is given by
$$\begin{array}{c}
\left\{
\begin{array}{ccl}
\dz^0(a_{1}) &=& f_{11}(x,\bar x)w_{1},\\
\cdots   &  & \qquad\cdots \\
\dz^0(a_{n}) &=& f_{n1}(x,\bar x)w_{1}
+\cdots+f_{nn}(x,\bar x)w_{n},
\end{array} \right.
\end{array}\eqno {(3.3\mbox{-}2)}$$
where $w_i$ is obtained by base changes;
$f_{ii}(x,\bar x)=\alpha_{ii}(x)h_{ii}(x,\bar x)\beta_{ii}(\bar x)$ by Formula (3.3-1),
such that $h_{ii}(x,x)$ is invertible in $k[x,\phi(x)^{-1}]$ for $1\leqslant i\leqslant n$;
and there is some minimal $1\leqslant s\leqslant n$, such that
$f_{ss}(x,\bar x)$ is non-invertible in $k[x,\bar x,\phi(x)^{-1},\phi(\bar x)^{-1}]$.

\smallskip

Suppose there exists some $1\leqslant n_1\leqslant n$, such that:
\begin{equation*}
\left\{
\begin{array}{cllll}
\dz^0(a_{1}) &=& f_{11}(x,\bar x)w_{1},\\
\cdots   &  & \qquad\cdots &\\
\dz^0(a_{n_1-1}) &=& f_{n_1-1,1}(x,\bar x)w_{1}
&+\cdots&+f_{n_1-1,1}(x,\bar x)w_{n_1-1},\\
\dz^0(a_{n_1}) &=& f_{n_1,1}(x,\bar x)w_{1}
&+\cdots& +f_{n_1,n_1-1}(x,\bar x)w_{n_1-1}+f_{n_1,n_1}(x,\bar x)\bar w,
\end{array} \right. \eqno {(3.3\mbox{-}3)}
\end{equation*}
where $f_{ii}(x,x)$ for $1\leqslant i<n_1$ are invertible in $k[x,\phi(x)^{-1}]$;
$\bar w=0$, or $\bar w\ne 0$ but $f_{n_1,n_1}(x,x)=0$.
Denote by $x_1$ the solid arrow $a_{n_1}$, there exists a polynomial $\psi(x,x_1)$
being divided by $\phi(x)$. Write $\delta^1$ the part of differential $\delta$
by deleting all the monomials involving any solid arrow except $x,x_1$. Suppose the further
unraveling for $x$ is restricted to $x\mapsto (\lambda)$ with $\psi(\lambda,x_1)\ne 0$. Then
$$\left\{
\begin{array}{llll}
\dz^1(a_{n_1+1}) &=K_{n_1+1}&+f_{n_1+1,n_1+1}(x,x_1,\bar x_1)w_{n_1+1},&\\
&\cdots \cdots&\cdots \\
\dz^1(a_n) &=K_n&+f_{n,n_1+1}(x,x_1,\bar x_1)w_{n_1+1}&+ \cdots+
f_{nn}(x,x_1,\bar x_1)w_{n},
\end{array} \right.\eqno {(3.3\mbox{-}4)}$$
where $K_{i}=\sum_{j=1}^{n_1-1}f_{ij}(x,x_1,\bar x_1)w_j$,
$w_i$ are given by the base changes described below Formula $(\ast)$ inductively, and
$f_{ii}(x,x_1,\bar x_1)$ are invertible in $k[x,x_1,\bar x_1,
\psi(x,x_1)^{-1},\psi(x,\bar x_1)^{-1}]$ for $n_1<i\leqslant n$.

\smallskip

{\bf MW4}\, $\bar w=0$, or $\bar w\ne 0$ but $(x-\bar x)^2\mid f_{n_1n_1}(x,\bar x)$ in Formula (3.3-3).

\smallskip

{\bf MW5}\, $\bar w\ne 0$ and $(x-\bar x)^2\nmid f_{n_1n_1}(x,\bar x)$ in Formula (3.3-3).\hfill$\Box$

\medskip

The proof of Classification 3.3.2 depends on Classification 3.3.5 of local bocses at the end of the subsection,
while the proof of 3.3.5 is based on formulae (3.3-2)-(3.3-9) and Lemma 3.3.3--3.3.4 below.

Let $\mf B$ be a local bocs having a layer
$L=(R;\omega; a_1, \cdots, a_n; v_1,\cdots,v_m)$. If $R=k1_X$ is trivial,
then the differentials of the solid arrows have two possibilities. First,
\begin{equation*}
\left\{
\begin{array}{ccl}
\dz^0(a_1) &=& f_{11}w_1,\\
\cdots   &  & \qquad\cdots \\
\dz^0(a_{n}) &=& f_{n1}w_1 +\cdots+f_{nn}w_{n},
\end{array} \right. \eqno {(3.3\mbox{-}5)}
\end{equation*}
where $f_{ij}\in k,f_{ii}\ne 0$ for $1\leqslant i\leqslant n$. Second, there exists
some $1\leqslant n_0\leqslant n$, such that:
\begin{equation*}
\left\{
\begin{array}{ccl}
\dz^0(a_1) &=& f_{11}w_1,\\
\cdots   &  & \qquad\cdots \\
\dz^0(a_{n_0-1}) &=& f_{n_0-1,1}w_1 +\cdots +f_{n_0-1,n_0-1}w_{n_0-1},\\
\dz^0(a_{n_0}) &=& f_{n_01}w_1 \quad +\cdots +f_{n_0,n_0-1}w_{n_0-1},
\end{array} \right. \eqno {(3.3\mbox{-}6)}
\end{equation*}
where $f_{ij}\in k,f_{ii}\ne 0$ for $1\leqslant i<n_0$.
Set $a_i\mapsto \emptyset,i=1,\cdots, n_0-1$,
by a series of regularization, then $a_{n_0}\mapsto (x)$ by a loop mutation,
an induced local bocs $\mf B'$ is obtained.

\medskip

Without loss of generality, the bocs $\mf B'$ may still be denoted by $\mf B$ with a layer $L$, but $R=k[x]$.
The differentials $\dz^0$
have again two possibilities. The first one is given by Formula (3.3-2),
such that $h_{ii}(x,x)\ne 0$, i.e. $(x-\bar x)\nmid f_{ii}(x,\bar x)$
for $i=1,\cdots,n$. Define a polynomial:
$$\begin{array}{c}\phi(x)=\prod_{i=1}^{n}c_i(x)h_{ii}(x,x),\end{array}\eqno {(3.3\mbox{-}7)}$$
where $c_i(x)$ appears at the localization in order to do a base
change before the $i$-th step of a regularization.

\medskip

{\bf Lemma 3.3.3}\, Let $\mf B$ be a bocs given by Formula (3.3-2) with a polynomial $\phi(x)$ given by
Formula (3.3-7). There exist two cases:

(i)\, $f_{ii}(x,\bar x)$ are invertible in $k[x,\bar x,\phi(x)^{-1}\phi(\bar x)^{-1}]$
for $1\leqslant i\leqslant n$;

(ii)\, $f_{ss}(x,\bar x)$ is not invertible in $k[x,\bar x,\phi(x)^{-1}\phi(\bar x)^{-1}]$
for some minimal $1\leqslant s\leqslant n$.\hfill$\Box$

\medskip

The second possibility of the differential $\dz^0$ in the case of $R=k[x$] is given by Formula (3.3-3)
for some fixed $1\leqslant n_1\leqslant n$,
such that $f_{ii}(x,x)\ne 0$ for $1\leqslant i<n_1$, and $\bar w=0$, or
$\bar w\ne 0$ but $f_{n_1,n_1}(x,x)=0$. Let
$$\phi(x)=\left\{\begin{array}{ll}\prod_{i=1}^{n_1-1}c_i(x)f_{ii}(x,x),&\bar w=0;\\[0.5mm]
c_{n_1}(x)\prod_{i=1}^{n_1-1}c_i(x)f_{ii}(x,x),&\bar w\ne 0,\end{array}\right.\eqno {(3.3\mbox{-}8)}$$
Thus, under the restriction $x\mapsto (\lambda),\phi(\lambda)\ne 0$, an induced bocs given by regularazations
with the first arrow $a_{n_1}$ is obtained.
There are two possibilities in the further reductions. The first
possibility is given by Formula (3.3-4),
such that $f_{ii}(x,x_1,x_1)\ne 0$ for $n_1<i\leqslant n$. There is a sequence
of localizations given by the polynomials $c_{i}(x,x_1)$
appeared before the base changes in order to do regularizations. Let
$$\begin{array}{c}\psi(x,x_1)=\phi(x)
\prod_{i=n_1+1}^n c_i(x,x_1)f_{ii}(x,x_1,x_1),\end{array}\eqno {(3.3\mbox{-}9)}$$

{\bf Lemma 3.3.4}\, Let the differentials in the bocs $\mf B$ be given by Formulae (3.3-3)--(3.3-4)
with polynomials $\phi(x)$ in (3.3-8), and $\psi(x,x_1)$ in (3.3-9). There exist two cases.

(i)\, There exists some $\lambda\in k$ with $\psi(\lambda,x_1)\ne0$, and
a minimal $n_1+1\leqslant s\leqslant n$, such that
$f_{ss}(\lambda,x_1,\bar x_1)$ is non-invertible in $k[x_1,\bar x_1,
\psi(\lambda,x_1)^{-1}\psi(\lambda,\bar x_1)^{-1}]$, i.e.,
after making an unraveling $x\mapsto (\lambda)$, followed by a series of
regularizations $a_i\mapsto\emptyset,w_i=0$
for $i=1,\cdots,n_1-1$, the induced local bocs $\mf B_{(\lambda)}$
with $R_{(\lambda)}=k[x_1,\psi(\lambda,x_1)^{-1}]$ is in case (ii) of Lemma 3.3.3.

(ii)\, For any $\lambda\in k$ with $\psi(\lambda,x_1)\ne0$,
$f_{ii}(\lambda,x_1,\bar x_1)$ are invertible for $n_1<i\leqslant n$ in
$k[x_1,\bar x_1,\psi(\lambda,x_1)^{-1}\psi(\lambda,\bar x_1)^{-1}]$,
i.e., the induced bocs $\mf B_{(\lambda)}$
with $R_{(\lambda)}=k[x_1,\psi(\lambda,x_1)^{-1}]$ is in case (i) of 3.3.3.

Case (ii) is equivalent to (ii)$'$: $f_{ii}(x,x_1,\bar x_1)$ are invertible
in $k[x,x_1,\bar x_1,\psi(x,x_1)^{-1},\psi(x,\bar x_1)^{-1}]$ for $n_1<i\leqslant n$.

\smallskip

{\bf Proof}\, It is only need to prove the equivalence of (ii) and (ii)$'$.

(ii)$\Longrightarrow$(ii)$'$\, If there exists some $n_1< s\leqslant n$ with
$f_{ss}(x,x_1,\bar x_1)$ non-invertible, then it contains a non-trivial factor
$g(x,x_1,\bar x_1)$ coprime to $\psi(x,x_1)\psi(x,\bar x_1)$. Consider the
variety $V=\{(\alpha,\beta,\gamma)\in k^3\mid g(\alpha,\beta,\gamma)=0,
\psi(\alpha,\beta)\psi(\alpha,\gamma)=0\}$.
Since dim$(V)\leqslant 1$, there exists a co-finite subset $\mathscr L\subset k$, such that
$\forall\,\lambda\in\mathscr L$, the plane $x=\lambda$
of $k^3$ intersects $V$ at only a finite number of points. Thus $g(\lambda,x,\bar x_1)$
and $\psi(\lambda,x_1)\psi(\lambda,\bar x_1)$ are coprime. Consequently $g(\lambda,x,\bar x_1)$,
thus $f_{ss}(\lambda,x_1,\bar x_1)$ is not invertible in $k[x_1,\bar x_1,
\psi(\lambda,x_1)^{-1}\psi(\lambda,\bar x_1)^{-1}]$.

(ii)$'\Longrightarrow$(ii)\, If $f_{ii}(x,x_1,\bar x_1)$ is invertible in $k[x,x_1,\bar x_1,
\psi(x,x_1)^{-1}\psi(x,\bar x_1)^{-1}]$, then for any
$\lambda\in k$ with $\psi(\lambda,x_1)\ne 0$, $f_{ii}(\lambda,x_1,\bar x_1)$ is
invertible in $k[x_1,\bar x_1,\psi(\lambda,x_1)^{-1}\psi(\lambda,\bar x_1)^{-1}]$. \hfill$\Box$

\medskip

The second possibility of $\dz^1$ is: there exists some $n_2$ with $n_1<n_2\leqslant n$,
such that
\begin{equation*}
\left\{\begin{array}{l}
\dz^1(a_{n_1+1})= K_{n_1+1}+f_{n_1+1,n_1+1}(x,x_1,\bar x_1)w_{n_1+1},\\
\qquad \cdots\quad \cdots \\
\dz^1(a_{n_2-1}) =K_{n_2-1}
+\cdots\cdots+f_{n_2-1,n_2-1}(x,x_1,\bar x_1)w_{n_2-1},\\
\dz^1(a_{n_2})\quad=\,\,\,K_{n_2}\,\,\,+\cdots\cdots\,\,\,\, +f_{n_2,n_2-1}(x,x_1,\bar x_1)w_{n_2-1}
+f_{n_2,n_2}(x,x_1,\bar x_1)\bar w',\end{array} \right. \eqno {(3.3\mbox{-}10)}
\end{equation*}
where $K_{i}=\sum_{j=1}^{n_1-1}f_{ij}(x,x_1,\bar x_1)w_j$ for $n_1< i\leqslant n_2$,
$f_{ii}(x,x_1,x_1)\ne 0$ for $n_1< i<n_2$, and $\bar w'=0$, or
$\bar w'\ne 0$ but $f_{n_2,n_2}(x,x_1,x_1)=0$. Define a polynomial
$$\psi_1(x,x_1)=\left\{\begin{array}{ll}\phi(x)
\prod_{i=n_1+1}^{n_2-1}c_i(x,x_1)f_{ii}(x,x_1,x_1),
&\mbox{if}\,\,\bar w'=0;\\[0.5mm]
c_{n_2}(x,x_1)\phi(x)\prod_{i=n_1+1}^{n_2-1}c_i(x,x_1)f_{ii}(x,x_1,x_1),
&\mbox{if}\,\,\bar w'\ne 0.\end{array}\right.\eqno {(3.3\mbox{-}11)}$$

Suppose a bocs $\mf B$ is med, the differential of which is given by Formula (3.3-3) and (3.3-10)
with a polynomial $\psi_1(x,x_1)$ of (3.3-11).
Fix any $\lambda_0\in k$
with $\psi_1(\lambda_0,x_1)\ne 0$, there is an induced
bocs $\mf B_{(\lambda_0)}$ with $R_{(\lambda_0)}=k[x_1,\psi_1(\lambda_0,x_1)^{-1}]$
given by an unraveling $x\mapsto(\lambda_0)$, and then
a series of regularization $a_i\mapsto\emptyset,w_i=0$ for $i=1,\cdots,n_1-1$.
There exist three cases:

1) $\mf B_{(\lambda_0)}$ is in case  (i) of
Lemma 3.3.4, then there exists some $\lambda_1$ with $\psi(\lambda_0,\lambda_1)\ne 0$,
such that after sending $x_1\mapsto (\lambda_1)$ by an unraveling, followed by a series of regularizations, the induced bocs
$\mf B_{(\lambda_0,\lambda_1)}$ satisfies Lemma 3.3.3 (ii);

2) $\mf B_{(\lambda_0)}$ is in case (ii)$'$ of Lemma 3.3.4;

3) $\mf B_{(\lambda_0)}$ is in the case
of Formulae (3.3-3) and (3.3-10).

\noindent In case 3), the above procedure is repeated once again
for $\mf B_{(\lambda_0)}$. By induction on the number of the finitely many solid arrows,
the case 1) or case 2) is finally reached.

\medskip

{\bf Classification 3.3.5}\, Let $\mf B$ be a local bocs with $R$ trivial, there exist four cases:

(i)  $\mf B$ has Formula (3.3-5).

(ii)  $\mf B$ has Formula (3.3-6). And by a series of regularizations
$a_i\mapsto\emptyset,i=1,\cdots,n_0-1$, followed by a loop mutation $a_{n_0}\mapsto (x)$, the induced bocs
$\mf B'$ with a polynomial (3.3-7) is in case (i) of Lemma 3.3.3.

(iii) $\mf B$ has an induced local bocs
$\mf B_{(\lambda_0,\lambda_1,\cdots,\lambda_{l})}$ for some $l<n$ in case (ii) of Lemma 3.3.3.

(iv) $\mf B$ has an induced local bocs
$\mf B_{(\lambda_0,\lambda_1,\cdots,\lambda_{l-1})}$ for some $l<n$ in case (ii) of Lemma 3.3.4.

\medskip

{\bf The proof of Classification 3.3.2}\, 1)\, Suppose a two-point wild bocs is med,
if $\mf B_X$ or $\mf B_Y$ is in case (iii) or (iv) of Classification 3.3.5,
the induced local bocs given by deleting $Y$ or $X$ may be considered, which is wild type.
Therefore it is assumed that one of bocs $\mf B_X$ or $\mf B_Y$ has Formula (3.3-5) and
another is in case (i) of Lemma 3.3.3, or both of $\mf B_X$ and $\mf B_Y$ are in case (i) of Lemma 3.3.3.
MW1 or MW2 follows.

2)\, If a local wild bocs in case  (iii) of Classification 3.3.5 is med,
then there is an induced bocs in case  (ii) of Lemma 3.3.3. MW3 is reached.

3)\, If a local wild bocs in case (iv) of Classification 3.3.5 is med,
then there is an induced bocs in case (ii)$'$ of Lemma 3.3.4. MW4 or MW5 is reached.
\hfill$\Box$

\bigskip
\bigskip
\noindent{\large\bf 3.4 Non-homogeneity in the cases of MW1-4}
\bigskip

Throughout the subsection, $(\mf A^0,\mf B^0)$ denotes any
pair of matrix bimodule problem and its associated bocs.

\medskip

{\bf Proposition 3.4.1} \cite{B1}\, If $\mf B^0$ has an induced
bocs $\mf{B}$ in the case of MW1,
then $\mf B^0$ is non-homogeneous.

\smallskip

{\bf Proof}\,  1)\, Let $\mf{B}_X$ be an induced local bocs of $\mf B$. Suppose
$\vartheta_1: R(\mf{B}_X)\rightarrow R(\mf{B})$,
$\vartheta_2:R(\mf{B})\rightarrow R(\mf{B}^0)$ are two induced functors, and
$\vartheta=\vartheta_2\vartheta_1:R(\mf{B}_X)\rightarrow R(\mf{B}^0)$.
For any $\lambda\in k$, $\phi(\lambda)\ne 0$, a
representation $S'_{\lambda}\in R(\mf{B}_X)$ given by
$(S'_{\lambda})_X=k$, $S'_{\lambda}(x)=(\lambda)$ is defined.
If $\mf B^0$ is homogeneous, then there is a co-finite subset
$\mathscr L\subseteq k\setminus\{\lambda\mid\phi(\lambda)\ne 0\}$, such that $\{\vartheta(S'_{\lambda})
\mid \lambda\in\mathscr L\}$ is a family of homogeneous iso-classes of
$R(\mf{B}^0)$. By Corollary 3.2.4, there is an almost split conflation
$(e'_{\lambda}):S'_{\lambda}\stackrel{\iota'}{\longrightarrow}
E'_{\lambda}\stackrel{\pi'}{\longrightarrow}
S'_{\lambda}$ in $R(\mf B_X)$ with $E'(x)=J_2(\lambda)$, such that $\vartheta(e'_\lambda)$
is an almost split conflation in $R(\mf B^0)$.
Fix any $\lambda\in \mathscr L$, and the conflation
$(e_{\lambda})=\vartheta_1(e'_\lambda): S_{\lambda}\stackrel{\iota}{\longrightarrow}
E_{\lambda}\stackrel{\pi}{\longrightarrow}
S_{\lambda}$ in $R(\mf{B})$, then $(S_\lambda)_X=k,(S_\lambda)_Y=0,S_\lambda(x)=(\lambda),
S_\lambda(a_i)=0$ for all solid arrow $a_i$ of $\mf B$, and $(E_\lambda)_X=k^2,(E_\lambda)_Y=0,E_\lambda(x)=J_2(\lambda),
E_\lambda(a_i)=0$. Since $\vartheta_2(e_\lambda)=\vartheta(e_\lambda')$ is almost split in $R(\mf B^0)$,
so is $(e_\lambda)$ in $R(\mf B)$ by Lemma 3.2.2 (ii) inductively.

2)\, Let $L\in R(\mf{B})$ be an object given by $L_X=L_Y=k$,
$L(x)=(\lambda)$, $L(a_1)=(1)$ and $L(a_i)=0$ for $i>1$. Let $g: L\rightarrow S_{\lz}$ be a morphism
with $g_{_X}=(1)$, $g_{_Y}=(0)$ and
$g(v)=0$ for all dotted arrows $v$ of $\mf B$. It is asserted that $g$ is not a retraction.
Otherwise, if there is a morphism $h:S_{\lz}\rightarrow L$ such
that $hg=id_{S_{\lz}}$, then $h_X=(1)$ and $h_Y=(0)$.
But $h$ being a morphism implies that
$(1)(1)=h_XL(a)=S_{\lz}(a)h_Y =(0)(0)$, a contradiction.

3)\, There exists a lifting $\tilde{g}: L\rightarrow E_{\lambda}$ with
$\tilde{g}\pi=g$. If $\tilde{g}_{_{X}}=(a,b)$, then
$\tilde{g}_{_{X}}\pi_{_{X}}=g_{_X}$, i.e., $(a,b){1\choose 0}=(1),a=1$.
But $\tilde{g}$ being a morphism implies: $\tilde{g}_{_X} E_{\lz}(x)=L
(x)\tilde{g}_{_X}$, i.e., $(1,b){{\lz\,\, 1}\choose {0\,\,
\lz}}=(\lz)(1,b)$, $(\lambda,1+b\lambda)=(\lambda,\lambda b)$, a
contradiction. Thus $\mf{B}^0$ is not homogeneous.\hfill$\Box$

\medskip

{\bf Proposition 3.4.2} \cite{B1}\, If $\mf B^0$ has an induced
bocs $\mf{B}$ in the case of MW2, then $\mf B^0$ is
non-homogeneous.

\smallskip

{\bf Proof}\, Since $f(x,y)$ is non-invertible in $k[x,y,\phi_{_X}(x)^{-1}\phi_{_Y}(y)^{-1}]$,
after dividing the dotted arrows $v_j$ by some powers
of $\phi_{_X}(x)$ and $\phi_{_Y}(y)$, it may be assumed that $f(x,y)\in k[x,y]$.
There exist three cases on $f(x,y)=\alpha(x) h(x,y)\beta(y)$ as in Formula (3.3-1).

Case \ding{172} $h(x,y)\notin k^\ast$, then
$h(x,y)$ and $\phi_{_X}(x)\phi_{_Y}(y)$ are coprime. There is an infinite set
$$\mathscr L'=\{(\lambda,\mu)\in k\times k\mid h(\lz,\mu)=0,
\phi_{_X}(\lz)\phi_{_Y}(\mu)\ne 0\}.$$
by Bezout's theorem. Clearly, $\mathscr L_X=\{\lambda\in k\mid
(\lambda,\mu)\in\mathscr L'\}$ is an infinite set.

Case \ding{173} $h(x,y)\in k^\ast$, but
there is an irreducible factor $\beta'(y)$ of $\beta(y)$ coprime to $\phi_{_Y}(y)$.
If $\beta'(\mu)=0$, then there is an infinite set $\mathscr L_X=
\{\lambda\in k\mid \phi_{_X}(\lz)\phi_{_Y}(\mu)\ne 0,\beta(\mu)=0\}$.

Case \ding{174} $h(x,y)\in k^\ast$, and $\beta(y)\mid\phi(y)^e$ for some $e\in\mathbb Z^+$, then there must exist an irreducible
factor $\alpha'(x)$ of $\alpha(x)$ coprime to $\phi_{_X}(x)$. If $\alpha'(\lambda)=0$,
there is an infinite set $\mathscr L_Y=
\{\mu\in k\mid \phi_{_X}(\lz)\phi_{_Y}(\mu)\ne 0,\alpha(\lambda)=0\}$.

The cases \ding{172}--\ding{173} are dealt with first in the following statement 1)--3).

1)\, The discussion is carried out as
in proof 1) of Proposition 3.4.1, then an infinite
set $\mathscr L\subseteq \mathscr L_X$ is obtained.

2)\, Let $L\in R(\mf{B})$ be an object given by $L_X=k=L_Y$, $L(x)=(\lambda),\lambda
\in\mathscr L; L(y)=(\mu)$ with $(\lambda,\mu)\in \mathscr L'$;
$L(a_1)=(1); L(a_i)=0$ for $i>1$.
Let $g: L\rightarrow S_{\lz}$ be a morphism in
$R(\mf{B})$ with $g_{_X}=(1)$, $g_{_Y}=(0)$ and
$g(v)=0$ for all dotted arrows $v$, then $g$ is not a retraction.
Otherwise, if there is a morphism $h:S_{\lz}\rightarrow L$ such
that $hg=id_{S_{\lz}}$, then $h_{_X}=(1)$ and $h_{_Y}=(0)$.
But $h$ being a morphism implies that
$-1=S_{\lz}(a_1)h_{_Y}-h_{_X}L(a_1)=h(\dz(a_1))
=f(\lambda,\mu)h(v)=0$, a contradiction.

3)\, There exists a lifting $\tilde{g}: L\rightarrow E_{\lambda}$
with $\tilde{g}\pi=g$. A contradiction appears as the same as in the
proof 3) of Proposition 3.4.1.

If case \ding{174} appears, then a set of homogeneous iso-classes
$\{S_\mu\mid \mu\in\mathscr L\}$ is used. Let $L$ be the same as in 2), and a morphism $g:S_\mu\rightarrow L$ be not a section.
There is an extension $\tilde g:E_\mu\rightarrow L$, which leads to a contradiction.   \hfill$\Box$

\medskip

{\bf Proposition 3.4.3} \cite{B1}\, If $\mf B^0$ has an induced bocs $\mf{B}$
in the case of MW3 with $R=k[x,\phi(x)^{-1}]$, then $\mf B^0$ is
non-homogeneous.

\smallskip

{\bf Proof}\, After a series of regularizations $a_i\mapsto\emptyset,i=1,\cdots,s-1$, it may be assumed that $s=1$ in MW3.
Since $f_{11}(x,\bar x)$ is non-invertible in $k[x,\bar x,\phi_{_X}(x)^{-1}\phi_{_X}(\bar x)^{-1}]$,
by a similar discussion as the beginning of the proof of Proposition 3.4.2,
there are infinite sets:
$$\begin{array}{c}\mathscr L_{1}=\{\lambda\mid h_{11}(\lambda,\mu)=0,\phi(\lambda)\ne 0\};\\
\mathscr L_{2}
=\{\lambda\mid\beta_{11}(\mu)=0,\phi(\lambda)\beta_{11}(\lambda)\ne 0\};\,\,
\mathscr L_{3}=\{\mu\mid \alpha_{11}(\lambda)=0, \phi(\mu)\alpha_{11}(\mu)\ne 0\}\end{array}$$
according to the cases \ding{172}--\ding{174} respectively. It is easy to see that
$\lambda\ne \mu$ in $\mathscr L_1$--$\mathscr L_3$.

Define a polynomial $\psi(x)=\phi(x)\prod_{i=1}^n\alpha_i(x)\beta_i(x)$, there is an induced bocs $\mf B'$ of $\mf B$
with $R'=k[x,\psi(x)^{-1}]$ given by a localization. Note that $\mf B'$
is not necessarily minimal. Set the induced functors
$\vartheta_1: R(\mf{B}')\rightarrow R(\mf{B})$,
$\vartheta_2:R(\mf{B})\rightarrow R(\mf{B}^0)$, and
$\vartheta=\vartheta_2\vartheta_1:R(\mf{B}')\rightarrow R(\mf{B}^0)$.
The case of $\mathscr L_1$ or $\mathscr L_2$ is dealt with first in the following 1)--3).

1) Let $S_{\lz}\in R(\mf{B}'),\lambda\in\mathscr L_1$ or $\mathscr L_2$ be an object given by $(S_{\lz})_X=k$,
$S_{\lz}(x)=(\lambda),\psi(\lz)\ne 0$, then
$S_{\lz}(a_i)=(0)$ for $1\leqslant i\leqslant n$ by Lemma 3.2.6 (ii), since
$f_{ii}(\lambda,\lambda)\ne 0$.
If $\mf B^0$ is homogeneous, then there is a co-finite subset
$\mathscr L\subseteq\mathscr L_1$ or $\mathscr L_2$, such that $\{\vartheta(S_{\lambda})\in
R(\mf{B}^0)\mid \lambda\in\mathscr L\}$ is a family of homogeneous iso-classes
of $\mf{B}^0$. By Theorem 3.2.3 with respect to $x$ of $\mf B'$, there is an almost split conflation
$(e_{\lambda}'):S_{\lambda}\stackrel{\iota}{\longrightarrow}
E_{\lambda}\stackrel{\pi}{\longrightarrow}
S_{\lambda}$ in $R(\mf B')$, such that $\vartheta(e_\lambda')$
is an almost split conflation in $R(\mf B^0)$. By Lemma 3.1.3 (ii) it may be assumed that
$\iota=(0\,1),\pi=(1\,0)^T$, thus $E_\lambda(x)={{\lambda\,c}\choose{0\,\lambda}}$.
By Lemma 3.2.6 (ii) once again, $E_\lambda(a_i)=(0)$ for $i=1,\cdots,n$, therefore $E_\lambda(x)=J_2(\lambda)$.
In fact, if $E_\lambda(x)$ was $\lambda I_2$, $(e_\lambda')$ would be splitable.
Fix any $\lambda\in\mathscr L$, and consider the
conflation $(e_\lambda)=\vartheta_1(e_\lambda')$ in $R(\mf B)$, which is almost split by Lemma 3.2.2 (ii).

2) Let $\mu\in k$ with $\phi(\mu)\ne 0$, $h_{11}(\lambda,\mu)=0$ or $\beta_{11}(\mu)=0$.
Define $L\in R(\mf{B})$ with $L_X=k^2$, $L(x)={{\lambda\, 0}\choose{0\,\mu}}$,
$L(a_1)=J_2(0)$, and $L(a_i)=0$ for $2\leqslant i\leqslant n$. $L$ is well defined.
In fact, if $\eta=\{\eta_{_X},\eta(v_j)\}:L\mapsto L$ is an isomorphism, then $L(x)\eta_{_X}=\eta_{_X}L(x)$
forces $\eta_{_X}={{a\,0}\choose{0\,b}}$, and $L(a_1)\eta_{_X}-\eta_{_X}L(a_1)=f_{11}(L(x),L(x))\eta(v_1)$ implies
${{\,0\,b-a}\choose{\,0\quad 0}}={{f_{11}(\lambda,\lambda)v_{11}\,f_{11}(\lambda,\lambda)v_{12}}\choose
{f_{11}(\mu,\mu)v_{21}\,f_{11}(\mu,\mu)v_{22}}}$, which has a solution $v_{ij}=0,a=b$.
Let $g: L \rightarrow S_{\lz}$ be a
morphism with $g_{_X}={1\choose 0}$ and $g(v_j)={0\choose 0}$ for all $j$, then
$g$ is not a retraction. Otherwise, there is a morphism $h:S_{\lz}\rightarrow L $ with
$hg=id_{S_{\lz}}$:
$$\xymatrix{L\ar@(ul,ur)^{{\lambda\,0}\choose{0\,\mu}}
\ar@(ul,dl)_{{0\,1}\choose{0\,0}}\ar@/^/[rr]^{g}&&
{S_\lambda}\ar@(ul,ur)^{(\lambda)}
\ar@(ur,dr)^{(0)}\ar@/^/[ll]^{h}}.$$
Thus $h_X=(1,b)$. Set $h(v_1)=(c,d)$. Then
$$\begin{array}{c}S_{\lz}(a_1)h_X-h_XL(a_1)=f_{11}
(S_{\lz}(x), L(x))h(v_1)\Longrightarrow\\
0(1,b)-(1,b){{0\, 1}\choose{0\, 0}}=f_{11}\big(\lambda,
{{\lambda\,0}\choose{0\,\mu}}\big)(c,d)
=\big(f_{11}(\lambda,\lambda)c, f_{11}(\lambda,\mu)d\big),\end{array}$$
which leads to $-(0,1)=(\ast,0)$, a contradiction.

3) There is a lifting $\tilde{g}: L \rightarrow E_{\lambda}$ with
$\tilde{g}\pi=g$. Set $\tilde{g}_X={{a\,\, b}\choose {c\,\, d}}$,
$\tilde{g}_{_X}\pi_{_X}=g_{_X}$ yields
$\tilde{g}_{_X}={{1\,\, b}\choose {0\,\, d}}$. On the other hand,
$\tilde{g}:L\mapsto E_{\lz}$ being a morphism leads to
$\tilde{g}_{_X} E_{\lz}(x)=L (x)\tilde{g}_{_X}$, i.e.,
${{1\,\,b}\choose {0\,\,d}}{{\lz\,\,1}\choose {0\,\,\lambda}}
={{\lz\,\, 0}\choose {0\,\, \mu}}{{1\,\, b}\choose{0\,\, d}}$, then ${{\lz\,\, 1+\lz b}\choose {0\,\,\,\,\,\,\lz d}}={{\lz\,\,
\lz b}\choose {0\,\, \mu d}}$, a contradiction. Therefore $\mf{B}^0$
is not homogeneous.

If $\mathscr L_3$ appears, then a set of homogeneous iso-classes
$\{S_\mu\mid \mu\in\mathscr L\}$ is used. Let $L$ be the same as in 2), and a morphism $g:S_\mu\rightarrow L$ be not a section.
There is an extension $\tilde g:E_\mu\rightarrow L$, which leads to a contradiction. \hfill$\Box$

\medskip

{\bf Proposition 3.4.4}\, If $\mf B^0$ has an induced bocs $\mf{B}$
in the case of MW4 with two polynomials $\phi(x)$ and $\psi(x,x_1)$, then $\mf B^0$ is
non-homogeneous.

\smallskip

{\bf Proof}\, Fix some $\lambda\in k$ with $\psi(\lambda, x_1)\ne 0$,
$\mathscr L'=\{\mu\mid\psi(\lambda,\mu)\ne 0\}\subseteq k$ is a co-finite subset.

Since $\phi(\lambda)\ne 0$, a series of regularizations $a_i\mapsto\emptyset,w_i=0$ for $i=1,\cdots,n_1-1$
and a loop mutation $a_{n_1}\mapsto (x_1)$ can be made in Formula (3.3-3). Thus
an induced bocs $\mf B'$ of $\mf B$ is obtained.
Furthermore, since $f_{ii}(\lambda, x_1,\bar x_1)$ are invertible in $k[x_1,\psi(\lambda,x_1)^{-1}
\psi(\lambda,\bar x_1)^{-1}]$ for $n_1< i\leqslant n$,
after regularizations $a_i\mapsto\emptyset, w_i=0$ in Formula (3.3-4),
an induced minimal local bocs $\mf B_\lambda$ is obtained and an induced functor $\vartheta_1$ as well.
$$\mf B_\lambda:\xymatrix {X\ar@(ul,dl)[]_{x_1}},\quad
R_\lambda=k[x_1,\psi(\lambda,x_1)^{-1}],\quad \vartheta_1: R(\mf{B}_\lambda)\rightarrow R(\mf{B}).
$$
Set $\vartheta_2:R(\mf B)\rightarrow R(\mf{B}^0)$
and $\vartheta=\vartheta_2\vartheta_1:R(\mf B_X)\rightarrow R(\mf{B}^0)$.

1) Let $S'_\mu\in R(\mf{B}_\lambda)$ for any $\mu\in\mathscr L'$ be given by $(S'_\mu)_X=k$ and $S'_\mu(x_1)=(\mu)$.
If $\mf B^0$ is homogeneous, then there exists a co-finite subset
$\mathscr L\subseteq\mathscr L'$, such that $\{\vartheta(S'_{\mu})\in
R(\mf{B}^0)\mid \mu\in\mathscr L\}$ is a family of homogeneous iso-classes.
By Corollary 3.2.4, there is an almost split conflation
$(e'_{\mu})$ in $R(\mf B_\lambda)$, such that $\vartheta(e'_\mu)$
is almost split in $R(\mf B^0)$. Fix any $\mu\in\mathscr L$, and consider the
conflation $(e_\mu)=\vartheta_1(e_\mu'):S_\mu\stackrel{\iota}{\longrightarrow}
E_\mu\stackrel{\pi}{\longrightarrow} S_\mu$ in $R(\mf B)$, where
$(E_\mu)_X=k^2, E_\mu(x)=\lambda I_2, E_\mu(a_{n_1})=J_2(\mu),E_\mu(a_{i})=(0),i\ne n_1$.
Then $(e_\mu)$ is almost split by Lemma 3.2.2 (ii).

2) Define a representation $L$ of $\mf{B}$ given by $L_X=k^2$,
$L(x)=J_2(\lambda)$, $L(a_{n_1})=\mu I_2$ and $L(a_i)=0$ for $i\ne n_1$. $L$ is well defined,
in fact if we make an unraveling
for $x\mapsto J_2(\lambda)$, then by Lemma 3.2.6 (ii), after a sequence of
regularizations ${{a_{i11}\,a_{i12}}\choose{a_{i21}\,a_{i22}}}$ (splitting from $a_i$)
$\mapsto\emptyset$; and ${{w_{i11}\,w_{i12}}\choose{w_{i21}\,w_{i22}}}$ (splitting from $w_i$)
$=0$ for  $i=1,\cdots,n_1-1$, an induced bocs with
$\dz^0\big({{a_{_{n_1,11}}\,a_{_{n_1,12}}}\choose{a_{_{n_1,21}}\,a_{_{n_1,22}}}}\big)=0$ is obtained.
Let $g: L \rightarrow S_{\mu}$ be a morphism with $g_{_X}={1\choose 0}$ and
$g(v_j)={0\choose 0}$ for all possible $j$. It is obvious that $g$ is not a retraction.

3) There exists a lifting $\tilde{g}: L \rightarrow E_{\lambda}$ with
$\tilde{g}\pi=g$:
$$\xymatrix{L\ar@(ul,ur)^{J_2(\lambda)}
\ar@(ul,dl)_{\mu I_2}\ar[rr]^{\tilde g}&&
{E_\lambda}\ar@(ul,ur)^{\lambda I_2}
\ar@(ur,dr)^{J_2(\mu)}}.$$
Since $\tilde{g}_{_X}\pi_{_X}=g_{_X}$, $\tilde{g}_{_X}={{1\,b}\choose {0\,d}}$. On the other
hand, $(x-y)^2\mid f_{11}(x,y)$ for $\bar w\ne 0$ and
$(J_2(\lz)-\lz I_2)^2=0$, hence $f_{11}(J_2(\lambda), \lz
I_2)=0$. $\tilde{g}:L\rightarrow E_{\lambda}$ being a morphism implies that
$$\begin{array}{c}L(a_1)\tilde{g}_{_X}-\tilde{g}_{_X} E_{\lambda}(a_1)=
f_{11}(L(x),E_{\lambda}(x))\tilde{g}(\bar w)=0\,\,\mbox{for}\,\,\bar w\ne0,
\,\,\,\mbox{or}\,\,0\,\,\mbox{for}\,\,\bar w=0\\[0.5mm]
\Longrightarrow \mu{{1\,\, b}\choose{0\,\,d}}-{{1\,\, b}\choose{0\,\,d}}J_2(\mu)=0,
\quad \mbox{i.e.} -{{0\, 1}\choose{0\,0}}=0.\end{array}$$
The contradiction shows that $\mf B^0$ is non-homogeneous. \hfill$\Box$

\medskip

{\bf Proposition 3.4.5}\, Let the bocs $\mf B^0$ have induced bocses $\mf{B}^s,\mf B^\varepsilon,\mf B^\tau$
given in Formula (3.2-3) and satisfying the condition (i)--(iii) stated between Corollary 3.2.4 and 3.2.5. Then $\mf B^0$ is non-homogeneous.

\smallskip

{\bf Proof}\, Suppose $\mf B=\mf B^\varepsilon$, where $\T=\{X,Y\},R_X=k[x,\phi(x)^{-1}], R_{Y}=k1_Y$, the layer
$L=(R;\omega;a_1,\cdots,a_n;
v_1,\cdots,v_{m})$ and $a_1:X\rightarrow Y,\dz(a_1)=0$:
$$\begin{array}{c}
\mf B:\quad \xymatrix {X\ar[rr]^{a_1}\ar@(ul, dl)[]_{x}
&& Y}\end{array}.
$$
Making a reduction given by $a_1\mapsto (1)$ of Proposition 2.2.7,
an induced local bocs $\mf B'=\mf B^\tau$ is obtained by a sequence of
regularizations with $\T'=\{Z\},R'=k[z,\phi'(z)^{-1}]$. Set the induced functors $\vartheta_1:R(\mf B')\rightarrow R(\mf B)$,
$\vartheta_2:R(\mf B)\rightarrow R(\mf B^0)$, and
$\vartheta=\vartheta_2\vartheta_1:R(\mf B')\rightarrow R(\mf B^0)$.

1) For any $\lambda\in k,\phi'(\lambda)\ne 0$,
there is an object $S'_\lambda\in R(\mf B')$
with $(S'_\lambda)_Z=k,S'_\lambda(z)=(\lambda)$.
If $\mf B^0$ is homogeneous, then there exists a co-finite subset
$\mathscr L\subseteq k\setminus\{\mu\mid\phi'(\mu)=0\}$, such that $\{\vartheta(S'_{\lambda})\in
R(\mf{B}^0)\mid \lambda\in\mathscr L\}$ is a family of homogeneous iso-classes.
Fix any $\lambda\in\mathscr L$, there is an almost split conflation
$(e_{\lambda}'):S'_{\lambda}\stackrel{\iota}{\longrightarrow}
E'_{\lambda}\stackrel{\pi}{\longrightarrow}
S'_{\lambda}$ in $R(\mf B')$ with $E'_\lambda(z)=J_2(\lambda)$,
such that $\vartheta(e_\lambda')$ is an almost split conflation in $R(\mf B^0)$ by Corollary 3.2.5.
Let $(e_\lambda)=\vartheta_1(e_\lambda'): S_\lambda\rightarrow
E_\lambda\rightarrow S_\lambda$ is an almost split conflation in $R(\mf B)$ by lemma 3.2.2 (ii) inductively, where
$$S_\lambda:\,
\xymatrix {k\ar[rr]^{(1)}\ar@(ul,dl)[]_{(\lambda)} && k},\quad E_\lambda:\,
\xymatrix {k^2\ar[rr]^{I_2}\ar@(ul,dl)[]_{J_2(\lambda)} && k^2}.$$

2)\, Define an object $L\in R(\mf B)$ with $L_X=k^2,L_Y=0$,
$L(x)=\lambda I_2$, $L(a_i)=(0),1\leqslant i\leqslant n$. Define a morphism
$g: S_\lambda\rightarrow L$ with $g_{_X}=(0\,1)$, $g_{_Y}=0$, and
$g(v)=0$ for any dotted arrow $v$. It is claimed that $g$ is not a
retraction. Otherwise, if there is a morphism
$h:L\rightarrow S_{\lz}$ with $gh=id_{S_{\lz}}$,
then $h_{_X}={{c}\choose{1}}$. Since $h_X S_\lambda(a_1)=L(a_1)h_Y$, there is ${{c}\choose{1}}(1)=0$, a
contradiction.
$$
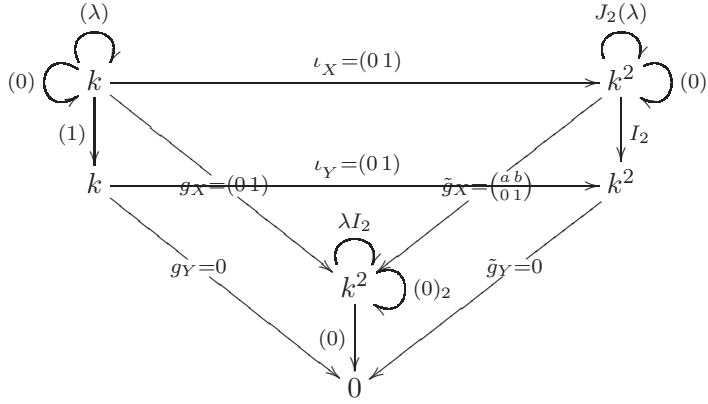
\xymatrix {k\ar@(ul,dl)[]_{(0)}\ar@(ul,ur)[]^{(\lambda)}
\ar[rrrrrr]^{\iota_{_X}={(0\,1)}}\ar[rrrdd]|(.5){g_X=
{(0\,1)}}\ar[d]_{(1)}
&&&&&&k^2\ar@(ul,ur)[]^{J_2(\lambda)}\ar@(ur,dr)[]^{(0)}
\ar[llldd]|(.5){\tilde{g}_X={{a\, b}\choose{0\,1}}}\ar[d]^{I_2}\\
k\ar[rrrrrr]^{\iota_{_Y}={(0\,1)}}\ar[rrrdd]|(.4){g_Y=0}
&&&&&& k^2\ar[llldd]|(.4){\tilde{g}_Y=0}
\\
&&& k^2\ar@(ul,ur)[]^{\lambda I_2}\ar@(ur,dr)[]^{(0)_2}\ar[d]_{(0)}&&&
\\&&&0&&&}$$

3) There exists an extension $\tilde
g: E_\lambda\rightarrow L$ with $\iota\tilde g=g$, which implies
$\tilde g_X={{a\,b}\choose{0\,1}}$. Since $\tilde g$ is a
morphism, there is $E(x)\tilde g_{_X}=\tilde g_{_X}
L(x)$, i.e. ${{0\, 1}\choose{0\, 0}}={{0\, 1}\choose{0\, 0}}{{a\, b}\choose {0\, 1}}
={{a\, b}\choose {0\, 1}}{{0\, 0}\choose{0\,
0}}$, a contradiction.
Therefore $\mf{B}^0$ is not homogeneous. \hfill$\Box$

\medskip

{\bf Remark 3.4.6 }\, For the sake of convenience, the following notation on MW5 is used.
Let $(\mf A,\mf B)$ be a pair with $R$ being trivial, and $({\mf A}^\nu,{\mf B}^\nu)$ be
an induced pair of $(\mf A,\mf B)$ given by a sequence of reductions in the sense of Lemma 2.3.2.
Suppose $(\mf A^\nu,\mf B^\nu)$
is local, and by calculating a series of
triangular formulae according to Subsection 3.3, an induced pair
$$\begin{array}{c}({\mf A}^\nu_{(\lambda_0,\lambda_1,\cdots,\lambda_{l-1})},
{\mf B}^\nu_{(\lambda_0,\lambda_1,\cdots,\lambda_{l-1})})\end{array}\eqno{(3.4\mbox{-}2)}$$
is obtained. Suppose in addition, the pair (3.4-2) is in the case of MW5, which satisfies Formula (3.3-3)--(3.3-4)
with the polynomials $\phi(x)$ and $\psi(x,x_1)$ given by Formulae (3.3-8) and (3.3-9).

Denote the pair (3.4-2) by $(\mf A^{s+1},\mf B^{s+1})$.
It means that $({\mf A}^{s},{\mf B}^{s})$, with the first arrow
$a^s_1$ and $\dz(a^s_1)=0$, is obtained by a sequence of reductions in the sense of Lemma 2.3.2 starting from
$({\mf A},{\mf B})$. After making a loop
mutation $a^s_1\mapsto (x)$, there is an induced pair $({\mf A}^{s+1},{\mf B}^{s+1})$ with $R^{s+1}=k[x]$.
Removing, by regularizations, the pairs of arrows $(a_1,w_1);\cdots;(a_{n_1-1},w_{n_1-1})$ from $\mf B^{s+1}$ in Formula
(3.3-3), an induced pair say $(\mf A^t,\mf B^t)$ under the assumption $x\mapsto(\lambda)$
is obtained. Denote the solid arrows of $\mf B^t$ by $a^t_{j+1},
j=0,\cdots,n-n_1$, then $a^t_1$ obtained from $a_{n_1}$, is the first arrow of $\mf B^t$. The bocs
$\mf B^t$ is also said to be {\it in the case of} MW5.

\bigskip
\bigskip
\bigskip

\centerline {\Large\bf 4. One-sided pairs}

\bigskip

In this section, some special quotient problems of matrix bimodule
problems are defined. The purpose is to deal with the most difficult part of the proof of the main theorem
in the subsections 5.4--5.5.

\bigskip
\noindent{\large\bf 4.1 Definition of one-sided pairs}
\bigskip

We give the definition of one-sided pairs; then consider the
induced pairs after some reductions via pseudo formal equations
in this subsection.

\medskip

Let $\mf A=(R,\K,\M,H=0)$ be a matrix bimodule problem
with $\T$ being trivial, let $\mf C,\mf B$ be the associated bi-comodule problem and
bocs respectively. Suppose a sequence of pairs
$$(\mf A,\mf B),(\mf A^1,\mf B^1),\cdots,(\mf A^r,\mf B^r)$$
is given by reductions in the sense of Lemma 2.3.2.
Assume that the leading position of
the first base matrix $A^r_1$ of $\M^r_1$ is $(p^r,q^r)$ contained
in the $(\textsf{p},\textsf{q})$-th leading block of of $A_1$ of $\M$ partitioned under $\T$.
It is further assumed that $d_1, \cdots, d_m$ are the first $m$ solid arrows of $\mf B^r$,
which locate at the $p^r$-th row of the formal product $\Theta^r$, such that $d_m$ is
sitting at the last column of the $(\textsf{p},\textsf{q})$-block, see the picture below:
\vspace{0mm}
\begin{center}
\setlength{\unitlength}{1mm}
\begin{picture}(40,25)
\put(0,0){\line(1,0){40}}  \put(0,0){\line(0,1){25}}
\put(0,25){\line(1,0){40}} \put(40,0){\line(0,1){25}}
\put(15,10){\line(1,0){25}} \put(15,10){\line(0,1){5}}
\put(15,15){\line(1,0){25}}
\put(20,10){\line(0,1){5}} \put(35,10){\line(0,1){5}}
\put(25,11){$\cdots$}\put(16, 11){$d_1$}
\put(35.5, 11){$d_m$}\put(86,12){$(4.1\mbox{-}1)$}
\end{picture}
\end{center}

Recalling the statement between Formula (1.2-3) and (1.2-4),
the quotient problem ${(\mf A^r)}^{[m]}=(R^r, \K^r,(\M^r)^{[m]},H^r)$ of $\mf A^r$;
the sub-bi-comodule problem $(\mf C^r)^{(m)}=(R^r, \C^r,(\N^r)^{(m)}$, $\partial\mid_{(\N^r)^{(m)}})$ of $\mf C^r$
with the quasi-basis $d_1,\cdots,d_m$ of $(\N^r)^{(m)}$, and the sub-bocs
$(\mf{B}^r)^{(m)}$ of $\mf B^r$, a {\it quotient-sub-pair $((\mf
A^r)^{[m]}, (\mf B^r)^{(m)})$} is obtained, and it is denoted by $(\bar{\mf A},\bar{\mf B})$ for simplicity.

Denote a set of integers by $\bar T=\bar T_R\times \bar T_C\subseteq T^r$, where
$\bar T_R=\{0\}$ and $\bar T_C=\{1,,2\cdots,m\}$ are the row and column
indices of $(d_1,d_2,\cdots,d_m)$ respectively.
A representation of size vector $\n=(n_0;n_1,\cdots,n_m)$ over $\bar{\mf A}$ can be written as
$\bar M=(\bar M_1,\cdots,\bar M_m)$, where $\bar M_i$ is an $n_0\times n_i$-matrix over $k$.
But a morphism between two representations must be discussed returning back to the category $R(\mf A^r)$.
Moreover, if $\bar{\mf A}'$ is any induced pair of $\bar{\mf A}$, a pseudo functor
$\bar\vartheta$ can be considered acting on the
objects over $\bar{\mf A}'$. Recall the formal equation of the pair $(\mf A^r,\mf B^r)$:
$$\begin{array}{c}(\sum_{Y\in \T^r}e_{_{Y}}^r\ast E^r_{Y})(\sum_{i}a^r_i\ast A^r_i)=
\big(H^r(\sum_{j}v^r_j\ast V^r_j)- (\sum_{j}v^r_j\ast V^r_j)H^r\big)\\[1.5mm]
+(\sum_{i}a^r_i\ast A_i^r)(\sum_{Y\in \T^r}e_{_{Y}}^r\ast E^r_{Y}+\sum_{j}v^r_j\ast V^r_j)-(\sum_{j}v^r_j\ast V^r_j)
(\sum_{i}a^r_i\ast A_i^r).\end{array}$$
Let $d_i: X\rightarrow Y_i$ (possibly $Y_i=X$, or $Y_i=Y_j$ for $i\ne j$). Then the $(p^r,q^r),\cdots,(p^r,q^r+m-1)$-th
equations of the formal equation of $(\mf A^r,\mf B^r)$ can be
rewritten as:
$$e_{_{X}}(d_1,d_2,\cdots,d_m)= (w_1,w_2,\cdots,w_m)+(d_1,d_2,\cdots,d_m)
\left(\begin{array}{cccc}e_{_{Y_1}}&w_{12}&\cdots & w_{1m}\\
&e_{_{Y_2}}&\cdots& w_{2m}\\&&\ddots &\vdots\\&&&e_{_{Y_m}}
\end{array}\right). \eqno (4.1\mbox{-}2)$$

{\bf Remark 4.1.1}\, We give some explanation on the notations of Formula (4.1-2).

(i)\, $e_{_{X}}$ is the $(p^r, p^r)$-th entry of the formal
product $\sum_{Y\in {\T^r}}e^r_{_{Y}}\ast E^r_{Y}$; and $e_{_{Y_\xi}}$ the $(q^r+\xi-1,q^r+\xi-1)$-th entry
of that for $\xi=1,\cdots,m$.

(ii)\, For $\xi=1,\cdots,m$, $w_\xi$ is the $(p^r,q^r+\xi-1)$-th entry of
$H^r(\sum_{j}v^r_j\ast V^r_j)-(\sum_{j}v^r_j\ast V^r_j)H^r$. In fact,
$w_\xi=\sum_j\alpha_{\xi}^jv^r_j$, where $s(v^r_j)\ni p^r,t(v^r_j)\ni q^r+\xi-1,\alpha_{\xi}^j\in k$.

(iii)\,  For $1\leqslant\eta<\xi\leqslant m$,
$w_{\eta\xi}$ is the $(q^r+\eta-1,q^r+\xi-1)$-th entry of $\sum_{j}v^r_j\ast V^r_j$. In fact,
$w_{\eta\xi}=\sum_j\beta_{\eta\xi}^jv^r_j$ where $s(v^r_j)\ni q^r+\eta-1,t(v^r_j)\ni q^r+\xi-1,\beta_{\eta\xi}^j\in k$.

\smallskip

The differential of $d_1,\cdots,d_m$ can be read off from Formula (4.1-2).
Note that in each monomial of the differential, there exists at most one solid arrow multiplying
the the dotted one from the left:
$$\begin{array}{c}
-\delta(d_i)=w_{i} +\sum_{j<i} d_jw_{ij},\quad 1\leqslant i\leqslant
m.\end{array}\eqno{(4.1\mbox{-}3)}
$$

{\bf Definition 4.1.2}\, A bocs $\mf{B}$  with a layer $L=(R;\omega;
a_1,\ldots,a_m,b_1,\cdots,b_n;\underline u_j,\underline v_j,\bar u_j,\bar v_j)$
is said to be {\it one-sided}, provided that $R$ is trivial and $\T=\{X,Y_1,\cdots,Y_h\}$;
the solid arrows $a:X\rightarrow Y,b:X\rightarrow X$,
and the dotted arrows are divided into four classes: $\bar u:X\rightarrow X$, $\underline u:Y\rightarrow Y'$,
$\bar v:Y\rightarrow X,\underline v:X\rightarrow Y$; the differentials of the solid arrows are given by
$$\begin{array}{c}\dz(a_i)=\sum_j\alpha_{ij}\underline v_j+\sum_{i'<i,j}\beta_{ii'j}a_{i'}\underline u_j
+\sum_{b_{i'}\prec a_i,j}\gamma_{ii'j}b_{i'}\underline v_j,\\
\dz(b_i)=\sum_j\lambda_{ij}\bar u_j+\sum_{a_{i'}\prec b_i,j}\mu_{ii'j}a_{i'}\bar v_j
+\sum_{i'<i,j}\nu_{ii'j}b_{i'}\bar u_j,\end{array}$$
with the coefficients $\alpha_{ij},\beta_{ii'j},\gamma_{ii'j},\lambda_{ij},\mu_{ii'j},\nu_{ii'j}\in k$.

The associated bocs $\bar{\mf B}$ of $\bar{\mf A}$ is one-sided by
Formula (4.1-3), and $(\bar{\mf A},\bar{\mf B})$ is called a {\it one-sided pair}.

\vskip -6mm{\unitlength=1mm
$$\begin{array}{c}
\begin{picture}(80,45)
\put(29,39){\circle{4}} \put(30.8,39){\vector(0,1){0.5}}
\put(32,39){\circle*{1}} \qbezier[10](32,39)(37,42)(38,39)
\qbezier[10](32,39)(37,36)(38,39)

\put(32,38){\vector(-1,-3){6.5}} \put(31,37.5){\vector(-1,-1){18}}
\put(33,37.5){\vector(1,-1){18}} \put(25,17){\circle*{1}}
\put(12,17){\circle*{1}}  \put(52,17){\circle*{1}}
\put(10,12){$Y_1$}  \put(50,12){$Y_h$} \put(25,12){$Y_2$}
\put(30,42){$X$}

\qbezier[10](26,17)(31,20)(32,17) \qbezier[10](26,17)(31,14)(32,17)

\mput(13,17)(2,0){6}{\line(1,0){1}}
\mput(28,36)(-1,-1){14}{\circle*{0.5}}
\put(15,23){\vector(-1,-1){2}}

\mput(36,36)(1,-1){16}{\circle*{0.5}} \put(37,35){\vector(-1,1){2}}

\put(39,39){\makebox{$\bar u$}} \put(18,29){\makebox{$\underline v$}}
\put(44,29){\makebox{$\bar v$}} \put(19,19){\makebox{$\underline u$}}
\put(24,39){\makebox{$b$}} \put(30,27){\makebox{$a$}}
\put(38,15.5){\makebox{$\cdots$}}
\end{picture}
\end{array}\eqno{(4.1\mbox{-}4)}
$$}\vskip -12mm

Let $(\mf A,\mf B)$ be any pair,
$(\mf A^{[h]},\mf B^{(h)}), h\geqslant 1$, be a quotient-sub-pair.
Note that if the reduction on $(\mf A,\mf B)$ is made with respect to
an admissible $R'$-$\bar R$-bimodule $L$ by Proposition 2.2.1--2.2.7,
then $L$ is also the admissible bimodule of $(\mf A^{[h]},\mf B^{(h)})$.
Thus there are two sequences of reductions as below, such that $(\bar{\mf A}^{i+1},\bar{\mf B}^{i+1})$
and $({\mf A}^{r+i+1},{\mf B}^{r+i+1})$ are obtained respectively
from $(\bar{\mf A}^i,\bar{\mf B}^i)$ and $({\mf A}^{r+i},{\mf B}^{r+i})$
by the same admissible bimodule, or by the same
regularization 2.1.8 for $0\leqslant i<s$:
$$\begin{array}{ccccccc}&(\bar{\mf A},\bar{\mf B}),&
(\bar{\mf A}^1,\bar{\mf B}^1),&\cdots,&(\bar{\mf A}^i,\bar{\mf B}^i)
&\cdots,&(\bar{\mf A}^s,\bar{\mf B}^s);\\[1mm]
&(\mf A^r,{\mf B}^r), &({\mf A}^{r+1},{\mf B}^{r+1}),&\cdots,&({\mf A}^{r+i},
{\mf B}^{r+i}),&\cdots,&({\mf A}^{r+s},
{\mf B}^{r+s}).\end{array} \eqno{(4.1\mbox{-}5)}$$
Let $\bar{\mf A}^i=(\mf A^{r+i})^{[m^i]}=(R^{r+i},\K^{r+i},\bar\M^i,F^i)$ be a
quotient problem with $m^i$ being the number of the base matrices of $\bar\M^i$, and
the associated sub-bi-comodule problem
$\bar{\mf C}^i=(\mf C^{r+i})^{(m^i)}=(R^{r+i},\C^{r+i},\bar\N^i,\bar\partial^i)$.
A formula of the pair $(\bar{\mf A}^i,\bar{\mf B}^i)$ is written as follows for $i=0,1,\cdots,s$:
$${e}^i_{_X}(F^i+\bar\Theta^i)=(W_1,\cdots, W_m)+(F^i+\bar\Theta^i)
\left( \begin{array}{cccc}e_{_{Y_1}}^i&W_{12}&\cdots & W_{1m}\\
& e_{_{Y_2}}^i&\cdots&W_{2m}\\&&\ddots &\vdots\\&&& e_{_{Y_m}}^i
\end{array}\right)\eqno{(4.1\mbox{-}6)}$$
where the diagonal parts of $e^i_{_X}$ and the most right matrix in Formula (4.1-6)
can be viewed as a reduced formal product $\bar\Upsilon^i$ of $(\K_0,\C_0)$.
$\bar\Theta^i$ is the formal product of $(\bar\M^i_1,\bar\N^i_1)$ containing
the solid arrows splitting from $d_1,\cdots,d_m$; $F^i+\bar\Theta^i$ is a $(1\times m)$-partitioned
matrix under $\bar \T$ with a size vector $\n^i=(n^i_0;n^i_1,\cdots,n^i_m)$,
and $F^i$ is sitting in the blank part in the picture:
\unitlength=1mm
$$F^i+\bar\Theta^i=\begin{picture}(85,10)
\put(0,-10){\line(1,0){20}} \put(0,-10){\line(0,1){20}}
\put(0,10){\line(1,0){20}} \put(20,-10){\line(1,0){65}}
\put(15,-10){\line(0,1){20}} \put(20,10){\line(1,0){65}}
\put(85,-10){\line(0,1){20}} \put(24,-10){\line(0,1){20}}
\put(25,7){$d_q^i\;\cdots$} \put(24,6){\line(1,0){61}}
\put(39,-10){\line(0,1){20}} \put(40,3){$d_p^i\; d_{p+1}^i \;
\cdots\;  $} \put(42,7){$\cdots\quad \cdots$} \put(50,-1){$d_1^i\;
d_2^i\cdots$} \put(39,2){\line(1,0){46}} \put(49,-2){\line(0,1){4}}
\put(49,-2){\line(1,0){36}} \put(63,-10){\line(0,1){20}}
\put(65,7){$\cdots$} \put(65,-1){$\cdots$} \put(65,3){$\cdots $}
\put(71,-10){\line(0,1){20}}\put(72,7){$\cdots\; \cdots$}
\put(72,3){$\cdots\; \cdots$} \put(72,-1){$\cdots\;d_{p-1}^i$}
\put(77,-7){$ $} \put(50,-7){$ $} 
\end{picture}
$$
\vskip 9mm

Reductions are performed according to Theorem 2.4.4.
More precisely, the system $\bar{\IF}^{ri}$ of Formula (2.4-5) for the pair
$(\mf A^{r+i},\mf B^{r+i})$ can be written as $\bar{\mathcal F}^i$ below,
which is said to be the {\it reduced defining system of the pair $(\bar{\mf A}^i,\bar{\mf B}^i)$}.
$$\begin{array}{c}{\bar{\mathcal F}}^i: \bar\Psi_{\n^i}^{l}F^i
\equiv_{\prec(\bar p^i,\bar q^i)} \bar\Psi_{\n^i}^{m}+F^i\bar\Psi_{\n^i}^{r},
\end{array} \eqno{(4.1\mbox{-}7)}$$
where the upper indices $l,m,r$ are used to show
the left, middle and right parts of the variable matrix $\bar\Psi_{\n^i}$;
$\bar\Psi_{\n^i}^{l}$ is the strict upper triangular part of $e^i_{_X}$;
$\bar\Psi_{\n^i}^{r}$ is that of the most right matrix in Formula (4.1-6); and
$\bar\Psi_{\n_i}^{m}=(W_1,\cdots,W_m)$. Namely, $\bar{\mathcal F}^i$ is obtained from
Formula (4.1-6) by removing the term $\bar\Theta^i$ and the diagonal parts
of the most left and right matrices. The strict upper triangular parts of
$e^i_{_X}$ and $e^i_{_{Y_j}}$ for $1\leqslant j\leqslant m$
are constructed inductively in Proof 2) of Theorem 2.4.4; while
$W_h$ and $W_{hl}$ are the splitting of $w_h,w_{hl}$.

Since it is difficult to calculate the dotted arrows after a reduction,
the linear relation of the dotted elements of $\bar\Psi$
appearing in the reduction will be described instead.

On the other hand,
$\bar\Pi^i=\bar\Psi_{\n^i}$ may be said to be a {\it pseudo formal product}
of $(\K_1^i,\C_1^i)$; and Formula (4.1-6) a {\it pseudo formal equation of the pair $(\bar{\mf A}^i,\bar{\mf B}^i)$},
since the entries of $\bar\Pi^i$ are dotted elements with some linear relations.
The reason that we borrowed the concept of ``formal product" was that
it is possible to read off the differential of the solid arrows from Formula (4.1-6) according to Theorem 1.4.2.
For example for the first arrow $d^i_{l\bar p\bar q}$:
$$\begin{array}{c}-\dz(d^i_{l\bar p\bar q})=w^i_{lpq}+\sum_{j<l,q}d^{i,0}_{j\bar pq}w^i_{jl,q\bar q}+
\sum_{q<\bar q}d^{i,0}_{l\bar pq}w^i_{Yq\bar q}-\sum_{q>\bar p}w^i_{X\bar pq}d^{i,0}_{lq\bar q},
\end{array}\eqno{(4.1\mbox{-}8)}$$
where $d^i_{l\bar p\bar q}$ is split from $d_l:X\rightarrow Y$, $d^{i,0}_{jpq}$ is the $(p,q)$-th
entry obtained from $d_j$ in $F^i$.


\medskip

{\bf Remark 4.1.3}\, Let $(\bar{\mf A}',\bar{\mf B}')$ be any induced pair of $(\bar{\mf A},\bar{\mf B})$
after several reductions in the sense of Lemma 2.3.2. And $(\bar{\mf A}'',\bar{\mf B}'')$
is an induced pair of $(\bar{\mf A}',\bar{\mf B}')$ given by one of three reductions of 2.3.2.

(i) If there is a linear relation of dotted elements $\sum_ju_j=0$ in ${\bar\Psi}_{\n'}$,
then $\sum_j\bar u_j=0$ in ${\bar\Psi}_{\n''}$ with $\bar u_j$ being the split of $u_j$.

(ii) Suppose $a_1'$ is the first arrow of $\mf B'$, and $\dz(a_1')=v+\sum_j\alpha_ju_j$, where
$v, u_j$ are dotted elements of ${\bar\Psi}_{\n'}$. If $v$ is a dotted arrow, and $v\notin\{u_j\}$, then $\dz(a_1')\ne 0$.

(iii) Set $a_1'\mapsto\emptyset$ in (ii) by a regularization, it is said that $v$ is {\it replaced by $-\sum_ju_j$
in $\bar\Psi_{\n''}$}.

(iv) If we are able to determine that a dotted element $v$ of $\bar\Psi_{\n''}$ is
linearly independent of all the others, then $v$ is said to be a {\it dotted arrow preserved in
$\bar{\mf B''}$}.

\bigskip
\bigskip
\noindent{\large\bf 4.2\, Differentials in one-sided pairs}
\bigskip

In this subsection, the classification of local one-sided bocses is discussed;
and the differentials of the solid arrows of non-local bocses are calculated.

\medskip

For the sake of simplicity,
the black letter $\boldsymbol\delta$ is used instead of
the symbol ``$-\dz$", see Formulae (4.1-3) and (4.1-8), up to the end of the whole section 4.
First of all, Classification 3.3.5 gives the following Classification on local one-sided bocses.
The forms in one-sided case are much simpler than those in general case.

\medskip

{\bf Classification 4.2.1}\, Let $\bar{\mf B}$ be a local one-sided bocs with a layer (letter $a$ is changed to the letter $b$):
$$\begin{array}{c}L=(R;\omega;b_1,\cdots,b_n;v_1,\cdots,v_m).\end{array}$$

(i)\, $\bar{\mf B}$ with $R=k1_X$ has differentials by Formula (3.3-5) after some base changes:
$$\begin{array}{c}{\boldsymbol\delta}^0(b_1)= \bar u_1,{\boldsymbol\delta}^0(b_2)= \bar u_2,
\cdots,{\boldsymbol\delta}^0(b_{n})= \bar u_{n}.\end{array}$$

(ii)\, $\bar{\mf B}$ with $R=k1_X$ has some integer $1\leqslant n_0<n$, such that the differentials are
$$\begin{array}{c}{\boldsymbol\delta}^0(b_i)= \bar u_i,\,\, 1\leqslant i<n_0;\,\,\,{\boldsymbol\dz}^0(b_{n_0})=\sum_{i=1}^{n_0-1}f_{n_0i}\bar u_i,\,f_{n_0i}\in k;\\
{\boldsymbol\dz}^1(b_{i})=\sum_{j\ne n_0,j=1}^{i-1}f_{ij}(b_{n_0})\bar u_j+f_{ii}\bar u_i,\quad\mbox{with respect to}\,\,b_{n_0},\\
\mbox{where},\,\,n_0<i\leqslant n,\,\,f_{ij}(b_{n_0})=\beta_{i,j}^0+\beta_{i,j}^1b_{n_0}\in k[b_{n_0}]\,\,\mbox{for}\,\, i<j, \,\,
f_{ii}\in k^\ast.\end{array}
\eqno{(4.2\mbox{-}1)}$$

$\bar{\mf B}$ with $R=k[x]$ is given by a sequence of regularizations, followed by $a_{n_0}\mapsto x$
in Formula (3.3-6). For the sake of simplicity,
Formulae (3.3-2) and (3.3-3) are written in a unified form:
$$\left\{\begin{array}{ll}{\boldsymbol\dz}^0(b_{1})&=f_{11}(x)\bar u_{1},\\
{\boldsymbol\dz}^0(b_{2})&=f_{21}(x)\bar u_{1}+f_{22}(x)\bar u_{2},\\
&\cdots\quad \cdots\\
{\boldsymbol\dz}^0(b_{t})&=f_{t1}(x)\bar u_{1}
+f_{t2}(x)\bar u_{2}+\cdots+f_{t,t-1}(x)\bar u_{t-1}+f_{tt}(x)\bar u_{t},
\end{array}\right.\eqno{(4.2\mbox{-}2)}$$
where $f_{ij}(x)=\alpha_{i,j}^0+\alpha_{i,j}^1x\in k[x]$ for
$1\leqslant i\leqslant j\leqslant t$, and $f_{ii}(x)\ne 0$ for $1\leqslant i<t$.

(iii)\, $\bar{\mf B}$ has an induced bocs $\bar{\mf B}_{(\lambda_0,\cdots,\lambda_{l})}$ with
Formula (4.2-2), where $t=n, f_{nn}(x)\ne 0$; $\phi(x)=1$ in Formula (3.3-7); and
there is some minimal $1\leqslant s\leqslant n$, such that
$f_{ss}(x)\in k[x]\setminus k$ is non-invertible in $k[x]$. It is in the case of MW3.

(iv)\, $\bar{\mf B}$ has an induced bocs $\bar{\mf B}_{(\lambda_0,\cdots,\lambda_{l-1})}$ with
Formula (4.2-2), where $t=n_1<n$, $f_{n_1,n_1}(x)=0$ and $\phi(x)=\prod_{i=1}^{n_1-1}f_{ii}(x)$
in Formula (3.3-8). Denoting $b_{n_1}$ by $x_1$, Formula (3.3-4) shows:
$$\left\{\begin{array}{l}{\boldsymbol\dz}^1(b_{n_1+1})=K_{n_1+1}+f_{n_1+1,n_1+1}(x)\bar u_{n_1+1},\\
{\boldsymbol\dz}^1(b_{n_1+2})=K_{n_1+2}+f_{n_1+1,n_1+1}(x,x_1)\bar u_{n_1+1}+f_{n_1+2,n_1+2}(x)\bar u_{n_1+2},\\
\qquad \cdots\qquad \cdots\\
{\boldsymbol\dz}^1(b_n)=K_{n}+h_{n,n_1+1}(x,x_1)\bar u_{n_1+1}
+\cdots+ f_{n,n-1}(x,x_1)\bar u_{n-1}+f_{nn}(x)\bar u_n,\end{array}\right.\eqno{(4.2\mbox{-}3)}$$
where $f_{ij}(x,x_1)\in k[x,x_1]$ for $i<j$, $f_{ii}(x)\ne 0$; and the polynomial
$\psi(x)=\phi(x)\prod_{i=n_1+1}^nc_i(x)f_{ii}(x)$ is given by Formula (3.3-9).
It is in the case of MW4.

Note in particular, that MW5 never occurs in one-sided case.

\smallskip

{\bf Proof} The one-sided bocs of (i) is finite type, and that of (ii) is tame infinite.

(iii) There is no localization needed in Formula (4.2-2), thus $c_i(x)=1$. It is clear that
$h_{ii}(x, x)=1$ in Formula (3.3-1), therefore $\phi(x)=1$. Finally,
any non-zero and non-invertible polynomial of $k[x]$ belongs to $k[x]\setminus k$.

(iv) Suppose the last terms of the formulae in (4.2-3) are $0\ne f_{ii}'(x,x_1)=f_{ii}(x)h_{ii}(x,x_1)$ with
$h_{ii}(x,x_1)\in k[x,x_1,\psi(x)^{-1}]\setminus k[x]$ or $h_{ii}(x,x_1)=1$. If there exists
a minimal integer $n_1<s\leqslant n$, such that $k_{ss}(x,x_1)\notin k[x]$,
then for any $\lambda\in k$ with $\psi(\lambda)\ne 0$, the induced bocs $\mf B_{\lambda}$ given by
$x\mapsto(\lambda)$ returns to case (iii). Therefore $f_{ii}'(x,x_1)=f_{ii}(x)\in k[x]$
for $n_1<i\leqslant n$. \hfill$\Box$

\medskip

From now on, a general one-sided pair $(\bar{\mf A},\bar{\mf B})$ with $|\bar\T|>1$ is dealt with. If $\bar{\mf B}_X$,
the induced local bocs of $\bar{\mf B}$, is in case (iii) or (iv) of Classification 4.2.1,
then $\bar{\mf B}_X$ is wild and non-homogeneous, so is $\bar{\mf B}$. Since
$\bar{\mf B}_X$ in case (i) of 4.2.1 is relatively simple, the discussion below
is concentrated on $\bar{\mf B}_X$ given by Formula (4.2-1) in case (ii) of 4.2.1. Denote the solid edges of $\bar{\mf B}$
before $b_{n_0}$ by $a_1,\cdots a_h$. The differential ${\boldsymbol\dz}^0$ acting on $a$'s has two possible expressions. First,
$$\begin{array}{c}
{\boldsymbol\delta}^0(a_1)=\underline v_1,
\,\,\cdots\cdots,\,\,
{\boldsymbol\delta}^0(a_h)=\underline v_{h}.
\end{array}\eqno(4.2\mbox{-}4)$$
Second, there exists some $1\leqslant h_1<h$, such that ${\boldsymbol\delta}^0(a_1)=\underline v_1,
\cdots,{\boldsymbol\dz}^0(a_{h_1-1})=\underline v_{h_1-1}$, but
${\boldsymbol\delta}^0(a_{h_1})=\sum_{j=1}^{h_1-1}\bar\alpha_{h_1,j}\underline v_{j}$. Inductively,
there are two subsets $\{h_1,\cdots,h_s\}\subseteq\{1,\cdots,h\}$,
and $\Lambda=\{1,\cdots,h\}\setminus\{h_1,\cdots,h_s\}$, such that
$$\left\{\begin{array}{llc}{\boldsymbol\dz}^0(a_i)&=\quad\quad\underline v_i,&
i\in\Lambda;\\[1mm]
{\boldsymbol\dz}^0(a_{h_l})&=\sum_{j\in\Lambda,j<h_l}\bar\alpha_{h_l,j}\underline v_j,&l=1,\cdots,s.
\end{array}\right. \eqno(4.2\mbox{-}5)$$

\medskip

{\bf Convention 4.2.2}\, Suppose $\bar{\mf B}$ is
a one-sided bocs with $\bar{\mf B}_X$ given by Formula (4.2-1).
All the loops $b_1,\cdots,b_n$ at $X$ are called {\it $b$-class arrows},
where the loop $\bar b=b_{n_0}$ is said to be {\it effective} or $\bar b$-class, the others
are {\it non-effective}.
The edges $a_1,\cdots,a_h$ before $\bar b$ are called {\it $a$-class arrows},
where $\{\bar a_i=a_{h_i}\mid 1\leqslant i\leqslant s\}$
are said to be {\it effective} or $\bar a$-class,
the others are {\it non-effective}.
Let $c_1,c_2,\ldots,c_t$ be the solid edges after $\bar b$, which are called
$c$-class arrows, and they are effective.

A solid arrow splitting from one of the classes $a,\bar a,b,\bar b,c$,
or a dotted element splitting from a dotted arrow of $\underline u,
\underline v,\bar u,\bar v$-classes in Picture (4.1-4) are said to be in the same class.

\medskip

A special case of the differential ${\boldsymbol\dz}^1$ with respect to $\bar b$ on $c$-class arrows is given by
$$\left\{\begin{array}{l}{\boldsymbol\dz}^1(c_1)=\sum_{j\in\Lambda}\gamma_{1j}(\bar b)
\underline v_j+\gamma_{1,h+1}(\bar b)\underline v_{h+1},\\
\quad\cdots\quad\cdots\quad\cdots\\
{\boldsymbol\dz}^1(c_t)=\sum_{j\in\Lambda}\gamma_{tj}(\bar b)
\underline v_j+\gamma_{t,h+1}(\bar b)\underline v_{h+1}+\cdots
+\gamma_{t,h+t}(\bar b)\underline v_{h+t},\\
\end{array}\right.\eqno{(4.2\mbox{-}6)}$$
where $\{\underline v_j\}_{j\in \Lambda}\cup\{\underline v_{h+j}\}_{1\leqslant j\leqslant t}$
are dotted arrows, $\gamma_{i,h+i}(\bar b)\ne 0$ for $1\leqslant i\leqslant t$.

\medskip

{\bf Lemma 4.2.3}\, Let $\bar{\mf B}$ be a one-sided bocs with $\bar{\mf B}_X$ being
in case  (ii) of Classification 4.2.1. If Formula (4.2-6) fails,
i.e. there exists some minimal $1\leqslant l\leqslant t$
with $\gamma_{l,h+l}(\bar b)=0$, then $\bar{\mf B}$ is non-homogeneous.

\smallskip

{\bf Proof}\, It is proceed with a sequence of regularizations:
$a_j\mapsto\emptyset,\underline v_j=0$ for
$j\in\Lambda$, $b_j\mapsto\emptyset,\bar u_j=0$ for $1\leqslant j<n_0$,
and edge reductions $\bar a_{h_i}\mapsto (0)$ for $i=1,\cdots,s$. Then after a loop
mutation $\bar b\mapsto (x)$, and defining a
polynomial $\phi(x)=\prod_{i=1}^{l-1}\gamma_{i,h+i}(x)$ by Formula (4.2-6),
an induced pair $(\mf A',\mf B')$ is obtained, such that $R'_X=k[x,\phi(x)^{-1}]$ and $\mf B'_X$ is minimal. Without loss of generality,
suppose $\T'=\{X,Y\}$ with $R_Y=k$. Making regularizations $c_i\mapsto\emptyset,\underline v'_{h+i}=0$ for $1\leqslant i<l$,
an induced bocs $\xymatrix{\cdot\ar@(ul,dl)[]|{x}\ar[r]^{c_{l}}&\cdot}$ of
two vertices with ${\boldsymbol\dz}(c_{l})=0$ follows, which is in the case of MW1.
Thus $\bar{\mf B}'$, consequently $\bar{\mf B}$, are wild and non-homogeneous.\hfill$\Box$

\medskip

Let $\bar{\boldsymbol\dz}$ be obtained from the differential ${\boldsymbol\dz}$ by removing all
the monomial involving any non-effective $a,b$-class solid arrows.
Now $\bar{\boldsymbol\dz}$ acting on all $a,b,c$-class arrows is written in the following
three formulae:
$$\begin{array}{c}\bar{\boldsymbol{\boldsymbol\delta}}(a_i)=\underline v_i
+\sum_{h_l<i}\bar a_{l}(\sum_{j}\epsilon_{ilj}\underline u_{j}),
\quad i\in\Lambda;\\[1mm]
\bar{\boldsymbol\dz}(\bar a_\tau)=\sum_{j<h_\tau}\bar\alpha_{\tau j}\underline v_j
+\sum_{l<\tau}\bar a_{l}(\sum_{j}\bar\epsilon_{\tau lj}\underline u_{j}),
\quad 1\leqslant \tau\leqslant s.\end{array}\eqno{(4.2\mbox{-}7)}$$
$$\begin{array}{c}\bar{\boldsymbol\dz}(b_i)=\bar u_i
+\sum_{\bar a_{l}\prec b_i}\bar a_{l}(\sum_{j}\varepsilon_{ilj}\bar v_{j}),\quad
 i<n_0;\\[1mm]
\bar{\boldsymbol\dz}(\bar b)=\sum_{j=1}^{n_0-1}\bar\beta_{j} \bar u_{j}
+\sum_{l=1}^s\bar a_{l}(\sum_{j}\bar\varepsilon_{lj}\bar v_{j}),\quad
i=n_0;\\[1mm]
\bar{\boldsymbol\dz}(b_i)=\bar u_i+\sum_{j=1,j\ne n_0}^{i-1}\beta_{ij}(\bar b)\bar
u_{j}+\sum_{l=1}^s\bar a_l(\sum_{j}\varepsilon_{ilj}\bar v_j)
+\sum_{c_l\prec b_i}c_l(\sum_j\varepsilon'_{ilj}\bar
v_{j}),\,\, i>n_0.\end{array}\eqno{(4.2\mbox{-}8)}$$
$$\begin{array}{c}\bar{\boldsymbol\dz}(c_\tau)=\sum_{l=1}^s\bar a_{l}
(\sum_{j}\zeta_{\tau lj}\underline u_{j})+
\sum_{j\leqslant h+\tau}\gamma_{\tau j}(\bar b)\underline v_j+
\sum_{l=1}^{i-1}c_l(\sum_j\xi_{\tau lj}
\underline u_{j}),\,\, 1\leqslant \tau\leqslant t.
\end{array}\eqno{(4.2\mbox{-}9)}$$
where all the coefficients $\epsilon,\bar\epsilon,\bar\alpha,\varepsilon,
\varepsilon',\bar\varepsilon,\bar\beta,\zeta,\xi\in k$, and
$\beta(\bar b)=\beta^0+\beta^1\bar b,\gamma(\bar b)=\gamma^0+\gamma^1\bar b\in k[\bar b]$.

\bigskip
\bigskip
\noindent {\large\bf 4.3. Reduction sequences of one-sided pairs}
\bigskip

The purpose of this subsection is tow folds: 1) present a condition
on formal products, which can be preserved
after some edge reductions;
2) construct a reduction sequence based on 1)
starting from a non-local one-sided pair $(\bar{\mf A},\bar{\mf B})$
with $\bar{\mf B}_X$ given by Formula (4.2-1).

\medskip

{\bf Condition 4.3.1 (BRC)}\, Let $({\mf A},{\mf B})$ be any pair with trivial $\T$.

(i)\, Suppose the solid arrows $\mathcal D=\{d_1,\cdots,d_{q}\}$ and
$\mathcal E=\{e_1,\cdots,e_{p}\}$
locate in the lowest non-zero row of $\Theta$ and form the first
$p+q$ arrows of $\mf B$ (not necessarily fulfilling of the whole row), such that $e_1,\cdots, e_{p-1}$ are edges starting from $X$,
$e=e_{p}$ is a loop at $X$, and $d_i\prec e_p,1\leqslant i\leqslant q$.
There exists a set of dotted arrows $\mathcal U=\{u_1,\cdots,u_q\}$, whose complement in $\V^\ast$
is $\mathcal W=\{w_1,\cdots w_t\}$.

(ii)\, Denote by $\bar{\boldsymbol\dz}$ the part of the differential of a solid arrow
in $\mathcal D\cup\mathcal E$ by removing all the monomials containing any solid arrow in $\mathcal D$.
$$\begin{array}{c}\bar{\boldsymbol\dz}(d_i)=u_i+
\sum_{j=1}^t(\sum_{e_l\prec d_i}\lambda_{ijl}e_l)w_j,\quad
1\leqslant i\leqslant q;\\[1.5mm]
\bar{\boldsymbol\dz}(e_i)=\sum_{d_j\prec e_i}\alpha_{ij}
u_j+\sum_{j=1}^{t}(\sum_{l=1}^{i-1}\mu_{ijl}
e_l)w_j,\quad
1\leqslant i\leqslant p,\end{array}$$
where the coefficients $\lambda_{ijl},\alpha_{ij},\mu_{ijl}\in k$.
Then it is said that $({\mf A},{\mf B})$ satisfies the
{\it bottom row condition} with respect to $(\mathcal D,\mathcal U)$
and $(\mathcal E,\mathcal W)$, or (BRC) for short.

\medskip

From now on, we use $G(k)$ instead of the reduction block $G$
given below Lemma 2.3.2 for the sake of simplicity
up to the end of Section 4, if there is not any confusion to be caused.
Suppose the pair $({\mf A},{\mf B})$ satisfies (BRC) with $p,q>1$, and
the first arrow of $\mf B$ is $a_1:X\rightarrow Y$. Now the condition (BRC) on the
induced pair $({\mf A}',{\mf B}')$ is discussed.

Case (i)\, $a_1=d_1$, and $d_1\mapsto\emptyset$, let $\mathcal D'=\{d_2,\cdots,d_q\},
\mathcal U'=\{u_2,\cdots,u_q\}$ and $\mathcal E'=\mathcal E,\mathcal W'=\mathcal W$.

Case (ii)\, $a_1=e_1$, and $e_1\mapsto(0)$, let $\mathcal D'=\mathcal D,\mathcal U'=\mathcal U$;
$\mathcal E'=\{e_2,\cdots,e_p\},\mathcal W'=\mathcal W$.

Denote by $\bullet_i$ a solid or dotted arrow of $\mf B$. Suppose after a reduction
below, $\bullet_i$ splits into a $2\times 1$ or $2\times 2$ matrix in $\mf B'$.
Then the arrows at the second row of the matrix are denoted by $\bullet_{i2}$, or $\bullet_{i21},\bullet_{i22}$.

Case (iii)\, $a_1=e_1$, and $e_1\mapsto{1\choose 0}$, let
$\mathcal D'=\{d_{i21},d_{i22} \mid\, t(d_i)=X\}\cup\{d_{j2} \mid\,
t(d_j)\ne X\}$, $\mathcal U'=\{u_{i21},$ $u_{i22}\mid\, t(d_i)=X\}\cup\{u_{j2} \mid\,
t(d_j)\ne X\}$; and $\mathcal E'=\{e_{22},\cdots,e_{p-1,2},e_{p21},e_{p22}\}$, $\mathcal W'$ is
obtained by the split of $\mathcal W$ and an additional dotted arrow
$\nu_1^\ast(\textsf e_{_{Z_{(X,1)}}}\otimes_k\textsf f_{_{Z_{(X,2)}}})$ defined in Proof 2)
of Proposition 2.1.5.

Case (iv)\, $a_1=e_1$ and $e_1\mapsto{{0\, 1}\choose {0\, 0}}$, let
$\mathcal D'=\{d_{i21},d_{i22} \mid\, t(d_i)=X\,$ or $Y\}\cup\{d_{j2} \mid\,
t(d_j)\ne X,Y\}$, $\mathcal U'=\{u_{i21},u_{i22} \mid\, t(d_i)=X\,$ or $Y\}\cup\{u_{j2} \mid\,
t(d_j)\ne X,Y\}$; and $\mathcal E'=\{e_{i21},e_{i22} \mid\, t(e_i)=Y\}\cup\{e_{j2} \mid\,
t(e_j)\ne Y\}\cup\{ e_{p21}, e_{p22}\}$, $\mathcal W$ is
obtained by the split of $\mathcal W$ and two additional dotted arrows
$\nu_1^\ast(\textsf e_{_{Z_{(X,1)}}}\otimes_k\textsf f_{_{Z_{(X,2)}}})$,
$\nu_1^\ast(\textsf e_{_{Z_{(Y,1)}}}\otimes_k\textsf f_{_{Z_{(Y,2)}}})$ given in Proof 2) of 2.1.5.

\medskip

{\bf Lemma 4.3.2}\, Suppose the pair $({\mf A},{\mf B})$ satisfies (BRC) with $p,q>1$, and
the first arrow of $\mf B$ is $a_1:X\rightarrow Y$. Then after making a reduction $a_1\mapsto G$ as above (i)--(iv),
the induced pair $({\mf A}',{\mf B}')$ satisfies (BRC) with respect to
$(\mathcal D',\mathcal U')$ and $(\mathcal E',\mathcal W')$.

\smallskip

{\bf Proof}\, The cases (i) and (ii) are trivial.

Suppose the reduction block $G={{1}\choose {0}}$ in (iii), resp. ${0\,1}\choose {0\, 0}$ in (iv),
then $X,Y$ split into two vertices $X', Y'$, resp. three vertices $X', Y', Y''$:
$$e_{_X}\mapsto e'_{_X}=\left (\begin{array}{cc}e_{_{Y'}}&w\\
&e_{_{X'}}\end{array}\right ),\quad e_{_Y}\mapsto e_{_Y}'=e{_{_{Y'}}},\,\,\, \mbox{or}\,\,\,
e_{_Y}\mapsto e_{_Y}'=\left (\begin{array}{cc}e_{_{Y''}}&w'\\
&e_{_{Y'}}\end{array}\right ).$$
In above-mentioned two cases, by condition (i) of (BRC) on $({\mf A},{\mf B})$, $e_{i2}$ or $e_{i21},e_{i22}$ for $1<i<p$
start at $X'$, but do not end at $X'$, because $t(e_i)\ne X$; the edge
$e_{p21}: X'\rightarrow Y'$, and the loop $e_{p22}:X'\rightarrow X'$. Therefore the pair $(\mf A',\mf B')$
still satisfies (i) of (BRC). By condition (ii) of (BRC) on $({\mf A},{\mf B})$, we have
$$\begin{array}{c}\bar {\boldsymbol\dz}(D_i)=U_i+\sum_{j=1}^t\lambda_{ij1}GW_j
+\sum_{j=1}^t(\sum_{e_l\prec d_i}\lambda_{ijl}E_l)w_j,\\[1.5mm]
\bar{\boldsymbol\dz}(E_i)=\sum_{d_j\prec e_i}\alpha_{ij}U_j+\sum_{j=1}^t\mu_{ij1}GW_j
+\sum_{j=1}^t(\sum_{l=2}^{i-1}\mu_{ijl}E_l)W_j+{{0\,w}\choose{0\,0}}E_i
-\big(0\,\mbox{or}\,E_i{{0\, w'}\choose{0\,0\,}}\big),\end{array}$$
where $\bar{\boldsymbol\dz}(M)=(\bar{\boldsymbol\dz}(a_{ij}))$ for $M=(a_{ij})$. Since the bottom row of $G$ is $(0)$
or $(0\,0)$, $({\mf A}',{\mf B}')$ still satisfies condition (ii) of (BRC). \hfill$\Box$

\medskip

{\bf Lemma 4.3.3}\, Let $(\bar{\mf A},\bar{\mf B})$ be a one-sided pair with $\bar{\mf B}_X$ given by Formula (4.2-1),
$\bar\T$ being trivial and $|\bar\T|>1$.
Then $(\bar{\mf A},\bar{\mf B})$ satisfies (BRC) with respect to the sets:
$$\begin{array}{lll}&\mathcal D=\{a_i,i\mid\Lambda\}\cup\{b_j\mid j<n_0\},
&\mathcal U=\{\underline v_i\mid i\in\Lambda\}\cup\{\bar u_j\mid j<n_0\};\\[1.5mm]
&\mathcal E=\{\bar a_\tau,1\leqslant\tau\leqslant s\}\cup\{\bar b\},
&\mathcal W=\{\underline v_i\mid i\notin\Lambda\}\cup\{\bar u_j\mid j>n_0\}\cup \{\underline u,\bar v\mbox{-class arrows}\}.\end{array}$$


{\bf Theorem 4.3.4}\, Let $(\bar{\mf A},\bar{\mf B})$ be a one-sided pair with
$\bar{\mf B}_X$ given by Formula (4.2-1), $\bar\T$ being trivial
and $|\bar\T|>1$.
Then there exists a sequence of reductions in the sense of Lemma 2.3.2 as the first part
 of a sequence towards a pair $(\bar{\mf A}^t,\bar{\mf B}^t)$
in the case of MW5 given by Remark 3.4.6:
$$(\bar{\mf A},\bar{\mf B})=(\bar{\mf A^0},\bar{\mf B}^0),(\bar{\mf A}^1,\bar{\mf B}^1),\cdots,
(\bar{\mf A}^\gamma,\bar{\mf B}^\gamma),(\bar{\mf A}^{\gamma+1},\bar{\mf B}^{\gamma+1}),
\cdots,(\bar{\mf A}^{\kappa-1},\bar{\mf B}^{\kappa-1}),
(\bar{\mf A}^\kappa,\bar{\mf B}^\kappa)
\eqno{(4.3\mbox{-}1)}$$ where $\kappa$ is the minimal index,
such that the pair $(\bar{\mf A}^\kappa,\bar{\mf B}^\kappa)$ satisfies the following condition (B).

Condition (B).\, If a row of the formal product $\bar\Theta^\beta$ of $(\M_1^\beta,\N_1^\beta)$
contains some $\bar b$-class arrows of $\bar {\mf B}^\beta$, then the same row of $F^\beta$ contains
one and only one nonzero entry which is a link in some reduction block $G^i_\beta$ of $H^\beta$ obtained by an edge reduction.

\begin{itemize}
\item[(i)] For $i=0, 1,\cdots, (\kappa-2)$,  the reduction from $\bar{\mf A}^{i}$ to
$\bar{\mf A}^{i+1}$  is a composition of a series of reductions $\bar{\mf
A}^{i}=\bar{\mf A}^{i,0}, \bar{\mf A}^{i,1}, \cdots,\bar{\mf A}^{i,r_{i}}$,
$\bar{\mf A}^{i,r_{i}+1}=\bar{\mf A}^{i+1}$:
\begin{itemize}
\item[\ding{172}] For $0\leqslant j< r_{i}-1$, the reduction from $\bar{\mf A}^{i,j}$ to $\bar{\mf
A}^{i,j+1}$ is a sequence of regularizations for non-effective
$a,b$-class arrows, and finally an edge reduction of the form $(0)$
for an effective $a$ or $b$-class arrow. The reduction form
$\bar{\mf A}^{i,r_{i}-1}$ to $\bar{\mf A}^{i,r_i}$ is a sequence
of regularizatiosn for non-effective $a,b$-class arrows.
\item[\ding{173}] The first arrow $a_1^i:X^i\mapsto Y^i$ of $\bar{\mf B}^{i,r_i}$
is an effective $a$ or $b$-class edge with ${\boldsymbol\dz}(a_1^{i})=0$. Making an edge reduction
$a_1^{i}\mapsto{{1}\choose {0}}$ or ${0\, 1}\choose {0\,0}$, the last term
$\bar{\mf A}^{i,r_{i}+1}=\bar{\mf A}^{i+1}$ is obtain.
\end{itemize}
\item[(ii)] It is possible that there exist a minimal integer $\gamma$, and an
index $1\leqslant j\leqslant r_{\gamma}+1$,
such that the first arrow of $\bar{\mf B}^{\gamma,j}$ locates outside the matrix block coming from $\bar b$
under $\bar\T$,
but the first arrow of $\bar{\mf B}^{\gamma,j+1}$ locates at the first column of
the block.
\item[(iii)] The reduction from $\bar{\mf A}^{\kappa-1,0}$ to
$\bar{\mf A}^{\kappa-1,r_{\kappa-1}}$  is a composition of a series of reductions
given by \ding{172} of (i). There are two possibilities.
\begin{itemize}
\item[\ding{172}] The first arrow
$a^{\kappa-1}_1$ of $\bar{\mf B}^{\kappa-1,r_{\kappa-1}}$ is an effective $a$ or $b$-class solid
edge with ${\boldsymbol\dz}(a^{\kappa-1}_1)=0$. Making an edge reduction
$a_1^{\kappa-1}\mapsto(1)$ or $(0\; 1)$, the last term
$\bar{\mf A}^{\kappa-1,r_{\kappa-1}+1}=\bar{\mf A}^\kappa$ is obtained.
\item[\ding{173}] The first arrow
$a^{\kappa-1}_1$ is an effective $b$-class loop at the
down-right corner of the matrix block
coming from $\bar b$ under $\bar\T$ with ${\boldsymbol\dz}(a^{\kappa-1}_1)=0$. Making a loop reduction
$a_1^{\kappa-1}\mapsto W$, a Weyr matrix over $k$,
$\bar{\mf A}^{\kappa-1,r_{\kappa-1}+1}=\bar{\mf A}^\kappa$ is obtained.
\end{itemize}
\end{itemize}

{\bf Proof}\, If there is not any $\bar a$-class edges, i.e. $s=0$, then after a series of regularizations,
the unique effective loop in the induced bocs becomes the first arrow with ${\boldsymbol\dz}(\bar b)=0$.
Since the induced pair is not local,
but the parameter $x$ appears only in a local pair by Remark 3.4.6, set $\bar b\mapsto W$ by a loop reduction of Lemma 2.3.2,
the final pair $(\bar{\mf A}^1,\bar{\mf B}^1)$ satisfies \ding{173} of (iii) with $\kappa=1$.
Suppose $s>0$, regularizations are made on $a_i,b_j$ before $\bar a_1$,
the corresponding $\underline v_i,\bar u_j=0$. Thus ${\boldsymbol\dz}(\bar a_1)=0$ by Formula
(4.2-7), if $\bar a_1\mapsto (0)$, $\bar{\mf A}^{0,1}$ given by \ding{172} of (i) is obtain.
If $r_0>1$, repeating the procedure in \ding{172} of (i), $\bar{\mf A}^{0,r_0}$ is finally reached
with the first arrow $a^{0}_1$ and ${\boldsymbol\dz}(a^{0}_1)=0$.
If $a^{0}_1$ is $\bar a$-class and $a^{0}_1\mapsto(1),(0\, 1)$, \ding{172} of (iii) is obtained;
if $a^{0}_1=\bar b$ then $\bar b\mapsto W$, \ding{173} of (iii) is obtained, and $\kappa=1$ in both cases.

Otherwise, if $a^{0}_1\mapsto{{1}\choose {0}}$ or ${0\, 1}\choose {0\,0}$ in case \ding{173} of (i),
the induced pair $(\bar{\mf A}^1,\bar{\mf B}^1)$ is obtained, which satisfies (BRC) by Lemma 4.3.2--4.3.3.

Suppose $(\bar{\mf A}^{i},\bar{\mf B}^{i})$ for some $i<\kappa-1$ given in (i) has been obtained. Now
we continue the reductions up to the induced pair $(\bar{\mf A}^{i+1},\bar{\mf B}^{i+1})$.
$(\bar{\mf A}^{i},\bar{\mf B}^{i})$ satisfies (BRC) by Lemma 4.3.2--4.3.3 inductively.
Suppose the first arrow of $\bar{\mf B}^{i,0}$ is $a_1^{i,0}=a_{\tau n^iq}$ or $b_{\tau n^iq}$ splitting from a non-effective
arrow $a_\tau$ or $b_{\tau}$ with $n^i$ being the index of the bottom row of $\bar\Theta^i$, and $q$
being the column index inside the splitting block partitioned under $\bar\T$,
then ${\boldsymbol\dz}(a_1^{i,0})=\underline v_{\tau n^iq}$,
or $\bar u_{\tau n^iq}$ by Formulae (4.2-7)--(4.2-8), (4.1-8) and Remark 4.1.3 (i).
Thus $a_1^{i,0}\mapsto\emptyset,\underline v_{\tau n^iq}=0$ or
$\bar u_{\tau n^iq}=0$ by Remark 4.1.3 (ii). The regularizations are continue
for the non-effective arrows inductively, and finally
an effective one is sent to $(0)$, then $\bar{\mf A}^{i,1}$ is
obtained by \ding{172} of (i). With a similar argument as above,
$\bar{\mf A}^{i,r_i}$ is reached, the first arrow of $\bar{\mf B}^{i,r_i}$ has the differential
${\boldsymbol\dz}(a^{i}_1)=0$. Let $a^i_1\mapsto{{1}\choose{0}}$ or
${0\,1}\choose{0\,0}$, the $(i+1)$-th pair is obtained.

If the procedure in (i) was continued without stop, the reduction
sequence would have been infinite. Meanwhile, the $\bar b$-class loop at the
down-right corner of the matrix block splitting from $\bar b$ partitioned under $\bar\T$ has never been reached
in (i). Therefore the procedure of (iii) must occur at some stage, say at the stage $\kappa-1$.

If the first arrow $a^{\kappa-1}_1$ of $\bar{\mf B}^{\kappa-1,r_{\kappa-1}}$ is an
edge, then $a_1^{\kappa-1}\mapsto(1)$ or $(0\, 1)$
gives case \ding{172} of (iii). If $a^{\kappa-1}_1$ is a loop,
then $a_1^{\kappa-1}\mapsto W$ by a loop reduction of Lemma 2.3.2
gives case \ding{173} of (iii), since $\bar{\mf A}^{\kappa-1,r_{\kappa-1}}$ is not local
from $\bar{\mf A}^{\kappa-1,0}$ being non-local by (BRC).
In both cases, the induced pair $(\bar{\mf A}^\kappa,\bar{\mf B^\kappa})$
has the minimal index $\kappa$ satisfying Condition (B).\hfill$\Box$

\medskip

Suppose $s(a_1^{i-1})=X^{i-1}$ in case (i) \ding{173} of Theorem 4.3.4,
the reduction on $a^{i-1}_1$ gives $e_{_{X^{i-1}}}\mapsto{{e_{_{Y^{i}}}\,\,
\bar w^{i}\,\,}\choose{\,\,0\,\,\,\, e_{_{X^{i}}}}}$ for $1\leqslant i<\kappa$.
Denote by $\bar W^i_\kappa$ the split of $\bar w^i$ in $e_{_{X}}^\kappa$
for $1\leqslant i<\kappa$. Then $\bar W^i_\kappa$ can be divided into $(\kappa-i)$ blocks $\bar W^i_{\kappa i},
\cdots,\bar W^i_{\kappa,\kappa-1}$. Denote by $n^\kappa_i$
the size of $e_{_{Y^i}}^\kappa$, and by $n^\kappa_\kappa$ that
of $e^\kappa_{_{X^\kappa}}$, which is $1$ in case  (iii) \ding{172} of Theorem 4.3.4,
or is the same as that of $W$ in (iii) \ding{173}.
Thus $\bar W^i_{\kappa j}$ has the size $n^\kappa_i\times n^\kappa_{j+1}$.
Write $n^\kappa=\sum_{i=1}^\kappa n^\kappa_i$, which is the size of $e_{_X}^\kappa$. We have
$$ e_{_X}^\kappa=\left (\begin{array}{ccccc} e^{\kappa}_{_{Y^1}}&\bar W^1_{\kappa1}
&\cdots &\cdots &\bar W^1_{\kappa,\kappa-1}\\
& e^{\kappa}_{_{Y^2}}&\cdots &\cdots &\bar W^2_{\kappa,\kappa-1}\\
&&\cdots&\cdots &\\
&&& e^{\kappa}_{_{Y^{\kappa-1}}}&\bar W^{\kappa-1}_{\kappa,\kappa-1}\\
&&&&e^\kappa_{_{X^\kappa}}\end{array}\right
)\eqno{(4.3\mbox{-}2)}$$

{\bf Corollary 4.3.5}\, The elements in $\bar W^i_\kappa$ for $1\leqslant i<\kappa$
are dotted arrows of $\mf B^\kappa$.

\smallskip

{\bf Proof} The assertion is already implied in the proof of Theorem 4.3.4.

\medskip

When we make an edge or a loop reduction of Lemma 2.3.2, the dotted arrows $\{{F'_i}^\ast\mid i=1,\cdots,l\}$,
given in proof 2) of Proposition 2.1.5 are said to be {\it $w$-class arrows},
where the dotted arrows in $\bar W^i_{\kappa}$ for $1\leqslant i<\kappa$ of Formula (4.3-2)
are specially said to be {\it $\bar w$-class}. Furthermore the elements splitting from $w$
or $\bar w$-class arrows are still said to be in the same class.

\bigskip
\bigskip
\noindent {\large\bf 4.4\, Major pairs}
\bigskip

We prove in this subsection that under
some further assumption, a one sided pair
$(\bar{\mf A},\bar{\mf B})$ with $\bar{\mf B}_X$
given by Formula (4.2-1) is not homogeneous.

\medskip

Let $(\bar{\mf A},\bar{\mf B})$ be a one-sided pair, where
$\bar{\mf B}_X$ is given by Formula (4.2-1), and $\bar{\mf B}$
has $s$ $\bar a$-class arrows with $s\geqslant 1$.
According to the coefficients of the first two Formulae of (4.2-8),
$s$ linear combinations of the $\bar v$-class arrows are define
in $\bar{\mf B}$:
$$\begin{array}{c}\hat{v}_\tau=\sum_{j}(\bar\varepsilon_{\tau j}
-\sum_{\bar a_\tau\prec b_i\prec\bar b}\bar\beta_i\varepsilon_{i\tau j})
\bar v_{j},\quad \tau=1,\cdots,s.\end{array}\eqno{(4.4\mbox{-}1)}$$
Fix any $1\leqslant\tau\leqslant s$, making reductions according to (iii) \ding{172}
of Theorem 4.3.4 for $\kappa=1$, such that $a^0_1=\bar a_\tau\mapsto (1)$,
the induced pair $(\bar{\mf A}^1,\bar{\mf B}^1)$ is reached. Then we continue to do further reductions
based on Formulae (4.2-7)--(4.2-8) inductively.
For $\bar a_{\eta}\prec a_i,b_i\prec\bar a_{\eta+1}, \tau\leqslant\eta<s$ and
$\bar a_s\prec a_i,b_i\prec\bar b$, by Remark 4.1.3 (ii):
$$\begin{array}{c} a_i\mapsto\emptyset,\, \underline v_i                                                                                                                                       +1\sum_j\epsilon_{i\tau j}\underline u_j=0,\quad
b_i\mapsto\emptyset,\, \bar u_i+1\sum_j\varepsilon_{i\tau j}\bar v_j=0.
\end{array}\eqno{(4.4\mbox{-}2)}$$
On the other hand, $\bar a_\eta\mapsto\emptyset$ or $(0)$ for $\tau<\eta\leqslant s$ corresponding to
${\boldsymbol\dz}(\bar a_{\tau+\eta})\ne 0$ or $=0$. the dotted element $\bar u_i$ is replaced by the linear
composition of $\bar v$-class arrows inductively by Remark 4.1.3 (iii),
the second formula of (4.2-8) shows the formula below in some induced pair:
$$\begin{array}{ll}{\boldsymbol\dz}(\bar b)&=\sum_{\bar a_\tau\prec b_i\prec\bar b}\bar\beta_{i} \bar u_{i}
+1(\sum_{j}\bar\varepsilon_{\tau j}\bar v_{j})\\[1mm]
&=\sum_{j}(\bar \varepsilon_{\tau j}-\sum_{\bar a_\tau\prec b_i\prec\bar b}
\bar\beta_{i}\varepsilon_{i\tau j})\bar v_{j}=\hat v_\tau.
\end{array}\eqno{(4.4\mbox{-}3)}$$

{\bf Lemma 4.4.1}\, Let $(\bar{\mf A},\bar{\mf B})$ be a one-sided pair with $\bar\T$ being trivial,
$s\geqslant 1$, and $\bar{\mf B}_X$ given by Formula (4.2-1). If there exists some $1\leqslant\tau\leqslant s$,
with $\hat v_\tau=0$ in Formula (4.4-1), then $\bar{\mf B}$ is wild and non-homogeneous.

\smallskip

{\bf Proof}\, If $\bar a_\tau:X\rightarrow Y$, it may be assumed that $\T=\{X,Y\}$.

1)\, Since $\bar{\mf B}_X$ is minimal with $R_X=k[x,\phi(x)^{-1}]$,
there is an almost split conflation $(e'_\lambda):S'_\lambda\rightarrow E'_\lambda\rightarrow S'_\lambda$
for any $\lambda\in\mathscr L'=k\setminus\{$roots of $\phi(x)\}$ in $R(\mf{\bar B}_X)$.
Let $\vartheta:R(\mf{\bar B}_X)\rightarrow R(\mf{\bar B})$ be the induced functor. If
$\mf{\bar B}$ is homogeneous, then there is a co-finite subset $\mathscr L\subseteq\mathscr L'$,
and a set of almost split conflations $\{(e_{\lambda})=\vartheta(e_\lambda'):
S_\lambda\rightarrow E_\lambda\rightarrow S_\lambda\mid\lambda\in\mathscr L\}$
by Corollary 3.2.4.

2)\, According to Formula (4.4-1)-(4.4-3), an induced pair $(\bar{\mf A}',\bar{\mf B}')$
is obtained with ${\boldsymbol\dz}(\bar b)=0$. Thus it is possible to
construct an object $L\in R(\bar{\mf B})$ with $L_X=k,L_Y=k,L(\bar a_\tau)=(1),L(\bar b)=(\lambda)$
and $L(b_i)=0,i>n_0,L(c_i)=0,i=1,\cdots,t$. The same argument given in 2)--3) of the proof of Lemma 3.4.1
shows that $\bar{\mf B}$ is non-homogeneous.
And $\bar{\mf B}$ is wild by \cite[Theorem A]{CB1}. \hfill$\Box$

\medskip

{\bf Theorem 4.4.2}\, Let $(\bar{\mf A},\bar{\mf B})$ be a one-sided pair
with $\bar\T$ being trivial, $s>1$, and
$\bar{\mf B}_X$ given by Formula (4.2-1).
If the elements $\{\hat v_1,\hat v_2,\cdots,\hat v_s\}$ defined
in Formula (4.4-1) are linearly dependent, then $\bar{\mf B}$ is wild and non-homogeneous.

\smallskip

{\bf Proof}\, Without loss of generality, it may be assumed $\bar\T=\{X,Y\}$.
Suppose there is a minimal linearly dependent subset
$\{\hat v_{\tau_1}, \hat v_{\tau_2}, \cdots, \hat v_{\tau_l}\}$ with $l$ vectors.
Since the case of $l=1$ has been treated in Lemma 4.4.1, it is assumed here that $l>1$. Suppose
$\tau_1<\tau_2<\cdots<\tau_l$ and
$$\begin{array}{c}\hat v_{\tau_1}=\beta_2\hat v_{\tau_2}+\cdots+\beta_l\hat v_{\tau_l}, \quad \beta_2,
\cdots,\beta_l\in k^\ast.
\end{array}\eqno{(4.4\mbox{-}5)}$$

1) Making reductions according to Theorem 4.3.4 (i) and (iii) \ding{172} for $\kappa=l$,
such that $a^{p-1}_1\mapsto{{1}\choose{0}}$
for $1\leqslant p<l$, and $a_1^l\mapsto(1)$, an
induced pair $(\bar{\mf A}^l,\bar{\mf B}^l)$ is obtained. The sum $F^l+\bar\Theta^l$
looks like (with only $\bar a,\bar b, c$-class arrows):
$$
\begin{array}{|c|c|c|c|c|c|c|c|c|c|c|} \hline
& 1& \cdots& \bar{a}_{\tau_2,1}&\cdots& \bar{a}_{\tau_3,1}&\cdots &
\bar{a}_{\tau_l,1}&\bar{a}_{\tau_l+1,1}\cdots \bar{a}_{s1}&
\bar b_{11}\ \bar b_{12}\ \cdots \ \bar b_{1l} & c_{11}\ \cdots\ c_{t1}\\
\cline{3-11}
&&&1&\cdots& \bar{a}_{\tau_3,2}&\cdots &
\bar{a}_{\tau_l,2}&\bar{a}_{\tau_l+1,2}\cdots \bar{a}_{s2}&\bar b_{21}\
\bar b_{22}\ \cdots \ \bar b_{2l} & c_{12}\ \cdots\ c_{t2}\\
\cline{5-11}
0&0&&0&&1&\cdots & \bar{a}_{\tau_l,3}&\bar{a}_{\tau_l+1,3}\cdots
\bar{a}_{s3}&\bar b_{31}\
\bar b_{32}\ \cdots \ \bar b_{3l} & c_{13}\ \cdots\ c_{t3}\\
\cline{7-11}
&&&&&&\cdots && &\cdots&\\ \cline{8-11}
 &&&&&0&&1&\bar{a}_{\tau_l+1,l}\cdots \bar{a}_{sl}&\bar b_{l1}\
\bar b_{l2}\ \cdots \ \bar b_{ll} & c_{1l}\ \cdots\ c_{tl}\\ \hline
\end{array}\eqno{(4.4\mbox{-}6)}$$
Since $\bar{\mf A}$ has two vertices, the dimension of
$\vartheta^{0l}(F^l(k))$ in $R(\bar{\mf A})$ is $l+1$, and the number of links
of $F^l$ is $l$, the pair $(\bar{\mf A}^l,\bar{\mf B}^l)$ is local by the assertion below Formula (2.3-7).

2) We make further reductions from $\bar{\mf B}^l$ inductively
for the $\bar p$-th row ordered by $\bar p=l,l-1,\cdots,2$ in the reduced formal product $\Theta^l$.
For $\bar p=l$, similar to Formulae (4.4-2)--(4.4-3):
$a_{il}\mapsto\emptyset,i\in\Lambda$; note that $\bar v_j: Y\rightarrow X$,
the matrix splitting from $\bar v_j$ in $\bar{\Psi}_{\n^l}$ is $(\bar v_{j1},\cdots,\bar v_{jl})$ of size $1\times l$,
$b_{ilq}\mapsto\emptyset,\bar u_{ilq}+1\sum_j
\varepsilon_{i\tau_l j}\bar v_{jq}=0,i<n_0; \bar a_{\eta l}\mapsto(0)$ or $\emptyset,
\tau_l<\eta\leqslant s$;
$$\begin{array}{ll}{\boldsymbol\dz}^0(\bar b_{lq})&=\sum_{\bar a_{\tau_l}\prec b_i\prec\bar b}\bar\beta_{i} \bar u_{ilq}
+1(\sum_{j}\bar\varepsilon_{\tau_l j}\bar v_{jq})\\[1mm]
&=\sum_{j}(\bar \varepsilon_{\tau_l j}-\sum_{\bar a_\tau\prec b_i\prec\bar b}
\bar\beta_{i}\varepsilon_{i\tau_l j})\bar v_{jq}=\hat v_{\tau_l q},
\end{array}$$
thus $\bar b_{lq}\mapsto\emptyset,
\hat v_{\tau_lq}=0$ for $q=1,\cdots,l$ inductively.
Next, $b_{ilq}\mapsto\emptyset$
for $i>n_0, 1\leqslant q\leqslant l$, by Remark 4.1.3 (ii); and
$c_{il}\mapsto (0)$ or $\emptyset$. The dotted arrows $\underline v_{ip},\bar u_{ipq}$
for all $i$ and $p<l,1\leqslant q\leqslant l$,
$\hat v_{iq}$ for $i<l$ and $1\leqslant q\leqslant l$ are preserved by 4.1.3 (iv).
The induced bocs $\bar{\mf B}^{l+1}$ follows.

3)\, Suppose an induced bocs $\bar{\mf B}^{2l-\bar p}$ is reached for some
$\bar p<l$. Denote the entries of $F^{2l-\bar p}$, which are not the entries of $G_{2l-\bar p}^j$
for $j=1,\cdots,l$, by $\bullet^0$ coming from $\bullet$, one of the $a,b,c$-class solid arrows.
Then $\bar{\mf B}^{2l-\bar p}$ satisfies the following two conditions.

\subitem\hspace{-3mm}\ding{172} for any $p>\bar p$, $a_{ip}^0=\emptyset,i\in\Lambda;\bar a_{ip}^0=\emptyset$ or $(0)$; $b_{ipq}^0=\emptyset,i\ne n_0;
\bar b_{pq}=\emptyset;c_{ip}=\emptyset$ or $(0)$.

\subitem\hspace{-3mm}\ding{173} The dotted arrows $\underline v_{ip}$,
$\underline u_{ipq}$ for all $i,q$ and $p\leqslant\bar p$; and $\hat v_{iq}$ for $i\leqslant \bar p$ and $1\leqslant q\leqslant l$
are preserved.

Now we continue with reductions on the solid arrows at the $\bar p$-th row of
$\Theta^{2l-\bar p}$. According to assumption \ding{172}--\ding{173} and Remark 4.1.3 (ii),
with a similar discussion as in 2):
$a_{i\bar p}^0\mapsto\emptyset,i\in\Lambda$, $b_{i\bar pq}\mapsto\emptyset,\bar u_{i\bar pq}
+1\sum_j\varepsilon_{i\tau_{\bar p} j}\bar v_{jq}=0,i<n_0$;
$\bar a_{\eta\bar p}\mapsto(0)$ or $\emptyset,\tau_{\bar p}<\eta\leqslant s$;
$\bar b_{\bar pq}\mapsto\emptyset,\hat v_{\tau_{\bar p}q}=0$; $b_{i\bar pq}\mapsto\emptyset,i>n_0$,
$c_{i\bar p}\mapsto(0)$ or $\emptyset$, an induced pair $(\bar{\mf A}^{2l-\bar p+1},\bar{\mf B}^{2l-\bar p+1})$
satisfying assumption \ding{172}--\ding{173} is reached. Finally a pair $(\bar{\mf A}^{2l-1},\bar{\mf B}^{2l-1})$ is obtained
by induction on $\bar p$.

4) Making reductions on the arrows before $\bar b_{11}$ in the first row of (4.4-6):
$a_{i1}\mapsto\emptyset,i\in\Lambda$; $b_{j1q}\mapsto\emptyset,j<n_0$;
$\bar a_{i1}\mapsto\emptyset$ or $(0)$ for $\tau_1<i\leqslant s$,
an induced bocs $\bar{\mf B}^{2l}$ is obtained. Formula (4.4-5) gives
$\hat v_{\tau_1q}=\beta_2\hat v_{\tau_2 q}+\cdots+\beta_r\hat v_{\tau_lq}=0$,
thus ${\boldsymbol\dz}^0(\bar b_{1q})=0$ for $q=1,\cdots,l$. Since $l\geqslant 2$, the bocs
$\bar{\mf B}^{2l}$ is wild and non-homogeneous
by (iii) or (iv) of Classification 4.2.1. And hence so is $\bar{\mf B}$.\hfill$\Box$

\medskip

{\bf Definition 4.4.3}\, A one-sided pair $(\bar{\mf A},\bar{\mf B})$ with $\bar{\mf
B}_X$ given by Formula (4.2-1) is said to be a {\it
major pair}, provided that $\{\hat v_1,\hat v_2,\cdots,\hat v_s\}$
of Formula (4.4-1) are linearly independent.

\bigskip
\bigskip
\noindent {\large\bf 4.5 Further reductions}
\bigskip

Throughout the subsection let $(\bar{\mf A},\bar{\mf B})$ be a one-sided pair, such that
$\bar{\mf B}_X$ is given by Formula (4.2-1); $\bar\T$ is trivial and $|\bar\T|>1$;
the pair is major; and the $c$-class arrows satisfy Formula (4.2-6). Suppose
$(\bar{\mf A}^\kappa,\bar{\mf B}^\kappa)$ is an induced pair
given by Theorem 4.3.4 (iii), we continue to construct a sequence of reductions of Lemma 2.3.2,
which is still the first part of a sequence towards a pair $(\bar{\mf A}^t,\bar{\mf B}^t)$
in the case of MW5. Let $(\bar{\mf A}^\varsigma,\bar{\mf B}^\varsigma)$
be an induced pair of $(\bar{\mf A}^\kappa,\bar{\mf B}^\kappa)$ in the sequence for some
$\varsigma\geqslant\kappa$. The explicit linear
relation of dotted elements in the induced bocs $\bar{\mf B}^{\varsigma+1}$ will be
expressed via calculating the differential of the first arrow of $\mf B^\varsigma$.

\medskip

Concern the formal product $\bar\Theta^\varsigma$ and the pseudo formal product $\bar\Pi^\varsigma$ of
the pair $(\bar{\mf A}^\varsigma,\bar{\mf B}^\varsigma)$.
Put a solid or dotted arrow in a square box; and
a matrix block in a rectangular box with four boundaries. The reduction block $G^j_\varsigma$ of $H^\varsigma$ for
$1\leqslant j\leqslant\kappa$ is defined before Formula (2.3-7), whose
upper boundary is that of $I^j_\varsigma$ and denoted by $m_\varsigma^{j-1}$,
the lower boundary is that of $F^\varsigma$; the left and right boundaries are
given by the dotted lines $l^j_\varsigma$ and $r^j_\varsigma$ respectively, see Picture (4.5-1) below.
The set of entries and solid arrows in $(F^\varsigma+\Theta^\varsigma)$ but not in $G^j_\varsigma$
for $1\leqslant j\leqslant\kappa$ is divided into $\kappa$ blocks, where the $j$-th block has the upper boundary
$m_\varsigma^{j-1}$ and the lower one $m_\varsigma^{j}$,
the left boundary $r^j$ and the right one being that of $F^\varsigma+\Theta^\varsigma$.
Denote by $A^j_\varsigma,B^j_\varsigma,C^j_\varsigma$ in the $j$-block containing  the sets of entries in $F^\varsigma$
obtained by reductions for $a,b,c$-class arrows and the sets of $a,b,c$-class solid arrows
in $\Theta^\varsigma$.
\vspace{-5mm}
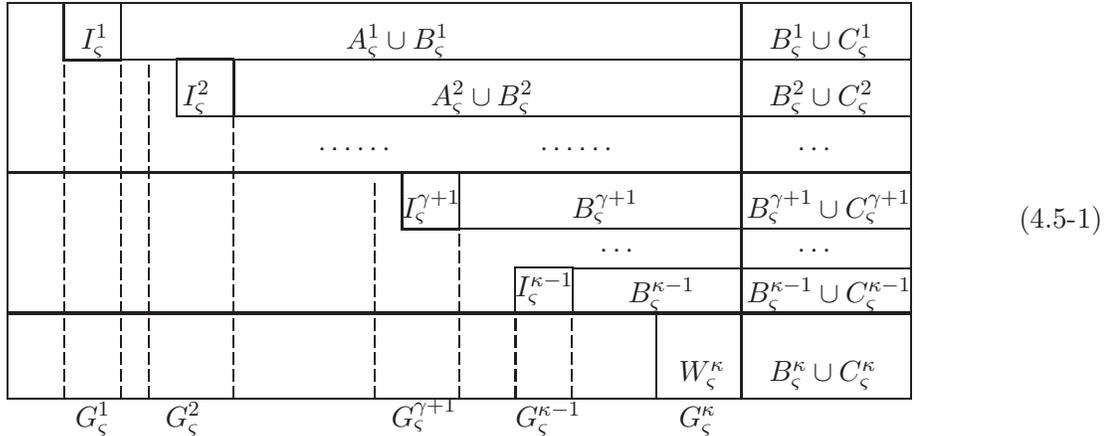
\begin{center}
\begin{equation*}
\begin{array}{c} {\unitlength=0.75mm
\begin{picture}(180,75)
\put(10,5){\framebox(160,70)}
\put(20,65){\framebox(10,10)}\multiput(20,5)(0,3){20}{\line(0,1){2}}
\multiput(30,5)(0,3){20}{\line(0,1){2}}\put(20,65){\line(1,0){150}}
\put(23,67){$I^1_\varsigma$}\put(70,67){$A^1_\varsigma\cup B^1_\varsigma$}
\put(145,67){$B^1_\varsigma\cup C^1_\varsigma$}\put(22,0){$G^1_\varsigma$}
\put(40,55){\framebox(10,10)}\multiput(35,5)(0,3){20}{\line(0,1){2}}
\multiput(50,5)(0,3){17}{\line(0,1){2}}\put(40,55){\line(1,0){130}}
\put(41,57){$I^2_\varsigma$}\put(85,57){$A^2_\varsigma\cup B^2_\varsigma$}
\put(145,57){$B^2_\varsigma\cup C^2_\varsigma$}\put(38,0){$G^2_\varsigma$}
\put(80,35){\framebox(10,10)}\multiput(75,5)(0,3){13}{\line(0,1){2}}
\multiput(90,5)(0,3){11}{\line(0,1){2}}\put(90,35){\line(1,0){80}}
\put(80,37){$I^{\gamma+1}_\varsigma$}\put(110,37){$B^{\gamma+1}_\varsigma$}
\put(141,37){$B^{\gamma+1}_\varsigma\cup C^{\gamma+1}_\varsigma$}\put(78,0){$G^{\gamma+1}_\varsigma$}
\put(100,20){\framebox(10,8)}\multiput(100,5)(0,3){5}{\line(0,1){2}}
\multiput(110,5)(0,3){5}{\line(0,1){2}}\put(110,20){\line(1,0){60}}\put(110,28){\line (1,0){60}}
\put(100,23){$I^{\kappa-1}_\varsigma$}\put(120,22){$B^{\kappa-1}_\varsigma$}
\put(141,22){$B^{\kappa-1}_\varsigma\cup C^{\kappa-1}_\varsigma$}\put(100,0){$G^{\kappa-1}_\varsigma$}
{\thicklines\put(10,45){\line(1,0){160}}}{\thicklines\put(10,20){\line(1,0){160}}}
\put(125,5){\line(0,1){15}}\put(140,5){\line(0,1){70}}
\put(129, 8){$W^{\kappa}_\varsigma$}\put(145,8){$B^{\kappa}_\varsigma\cup C^{\kappa}_\varsigma$}
\put(129,0){$G^{\kappa}_\varsigma$}
\put(65,48){$\cdots\cdots\qquad\qquad\quad\cdots\cdots$}\put(150,48){$\cdots$}
\put(115,30){$\cdots$}\put(150,30){$\cdots$}
\end{picture}}\end{array}\eqno{(4.5\mbox{-}1)}\end{equation*}\end{center}

In $A^j_{\varsigma},B^j_{\varsigma},C^j_{\varsigma}$, a solid arrow is denoted by
$\bullet_{ipq}^j$ splitting from $a_i,b_i,c_i$
respectively, an entry of $F^\varsigma$ by $\bullet_{\varsigma i,pq}^{j,0}$,
where $(p,q)$ is the index in the
$n^\varsigma\times n^\varsigma_{t(a_i)}$,$n^\varsigma\times n^\varsigma$ or
$n^\varsigma\times n^\varsigma_{t(c_i)}$-block matrices.
For the sake of convenience, $\Phi_{\m^\varsigma}^{m}$ of Formula (4.1-7) is also partitioned
by the lines $m^\varsigma_j,l^\varsigma_j,r^\varsigma_j$ in the same way as in $F^\varsigma+\Theta^\varsigma$.

\medskip

{\bf Remark}\, From now on, the pseudo formal equation (4.1-7)
at the $\varsigma$-th step for presenting the differential of the first arrow is considered.
After a loop or an edge reduction, some $w$-class dotted arrows may be added into the induced bocs,
but the linear relations among the splits of the dotted elements in the induced pair can be
obtained completely by Remark 4.1.3 (i). Therefore
the new relation of $\bar u,\bar v, \underline u,\underline v,\bar w,w$-class
elements during the regularization from $\bar{\mf B}^\varsigma$ to $\bar{\mf B}^{\varsigma+1}$ will be concentrated on.

Throughout the subsection, suppose the first arrow
$a^\varsigma_1=\bullet_{\tau\bar p\bar q}^\iota$ of $\bar{\mf B}^\varsigma$ belongs to $A^\iota_\varsigma
\cup B^\iota_\varsigma\cup C^\iota_\varsigma$ in Picture (4.5-1).
Denote by $\bar u^j_{ipq}$ or $\underline v^j_{ipq}$ the dotted element
in $\bar\Psi^m_{\n^\varsigma}$, which
corresponds to the entry $b^{j,0}_{ipq}$ or $a^{j,0}_{ipq}$ in $F^\varsigma$,
thus $j\geqslant\iota$, $p>\bar p$, or $p=\bar p$ but $q<\bar q$. On the other hand, denote by $\bar u^{j'}_{i'p'q'}$
or  $\underline v^{j'}_{i'p'q'}$ the dotted element corresponding to the solid arrow
$b^{j'}_{i'p'q'}$ or $a^{j'}_{i'p'q'}$ in $\Theta^\varsigma$,
thus $j'\leqslant\iota$, and $p'<\bar p$ or $p'=\bar p,q'\geqslant\bar q$.
All of them are coming from $(\mathcal D,\mathcal U)$ in Lemma 4.3.3,

In the following Lemmas, the index $n_0$ is defined in Formula (4.2-1), and the number
$h$ and the set $\Lambda$ are defined in Formula (4.2-5).

\medskip

{\bf Lemma 4.5.1}\, Let $(\mf A^\varsigma,\mf B^\varsigma)$ be a pair induced from
$(\mf A^{\kappa},\mf B^{\kappa})$ given by Theorem 4.3.4 (iii) \ding{173}.
Suppose the first arrow of $\mf B^\varsigma$, $a^\varsigma_1=\bullet_{\tau\bar p\bar q}^{\kappa}\in
B^{\kappa}_\varsigma\cup C^{\kappa}_\varsigma$, (see the second thick line below in Picture
(4.5-1) for example). Assume that

(i)\, all $b^{\kappa,0}_{ipq}=\emptyset,i>n_0$, and the corresponding dotted element
$\bar u_{ipq}^{\kappa}$ is replaced by a linear combination of some $\bar v$-class arrows in $\bar{\mf B}^\varsigma$;
while the dotted arrows $\bar u_{i'p'q'}^{j'}$
are preserved;

(ii)\, if $c^{\kappa,0}_{ipq}=\emptyset$, there is a linear relation among some elements
$\underline v_{i_1p_1q}^\kappa,h<i_1\leqslant h+i$, $p_1\geqslant p$ and some $\underline u,w$-class
arrows in $\bar{\mf B}^\varsigma$; while all the dotted arrows $\underline v_{i'p'q'}^{j'}$ are preserved.

Then after a regularization, the induced pair $(\mf A^{\varsigma+1},\mf B^{\varsigma+1})$
still satisfies (i)-(ii). In particular all the dotted arrows $\underline v_{i'p'q'}^{j'}$
except $\underline v_{\tau\bar p\bar q}^{\kappa}$ in case  (ii); $\bar u^{j'}_{i'p'q'}$,
except $\bar u^{\kappa}_{\tau\bar p\bar q}$ in case  (i);
and all the $\bar w,\bar v$-class arrows are preserved.

\smallskip

{\bf Proof}\, The assumption (i)--(ii) are valid for $\varsigma=\kappa$ by Theorem 4.3.4 and Corollary 4.3.5.

(i)\, If $a^\varsigma_1=b^\kappa_{\tau\bar  p\bar q},\tau>n_0$,
then according to the third formula of (4.2-8) and Formula (4.1.8),
$$\begin{array}{l}
{\boldsymbol\dz}(b_{\tau\bar p\bar q}^\kappa)=\bar u_{\tau\bar p\bar q}^\kappa
+\sum_{n_0<i<\tau}(\beta^0_{\tau i}\bar u_{i\bar p\bar q}^\kappa
+\sum_{q}\beta^1_{\tau i}\bar b^{\kappa,0}_{\bar pq}\bar u_{iq\bar q}^\kappa)
+\sum_{c_i\prec b_\tau,q}c^{\kappa,0}_{i\bar pq}(\sum_{l}\varepsilon'_{\tau il}
\bar v_{lq\bar q}).\end{array}$$
Since $W^{\kappa}_{\zeta}$ is upper triangular, the index $\bar p\leqslant q$ in $\bar b^{\kappa,0}_{\bar p q}$.
By assumption (i), $\bar u_{\tau\bar p\bar q}^\kappa$ is a dotted arrow,
thus $b_{\tau\bar p\bar q}^\kappa\mapsto\emptyset$,
$\bar u_{\tau\bar p\bar q}^\kappa$ is replaced by a linear combination of some $\bar v$-class arrows by 4.1.3 (ii)--(iii),
since $\bar u_{i\bar p\bar q}^\kappa,\bar u_{iq\bar q}^\kappa$ are already replaced by those arrows still by assumption (i).

(ii)\, If $a^\varsigma_1=c_{\tau\bar p\bar q}^\kappa$, then according to Formula (4.2-9) and (4.1-8),
$$\small\begin{array}{l}{\boldsymbol\dz}(c_{\tau\bar p\bar q}^\kappa)=
\sum_{h<i\leqslant h+\tau}(\gamma^0_{\tau i}\underline v_{i\bar p\bar q}^\kappa
+\sum_{q}\gamma_{\tau i}^1\bar b^{\kappa,0}_{\bar pq}\underline v_{iq\bar q}^\kappa)
+\sum_{i,q}
c^{\kappa,0}_{i\bar pq}(\sum_{l}\xi_{\tau il}\underline u_{lq\bar q})
+\sum_{q<\bar q}c^{\kappa,0}_{\tau\bar pq} w_{q\bar q}-\sum_{p>\bar p} w_{\bar pp}c^0_{\tau\bar pq}.
\end{array}$$
In the case of ${\boldsymbol\dz}(c_{\tau\bar p\bar q})\ne 0$,
$c_{\tau\bar p\bar q}\mapsto\emptyset$, and a linear relation
among elements $\underline v_{i\bar p\bar q},\underline v_{iq\bar q},h<i\leqslant h+\tau, q\geqslant\bar p$,
and some $\underline u,w$-class elements is added, which is given by the right-hand
side of the above formula being equal to $0$.

The required $\underline v,\bar u$-class and all the
$\bar v,\bar w$-class dotted arrows are preserved,
the pair $(\mf A^{\varsigma+1},\mf B^{\varsigma+1})$ still satisfies assumption
(i)-(ii).\hfill$\Box$

\medskip

If $\gamma$ exists in Theorem 4.3.4 (ii), suppose
$\gamma<j\leqslant\kappa$ in case  (iii) \ding{172} of 4.3.4, or $\gamma<j<\kappa$, in case  (iii) \ding{173},
see what between the two thick lines of Picture (4.5-1) for example.

A $\bar w$-class elements of $e_{_X}^\varsigma$ is described by two indices $\bar w_{pq}$,
where $p$ is the row index and $q$ the column index of $e_{_X}^\varsigma$.
Suppose the identity matrix $I^j_\varsigma$
intersects the $p$-th row of $F^\varsigma$ at the $q_p^j$-th column with $\bar b^{j,0}_{pq_p^j}=(1)$.
Denote by $\bar w_{q_{p}^{j} q}$ for any possible $q$ the dotted element with the row index $q_p^j$
in $e_{_X}^\varsigma$. Let $\gamma<\iota\leqslant\kappa$ (or $<\kappa$) be
given above. If $p>\bar p$, then $\bar w_{q_p^jq}$ is sitting below $\bar w_{q_{\bar p}^{\iota}q}$.
Note that $\bar W^j_\varsigma$, the $j$-block in $e^\varsigma_{_X}$, defined by Formula (4.3-2) is also sitting in the
$(n_0,n_0)$-th block of $\bar\Psi^r_{\n^\varsigma}$ partitioned under $\bar\T$. Clearly, the row indices of $\bar W^j_\varsigma$ coincide
with the column indices of $I^j_\varsigma$, see $(I^3,\bar W_3),(I^4,\bar W_4)$ in Picture (4.5-2) at the end of the subsection.

\medskip

{\bf Lemma 4.5.2}\, Let $(\bar{\mf A}^\varsigma,\bar{\mf B}^\varsigma)$ be an induced pair of
$(\bar{\mf A}^\kappa,\bar{\mf B}^\kappa)$ with $\gamma$ existing in Theorem 4.3.4 (ii).
Suppose the first arrow of $\mf B^\varsigma$, $a^\varsigma_1=\bullet_{\tau\bar p\bar q}^\iota\in
B^\iota_\varsigma\cup C^\iota_\varsigma$ with $\gamma<\iota\leqslant\kappa$ (or $<\kappa$). Assume that

(i)\, all $\bar b^{j,0}_{pq}=\emptyset$, the
corresponding dotted element $\bar w_{q_{p}^{j} q}$ is replaced by a linear combination of some
$\bar w_{pq}$ for $p>q_p^j$, and some $w$-class elements in $\bar{\mf B}^\varsigma$;
while the dotted arrows $\bar w_{p'q'}$ for $p'\leqslant q_{\bar p}^j$
are preserved;

(ii)\, all $b^{j,0}_{ipq}=\emptyset,i>n_0$, the corresponding
$\bar u_{ipq}^{j}$ is replaced by a linear combination of
an element $\bar u_{i_{1}p_1q}^{j_1}$, with $n_0<i_1<i,j\leqslant j_1<\gamma,p_1< p$,
and some $\bar v$-class elements in $\bar{\mf B}^\varsigma$; while the dotted arrows $\bar u_{i'p'q'}^{j'}$
are preserved;

(iii)\, if $c^{j,0}_{ipq}=\emptyset$, there is a linear relation among some elements
$\underline v_{i_1p_1q}^{j_1},h<i_1\leqslant h+i,j\leqslant j_1<\gamma,p_1\leqslant p$,
and some $\underline u,w$-class elements in $\bar{\mf B}^\varsigma$;
while the dotted arrow $\underline v_{i'p'q'}^{j'},i'\in\Lambda$,
are preserved.

Then after a regularization, the induced pair $(\mf A^{\varsigma+1},\mf B^{\varsigma+1})$
still satisfies (i)-(iii). In particular all the dotted arrows $\underline v_{i'p'q'}^{j'},i'\in\Lambda$ in case  (iii);
$\bar u^{j'}_{i'p'q'}$ except $\bar u^{j}_{\tau\bar p\bar q}$ in case  (ii);
$\bar w_{p'q'},p'<q_{\bar p}^j$, in case  (i); and all the $\bar v$-class dotted arrows are preserved.

\smallskip

{\bf Proof}\, The assumption (i)--(iii) are valid in the following two cases. First,
the box of $a_1^{\varsigma}$ has the bottom and right boundaries $(m_{\varsigma}^{\kappa}, r_{\varsigma}^{\kappa})$
in case  (iii) \ding{172} of Theorem 4.3.4. Second, $a_1^{\varsigma}$ has those boundaries
$(m_{\varsigma}^{\kappa-1}, r_{\varsigma}^{\kappa-1})$
in case  (iii) \ding{173} of 4.3.4 according to Lemma 4.5.1.

(i)\, If $a^\varsigma_1=\bar b_{\bar p\bar q}^\iota$, then
$\bar b^{\iota,0}_{\bar pq_{\bar p}^\iota}=(1)$,
$\bar b^{\iota,0}_{\bar pq}=\emptyset,q_{\bar p}^\iota<q<\bar q,$ by assumption (i), Formula (4.2-8) shows
$$\begin{array}{ll}{\boldsymbol\dz}(\bar b^\iota_{\bar p\bar q})&=
1\bar w_{q^{\iota}_{\bar p} \bar q}
-\sum_{j\geqslant\iota,q>\bar p}w_{\bar pq}\bar b^{j,0}_{q\bar q},
\end{array}$$
where $w_{\bar pp}$ is $\bar w$ or $w$-class.
Since $\bar w_{q^{\iota}_{\bar p} \bar q}^\iota$ is a dotted arrow still by (i),
$\bar b^\iota_{\bar p\bar q}\mapsto\emptyset$ and $\bar w_{q^{\iota}_{\bar p} \bar q}$
is replaced by a linear combination of some $\bar w$-class elements below the $q^{\iota}_{\bar p}$-th row
and some $w$-class elements in the pair $(\mf A^{\varsigma+1},\mf B^{\varsigma+1})$.

(ii)\, If $a^\varsigma_1=b^\iota_{\tau\bar p\bar q},\tau>n_0$, then
$\bar b^{\iota,0}_{\bar pq_{\bar p}^\iota}=(1)$,
$\bar b^{\iota,0}_{\bar pq}=\emptyset,\forall\, q_{\bar p}^\iota<q<\bar q$, Formula (4.2-8) gives
$$\begin{array}{l}
{\boldsymbol\dz}(b_{\tau\bar p\bar q}^\iota)=\bar u_{\tau\bar p\bar q}^\iota
+\sum_{n_0<i<\tau}(\beta^0_{\tau i}\bar u_{i\bar p\bar q}^\iota
+\beta^1_{\tau i}\bar b^{\iota,0}_{\bar pq^\iota_{\bar p}}\bar u_{iq^\iota_{\bar p}\bar q}^{\iota_1})
+\sum_{c_i\prec b_\tau,q}c^0_{i\bar pq}(\sum_{l}\varepsilon'_{\tau il}
\bar v_{lq\bar q}),\end{array}$$
where $\iota\leqslant\iota_1<\gamma, q^\iota_{\bar p}< \bar p$.
Since $\bar u_{\tau\bar p\bar q}^\iota$ is a dotted arrow
by assumption (ii), $b_{\tau\bar p\bar q}^\iota\mapsto\emptyset$,  and $\bar u_{\tau\bar p\bar q}^\iota$
is replaced by a $\bar u$-class element and some $\bar v$-class
elements by Remark 4.1.3 (ii)--(iii).

(iii) If $a^\varsigma_1=c^\iota_{\tau\bar p\bar q}$, then according to Formula (4.2-9):
$$\begin{array}{ll}{\boldsymbol\dz}(c_{\tau\bar p\bar q}^\iota)&=
\sum_{h<i\leqslant h+\tau,\iota^1<\iota}(\gamma^0_{\tau i}\underline v_{i\bar p\bar q}^\iota
+\gamma_{\tau i}^1\bar b^{\iota,0}_{\bar pq^\iota_{\bar p}}\underline v_{iq^\iota_{\bar p}\bar q}^{\iota_1})\\[0.5mm]
&+\sum_{i,q}c^{\iota,0}_{i\bar pq}(\sum_{l}\xi_{\tau il}\underline u_{lq\bar q})
+\sum_{q<\bar q}c^{\iota,0}_{\tau\bar pq} w_{q\bar q}-\sum_{p>\bar p} w_{\bar pp}c^{j,0}_{\tau\bar pq},
\end{array}$$
where $\iota\leqslant\iota_1<\gamma, q^\iota_{\bar p}< \bar p$.
If $c^\iota_{\tau\bar p\bar q}\mapsto\emptyset$, a linear relation among some $\underline u,w$-class elements
and some $\underline v$-class elements with the subscripts being bigger than $h$
is added.

The required $\underline v,\bar u,\bar w$-class and all the
$\bar v$-class dotted arrows are preserved, and
the pair $(\mf A^{\varsigma+1},\mf B^{\varsigma+1})$ still satisfies assumption
(i)-(iii).\hfill$\Box$

\medskip

Suppose $j\leqslant\gamma$ if $\gamma$ exists, otherwise
$j\leqslant\kappa$ in case  (iii) \ding{172}
of Theorem 4.3.4, or $j<\kappa$ in (iii) \ding{173}, see the first thick line above in Picture (4.5-1) for example.

Assume $I^j_\varsigma$ intersects the $p$-th row of $F^\varsigma$ at the $q_p^j$-th column
in the $n^\varsigma_{_X}\times n^\varsigma_{t(a_{i^j})}$-block
obtained from $\bar a_{i^j}$ partitioned under $\bar\T$ with $\bar a^{j,0}_{i^jpq_p^j}=(1)$.
Denote by $\hat v_{i^jq_{p}^{j} q}$ the $(q_{p}^{j}, q)$-element in
the block of size $n^\varsigma_{t(a_{i^j})}\times n^\varsigma_{_X}$ splitting from $\hat v_{i^j}$.
Denote by $\hat V^j_\varsigma$ the block in the $n_0$-th block-column of $\bar\Psi^{r}_{\n^\varsigma}$
partitioned under $\bar\T$, such that the row indices of $\hat V^j_\varsigma$ coincide
with the column indices of $I^j_\varsigma$, see $(I^1,\hat V_1),(I^2,\hat V_2)$ in Picture (4.5-2) below.

\medskip

{\bf Lemma 4.5.3}\, Let $(\bar{\mf A}^\varsigma,\bar{\mf B}^\varsigma)$ be an induced pair of
$(\bar{\mf A}^\kappa,\bar{\mf B}^\kappa)$. Suppose
the first arrow of $\bar{\mf B}^\varsigma$, $a^\varsigma_1=\bullet_{\tau\bar p\bar q}^\iota\in
A^\iota_\varsigma\cup B^\iota_\varsigma\cup C^\iota_\varsigma$, where $\iota\leqslant\gamma$ if $\gamma$ exists, otherwise
$\iota\leqslant\kappa$ in case  (iii) \ding{172} of Theorem 4.3.4, or $\iota<\kappa$ in (iii) \ding{173}. Assume that

(i)\, all $a^{j,0}_{ipq}=\emptyset, i\in\Lambda$,
the corresponding $\underline v_{ipq}^j$ is replaced by a linear
combination of some $\underline u$-class elements in
$\bar{\mf B}^\varsigma$; while all the dotted arrows
$\underline v_{i'p'q'}^{j'}$ are preserved;

(ii)\, if $\bar a^{j,0}_{ipq}=\emptyset$, there is a linear relation among
$\underline u,\bar w,w$-class elements in $\bar{\mf B}^{\varsigma}$; while all the
dotted arrows $\underline v_{i'p'q'}^{j'}, i'\in\Lambda$, are preserved;

(iii)\, all $b^{j,0}_{ipq}=\emptyset, i<n_0$, the corresponding $\bar u_{ipq}^j$ is replaced by a
linear combination of some $\bar v$-class elements in $\bar{\mf B}^\varsigma$; while all the dotted
arrows $\bar u_{i'p'q'}^{j'}$ are preserved;

(iv)\, all $\bar b^{j,0}_{pq}=\emptyset$, $\hat v_{i^jq_{p}^{j} q}$ corresponding to
$\bar a^{j,0}_{i^jpq_p^j}$ is replaced by a linear combination of some $\bar v$-class elements
below and some $\bar w, w$-class in $\bar{\mf B}^\varsigma$;
while the dotted arrows $\hat v_{i'p'q'}, p'\leqslant q_{\bar p}^{\iota}$, are preserved;

(v)\, all $b^{j,0}_{ipq}=\emptyset, i>n_0$, the corresponding element $\bar u_{ipq}^j$
is replaced by a linear combination of some $\bar v$-class elements in
$\bar{\mf B}^\varsigma$; while all the dotted arrows $\bar u_{i'p'q'}^{j'}$ are preserved;

(vi)\, if $c^{j,0}_{ipq}=\emptyset$, there is a linear relation among some elements
$\underline v_{i_1p_1q}^{j_1}, h<i_1<h+\tau,p_1=p$,
and some $\underline u,w,\bar w$-class elements in $\bar{\mf B}^\varsigma$; while all the dotted arrows
$\underline v_{i'p'q'}^{j'}, i'\in\Lambda$, are preserved.

Then after a regularization, the induced pair $(\mf A^{\varsigma+1},\mf B^{\varsigma+1})$
still satisfies (i)-(vi). In particular, all the dotted arrows $\underline v_{i'p'q'}^{j'}, i'\in\Lambda$,
except $\underline v_{\tau\bar p\bar q}^{j}$ in cases (i), (ii), (vi); $\bar u^{j'}_{i'p'q'}$ except $\bar u^{j}_{\tau\bar p\bar q}$
in case (iii) or (v) and $\hat v^{j'}_{p'q'}$ for $p'<q_{\bar p}^\iota$ in case  (iv), are preserved.

\smallskip

{\bf Proof}\, We claim first, that if $\gamma$ exists, than
$\bar u_{i_{1}p_1q}^{j_1}\succ\bar u_{ipq}^{j}$ given in case  (ii) of Lemma 4.5.2
can be replaced inductively by some $\bar v$-class arrows, when the reductions inside the
$(\gamma+1)$-th block in Picture (4.5-1) are finished. Thus the assumption (i)--(vi) are valid, if $a_1^\varsigma$ has
the bottom and right boundaries $(m_\varsigma^\gamma,r_{\varsigma}^\gamma)$,
when $\gamma$ exists, according to Lemma 4.5.2; otherwise $a^\varsigma_1$ has those
$(m_\varsigma^\kappa,r_{\varsigma}^\kappa)$ in
case (iii) \ding{172} of Theorem 4.3.4; or $(m_\varsigma^{\kappa-1},r_{\varsigma}^{\kappa-1})$
in (iii) \ding{173} by Lemma 4.5.1.

(i)\, If $a^\varsigma_1=a^\iota_{\tau\bar p\bar q}, \tau\in\Lambda$,
$\underline v_{\tau\bar p\bar q}$ is a dotted arrow by assumption (i). Formula (4.2-7) tells
$$\begin{array}{c}{\boldsymbol\dz}(a^\iota_{\tau\bar p\bar q})=\underline v^\iota_{\tau\bar p\bar q}+
\sum_{i,q}\bar a^{\iota,0}_{i\bar pq}(\sum_{l}\epsilon_{\tau il}\underline u_{lq\bar q})\,\,\Longrightarrow
\,\,a_{\tau\bar p\bar q}^\iota\mapsto\emptyset,\,\, v_{\tau\bar p\bar q}^\iota=
-\sum_{i,q}\bar a^{\iota,0}_{i\bar pq}
(\sum_{l}\epsilon_{\tau il}\underline u_{lq\bar q}).\end{array}$$

(ii)\, If $a^\varsigma_1=\bar a^\iota_{\tau pq}$ is effective, then by
substituting $\underline v^\iota_{i'\bar p\bar q}$ given by the formula above,
$$\begin{array}{c}{\boldsymbol\dz}(\bar a_{\tau\bar p\bar q}^\iota)=\sum_{\bar a_{\tau^\iota}\prec a_{i'}
\prec\bar a_{\tau}}\bar\alpha_{\tau i'}\underline v^\iota_{i'\bar p\bar q}
+\sum_{i,q}\bar a^{\iota,0}_{i\bar pq}(\sum_{l}\bar\epsilon_{\tau il}\underline u_{lq\bar q})
+\sum_{q<\bar q}\bar a^{\iota,0}_{\tau\bar pq} w_{q\bar q}
-\sum_{j\geqslant\iota,p>\bar p} w_{\bar pp}\bar a^{j,0}_{\tau\bar pq}\\[0.5mm]
=\sum_{i,q}\bar a^{\iota,0}_{i\bar pq}\big(\sum_{l}(\bar\epsilon_{\tau il}-
\sum_{\bar a_{\tau^\iota}\prec a_{i'}\prec\bar a_i\preccurlyeq\bar a_{\tau}}
\bar\alpha_{\tau i'}\epsilon_{i'il})\big)\underline u_{lq\bar q}
+\sum_{q<\bar q}\bar a^{\iota,0}_{\tau\bar pq} w_{q\bar q}-\sum_{j\geqslant\iota,p>\bar p}
w_{\bar pp}\bar a^{j,0}_{\tau\bar pq}.\end{array}$$
If $\bar a^\iota_{\tau\bar p\bar q}\mapsto\emptyset$, then
a linear relation among some $\underline u$, $w$ and $\bar w$-class elements is added.

(iii)\, If $a^\varsigma_1=b^\iota_{\tau\bar p\bar q},\tau<n_0$,
$\bar u_{\tau\bar p\bar q}^\iota$ is a dotted arrow by assumption (iii):
$$\begin{array}{c}{\boldsymbol\dz}(b^\iota_{\tau\bar p\bar q})=\bar u^\iota_{\tau\bar p\bar q}+
\sum_{i,q} \bar a^{\iota,0}_{i\bar pq}(\sum_{l}\varepsilon_{\tau il}\bar v_{lq\bar q})\,\,\Longrightarrow\,\,
b^\iota_{\tau\bar p\bar q}\mapsto\emptyset,\,\, \bar u^\iota_{\tau\bar p\bar q}=-
\sum_{i,q}\bar a^{\iota,0}_{i\bar pq}(\sum_{l}\varepsilon_{\tau il}\bar v_{lq\bar q}).\end{array}$$

(iv)\, If $a^\varsigma_1=\bar b^\iota_{\bar p\bar q}$ is effective,
by substituting $\bar u^\iota_{i'\bar p\bar q}$ given in (iii) and Formula (4.4-1),
$$\begin{array}{ll}\bar{\boldsymbol\dz}(\bar b^\iota_{\bar p\bar q})&=\sum_{i'<n_0}
\bar\beta_{i'} \bar u_{i'\bar p\bar q}^\iota+\sum_{i,q}\bar a^{\iota,0}_{i\bar pq}
(\sum_{l}\bar\varepsilon_{il}\bar v_{lq\bar q})+\sum_{q<\bar q}\bar b^{\iota,0}_{\bar pq}
w_{q\bar q}-\sum_{j\geqslant\iota,p>\bar p} w_{\bar pp}\bar b^{j,0}_{p\bar q}\\[0.5mm]
&=\sum_{i,q}\bar a_{i\bar pq}^{\iota,0}\big(\sum_{l}
(\bar\varepsilon_{il}-\sum_{\bar a_{\tau^\iota}
\preccurlyeq\bar a_{i'}\prec b_i\prec\bar b}\bar\beta_{i'}\varepsilon_{i'il})\bar v_{lq\bar q}\big)
+\sum_{q<\bar q}\bar b^{\iota,0}_{\bar pq}
w_{q\bar q}-\sum_{j\geqslant\iota,p>\bar p} w_{\bar pp}\bar b^{j,0}_{p\bar q}\\[0.5mm]
&=\hat v_{\tau^\iota q^\iota_{\bar p}\bar q}+\sum_{\bar a_{i\bar pq}^{\iota,0}
\succ\bar a^{\iota,0}_{\tau^\iota\bar p q^\iota_{\bar p}};q}\bar a_{i\bar pq}^{\iota,0}
\hat v_{iq\bar q}+\sum_{q<\bar q}\bar b^{\iota,0}_{\bar pq}
w_{q\bar q}-\sum_{j\geqslant\iota,p>\bar p} w_{\bar pp}\bar b^{j,0}_{p\bar q}.\end{array}$$
Since $\bar a^{\iota,0}_{\tau^\iota\bar p q^\iota_{\bar p}}=1$, and
$\hat v_{\tau^\iota q^\iota_{\bar p}\bar q}$ is a dotted arrow by assumption (iv),
$\bar b^\iota_{\bar p\bar q}\mapsto\emptyset$, and
$\hat v_{\tau^\iota q^\iota_{\bar p}\bar q}$ is replaced by some $\hat v$-class elements
below the $q_{\bar p}^\iota$-th row, and some
$\bar w,w$-class elements.

(v)\, If $a^\varsigma_1=b^\iota_{\tau\bar p\bar q}, \tau>n_0$,
since $\bar b_{\bar pq}^0=\emptyset$ for all possible $q$ by (iv) above, Formula (4.2-8) shows
$$\begin{array}{l}{\boldsymbol\dz}(b^\iota_{\tau\bar p\bar q})=\bar u_{\tau\bar p\bar q}^\iota
+\sum_{i\ne n_0,i<\tau}\beta^0_{\tau i}\bar u^\iota_{i\bar p\bar q}
+\sum_{i,q}\bar a^{\iota,0}_{i\bar pq}
(\sum_{l}\varepsilon_{\tau il}\bar v_{lq\bar q})+\sum_{c_i\prec b_\tau;q}\bar c^{\iota,0}_{i\bar pq}
(\sum_{l}\varepsilon'_{\tau il}\bar v_{lq\bar q}).\end{array}$$
Since $\bar u^\iota_{\tau\bar p\bar q}$ is a dotted arrow by assumption (v),
$b_{\tau\bar p\bar q}^\iota\mapsto\emptyset$, $\bar u^\iota_{\tau\bar p\bar q}$ is replaced
by some $\bar v$-class elements according to the replacement given in (iii) and assumption (v).

(vi)\, If $a_1^\varsigma=c_{\tau\bar p\bar q}^\iota$, since $\bar b_{\bar pq}^0=\emptyset$
for all possible $q$ by (iv), Formula (4.2-9) gives
$$\begin{array}{ll}{\boldsymbol\dz}(c^\iota_{\tau\bar p\bar q})&=\sum_{i\leqslant h+\tau}\gamma^0_{\tau i}
\underline v_{i\bar p\bar q}^\iota+\sum_{i,q}\bar a^{\iota,0}_{i\bar pq}
(\sum_{l}\zeta_{\tau il}\underline u_{lq\bar q})\\[0.5mm]
&+\sum_{i<\tau;q}\bar c^{\iota,0}_{i\bar pq}(\sum_{l}\xi_{\tau il}\underline u_{lq\bar q})
+\sum_{q<\bar q}c^{\iota,0}_{\tau\bar pq} w_{q\bar q}
-\sum_{j\geqslant\iota,p>\bar p} w_{\bar pp}c^{j,0}_{\tau p\bar q}.\end{array}$$
If ${\boldsymbol\dz}(c^\iota_{\tau\bar p\bar q})\ne 0$, then $\bar{\mf B}^{\varsigma+1}$ is given by
$c^\iota_{\tau\bar p\bar q}\mapsto\emptyset$, and a linear relation is added among
$\underline v^\iota_{i\bar p\bar q}, h<i\leqslant h+\tau$ and some $\underline u,w,\bar w$-class elements,
because $\underline v_{i\bar p\bar q}^\iota$ for $i\leqslant h$ have already been replaced
by some $\underline u$-class arrows given in (i).

The required $\underline v,\bar u,\bar w,\bar v$-class dotted arrows are preserved,
the pair $(\bar{\mf A}^{\varsigma+1},\bar{\mf B}^{\varsigma+1})$ still satisfies assumption
(i)-(vi).\hfill$\Box$

\medskip

The following picture shows a pseudo formal equation $\bar\Theta^\varsigma$ of $(\bar{\mf
A}^\varsigma,\bar{\mf B}^\varsigma)$ for $\kappa=5, \gamma=2$ in case  (iii) \ding{173} of
Theorem 4.3.4 only  with effective arrows. From this, it is possible to see the correspondences
of $(\bar B^1_\varsigma,\hat V^1_\varsigma)$, $(\bar B^2_\varsigma,\hat V^2_\varsigma)$,
$(\bar B^3_\varsigma,\bar W^3_\varsigma)$, $(\bar B^4_\varsigma,\bar W^4_\varsigma)$ respectively.

\medskip

\vspace{7mm}
\begin{equation*}\begin{array}{c} {\unitlength=0.8mm
\begin{picture}(185,60)\thicklines
\put(0,16){\framebox(30,30){}} \put(0,40){\line(1,0){30}}
\put(6,34){\line(1,0){24}} \put(12,28){\line(1,0){18}}
\put(18,22){\line(1,0){12}}

\put(6,34){\line(0,1){12}} \put(12,28){\line(0,1){18}}
\put(18,22){\line(0,1){24}} \put(24,16){\line(0,1){30}}

\put(0,42){\makebox{\small$e^\varsigma_{_{Y^1}}$}} \put(6,36){\makebox{\small$e^\varsigma_{_{Y^2}}$}}
\put(12,30){\makebox{\small$e^\varsigma_{_{Y^3}}$}} \put(18,24){\makebox{\small$e^\varsigma_{_{Y^4}}$}}
\put(24,18){\makebox{\small$e^\varsigma_{_{X^5}}$}} \put(22,32){\makebox{$\bar W$}}

\thicklines \put(40,16){\framebox(70,30){}}
\put(52,40){\line(0,1){6}} \put(58,34){\line(0,1){6}}
\put(70,28){\line(0,1){6}} \put(82,22){\line(0,1){6}}
\put(94,16){\line(0,1){6}} \put(52,40){\line(1,0){6}}
\put(58,34){\line(1,0){12}}\put(70,28){\line(1,0){12}}
\put(82,22){\line(1,0){12}}

\thinlines \put(52,16){\line(0,1){30}} \put(64,16){\line(0,1){30}}
\put(94,16){\line(0,1){30}} \put(58,16){\line(0,1){30}}
\put(46,16){\line(0,1){30}}

\put(40,40){\line(1,0){70}} \put(52,34){\line(1,0){58}}
\put(64,28){\line(1,0){46}} \put(76,22){\line(1,0){34}}

\put(70,16){\line(0,1){18}} \put(76,16){\line(0,1){12}}
\put(82,16){\line(0,1){12}} \put(88,16){\line(0,1){6}}

\put(47,42){\makebox{$I^1$}}\put(42,42){\makebox{$0$}}
\put(53,36){\makebox{$I^2$}} \put(65,30){\makebox{$I^3$}}
\put(77,24){\makebox{$I^4$}}
\put(88,17.5){\makebox{$\scriptstyle{W}$}}
\put(56,41){\makebox{\small$\bar A^1_2$}} \put(58,35){\makebox{\small$\bar A_{2}^2$}}
\put(78,42){\makebox{\small$\bar B^1$}}
\put(78,36){\makebox{\small$\bar B^2$}}
\put(78,30){\makebox{\small$\bar B^3$}}
\put(87,24){\makebox{\small$\bar B^4$}} \put(95,41){\makebox{\small$C_{1}^1$}}
\put(95,35){\makebox{\small$C_{1}^2$}} \put(95,29){\makebox{\small$C_{1}^3$}}
\put(95,23){\makebox{\small$C_{1}^4$}} \put(95,17){\makebox{\small$C_{1}^5$}}
\put(103,41){\makebox{\small$C_{2}^1$}} \put(103,35){\makebox{\small$C_{2}^2$}}
\put(103,29){\makebox{\small$C_{2}^3$}} \put(103,23){\makebox{\small$C_{2}^4$}}
\put(102,17){\makebox{\small$C_{2}^5$}}

\mput(64,41.5)(0,1.5){4}{\line(1,0){30}}
\mput(65.5,40)(1.5,0){20}{\line(0,-1){6}}
\mput(70,29.5)(0,1.5){3}{\line(1,0){24}}
\mput(83.5,28)(1.5,0){8}{\line(0,-1){6}}
\put(102,16){\line(0,1){30}}

\thicklines \put(120,0){\framebox(70,70){}}
\mput(120,58)(0,6){2}{\line(1,0){70}} \put(132,52){\line(1,0){58}}
\put(132,46){\line(1,0){58}} \put(144,16){\line(1,0){46}}
\put(174,8){\line(1,0){16}} \put(132,70){\line(0,-1){24}}
\put(144,70){\line(0,-1){54}} \put(174,70){\line(0,-1){62}}
\put(182,70){\line(0,-1){70}}

\put(144,40){\line(1,0){30}} \put(150,34){\line(1,0){24}}
\put(156,28){\line(1,0){18}} \put(162,22){\line(1,0){12}}

\put(150,34){\line(0,1){12}} \put(156,28){\line(0,1){18}}
\put(162,22){\line(0,1){24}} \put(168,16){\line(0,1){30}}

\thinlines\mput(144,59.5)(0,1.5){4}{\line(1,0){30}}
\mput(145.5,52)(1.5,0){20}{\line(0,1){6}}
\mput(150,41.5)(0,1.5){4}{\line(1,0){24}}
\mput(163.5,28)(1.5,0){7}{\line(0,1){6}}

\put(126,70){\line(0,-1){12}} \put(138,58){\line(0,-1){12}}

\put(144,41){\makebox{\small$e^\varsigma_{_{Y^1}}$}} \put(150,35){\makebox{\small$e^\varsigma_{_{Y^2}}$}}
\put(156,29){\makebox{\small$e^\varsigma_{_{Y^3}}$}} \put(162,23){\makebox{\small$e^\varsigma_{_{Y^4}}$}}
\put(168,17){\makebox{\small$e^\varsigma_{_{X^5}}$}} \put(159,40){\makebox{\small$\bar W^3$}}
\put(165,28){\makebox{\small$\bar W^4$}}

\put(123,62){\makebox{$T^\varsigma_1$}} \put(135,50){\makebox{$T^\varsigma_2$}}
\put(176,10){\makebox{$T^\varsigma_3$}} \put(183,1){\makebox{$T^\varsigma_4$}}
\put(155.5,59){\makebox{\small$\hat V^1$}} \put(155.5,47.5){\makebox{$\,$}}
\put(156,53){\makebox{\small$\hat V^2$}} \put(134.5,62){\makebox{$\underline U$}}
\put(179.5,53){\makebox{$\underline U$}} \put(179.5,27){\makebox{$\underline V$}}
\put(185,10){\makebox{\small$\underline U$}}

\put(63,14){\makebox{$\underbrace{\hskip 25mm}$}}
\put(77,7){\makebox{$\,$}} \put(42,9){\makebox{$\,$}}
\put(55,9){\makebox{$\,$}} \put(95,9){\makebox{$\,$}}
\put(105,9){\makebox{$\,$}} \put(13,9){\makebox{$\,$}}
\put(153,-5){\makebox{$\,$}}
\put(-0.5,46.5){\makebox{$\overbrace{\hskip 24mm}$}}
\put(52,46.5){\makebox{$\overbrace{\hskip 9mm}$}}
\put(64,46.5){\makebox{$\overbrace{\hskip 24mm}$}}
\put(94.5,46.5){\makebox{$\overbrace{\hskip 12.5mm}$}}
\put(119.5,70.5){\makebox{$\overbrace{\hskip 56mm}$}}

\put(13,50){\makebox{$e_{_X}^\varsigma$}} \put(56,50){\makebox{$A^\varsigma$}}
\put(77,51){\makebox{$B^\varsigma$}} \put(101,50){\makebox{$C^\varsigma$}}
\put(153,74){\makebox{$T^\varsigma$}}

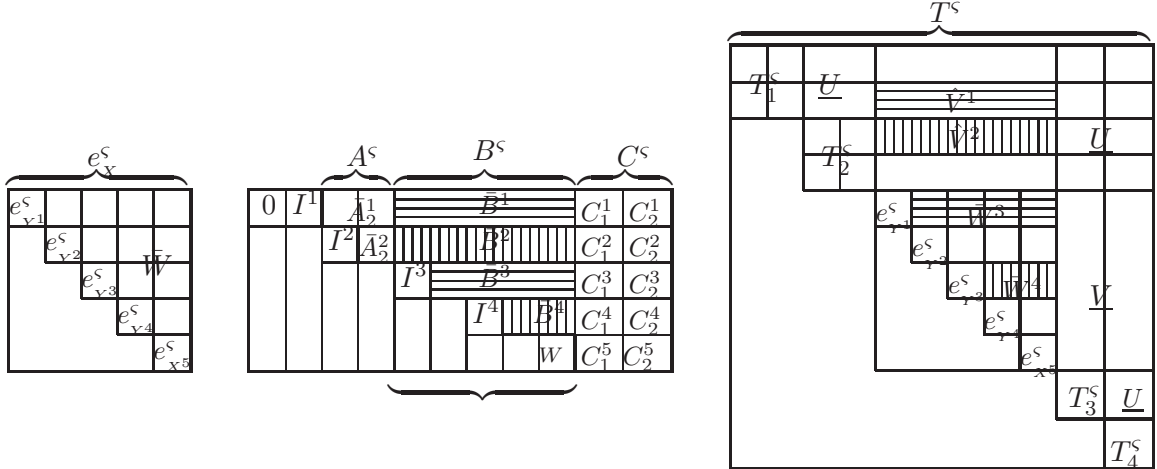
\put(77,-4){\makebox{$\mbox{Picture}\,(4.5\mbox{-}2)$}}
\end{picture}}
\end{array}
\end{equation*}

\medskip
\bigskip
\bigskip
\noindent {\large\bf 4.6 Regularizations on non-effective $a$ class and all $b$ class arrows}
\bigskip

Let $(\bar{\mf A},\bar{\mf B})$ be a one-sided pair, and the induced pair $(\bar{\mf A}^\kappa,\bar{\mf B}^\kappa)$
be given by Theorem 4.3.4. Using the notation of Remark 3.4.6, we may assume that
an induced local pair $(\bar{\mf A}^s,\bar{\mf A}^s)$ of $(\bar{\mf A}^\kappa,\bar{\mf B}^\kappa)$
is obtained by a sequence of reductions in the sense of Lemma 2.3.2.
Set $a^s_1\mapsto(x)$, the induced pair $(\bar{\mf A}^{s+1},\bar{\mf A}^{s+1})$ is obtained with $\bar R=k[x]$,
and $(\bar{\mf A}^t,\bar{\mf A}^t)$ in the case of MW5 is given by a series of regularizations.
We will show in the last subsection, that $x,a^t_1$ in $\bar{\mf B}^t$ can be only split from some
edges of $\bar{\mf B}$.

\medskip

It is clear by Lemma 4.5.1--4.5.3 that the
non-effective $a$-class and all the $b$-class solid arrows are regularized during the reductions.
Note that the conclusions of Lemmas 4.5.1--4.5.3 are still valid, if we deal with the linear relations
over the fractional field $k(x)$ of the polynomial ring $k[x]$, or over the field $k(x,x_1)$
of two indeterminants instead of the base field $k$.  Then the non-effective $a$-class
and all the $b$-class solid arrows are regularized,
which implies the following theorem.

\medskip

{\bf Theorem 4.6.1}\, Let $(\bar{\mf A},\bar{\mf B})$ be a one-sided pair
having at least two vertices, such that the induced local bocs $\bar{\mf B}_X$ is given by Formula (4.2-1),
the pair is major, and the $c$-class arrows satisfy Formula (4.2-6).
If $(\bar{\mf A}^\kappa,\bar{\mf B}^\kappa)$ given by Theorem 4.3.4 has an induced pair
$(\bar{\mf A}^t,\bar{\mf B}^t)$ in the case of MW5 defined by Remark 3.4.6, then the parameter
$x$ and the first arrow $a^t_1$ of $\bar{\mf B}^t$ belong to $\bar a$ or $c$-class.

\medskip

Finally, let $(\bar{\mf A},\bar{\mf B})$ be a one-sided pair
having at least two vertices, such that $\bar{\mf B}_X$ is in case (i) of Classification 4.2.1.
Then $\bar{\mf B}$ has only $a,b$-class solid arrows, where $b_1,\cdots,b_n$
are all non-effective, and Formulae (4.2-5) is also suitable for $a$-class arrows:
$a_i$ for $i\in\Lambda$ satisfying the first formula of
(4.2-5) are non-effective, while $\bar a_i=a_{h_i}$ for $1\leqslant i\leqslant s$ satisfying the second
one are effective. If there is an induced pair $(\bar{\mf A}^t,\bar{\mf B}^t)$
in the case of MW5 according to Remark 3.4.6, there are following observations.

1)\, Let $(\bar{\mf A},\bar{\mf B})$ be a pair with $\T$ being trivial.
The condition (BRC)$'$ is constructed parallel to (BRC) in Condition 4.3.1 as follows.

\subitem\hspace{-6mm}(i)\, Let $\mathcal D=\{d_1,\cdots,d_q\}$ be a
set of solid arrows, and $\mathcal E=\{e_1,\cdots,e_p\}$ be another set of solid edges without any loop,
such that $\mathcal D\cup\mathcal E$ forms the lowest non-zere row of the formal product $\Theta$.
And let $\mathcal U=\{u_1,\cdots,u_q\}$ be a set of dotted arrows, while $\mathcal W=\V\setminus\U$.

\subitem\hspace{-6mm}(ii)\, $\bar{\boldsymbol\dz}(d_i)$ and $\bar{\boldsymbol\dz}(e_i)$ satisfy the formulae of 4.3.1 (ii).

\noindent Then after a reduction given by Cases (i)--(iv) stated before Lemma 4.3.2, the induced pair still
satisfies (BRC)$'$. On the other hand, the original pair $(\bar{\mf A},\bar{\mf B})$ satisfies (BRC)$'$
parallel to Lemma 4.3.3. The proofs of the two facts are much easier than those of above two lemmas.

2)\, For constructing a reduction sequence from $(\bar{\mf A},\bar{\mf B})$ up to
$(\bar{\mf A}^\kappa,\bar{\mf B}^\kappa)$, what is needed is only the part (i) and (iii) \ding{172} of Theorem 4.3.4.
In fact, the reduction block $G^j$ is ${1}\choose{0}$ or ${0\,1}\choose{0\,0}$
for $j<\kappa$, and $G^\kappa=(1)$ or $(0\,1)$. 

3)\, For further reductions, what is needed is only Theorem 4.5.3 (i)-(iii),
then an induced pair $(\bar{\mf A}^\varsigma,\bar{\mf B}^\varsigma)$ is reached, where all the $b$-class,
and non-effective $a$-class arrows are regularized step by step.

\medskip

{\bf Corollary 4.6.2}\, Let $(\bar{\mf A},\bar{\mf B})$ be a one-sided pair
having at least two vertices, such that the induced local bocs $\bar{\mf B}_X$ is in case  (i) of
Classification 4.3.1. If $(\bar{\mf A}^\kappa,\bar{\mf B}^\kappa)$ given in 2) above has an induced pair $(\bar{\mf A}^t,\bar{\mf B}^t)$
in the case of MW5 defined by Remark 3.4.6, then the parameter $x$ and the first arrow $a^t_1$ of $\bar{\mf B}^t$
belong to $\bar a$-class.

\bigskip

\begin{center} {\Large\bf 5. Non-homogeneity of bipartite matrix\\[1mm] bimodule problems of wild type}\end{center}


This section is devoted to proving the non-homogeneous property for a wild
bipartite matrix bimodule problem satisfying RDCC condition in the case of MW5.
As a consequence, the main theorem 3 is proved in the last subsection.

\bigskip
\bigskip
\noindent{\large\bf 5.1  An inspiring example}
\bigskip

The purpose of this subsection is two folds: 1) classify the positions
of the first arrows of bocses in the case of MW5 at formal products; 2) define a notion of the
bordered matrix of a matrix, then prove a preliminary lemma on both matrices.

\medskip

Let $\mf A=(R,\K,\M,H=0)$ be a bipartite matrix bimodule problem, which has a trivial vertex set $\T=\T_{(1)}\dot{\cup}\T_{(2)}$
and satisfies RBDCC condition, see Remark 1.4.4. Suppose
$\mf A'=(R',\K',\M',$ $H')$ is an induced matrix bimodule problem in the case of MW5 defined by Remark
3.4.6. We show the classification of the position of the first arrow $a'_1$ in the sum $H'+\Theta'$:
$$\begin{array}{c}H'+\Theta'=H'+ \sum_{i=1}^{n'}
{a'}_i\ast{A'}_i.\end{array}\eqno{(5.1\mbox{-}1)}$$
Denote by $(p,q')$ the leading position of $A'_1$ over $\T'$, which
locates in the $(\textsf p,\textsf q)$-th leading block of some base
matrix of $\M_1$ partitioned under $\T$. By the bipartite property and RDCC condition, $\textsf q$ is the index of a main block
column, say $\textsf q=\textsf q_Z$ for some vertex $Z\in\T_{(2)}$.

\medskip

{\bf Classification 5.1.1}\, Suppose the bocs $\mf B'$ is in the case of MW5. There are two possible
relations between the row index $p$ and the row indices of the links of $H'$ defined below
Formula (2.3-7) in Formula (5.1-1):

{\bf case (I)} $p <$ the row indices of all the links in the
$(\textsf p,\textsf q)$-block of $H'$;

\smallskip

{\bf case (II)} $p \geqslant$ some row index of at least one link in the
$(\textsf p,\textsf q)$-block of $H'$.

\smallskip

It is clear that there is no link above the $(\textsf p,\textsf
q)$-th block, since $\mf A'$ is already local.

\medskip

{\bf Lemma 5.1.2} Let $p_x$ be the row index of $x$ in $H'$, then
$p_x>p$ in Classification 5.1.1.

\smallskip

{\bf Proof.}  Since $x$ appears before the first arrow $a_1'$ in $\mf B'$,
$p_x\geqslant p$ by the order of reductions according to matrix indices.
If $p_x=p$, then the parameter $x$ locates at the left side of $a_1'$
in $H'+\Theta'$, $\delta(a_1')$ contains only the terms of the form
$\alpha xv,\alpha\in k$, which contradicts
to the assumption that $\mf{B}'$ is in the case of MW5. Thus $p_x>p$.\hfill$\Box$

\medskip

{\bf Example 5.1.3}\, Let $(\mf A,\mf B)$ be a pair constructed by an algebra defined in Example 1.4.5.
There is a reduction sequence $(\mf A,\mf B),(\mf A^1,\mf B^1),(\mf A^2,\mf B^2),(\mf A^3,\mf B^3)$ given in
Examples 2.4.5, such that the bocs $\mf B^3$ is strongly
homogeneous in the case of MW5 described in Remark 3.1.7 (iii).
In order to prove that $(\mf A,\mf B)$ is not homogeneous,
another way different from the proof of MW1--MW4 must be found. More
precisely, we will reconstruct a new reduction sequence based on the matrix
$\tilde M$ over $k[x]$ with the size vector $\tilde{\m}=(2,2,2,2,2,3,3,3,3,3)$:
$$\tilde M
=\small\left(
\begin{array}{c|c|c|c|c}
\quad0\quad& \quad0\quad& \begin{array}{ccc}0& 1 &0\\
0&0&1\end{array} &
\begin{array}{ccc} 0&0& 1 \\ 0&0&0\end{array}&
0\\ \hline
& 0&\begin{array}{ccc}  0&0& 1 \\ 0&0&0\end{array}& 0&
\begin{array}{ccc} 0& 0&0 \\ 0&x&0\end{array}\\ \hline
&& 0&0& \begin{array}{ccc}0& 0&1 \\ 0&0&0\end{array}\\ \hline
&&&0&\begin{array}{ccc} 0&1 &0\\ 0&0&1\end{array}\\ \hline
&&&\begin{array}{ccc}&&\\&&\end{array}&0
\end{array}\right).
$$
There is a reduction sequence
$(\tilde{\mf{A}},\tilde{\mf B}),(\tilde{\mf{A}}^1,\tilde{\mf{B}}^1),(\tilde{\mf{A}}^2,\tilde{\mf{B}}^2),
(\tilde{\mf{A}}^{3},\tilde{\mf{B}}^3)$ corresponding to the steps (i)--(iii) of Example 2.4.5,
where the reduction from $\tilde{\mf B}$ to $\tilde{\mf B}^1$ is given by $a\mapsto (0\,1)$
in the sense of Lemma 2.3.2. Thus $b$ splits into $b_1, b_2$ in $\tilde{\mf B}^1$, and set
$b_1\mapsto(0), b_2\mapsto {{0\,1}\choose{0\,0}}$ from $\tilde{\mf{B}}^1$ to $\tilde{\mf{B}}^2$.
$\tilde{\mf B}^3$ is obtained from  $\tilde{\mf B}^2$
by an edge reduction $(0)$, followed by a loop mutation, then four regularizations:
$$\begin{array}{c}
\tilde{H}^3=\left(\begin{array}{ccc}0&1_{X}& 0\\
0&0&1_X\end{array}\right)\ast A +\left(\begin{array}{ccc}0&0& 1_X\\
0&0&0\end{array}\right)\ast B+\left(\begin{array}{ccc}\emptyset&\emptyset& \emptyset\\
0&x1_X&\emptyset\end{array}\right)\ast C.\end{array}$$
The $(1,5)$-th block partitioned under $\T$ in the formal equation of $(\tilde{\mf
A}^3,\tilde{\mf B}^3)$ is of the form:
$$\begin{array}{c}\tiny\left(\begin{array}{cc} e& v\\
0&e\end{array}\right)\left(\begin{array}{ccc} d_{10}&d_{11}& d_{12}\\
d_{20}&d_{21}&d_{22}\end{array}\right)+\left(\begin{array}{cc} u^1_{11}& u^1_{12}\\
u^1_{21}&u^1_{22}\end{array}\right)\left(\begin{array}{ccc}0& 0& 1\\
0&0&0\end{array}\right)+\left(\begin{array}{cc} u^2_{11}& u^2_{12}\\
xv&u^2_{22}\end{array}\right)\left(\begin{array}{ccc}0& 1& 0\\
0&0&1\end{array}\right)=\\[2mm]
\tiny\left(\begin{array}{ccc}0& 1&0 \\
0&0&1\end{array}\right)\left(\begin{array}{ccc}
v^{_2}_{_{00}}& v^{_2}_{_{01}}& v^{_2}_{_{02}}\\
v^{_2}_{_{10}}&v^{_2}_{_{11}}&v^{_2}_{_{12}}\\ 0&vx&v^{_2}_{_{11}}\end{array}\right)+
\left(\begin{array}{ccc}0& 0& 1\\
0&0&0\end{array}\right)\left(\begin{array}{ccc}
v^{_1}_{_{00}}& v^{_1}_{_{01}}& v^{_1}_{_{02}}\\
v^{_1}_{_{10}}& v^{_1}_{_{11}}& v^{_1}_{_{12}}\\
v^{_1}_{_{20}}&v^{_1}_{_{21}}&v^{_1}_{_{22}}\end{array}\right)
+\left(\begin{array}{ccc}d_{_{10}}& d_{_{11}}& d_{_{12}}\\
d_{_{20}}&d_{_{21}}&d_{_{22}}\end{array}\right)\left(\begin{array}{ccc}s_{00}&s_{01}&s_{02} \\0&e& v\\
0&0&e\end{array}\right)\end{array}$$
with $e=e_{_{X}}, s_{00}=e_{_{Y}}$. The
differentials of the solid arrows of $\tilde{\mf B}^3$ can be read off as follows:
$$d_{20}: X\mapsto Y,\,\,\delta(d_{20})=0;\quad
d_{21}:X\mapsto X,\,\,\dz(d_{21})=xv-vx-d_{20}s_{01},$$
and $\bar\dz(d_{22})=u_{21}^1-v^2_{11}-d_{20}s_{02}-d_{21}v,
\bar\dz(d_{10})=-v_{10}^2-v_{20}^1+vd_{20}, \bar\dz(d_{11})=u_{11}^2-v_{11}^2-v_{21}^1
-d_{10}s_{01}, \bar\dz(d_{12})=u_{11}^1+u_{11}^2-v_{12}^1-v_{22}^1$, where $\bar\dz$ is obtained
from $\dz$ by removing the monomials, which involve a solid arrow $d_{22},d_{10}$ or $d_{11}$.
It is clear that the bocs $\wt{\mf B}^3$ satisfies the hyperthesis of Proposition 3.4.5. In fact,
as $d_{20}\mapsto(1)$, the solid loops
$d_{21},d_{22},d_{10},d_{11},d_{12}$ will be regularized step by step, because
$s_{01},s_{02},v_{10}^2,u_{11}^2,u_{11}^1$ are pairwise
different dotted arrows.
Therefore $(\mf A,\mf B)$ is not homogeneous.

\medskip

Motivated by Example 5.1.3, the general cases are considered.
Since the example satisfies Case (I) of Classification 5.1.1,
we start from Case (I) in subsection 5.1--5.3.

Let $\mf{A}=(R,\K,\M,H=0)$ be a bipartite matrix bimodule problem
satisfying RDCC condition. Let $\mf{A}'$ be an induced matrix
bimodule problem with trivial $R'$, and let
$\vartheta: R(\mf A')\rightarrow R(\mf A)$ be the induced
functor. Suppose $ M= \vartheta (H'(k))=\sum_j{M}_j\ast A_j\in R(\mf A)$
with a size vector $\l\times \n$ over $\T$. Let $\textsf q=\textsf q_Z\in T_2$ for some $Z\in \T_2$.
Define a size vector $\l\times\tilde{\n}$ over $\T$, and construct a {\it bordered matrix}
$\tilde{M}=\sum_j \tilde{M}_j\ast A_j\in R(\mf{A})$ with $0$ a zero column as follows:
$$\tilde{n}_j=\left\{\begin{array}{ll}n_j, &\mbox{if}\,\,\,
j\not\in Z,\\n_j+1,& \mbox{if}\,\,\, j\in Z.\end{array}\right.\quad
\tilde{M}_j=
\left\{\begin{array}{ll}M_j,& \mbox{if}\,\,\, A_j1_Z=0,\\(0\,M_j),&
\mbox{if}\,\,\, A_j1_Z=A_j.\end{array}\right.\eqno(5.1\mbox{-}2)$$

Denote by $(p,q+1)$ the leading
position of $M$, such that $q+1$ is the index of the
first column of the $\textsf q$-th block-column of M partitioned under $\T$. Denote by
$\tilde q$ the added column index of $\tilde M$, the column is sitting in the
$\textsf{q}_Z$-th block-column as the first one. Applying Theorem 2.4.1, the defining
system $\IE$ of $\K'_0\oplus\K'_1$ given by Formula (2.4-2), and a matrix
equation $\tilde{\IE}$ are considered:
$${\IE}:\,\,
\Phi^1_{\l}M\equiv_{_{\prec(p, q+1)}}
M\Phi^2_{\n},\quad
\tilde{\IE}:\,\, \Phi^1_{\l}\tilde M\equiv_{_{\prec({p},
\tilde{q})}}\tilde{M}\Phi^2_{\wt{\n}},
\eqno{(5.1\mbox{-}3)}$$
where the upper scripts $1,2$ on $\Phi$ stand for the left and right parts of the
bipartite variable matrix $\Phi$, the two sets of variables in two parts do not intersect.
Since $\mf A$ satisfies RDCC condition, the main block-column
in $\Phi^2_{\n}$ determined by $Z\in \T$ can be written as $\Phi^2_{\n,Z}=
(\Phi_{1}^2,\cdots,\Phi_{n}^2)^T$, such that either $\Phi_{l}^2=0$ or
$\Phi^2_{l}=(z^l_{pq})\ne 0$, where $z^l_{pq}$ are variables over $k$. 

It is clear that the $\tilde q$-th column of $\Phi^1_{\l}\tilde M$ is a zero column, we may
define two new matrix equations respectively:
$$ {\IE}_\tau:\,\,0\equiv_{_{\prec(p, q+1)}}
M\Phi^2_{\n},\quad \tilde{\IE}_\tau:\,\,
0\equiv_{_{\prec({p},
\tilde{q})}}\tilde{M}\Phi^2_{\wt{\n}}.\eqno{(5.1\mbox{-}4)}$$
Taken any integer $p'\geqslant p$ and $1\leqslant j\leqslant
n_{_Z}$, the $(p',q+j)$-th entry of the right side of $\IE_\tau$ is:
$$\begin{array}{c}\sum_{\Phi_{\n,Z}^2\ne
0}\sum_{q'}\alpha^l_{p'q'}z_{q',q+j}^l.
\end{array}\eqno{(5.1\mbox{-}5)}$$ It is easy to see that
$z_{q',q+j}^l$ for all possible $j$ have the same coefficient
$\alpha_{p',q'}^l$, the $(p',q')$-th entry of $H(k)$.
In the picture below, $n_{_Z}=5, p'=p$, those five equations are indicated by
five circles, and the five variables at the same row of $\Phi^2_{\n,Z}$ have the same coefficients.

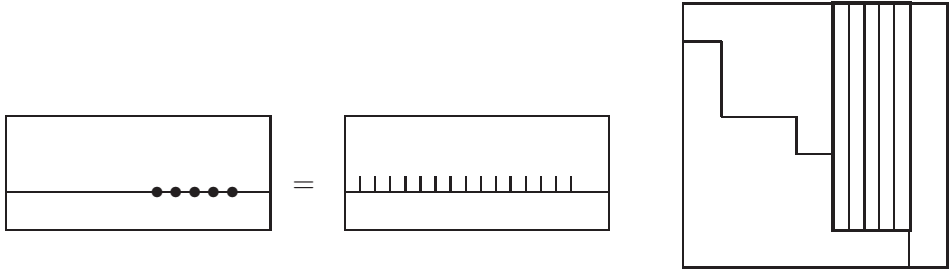
\begin{center}
\setlength{\unitlength}{1mm}
\begin{picture}(145,38)\thinlines{\linethickness{0.1mm}}
\put(30,15){\circle*{1.5}}\put(32.5,15){\circle*{1.5}}
\put(35,15){\circle*{1.5}}
\put(37.5,15){\circle*{1.5}}\put(40,15){\circle*{1.5}}
\put(10,10){\framebox(35,15)}\put(10,15){\line(1,0){35}}
\put(48, 15){$=$}

\put(55,10){\framebox(35,15)}\put(55,15){\line(1,0){35}}
\mput(55,15)(2,0){16}{\line(0,1){2}}

\put(100,5){\framebox(35,35)}
\put(100,35){\line(1,0){5}}
\put(105,25){\line(1,0){10}}\put(105,25){\line(0,1){10}}
\put(115,20){\line(1,0){5}}\put(115,20){\line(0,1){5}}
\mput(122,10)(2,0){4}{\line(0,1){30}}
\put(130,5){\line(0,1){5}}

\thicklines{\linethickness{1mm}}
\put(120,10){\framebox(10,30)}
\end{picture}
\end{center}

{\bf Lemma 5.1.4}\, With the notations as above.

(i)\, The $(p,q+j_1)$-th equation is a linear combination of
the previous equations in $\IE_\tau$ if and only if so is the
$(p,q+j_2)$-th equation. Similarly, the same result is valid in $\tilde\IE_\tau$.

(ii)\, The equations in the system $\tilde{\IE}$ (resp. $\tilde{\IE}_{\tau}$) and
those in ${\IE}$ (resp. ${\IE}_{\tau}$) are the same at the same positions of each
main block column, whenever the added $\tilde q$-th column has been dropped from
$\tilde{\IE}$ (resp. $\tilde{\IE}_{\tau}$).

(iii)\, A subsystem of $\tilde\IE$ (resp. $\tilde{\IE}_{\tau}$)
consisting of the $\tilde q$-column in both sides is
$\tilde{\IE}_{\tau\tilde q}: 0\equiv \tilde M\Phi_{\tilde\n,\tilde q}$,
where $\Phi_{\tilde\n,\tilde q}$ stands for the
$\tilde q$-th column of $\Phi_{\tilde\n}$. And $\tilde{\IE}_{\tau\tilde q}$ can be solved independently.

(iv)\, If the $(p,q+j)$-th equation is a linear combination of the previous
equations of $\IE$, then so is the $p$-th equation of $\tilde\IE_{\tau\tilde q}$.

\smallskip

{\bf Proof}\, (i) The assertion follows from Formula (5.1-5).

(ii) Recall that $\textsf q\in Z$, denote by $0$ the index of the first column
(row) in the $\textsf q'$-th block column (row) for any $\textsf q'\in Z$.
For any $X\in\T_{(1)},Y\in\T_{(2)}$, set
$1\leqslant\alpha\leqslant n_{_{X}}$ and $1\leqslant\beta\leqslant n_{_{Y}}$. We claim that
the $(\alpha,\beta)$-th equations in the $(\textsf h,\textsf q_{_Y})$-th block of $\IE$ and $\tilde \IE$
for any $\textsf h\in X$ are the same. In fact, the variable matrix $\Phi^1_{\l}$
in the two systems $\IE$ and $\tilde{\IE}$ is common;
and the $\beta$-th column of the $\textsf h$-th block row in $M$ and $\tilde M$ are the same
in the left side of two equations. Now consider the right side of two equations. Let $(M_1,\cdots,M_t)$
be the $\alpha$-row of the $\textsf h$-th block row in $M$ with
$M_j=(\lambda_{j1},\cdots,\lambda_{jn_j})$ and that in $\tilde M$
is $(\tilde M_1,\cdots,\tilde M_t)$. Then $\tilde M_j=M_j,\forall j\notin Z$; but
$\tilde M_j=(0,\lambda_{j1},\cdots,\lambda_{j n_j}),\forall j\in Z$. Let $(\Phi_1,\cdots,\Phi_t)^T$
be the $\beta$-column of the $\textsf q_{_Y}$-th block column in $\Phi^2_{\n}$ with
$\Phi_j=(x_{j1},\cdots,x_{jn_j})$ and that in $\tilde\Phi^2_{\tilde\n}$
is $(\tilde\Phi_1,\cdots,\tilde\Phi_t)^T$, then $\tilde\Phi_j=\Phi_j,\forall j\notin Z$; but
$\Phi_j=(x_{j0},x_{j1},\cdots,x_{j n_j})^T,\forall j\in Z$.
Thus the additional variables $x_{j0},\forall\,j\in Z$,
are killed by $0$ in the right side of the equation of $\tilde \IE$.

(iii) The $\tilde q$-column in the left side of the matrix equation $\tilde\IE$ is a zero column.
The variables of $\Phi_{\tilde\n,\tilde q}$  are different from those in $\Phi_{\l}$ and in the main block columns
of $\Phi_{\tilde\n}$ except the $\tilde q$-th column.

(iv) If there exists some $1\leqslant j\leqslant n_{_Z}$, such
that the $(p,q+j)$-th equation is a linear
combination of the previous equations in $\IE$, then so is the $(p,q+j)$-th equation in $\IE_\tau$ after deleting
$\Phi_{\l}$, since the variables of $\Phi_{\l}$ and $\Phi_{\n}$ are different.
Thus so is the $(p,\tilde q+j)$-th equation in $\tilde\IE_\tau$ by (ii),
and so is the $(p,\tilde q)$-th equation in $\tilde\IE_\tau$ by (i), finally, so is the $p$-th
equation in $\tilde\IE_{\tau\tilde q}$.

\bigskip
\bigskip
\noindent{\large\bf 5.2 Bordered matrices in bipartite case}
\bigskip

This subsection is devoted to constructing a reduction sequence
based on a given sequence and a bordered matrix, which generalizes Example 5.1.3.

\medskip

Let $\mf{A}=(R,\K,\M,H=0)$ be a bipartite matrix bimodule problem, which has a trivial $\T$
and satisfies RDCC condition. Let $\mf{A}^s=({R}^s,\K^s,\M^s,H^s)$ be an induced matrix bimodule problem with $R^s$
being local and trivial. Then there is a unique sequence of reductions
in the sense of Lemma 2.3.2 by Corollary 2.3.5:
$$\begin{array}{c}\mf{A}, \mf{A}^1, \cdots,  \mf{A}^i,\mf{A}^{i+1}, \cdots,
\mf{A}^s.\end{array} \eqno {({\ast})}$$
Write $M=\vartheta^{0s} (H^s(k))=\sum_{j} {M}_j\ast A_i\in
R(\mf A)$ with a size vector $\l\times \n$, where $M_j=G^{j+1}_s(k)$ are given by Formula (2.3-7).
Suppose the size vector $\l\times\tilde{\n}$ and the representation $\tilde{M}=\sum_{j}
\tilde{M}_j\ast A_j\in R(\mf{A})$ are defined by Formula (5.1-2).

\medskip

{\bf Theorem 5.2.1}\, There exists a unique reduction sequence
based on the sequence $(\ast)$:
$$\begin{array}{c}\mf{A},
\tilde{\mf{A}}^1, \cdots, \tilde{\mf{A}}^{i},
\tilde{\mf{A}}^{i+1}, \cdots, \tilde{\mf{A}}^s\end{array} \eqno
{(\tilde\ast)}
$$
where $\tilde{\mf{A}}^i=(\tilde{R}^i, \tilde{\K}^i,
\tilde{\M}^i,\tilde{H}^i)$, the reduction from $\tilde{\mf{A}}^{i}$
to $\tilde{\mf{A}}^{i+1}$ is a reduction or a composition of two
reductions in the sense of Lemma 2.3.2. Moreover, $\tilde\T^s$ has two
vertices, and $$\tilde{\vartheta}^{0s}(\tilde{H}^s(k))=\tilde{M}.$$

{\bf Proof}\,  We may assume that $\l\times \n$ is
sincere over $\T$. Otherwise, it is possible to
obtain an induced problem $\mf{A}'$ by a suitable deletion, such that $M$ has a sincere size vector
over $\mf{A}'$. In particular, $\mf A'$ is still bipartite and satisfies RDCC condition.

We will construct a sequence $(\tilde\ast)$ inductively. The original
term in the sequence is $\tilde{\mf{A}}=\mf{A}$. Suppose that
a sequence $\tilde{\mf{A}},\tilde{\mf{A}}^1,\cdots, \tilde{\mf{A}}^{i}$ for some $0\leqslant i<s$ has been constructed and $\tilde
\vartheta^{0i}: R(\tilde{\mf{A}}^{i})\mapsto R(\tilde{\mf A})$ is the induced functor, such
that there exists a representation
$$\begin{array}{c}\tilde M^{i}=\tilde H^{i}_{\tilde
\m^{i}}(k)+\sum_{j=1}^{n^i} \tilde M^{i}_j\ast \tilde A^{i}_j\in R(\tilde{\mf{A}}^{i}),\quad
\mbox{with}\quad\tilde\vartheta^{0i}(\tilde{M}^{i})\simeq \tilde M\in
R(\tilde{\mf A}^0).\end{array}\eqno{(5.2\mbox{-}1)}$$
Write $M^i=\vartheta^{is}({H^s}(k))=H^{i}_{\m^{i}}(k)+\sum_{j=1}^{n^{i}}M_j^{i}\ast
A_j^{i}\in R(\mf{A}^{i})$, where $M_1^{i}=G^{i+1}_s(k)$
by Formula (2.3-7), and is denoted by $B$ for simplicity. The
first column in the $\textsf q$-th main block-column of $M$ under the partition $\T$
is denoted by $\beta$. Now we are constructing $\tilde{\mf A}^{i+1}$.

\smallskip

{\bf Case 1}\, $\tilde{\T}^{i}=\T^{i}$ and $\tilde{\mf A}^i=\mf A^i$.

\smallskip

{\bf 1.1} $B\cap \beta$ is empty. Then
$\tilde{G}^{i+1}=G^{i+1}$, $\tilde{H}^{i+1}=H^{i+1}$ and
$\tilde{\mf{A}}^{i+1}={\mf A}^{i+1}$.

\smallskip

Before giving the following cases, we claim that if $B\cap \beta$ is non-empty, $B$ thus
$G^{i+1}$ can not be Weyr matrices. Otherwise, the first arrow $a_1^{i}$ of $\mf
B^{i}$ will be a loop. Since $\tilde{\T}^{i}=\T^{i}$, $\tilde{a}_1^{i}$
will also be a loop and hence the numbers of rows and columns of
$\tilde{B}$ are the same. When the matrix $B$ is enlarged by one
column, then $B$ is also enlarged by one row, which is a contradiction to the construction
of $\tilde M$. So $B$ is either $(\emptyset)$ from a regularization
or ${{0\,\, I_r}\choose{0\,\,\,0}}$ from an edge reduction.

\smallskip

{\bf 1.2} $B\cap \beta$ is non-empty, and $B=(\emptyset)$. Then $\tilde{B}=(\emptyset \, B)$
with $\emptyset$ being a distinguished zero column, $\tilde{H}^{i+1}=H^{i+1}$ and $\tilde{\mf{A}}^{i+1}=\mf{A}^{i+1}$
by a regularization.

\smallskip

{\bf 1.3} $B\cap \beta$ is non-empty, $B={{0\,\, I_r}\choose
{0\,\,\,0}}$ and $r<$ the number of columns of $B$. Then
$\tilde{B}=(0\, B)$ with $0$ being a zero column, $\tilde{H}^{i+1}=H^{i+1}$ and
$\tilde{\mf{A}}^{i+1}=\mf{A}^{i+1}$ by an edge reduction.

\smallskip

{\bf 1.4} $B\cap \beta$ is non-empty, $B={{I_r}\choose 0}$. Then
$\tilde{B}=(0\, B)$ with $0$ being a
zero column. Recall from Formula (2.3-5) and Theorem 2.3.3, the following is defined:
$$
\tilde{G}^{i+1}=\left\{\begin{array}{ll}  \left(\begin{array}{cc} 0&
1_{Z_2}\end{array}\right),&
\mbox{if } G^{i+1}=(1_{Z_2});\\[1.5mm]
\left(\begin{array}{cc} 0& 1_{Z_2}\\
0&0\end{array}\right),& \mbox{if } G^{i+1}=
\left(\begin{array}{c}  1_{Z_2}\\
0\end{array}\right).\end{array}\right.$$
Then $\tilde{H}^{i+1}=\sum_{X\in \T^{i}}I_{X}\ast H^{i}_{X}+\tilde{G}^{i+1}\ast A_1^{i}$.
Consequently $\tilde{\mf{A}}^{i+1}$ is induced from $\tilde{\mf A}^i$ by an edge
reduction in the sense of Lemma 2.3.2.

We stress, that after the edge reduction in the subcase 1.4, $\tilde{\T}^{i+1}=\T^{i+1}\cup
\{Y\}$, where $Y$ is an equivalent class consisting of the indices of the added columns
in the sum $\tilde H^{i+1}+\tilde\Theta^{i+1}$ of the pair $(\tilde{\mf A}^{i+1},\tilde{\mf B}^{i+1})$,
and $(\tilde{\mf{A}}^{i+1},\tilde{\mf B}^{i+1})\ne (\mf A^{i+1},\mf B^{i+1})$ from this stage.
The above $B$ is show in four cases as a small block in the corresponding leading block partitioned under $\T$ in $\tilde M$:

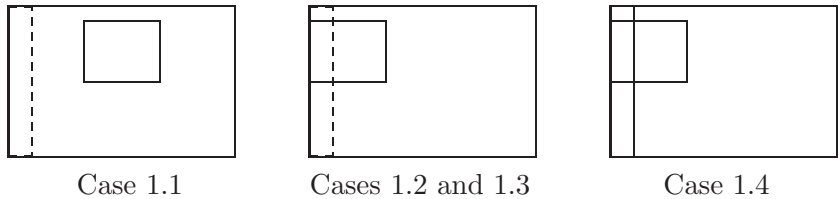
\bcen \unitlength 1mm
\begin{picture}(110,20)

\put(0,0){\framebox(30,20)[l]{\dashbox{1}(3,20)}}
\put(10,10){\framebox(10,8)}

\put(40,0){\framebox(30,20)[l]{\dashbox{1}(3,20)}}
\put(40,10){\framebox(10,8)}

\put(80,0){\framebox(30,20)} \put(83, 0){\line(0,1){20}}
\put(80,10){\framebox(10,8)}


\put(9, -5){Case 1.1} \put(40, -5){Cases 1.2 and 1.3} \put(87,
-5){Case 1.4}
\end{picture}

\end{center}
\vskip 3mm

\smallskip

{\bf Case 2.} $\tilde{\T}^{i}=\T^{i}\cup \{Y\}$.

\smallskip

{\bf 2.1} $B\cap \beta$ is empty. Then $\tilde{B}=B$,
$\tilde{G}^{i+1}=G^{i+1}$, and $\tilde{H}^{i+1}=\sum_{\tilde{X}\in
\tilde\T^{i}}\tilde{I}_{\tilde{X}}\ast\tilde H^{i}_{\tilde
{X}}+\tilde{G}^{i+1}\ast\tilde{A}_1$.

\smallskip

If $B\cap \beta$ is non-empty. Denote by $\tilde a^{i}_0$ and $\tilde
a^{i}_1$ the first and the second solid arrows of $\tilde {\mf B}^{i}$, which
locate at $(p^{i}, \tilde q_0^{i})$ and $(p^{i},\tilde q^{i}_0+1)$
in the formal product $\tilde\Theta^i$ respectively.

\smallskip

{\bf 2.2} $B\cap \beta$ is non-empty, and there exists some
$1\leqslant j\leqslant n^i_Z$, such that the
$(p^i,q^i+j)$-th equation is a linear combination of previous equations in
$\IE_\tau^i$. Then $\delta(\tilde{a}_0^{i})=0$ by Lemma 5.1.4 (ii) then (i), and Corollary 2.4.2.
Two reductions are made: the first one is an
edge reduction by $\tilde{a}_0^{i}\mapsto (0)$; and the second
one for $\tilde{a}_1^i$ is made in the same way as that for $a_1^i$ by Lemma 5.1.4 (ii).
Then an induced problem $\tilde{\mf{A}}^{i+1}$ is obtain, and $\tilde{B}=(0 \, B)$ with
$0$ being a zero column.

\smallskip

{\bf 2.3} $B\cap \beta$ is non-empty, and for all $1\leqslant
j\leqslant m^i_Z$, the $(p^i,q^i+j)$-th equation
is not a linear combination of previous equations in
$\IE^i_\tau$. Thus $\dz(\tilde a^i_0)\ne 0$ by Lemma 5.1.4 (i) and Corollary 2.4.2.
And $\dz^0(\tilde a^i_{j})\ne 0$ for any $1\leq j\leq n^i_Z$ by 5.1.4 (i)--(ii) and Corollary 2.4.2.
Then two regularizations $\tilde{a}_0^{i}\mapsto (\emptyset),\; \tilde{a}_1^{i}\mapsto (\emptyset)$
are made, and $\tilde B=(\emptyset\, B)$ with $\emptyset$ being a distinguished zero
column.

\smallskip

In the cases 2.2 and 2.3, there are two reduction blocks $\tilde{G}^{i+1,0}=(0)$ or $(\emptyset)$,
$\tilde{G}^{i+1,1}=G^{i+1}$, thus
$\tilde{H}^{i+1}=\sum_{\tilde X\in\tilde{\T}^{i}}
\tilde{I}_{\tilde X}\ast \tilde {H}^{i}_{{\tilde
X}}+\tilde{G}^{i+1,0}\ast\tilde{A}_0^{i}
+\tilde{G}^{i+1,1}\ast\tilde{A}_1^{i}.$

By summing up all the cases, an induced pair
$(\tilde{\mf{A}}^{i+1},\tilde{\mf{B}}^{i+1})$ and a representation
$\tilde M^{i+1}$ with $\tilde\vartheta^{0,i+1}(\tilde M^{i+1})\simeq\tilde M$ are obtained.
The theorem follows by induction.\hfill$\Box$

\medskip

{\bf Corollary 5.2.2} The main diagonal block $\tilde{e}^i_Z,Z\in\T,$ of $\tilde{\K}^i_0\oplus\tilde{\K}^i_1$
is of the form:
$$
\left(\begin{array}{ccccc} s_{00}& s_{01}& s_{02}&\cdots& s_{0m}\\
&s_{11} & s_{12}& \cdots & s_{1m}\\
&&s_{22}& \cdots &s_{2m}\\
&&&\ddots &\vdots\\
&&&&s_{mm}
\end{array}\right).
$$
where $m=n^i_Z$, $s_{01}, s_{02},\ldots,s_{0m}$ are dotted arrows of $\tilde{\mf B}^i$.

\smallskip

{\bf Proof}\, By the construction of $\tilde{H}^i$, the added
``$0$-column'' can be only $0$ or $\emptyset$. Therefore, except
$s_{00}$, the elements at the $0$-th row: $s_{01}, s_{02},
\ldots,s_{0m}$ do not appear in any equation of the defining system of
$\tilde{\mf{A}}^i$. Thus they are free. \hfill$\Box$

\bigskip
\bigskip
\noindent{\large\bf 5.3   Non-homogeneity in the case of MW5 and classification (I)}
\bigskip

The discussion of this subsection is two folds: 1) extend the reduction sequence $(\tilde\ast)$ of Theorem 5.2.1
into a sequence $(\tilde\ast')$, such that there is a parameter $x$ appearing from the $(s+1)$-th step;
2) prove that any bipartite pair with an induced minimal wild pair in the case of MW5 and
Classification 5.1.1 (I) is not homogeneous.

\medskip

Suppose we have a reduction sequence ending at $\mf A^t$ defined by Remark 3.4.6:
$$\mf{A}, \mf{A}^1,\cdots,\mf{A}^s, \mf{A}^{s+1},
\cdots,\mf{A}^\epsilon,\cdots,\mf A^t=\mf A', \eqno{(\ast')}$$
where the reduction from $\mf A^{i}$ to $\mf A^{i+1}$ is in the sense of
Lemma 2.3.2 for $1\leqslant i<s$; $\mf A^{s}$ is local with $\dz(a_1^s)=0$,
set $a^s_1\mapsto (x)$, $\mf A^{s+1}$ has a parameter $x$ locating
at the $(p_x, q_x)$-position of $H^{s+1}$ and $R^{s+1}=k[x]$; $\mf A^{i+1}$ is
obtained from $\mf A^{i}$ by a regularization for
$s<i<t$. The pair $(\mf A^t,\mf B^t)$ is in the case of MW5 given by Remark 3.4.6 and
satisfying Classification 5.1.1 (I).

Note that the set of integers $\T^i$ and its partition $\T^i$ are all the same for $i=s,\cdots,t$.
Suppose the first arrow $a^t_1$ of $\mf B^t$ locates at the $(p,q')$-th position
in the formal product $\Theta^t$ with $q'=q+j$ for
some $1\leqslant j\leqslant n^t_Z$; the first arrow
$a^\epsilon_1$ of ${\mf B}^\epsilon$ locates at the $(p,q+1)$-th position in $\Theta^\epsilon$,
where $q+1$ is the index of the
first column in the $\textsf q$-th block-column.
The picture below shows the position of the first solid arrows in
the formal products $\Theta^i$ of $({\mf A}^i,{\mf B}^i)$ for $i=s,
\epsilon, t$ (whenever the added $\tilde q_0$-th column is ignored):

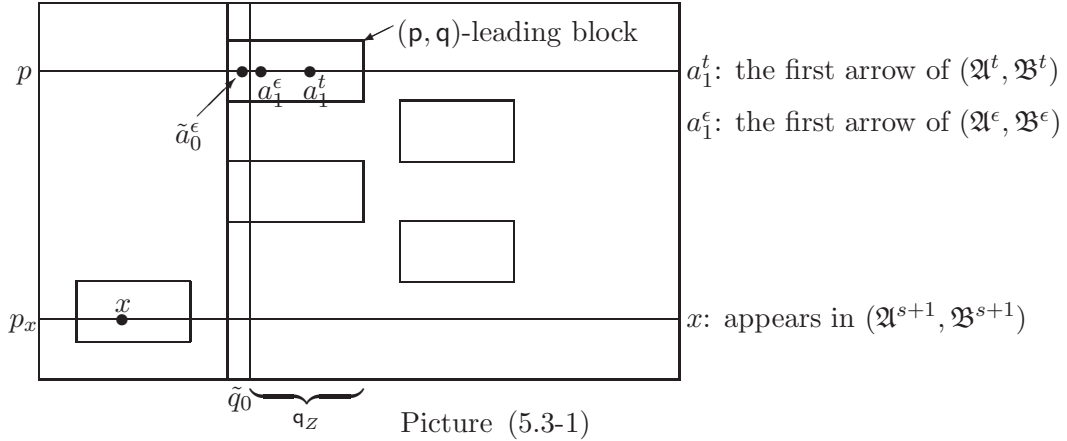
\begin{center}
\unitlength=1mm
\begin{picture}(120,50)
\put(0,0){\framebox(85,50){}} \put(5,5){\framebox(15,8){}}
\put(25,21){\framebox(18,8){}} \put(25,37){\framebox(18,8){}}
\put(48,13){\framebox(15,8){}} \put(48,29){\framebox(15,8){}}
\put(25,0){\line(0,1){50}} \put(28,0){\line(0,1){50}}
\put(0,8){\line(1,0){85}} \put(0,41){\line(1,0){85}}
\put(-4,7){$p_x$} \put(-3,40){$p$} \put(10,7){$\bullet$}
\put(26,40){$\bullet$}\put(35,40){$\bullet$}
\put(10,9){$x$}\put(28.5,40){$\bullet$}
\put(18,31.5){$\tilde{a}_0^\epsilon$}
\put(35,37.3){$a^t_1$}\put(29,37.3){${a}_1^\epsilon$}
\put(86,7){$x$: appears in $(\mf{A}^{s+1}, \mf{B}^{s+1})$}
\put(86,40){$a^t_1$: the first arrow of $(\mf{A}^t,
\mf{B}^t)$} \put(86,33){${a}_1^\epsilon$: the first arrow of
$(\mf{A}^\epsilon, \mf{B}^\epsilon)$} \put(21,35){\vector(1,1){5}}
\put(25,-3.5){$\tilde{q}_0$} \put(47,47){\vector(-2,-1){4}}
\put(47.5,45){$(\textsf p,\textsf q)$-leading block}
\put(28,-1){\makebox{$\underbrace{\hskip 15mm}_{\textsf q_Z}$}}
\put(48,-7){$\mbox{Picture}\,\,\,(5.3\mbox{-}1)$}
\end{picture}\vskip 9mm
\end{center}
Assume $R^i=k[x,\phi^i(x)^{-1}]$, and $H^i$ has
size vector $\l^i\times\n^i$ over $\T$ for $s<i\leqslant t$. Denote uniformly the same size vector $\l^i$
by $\l$ over $\T_{(1)}$ with size $l$, and $\n^i$ by $\n$ over $\T_{(2)}$ with size $n$.
Then $H^i\in\IM_{\l\times \n}(R^i)$ are all in the same
form but with different $\phi^i(x)$. Since $k(x)$ is an $R^i$-bimodule,
$$H^i\otimes_{R^i}1_{k(x)}\in \IM_{\l\times\n}(R^i)
\otimes_{R^i}k(x)\simeq\IM_{l\times n}(k(x)).$$

{\bf Remark 5.3.1}\, Lemma 5.1.4 (i)--(iv) are still valid if the matrices
$M^i$ and $\tilde M^i$ over $k(x)$ instead of over $k$ are considered.

\medskip

Recall the matrix equation defined in 5.1.4 (iii):
$$\begin{array}{c}
\tilde{\IE}_{\tau\tilde q}: 0
\equiv \tilde M\Phi_{\tilde\n,\tilde q}.\end{array} \eqno{(5.3\mbox{-}2)}$$

Suppose $x$ locates in the $(\textsf p_x,\textsf q_{Z'})$-th
main block partitioned under $\T$, see Picture (5.3-1) for$Z'\ne Z$; and Example 5.1.3 for $Z'=Z$.
Thus the equation system $\tilde\IE_{\tau\tilde q}^{(>p_x)}$, consisting of the
equations below the $p_x$-th row, is over the base field $k$.
Denote by $\tilde n$ the size of $\tilde\n$. The solution space
of the system $\tilde\IE_{\tau\tilde q}^{(>p_x)}$ contained in $\IM_{\tilde n\times 1}(k)$
is a direct sum of two subspaces by Theorem 5.2.1 and Lemma 5.1.4 (iii):

\subitem\hspace{-6mm} (i) the first subspace
is isomorphic to a space spanned by $E_Y$ in $\tilde\K^s_0$, it has a base matrix
with the $\tilde q$-th entry being $1_Y$ and others zero;

\subitem\hspace{-6mm} (ii) the second one is isomorphic to a subspace of  $\tilde\K^s_1$, and its base
matrices have non-zero entries only above the $\tilde q$-th entry.

\noindent  The second subspace is denoted by $\K_{\tau\tilde q}^{(> p_x)}$ with a polynomial
$d^{p_x}(x)=1$; a minimal algebra $R^{(>p_x)}_{\tau\tilde q}=\tilde R^s=k1_X\times k1_Y$; and a basis $\{U^0_1,\cdots,U^0_{\beta}\}$
with the dual basis $\{u^0_1,\cdots,u^0_{\beta}\}\subset\Hom_{}(\K_{\tau\tilde q}^{(> p_x)},R^{(>p_x)}_{\tau\tilde q})$.
From now on, we consider the equation $\tilde\IF^{(>p_x)}_{\tau\tilde q}:\,0\equiv \tilde M\Psi^{(> p_x)}_{\tau\tilde q}$ with
the variable matrix $\Psi^{(>p_x)}_{\tau\tilde q}=\sum_{\zeta=1}^{\beta}
u^0_\zeta\ast U^0_\zeta$ according to Theorem 2.4.4 and Formula (2.4-5).

Suppose the system $\tilde\IF^{(>h)}_{\tau\tilde q}$ for some $p< h\leqslant p_x$ have been solved, the solution space
$\K_{\tau\tilde q}^{(>h)}$ has a basis $\{U_1,\cdots,U_\kappa\}\subset \IM_{\tilde n\times 1}(R^{(>h)}_{\tau\tilde q})$,
where $R^{(>h)}_{\tau\tilde q}=k[x,(\gamma^{h+1}(x))^{-1}]1_X\times k1_Y$ is a minimal algebra with
$\gamma^{h+1}(x)=\prod_{\eta=p_x}^{h+1}d^{\eta}(x)$.
Let $\{u_1,\cdots,u_\kappa\}\subset\Hom(\K_{\tau\tilde q}^{(>h)},R^{(>h)}_{\tau\tilde q})$ be the dual basis of $\{U_j\}_j$.
According to Formula (2.4-5):
$$\begin{array}{c}\Psi^{(>h)}_{\tau\tilde q}=\sum_{\zeta=1}^{\kappa}
u_\zeta\ast U_\zeta,\quad\quad
\tilde\IF^{(>h)}_{\tau\tilde q}:\,0\equiv \tilde M\Psi^{(> h)}_{\tau\tilde q}.\end{array}\eqno{(5.3\mbox{-}3)}$$
The $h$-th equation of $\tilde\IF^{(>h)}_{\tau\tilde q}$ is $\sum_{\zeta=1}^{\kappa} f_\zeta(x)u_\zeta$
with $f_\zeta(x)\in\,R^{(>h)}_{\tau\tilde q}$. There are two possibilities.

(i)\, $f_\zeta(x)=0$ for $\zeta=1,\cdots,\kappa$. In this case
$\K^{(> h-1)}_{\tau\tilde q}=\K^{(>h)}_{\tau\tilde q}$, and the
quasi-basis of $\K^{(>h)}_{\tau\tilde q}$ are preserved in
$\K^{(> h-1)}_{\tau\tilde q}$. Let $d^{h}(x)=1$.

(ii)\, There exists some $f_{\zeta}(x)\ne 0$. Without loss of generality,
it may be assumed that $f_\kappa(x)\ne 0$. Choose a new basis in the dual space Hom$_{k(x)}(\K^{(>h)}_{\tau\tilde q}
\otimes_{R^{(>h)}_{\tau\tilde q}}k(x),k(x))$
at the first line of the formula below, the corresponding basis
of $\K^{(>h)}_{\tau\tilde q}$ is shown at the second line:
$$\left\{\begin{array}{l}u_\zeta'=u_\zeta,\\[1mm]
U'_\zeta=U_\zeta-f_\zeta(x)/f_\kappa(x)U_\kappa;
\end{array}\right.\,\,\mbox{for}\,\,1\leqslant\zeta<\kappa;\quad\left\{\begin{array}{l}
u'_\kappa=\sum_{\zeta=1}^\kappa f_\zeta(x)u_\zeta,\\[1mm]
U'_\kappa=1/f_\kappa(x)U_\kappa,\end{array}\right.
\eqno{(5.3\mbox{-}4)}$$
where $u_\kappa'=0$ is the solution of the $h$-th equation in the system (5.3-3).
Let $d^{h}(x)\in k[x]$ be the numerator
of $f_{\kappa}(x)$, and $\gamma^h(x)=d^h(x)\gamma^{h+1}(x)$, then $R^{(>h-1)}_{\tau\tilde q}=
k[x,(\gamma^h(x))^{-1}]1_X\times k1_Y$. Thus $\K^{(>h-1)}_{\tau\tilde q}$ has a quasi-basis
$\{U_\zeta'\mid \zeta=1,\cdots,\kappa-1\}$ over $R^{(>h-1)}_{\tau\tilde q}$.
The system $\tilde\IF^{(> p-1)}_{\tau\tilde q}$ with the solution space $\K_{\tau\tilde q}^{(>p-1)}$
and polynomial $\gamma^p(x)$ is finally reached by inverse-order induction.

Suppose $R^i=k[x,\phi^i(x)^{-1}]$, and the row index of the first arrow of $\mf B^i$
in the formal product $\Theta^i$ is $p^i, p_x\leqslant p^i\leqslant p$ for $s\leqslant i\leqslant \epsilon$.
Define
$$\begin{array}{c}\tilde{\phi}^{i}(x)=\phi^i(x)\gamma^{p^i}(x)\in k[x],\quad
\mbox{in particular}\quad \tilde{\phi}^{t}(x)=\phi^t(x)\gamma^{p}(x).
\end{array}\eqno{(5.3\mbox{-}5)}$$
Now we deal with representations of $\mf A^i$ over the field $k(x)$ instead of over $k$.
Suppose the matrix $M^i=\vartheta^{0i}(H^i(k[x,\tilde{\phi}^t(x)^{-1}])=\sum_j{M}_j^i\ast A_j$
has a size vector $\l\times \n$ over $\T$, and a matrix $\tilde M^{i}=\sum_{j}\tilde{M}_j^i\ast A_j$
of size vector $\l\times\tilde{\n}$ is defined by Formula (5.1-2). Returning to Theorem 2.4.1,
the matrix equations for $i\geqslant s$ are defined as follows:
$$\begin{array}{ll}\IE^i:
\Phi_{\l} M^i\equiv_{\prec (p^i,q^i)}
M^i\Phi_{{\n}},&\tilde\IE^i:
\Phi_{\l}\tilde M^i\equiv_{\prec (p^i,q^i)}
\tilde M^i\Phi_{{\tilde\n}};\\ \IE^i_\tau:0\equiv_{\prec(p^i,q^i)}
M^i\Phi_{{\n}},&\tilde\IE^i_\tau:0\equiv_{\prec(p^i,q^i)}
\tilde M^i\Phi_{{\tilde\n}}.\end{array}\eqno{(5.3\mbox{-}6)}$$

{\bf Theorem 5.3.2}\, There exists a unique reduction sequence
based on the sequence $(\ast')$:
$$\begin{array}{c}\mf{A},
\tilde{\mf{A}}^1, \cdots, \tilde{\mf{A}}^{s},
\tilde{\mf{A}}^{s+1},\cdots, \tilde{\mf{A}}^\varepsilon,\cdots,
\tilde{\mf{A}}^\epsilon,\cdots,
\tilde{\mf A}^t=\tilde{\mf A}',\end{array}\eqno{(\tilde\ast')}
$$
where the first part of the sequence up to $\tilde{\mf A}^s$ is given
by Theorem 5.2.1; the reduction from $\tilde{\mf A}^s$ to $\tilde{\mf
A}^{s+1}$ is given by a loop mutation $a^{s+1}_1\mapsto(x)$, or
an edge reduction $(0)$ followed by a loop mutation $(x)$; the
reduction from $\tilde{\mf A}^i$ to $\tilde{\mf A}^{i+1}$ for $s<
i<t$ is given by one regularization, or two regularizations,
or a reduction given by Lemma 2.2.6,  followed by a regularization.

\smallskip

{\bf Proof}\, If the first arrow $a^s_1$ of $\mf B^s$ does not locate at
the $(q+1)$-th column of the formal product $\Theta^{s}$, a loop mutation from $\tilde{\mf{A}}^{s}$ to
$\tilde{\mf{A}}^{s+1}$ is made. Otherwise an edge reduction is made by Remark 5.3.1, see 5.1.4 (iv)
and Corollary 2.4.2 for details, then followed by a loop mutation.

Now suppose we have an induced bimodule problem $\tilde{\mf A}^{i}$ for some $i>s$.
If the first arrow $a^{i}_1$ of ${\mf B}^{i}$ does not locate
at the $(q+1)$-th column of $\Theta^{i}$, a
regularization $\tilde a^{i}_1\mapsto \emptyset$ is made.
Otherwise, there are two possibilities.
\ding{172} There exists some $1\leqslant j\leqslant n_Z$, the $(p^i,q+j)$-th equation
is a linear combination of the previous equations in $\IE^i_\tau$,
then $\dz(\tilde a_0^i)=0$ by Remark 5.3.1 and Corollary 2.4.2. Set $\tilde a^i_0\mapsto(0)$ by Lemma 2.6.6, $\tilde
a^i_1\mapsto \emptyset$. \ding{173} Otherwise $\dz(\tilde a_0^i)\ne0$. Set $\tilde a^i_0\mapsto\emptyset$ and $\tilde a^i_1\mapsto
\emptyset$. The sequence $(\tilde\ast')$ is completed by induction as desired.\hfill$\Box$

\medskip

{\bf Corollary 5.3.3}\, If the bocs $\mf B^t$ in the sequence $(\ast')$ satisfies MW5 defined by Remark 3.4.6
and Classification 5.1.1 (I), then $\dz(\tilde a^\epsilon_{0})=0$ in $\tilde{\mf B}^\epsilon$
in the sequence $(\tilde \ast')$.

\smallskip

{\bf Proof}\, Since the first arrow $a_1^t$ of $\mf B^t$ locates at
the $(p,q+j)$-th position with $\dz(a^t_1)=0$, $\dz^0(a^\epsilon_{j})=0$ in $\mf B^\epsilon$.
Thus $\dz(\tilde a_0^\epsilon)=0$ in $\tilde{\mf B}^\epsilon$ by Remark 5.3.1 and Corollary 2.4.2.\hfill$\Box$

\medskip

{\bf Proposition 5.3.4}\, Let $(\mf A,\mf B)$ be a pair with $\T$ being trivial, such that
$\mf{A}=(R,\K,\M,H=0)$ is a bipartite
matrix bimodule problem satisfying RDCC condition. If there exists an induced
pair $(\mf A',\mf B')$ of $(\mf A,\mf B)$ in the case of MW5 defined by Remark 3.4.6, and the sum
$H'+\Theta'$ of $(\mf A',\mf B')$ satisfies Classification 5.1.1 (I), then $\mf{B}$ is not
homogeneous.

\smallskip

{\bf Proof}\, Suppose we have a sequence $(\ast')$ with $\mf B'=\mf
B^t$, then there is a sequence
$(\tilde\ast')$ based on $(\ast')$ by Theorem 5.3.2. Corollary 5.3.3 tells that the
first arrow $\tilde{a}_0^\epsilon$ of $\tilde{\mf B}^\epsilon$ is an edge
with $\dz(\tilde a_0^\epsilon)=0$, and hence
$\tilde{a}_0^\epsilon\mapsto (1)$ may be set according to Proposition 2.2.7.
The induced pair is obviously local. Thus it is possible to use the triangular formulae
of Subsection 3.3, and an induced pair $({\mf A}'',{\mf B}'')$ is obtain in one
of the cases (ii)-(iv) of Classification 3.3.5.

Case 1\, If 3.3.5 (ii) is met, then $\tilde{\mf
B}^\epsilon$ satisfies the hypothesis of Proposition 3.4.5. It is done.

Case 2\, If 3.3.5 (iii) is met, then $\tilde{\mf B}''$ satisfies MW3,
it is done by Proposition 3.4.3.

Case 3\,  If 3.3.5 (iv) is met, and $\mf B''$ satisfies MW4, it is done by Proposition 3.4.4.

Case 4\,  If 3.3.5 (iv) is met, and $\mf B''$ satisfies MW5, then there is an induced pair
$(\hat{\mf A}^1, \hat{\mf B}^1)$ in the case of MW5 defined by Remark 3.4.6.
In this case the pair $(\mf A^t,\mf B^t)$ is denoted
by $(\hat{\mf A},\hat{\mf B})$ in order to unify the notations.
Suppose the first arrow $\hat a^1_1$ of $\hat{\mf B}^1$
locates at the $p^1$-th row in the formal product $\hat\Theta^1$. We claim that
$p^1<p$. In fact, the solid arrows
$\tilde{a}^\epsilon_j$ for $j=1,\cdots,n_{_Z}$ at the $p$-th row of $\hat\Theta^\epsilon$ have
differentials $\delta^0 (\tilde a^\epsilon_j)=s_{0j}+\cdots$ according to Corollary 5.2.2, and hence
those arrows will be regularized step by step.

Repeating the above mentioned procedure for $(\hat{\mf A}^1, \hat{\mf B}^1)$, if
one of the cases 1--3 is met, the procedure stops. Otherwise, if the case 4 is
met repeatedly, there exist a sequence of local pairs and a decreasing sequence of the row indices:
$$\begin{array}{ccccc}(\hat{\mf A}, \hat{\mf B}),&
(\hat{\mf A}^1, \hat{\mf B}^1),& (\hat{\mf A}^2, \hat{\mf B}^2),&
\cdots, &(\hat{\mf A}^\beta, \hat{\mf B}^\beta),\\
\qquad\quad p\qquad>&\qquad\quad p^1\qquad>&\qquad\quad p^2\qquad>&\qquad\cdots,\quad>& p^\beta.\end{array}$$
Since the number of the rows of $\hat{H}^i$ for $i=1,\cdots,\beta$ is
fixed, the procedure must stop at some stage $\beta$, where one of the
cases 1--3 appears. The conclusion follows by induction.\hfill$\Box$

\bigskip
\bigskip
\noindent{\large\bf 5.4  Bordered matrices in one-sided case}
\bigskip

In this subsection, a notion of reduced defining systems of Formula (2.4-3)
is given for some induced pairs of a one-sided pair, which is different from Formula (4.1-7). Then some reduction
sequences are constructed starting from one-sided pairs based on bordered matrices.

\medskip

Let $(\mf{A},\mf B)$ be a bipartite pair satisfying RDCC condition, and
$(\mf{A}^r,\mf B^r)$ be an induced pair with $R^r$ being trivial. Suppose
$(\mf{A}^r,\mf B^r)$ has a quotient-sub pair
$((\mf A^r)^{[m]},(\mf B^r)^{(m)})$ denoted by $(\bar{\mf A},\bar{\mf B})$ given
in Formulae (4.1-1) and (4.1-2), where the vertex set $\bar T=\bar T_R\times \bar T_C\subseteq\T^r$,
and $\bar{\mf B}$ has a layer $L=(R; \omega; d_1,\cdots, d_m;
\bar u,\underline u,\bar v,\underline v)$ by Definition 4.1.2.

\medskip

{\bf Remark 5.4.1}\, Suppose $(\bar{\mf A}',\bar{\mf B}')$ with $\bar\T'$ being trivial is
an induced pair of $(\bar{\mf A},\bar{\mf B})$, then it is
a quotient-sub-pair of $(\mf A^{r+\prime},\mf B^{r+\prime})$ by Formula (4.1-5), where $\prime$ stands for some index.
Recall Theorem 2.4.1 and the defining system $\IF^{r,r+\prime}$, the variable matrices $\Psi_{\m^{r,r+\prime}}$,
$\Psi_{\m^{r,r+\prime}}^0$ given by Formula (2.4-3), we now
define its reduced form consisting of the $(p^r,q^r),\cdots(p^r,q^r+m-1)$-th blocks of $\IF^{r,r+\prime}$
according to Remark 4.1.1. Suppose there is a functor
$\bar\vartheta':R(\bar{\mf A}')\rightarrow R(\bar{\mf A})$ acting on objects, $F$
is defined below Formula (4.1-6), and $\bar\vartheta(F'(k))$ has a size vector
$\n=(n_0;n_1,\cdots,n_m)$ partitioned under $\bar\T$.

(i)\, Denote by $\bar Z_0$ the
$(p^r,p^r)$-the square block of $\Psi_{\m^{r,r+\prime}}$ of size $n_0$ with $p^r\in X^r$; by $\bar
Z_{\xi\xi}$ the $(q^r+\xi-1,q^r+\xi-1)$-th square block of size $n_\xi$
with $q^r+\xi-1\in Y^r$ for any $Y^r$. Set
$$\begin{array}{c}\bar Z_0=Z_{X^r}=(z^{X^r}_{pq})_{n_0\times n_0},
\quad \bar Z_{\xi\xi}=Z_{Y^r}=(z^{Y^r}_{pq})_{n_\xi\times n_\xi}.\end{array}
$$

(ii)\, Denote by $\bar Z_\xi$ the $(p^r,q^r+\xi-1)$-block of $\Psi_{\m^{r,r+\prime}}^0$ with size
$n_0\times n_\xi$, and by $\bar Z_{\eta\xi}$ the $(q^r+\eta-1,q^r+\xi-1)$-block for $\eta<\xi$ of $\Psi_{\m^{r,r+\prime}}$
with size $n_\eta\times n_\xi$. Write the matrices
$Z_j=(z^{j}_{pq})_{n_{s(v^r_j)}\times n_{t_{(v^r_j)}}}$ for all the dotted arrows of $\mf B^r$, where $\{z^{j}_{pq}\}_{(p,q),j}$
are different variables over $k$. Suppose
$$\begin{array}{lll} \bar Z_\xi=\sum_{j}\alpha_{\xi}^{j} Z_{j},&\mbox{where}
\,\, s(v^r_j)\ni p^r,\,\,t(v^r_{j})\ni q^r+\xi-1,& \alpha_\xi^j\in k;\\
\bar Z_{\eta\xi}=\sum_j\beta_{\eta\xi}^j Z_j,&\mbox{where}\,\,
s(v^r_j)\ni q^r+\eta-1,\,\, t(v^r_j)\ni q^r+\xi-1, &\beta_{\eta\xi}^j\in k.
\end{array}$$

Return to the pair $(\bar{\mf A}',\bar{\mf B}')$, some indices $\bar p,\bar q,\cdots,$ in the
formal product $\bar\Theta'$  will bee used in order to distinguish with indices $p,q,\cdots,$
in the formal product $\Theta^{r+\prime}$ of $(\mf A^{r+\prime},\mf B^{r+\prime})$ of Formula (5.1-1).
Fix an integer $l\in \{1, \cdots,m\}$ with $Y_l\ne X$ in Definition 4.1.2,
thus $d_l:X\rightarrow Y_l$ is a solid edge. Suppose
$(\bar q_l+1)$ is the index of the first column in the $l$-th block-column
of the formal product $\bar\Theta'$, such that
$(\bar p,\bar q_l+1)$ is the leading position of the first base matrix of $\bar\M'$.
Write $\bar M=\bar\vartheta(F(k))=(\bar M_1,\cdots,\bar M_m)\in R(\bar{\mf A})$
with the size vector $\n$ over $\bar\T$. Then
$$\begin{array}{c}{\mathcal F}:\,\,\bar Z_0\bar M
\equiv_{\prec (\bar p,\bar q_l+1)}(\bar Z_1,\cdots,\bar Z_m)+
\bar M(\bar Z_{\eta\xi})_{1\leqslant\eta\leqslant\xi\leqslant m}\end{array}\eqno{(5.4\mbox{-}1)}$$
is called a {\it reduced defining system} based on Theorem 2.4.1, which is different from $\bar{\mathcal F}$
given by Formula (4.1-7). Similarly as in Equation (5.1-4), there is an equation system:
$$\begin{array}{c}{\mathcal F}_\tau:\,\,0
\equiv_{\prec (\bar p,\bar q_l+1)}(\bar Z_1,\cdots,\bar Z_m)+
\bar M(\bar Z_{\eta\xi})_{1\leqslant\eta\leqslant\xi\leqslant m}.\end{array}\eqno{(5.4\mbox{-}2)}$$

Define a size vector $\tilde {\underline n}=(\tilde n_0;\tilde n_1, \cdots,\tilde n_m)$
over $\bar\T$ as follows: $\tilde n_j=n_j$ if $j\notin Y_l$; $\tilde
n_j=n_j+1$ if $j\in Y_l$. Construct a representation of $\bar{\mf A}$ based on
$\bar M$:
$$\tilde {M}=(\tilde{M}_1,\cdots,\tilde M_m)\in R(\bar{\mf A}),\quad
\tilde{M}_j=\left\{\begin{array}{ll}\bar{M}_j, & \mbox{if}\,\, j\notin\,Y_l;\\
(0\,\bar{M}_j),&\mbox{if}\,\,
j\in\, Y_l,\end{array}\right.\eqno{(5.4\mbox{-}3)}$$
where $0$ is a zero column. Write $\tilde Z_0, \tilde Z_{\xi\xi}$ the variable matrices of size
$\tilde n_0\times \tilde n_0, \tilde n_\xi\times \tilde n_\xi$;
and $\tilde Z_\xi=\sum_j\alpha_{\xi}^j \tilde Z_j$ of size $\tilde n_0\times\tilde n_\xi$,
$\tilde Z_{\eta\xi}=\sum_{\eta\xi}^j\beta_{\eta\xi}^j \tilde Z_j$ of size $\tilde n_\eta\times\tilde n_\xi$ respectively
according to Remark 5.4.1. Denote by $\tilde q_{l}$ the index of the
first column of the $l$-th block-column of $\tilde M$, we obtain the following two matrix equations:
$$\begin{array}{ccc}
\tilde{\mathcal F}:&\tilde Z_0\tilde M
&\equiv_{\prec (\bar p,\tilde q_l)}(\tilde Z_1,\cdots,\tilde Z_m)+
\tilde M(\tilde Z_{\eta\xi})_{1\leqslant\eta\leqslant\xi\leqslant m}\\[1mm]
\tilde{\mathcal F}_\tau:&0&\equiv_{\prec (\bar p,\tilde
q_{l})}(\tilde Z_1,\cdots,\tilde Z_m)+
\tilde M(\tilde Z_{\eta\xi})_{1\leqslant\eta\leqslant\xi\leqslant m}.\end{array}\eqno{(5.4\mbox{-}4)}$$

For any $l'\in Y_l$, denote by $\bar q_{l'}+1$ the index of the first column of the $l'$-th
block-column of $\bar M$. Set any integer $\bar p'\geqslant \bar p$ and  $1\leqslant
h\leqslant n_{_{Y_l}}$, the $(\bar p',\bar q_{l'}+h)$-th entry in the right side of
$\mathcal F_\tau$ of Formula (5.4-2) equals
$$\begin{array}{c} \sum_{\bar q}\gamma_{\bar p'\bar q}z^{Y_l}_{\bar q,\bar q_{l'}+h} +\sum_{\bar q}\nu_{\bar p'\bar q}z^{j}_{\bar q,
\bar q_{l'}+h},\quad \gamma_{\bar p'\bar q},\nu_{\bar p'\bar q}\in k,\, t(v_j)=Y_l.\end{array}\eqno{(5.4\mbox{-}5)}$$
The picture below shows four equations (abridged by four circles)
of $\mathcal F_\tau$. There are three solid
edges ending at $Y_l$, i.e. $|Y_l|=3$, and $n_{_{Y_l}}=4$. If $l'$ is the second index of $Y_l$, then the
equations at the $(\bar p',\bar q_{l'}+h)$-th positions have
the same coefficients for $h=1,2,3,4$.

\medskip

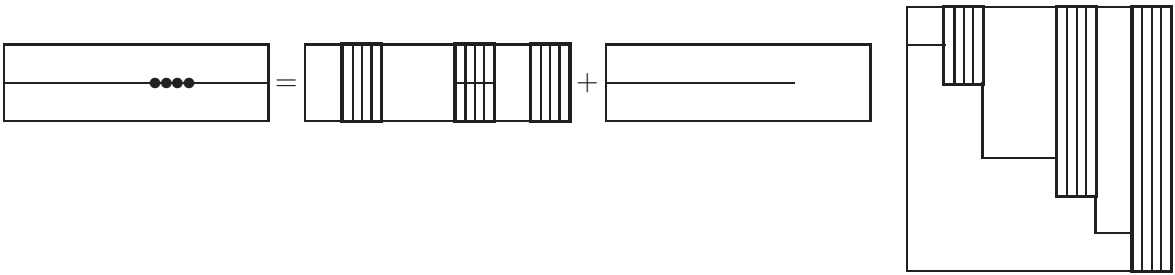
\begin{center}
\setlength{\unitlength}{1mm}
\begin{picture}(170,33)
\thinlines{\linethickness{0.1mm}}
\put(5,25){\framebox(35,10)}\put(5,30){\line(1,0){35}}
\put(25,30){\circle*{1.5}}\put(26.5,30){\circle*{1.5}}\put(28,30)
{\circle*{1.5}}\put(29.5,30){\circle*{1.5}}
\put(41,29){$=$}
\put(81,29){$+$}
\put(45,25){\framebox(35,10)}
\mput(51.25,25)(1.25,0){3}{\line(0,1){10}}\mput(66.25,25)(1.25,0){3}{\line(0,1){10}}
\mput(76.25,25)(1.25,0){3}{\line(0,1){10}}\put(65,30){\line(1,0){5}}

\put(85,25){\framebox(35,10)}\put(85,30){\line(1,0){25}}

\put(125,5){\framebox(35,35)}
\mput(131.25,30)(1.25,0){3}{\line(0,1){10}}\mput(146.25,15)(1.25,0){3}{\line(0,1){25}}
\mput(156.25,5)(1.25,0){3}{\line(0,1){35}}
\put(125,35){\line(1,0){5}}
\put(135,20){\line(1,0){10}}\put(135,20){\line(0,1){10}}
\put(150,10){\line(1,0){5}}\put(150,10){\line(0,1){5}}

\thicklines{\linethickness{1mm}}
\put(50,25){\framebox(5,10)}\put(65,25){\framebox(5,10)}\put(75,25){\framebox(5,10)}
\put(130,30){\framebox(5,10)}\put(145,15){\framebox(5,25)}\put(155,5){\framebox(5,35)}
\end{picture}
\end{center}

{\bf Lemma 5.4.2}\, Being parallel to Lemma 5.1.4, there are
following assertions.

(i)\, For any $1\leqslant h_1,h_2\leqslant n_l$ and any $l'\in Y$, the $(\bar p,\bar q_{l'}+h_1)$-th
equation is a linear combination of
the previous equations in $\mathcal F_\tau$, if and only if so is the
$(\bar p,\bar q_{l'}+h_2)$-th equation. Similarly, the same result is valid in $\tilde{\mathcal F}_\tau$.

(ii)\, The equations in the system $\tilde{\mathcal F}$ (resp. $\tilde{\mathcal F}_{\tau}$) and
those in $\mathcal F$ (resp. $\mathcal F_{\tau}$) are the same at the same positions of each
block column, whenever the added $\tilde q_{l'}$-th columns for all $l'\in Y_l$ have
been dropped from $\tilde{\mathcal F}$ (resp. $\tilde{\mathcal F}_{\tau}$).

(iii)\, A subsystem of $\tilde{\mathcal F}$ (resp. $\tilde{\mathcal F}_{\tau}$)
consisting of the $\tilde q_{l'}$-column in both sides is
$\tilde{\mathcal F}_{\tau\tilde q_{l'}}: 0\equiv \tilde M\Phi_{\tilde\n,\tilde q_{l'}}$, where $\Phi_{\tilde\n,\tilde q_{l'}}$
stands for the $\tilde q_{l'}$-th column of $\Phi_{\tilde\n}$. And
the system  $\{\tilde{\mathcal F}_{\tau\tilde q_{l'}}\mid
\forall\, \bar l'\in Y_l\}$ can be solved independently.

(iv)\, If the $(\bar p,\bar q_{l'}+h)$-th equation for some $1\leqslant h\leqslant n_{_Y}$ is a linear combination of the previous
equations in $\mathcal F$, then so is the $(\bar p,\tilde q_{l'})$-th equation in the system $\{\tilde{\mathcal F}_{\tau\tilde q_{l'}}\mid
\forall\, l'\in Y_l\}$.

\smallskip

{\bf Proof}\, (i)--(ii) See Proof (i)--(ii) of Lemma 5.1.4.

(iii) Note that $\forall\, l'\in Y_l$ the variables at the $\tilde q_{l'}$-th column
of $(\tilde Z_\xi)_{1\leqslant\xi\leqslant m}$ and
$(\tilde Z_{\eta\xi})_{1\leqslant\eta\leqslant\xi\leqslant m}$
are different from those at the $\bar h$-th column for all $\bar h\ne \tilde q_{l'}, \forall\, l'\in Y$,
and different from those in $\tilde Z_0$.

(iv) See proof (iv) of 5.1.4.

\medskip

Let $(\bar{\mf A}^s,\bar{\mf B}^s)$ be an induced pair of $(\bar{\mf A},\bar{\mf B})$
with $\bar R^s$ being trivial and local. Then
there are two sequences of reduction in the sense of Lemma 2.3.2 according to Formula (4.1-5):
$$\begin{array}{llllllll}\mf A, \mf A^1,\cdots,\mf A^{r-1},
&\mf A^r,&\mf A^{r+1},&\cdots,&\mf A^{r+i},&\mf A^{r+i+1},&\cdots,&
\mf{A}^{r+s};\\
& \bar{\mf A},&\bar{\mf{A}}^1,&\cdots,&\bar{\mf{A}}^{i},&\bar{\mf{A}}^{i+1},&\cdots,&
\bar{\mf{A}}^s.\end{array}\eqno {(\bar \ast)}
$$
Set $\bar M=\vartheta^{0s}(F^s(k))$ of size vector $\n$ over $\bar\T$, a bordered matrix
$\tilde M$ of size vector $\tilde{\n}$ can be constructed according to Formula (5.4-3).

\medskip

{\bf Theorem 5.4.3}\, Being parallel to Theorem 5.2.1, there exists a unique reduction sequence based
on the sequence $(\bar \ast)$, where ${\sideset{^\cdot}{}{\mathop{\tilde{\mf A}}}}$ stands
for $\tilde{\bar{\mf A}}$ in order to simplify the notation:
$$\begin{array}{llllllll}\tilde{\mf A}, \tilde{\mf A}^1,\cdots,\tilde{\mf A}^{r-1},
&\tilde{\mf A}^r,&\tilde{\mf A}^{r+1},&\cdots,&\tilde{\mf A}^{r+i},&\tilde{\mf
A}^{r+i+1},&\cdots,&\tilde{\mf A}^{r+s};\\
&{\sideset{^\cdot}{}{\mathop{\tilde{\mf A}}}},&{\sideset{^\cdot}{^1}{\mathop{\tilde{\mf A}}}},&\cdots,
&{\sideset{^\cdot}{^i}{\mathop{\tilde{\mf A}}}},
&{\sideset{^\cdot}{^{i+1}}{\mathop{\tilde{\mf A}}}},&\cdots,&
{\sideset{^\cdot}{^s}{\mathop{\tilde{\mf A}}}}.\end{array}\eqno {(\tilde{\bar\ast})}
$$

(i)\, $\tilde{\mf{A}}^i=\mf{A}^i$ for $i=0,1,\cdots,r$.

(ii)\, The reduction from ${\sideset{^\cdot}{^i}{\mathop{\tilde{\mf A}}}}$ to
${\sideset{^\cdot}{^{i+1}}{\mathop{\tilde{\mf A}}}}$ is a reduction or a composition of two
reductions in the sense of Lemma 2.3.2 for $i=0,\cdots,s-1$, such
that ${\sideset{^\cdot}{}{\mathop{\tilde{\T}}}}^s$ has two vertices, and
$\tilde{\vartheta}^{0s}(\tilde{F}^s)=\tilde{M}$.

(iii)\, The reduction from $\tilde{\mf{A}}^{r+i}$ to $\tilde{\mf{A}}^{r+i+1}$ is done in the same way as
that from ${\sideset{^\cdot}{^i}{\mathop{\tilde{\mf A}}}}$ to ${\sideset{^\cdot}{^{i+1}}{\mathop{\tilde{\mf A}}}}$.

(iv)\, The diagonal block $\tilde e_{_X}$ in ${\sideset{^\cdot}{^s}{\mathop{\tilde{\K}}}}$ of
${\sideset{^\cdot}{^s}{\mathop{\tilde{\mf A}}}}$ partitioned under $\bar\T$ is of the form of Corollary 5.2.2.

\smallskip

{\bf Proof}\, (i) is clear.

(ii) The proof is parallel to that of Theorem 5.2.1, the only difference appears
in the item 1.4 of Case 1. Suppose an edge reduction is made from $\bar{\mf A}^i$
to $\bar{\mf A}^{i+1}$ with the reduction block $G^{i+1}$
being at the $l^0$-th block column with $l^0\in Y_l$. Then $\tilde F^{i+1}$ has a size vector
$\l^{i+1}\times\tilde\n^{i+1}$ over $\bar\T$ with $\tilde n^{i+1}_l=n^{i+1}_l+1$, and
a zero column  is added into the $l'$-th block column from the left hand side for every $l'\in Y_l$ .

(iii)\, follows from Formula (4.1-5).

(iv)\, The proof is parallel to that of Corollary 5.2.2. \hfill$\Box$

\medskip

Being parallel to $(\ast')$ at the beginning of Subsection 5.3, there are following two sequences:
$$\begin{array}{llllllllll}\mf A,\mf A^1,\cdots,\mf A^{r-1},
&\mf{A}^{r}, &\mf{A}^{r+1},&\cdots,&{\mf{A}}^{r+s},
&{\mf{A}}^{r+s+1},  &\cdots, &{\mf{A}}^{r+\epsilon},&\cdots,&{\mf
A}^{r+t}, \\
&\bar{\mf{A}},&\bar{\mf{A}}^1,&\cdots,&\bar{\mf{A}}^s,
&\bar{\mf{A}}^{s+1},
&\cdots,&\bar{\mf{A}}^\epsilon,&\cdots,&\bar{\mf A}^t.\end{array} \eqno{(\bar\ast')}
$$
The reductions from $\bar{\mf A}$ to $\bar{\mf A}^{s}$
is given by $(\bar\ast)$; from $\bar{\mf A}^{s}$
to $\bar{\mf A}^{s+1}$ is a loop mutation and
a parameter $x$ appears; the reduction from $\bar{\mf A}^{i}$
to $\bar{\mf A}^{i+1}$ is a
regularization for $i=s+1,\cdots,t-1$. The first arrow of $\bar{\mf B}^t$ locates
at the $(\bar p,\bar q+j)$-th position in the formal product $\bar\Theta^t$ for some $1\leqslant j\leqslant n_l$, and that of
$\bar{\mf B}^\epsilon$ at $(\bar p,\bar q+1)$-th position in $\bar\Theta^\epsilon$. The pair $(\mf A^{r+t},\mf B^{r+t})$ is
minimal wild in the case of MW5 of Remark 3.4.6 and Classification 5.1.1 (II).

\medskip

{\bf Remark 5.4.4}\, (i)\, If the first arrow $a^t_1$ of $\bar{\mf B}^t$ is splitting from $d_l$ of
the one-sided bocs $\bar{\mf B}$,
then $d_l:X\mapsto Y_l$ is an edge by Theorem 4.6.1 and Corollary 4.6.2. Therefore
it is possible to apply Theorem 5.4.3 with respect to the vertex $Y_l\in\bar\T$ for the sequence
$(\bar\ast)$, and obtain the sequence $(\tilde{\bar\ast})$.

(ii)\, We will describe how to determine
$\mf A^r$, thus $\bar{\mf A}$, in the next subsection.

(iii)\, Being parallel to Formula (5.3-2), the equation
system $\{\tilde{\mathcal F}_{\tau\tilde q_{l'}}\mid
\forall\, l'\in Y_l\}$ given by Lemma 5.4.2 (iii) is considered. Thus
some polynomials $d^{jl'}(x)$ are obtained for $l'\in Y_l,j=\bar p_x,\cdots,\bar p$ inductively,
by an analogous discussion as in the subsection 5.3. If $\bar R^t=k[x,\phi^t(x)^{-1}]$, define a polynomial
similar to Formula (5.3-5):
$$\begin{array}{c}\tilde\phi^t(x)=\phi^t(x)\prod_{j=\bar p_x}^{\bar p}\prod_{l'\in Y_l}d^{jl'}(x).\end{array}$$

{\bf Theorem 5.4.5}\, Being parallel to Theorem 5.3.2, there exist two unique reduction sequences based
on the sequences $(\bar \ast')$:
$$\begin{array}{llllllllll}\tilde{\mf A},\tilde{\mf A}^1,\cdots,\tilde{\mf A}^{r-1}
&\tilde{\mf{A}}^{r}, &\tilde{\mf A}^{r+1},&\cdots,&\tilde{\mf{A}}^{r+s},
&\tilde{\mf{A}}^{r+s+1},  &\cdots,
&\tilde{\mf{A}}^{r+\epsilon},&\cdots,&\tilde{\mf
A}^{r+t}; \\
&{\sideset{^\cdot}{}{\mathop{\tilde{\mf A}}}},&{\sideset{^\cdot}{^1}{\mathop{\tilde{\mf A}}}},
&\cdots,&{\sideset{^\cdot}{^s}{\mathop{\tilde{\mf A}}}},
&{\sideset{^\cdot}{^{s+1}}{\mathop{\tilde{\mf A}}}},
&\cdots,&{\sideset{^\cdot}{^\epsilon}{\mathop{\tilde{\mf A}}}},&\cdots,
&{\sideset{^\cdot}{^t}{\mathop{\tilde{\mf A}}}}.\end{array} \eqno{(\tilde{\bar\ast}')}
$$

(i) The first parts of the two sequences up to $\tilde{\mf A}^{r+s}$ and ${\sideset{^\cdot}{^s}{\mathop{\tilde{\mf A}}}}$ respectively
are given by $(\tilde{\bar\ast})$.

(ii)\, ${\sideset{^\cdot}{^{s+1}}{\mathop{\tilde{\mf A}}}}$ is induced from
${\sideset{^\cdot}{^s}{\mathop{\tilde{\mf A}}}}$ by a loop mutation $a^{s+1}_1\mapsto(x)$,
or an edge reduction $(0)$ followed by a loop mutation $(x)$; the reduction from
${\sideset{^\cdot}{^{i}}{\mathop{\tilde{\mf A}}}}$ to ${\sideset{^\cdot}{^{i+1}}{\mathop{\tilde{\mf A}}}}$
for $i>s$ is given by a regularization, or two regularizations,
or a reduction given by Lemma 2.2.6 followed by a regularization.

(iii) The reduction from ${\sideset{^\cdot}{^{r+i}}{\mathop{\tilde{\mf A}}}}$ to
${\sideset{^\cdot}{^{r+i+1}}{\mathop{\tilde{\mf A}}}}$ is done in the same way as that from ${\sideset{^\cdot}{^{i}}{\mathop{\tilde{\mf A}}}}$ to
${\sideset{^\cdot}{^{i+1}}{\mathop{\tilde{\mf A}}}}$ for $s<i<t$.

(iv)  If the bocs $\bar{\mf B}^t$ in the sequence $(\bar\ast')$ satisfies MW5 of Remark 3.4.6
and Classification 5.1.1 (II), then $\dz(\tilde a^\epsilon_{0})=0$ for the first arrow $\tilde a^\epsilon_{0}$
 of ${\sideset{^\cdot}{^\epsilon}{\mathop{\tilde{\mf B}}}}$ in $(\tilde{\bar\ast}')$.

\smallskip

{\bf Proof}\, (i) is obvious. The proof of (ii) is parallel to that of Theorem 5.3.2. (iii)
follows from Formula (4.1-5). The proof of (iv) is parallel to that of Corollary 5.3.3 by Corollary 2.4.3. \hfill$\Box$

\bigskip
\bigskip
\noindent{\large\bf 5.5 Non-homogeneity in the case of MW5 and classification (II)}
\bigskip

Suppose a bipartite pair $(\mf A,\mf B)$ has an induced pair $(\mf A',\mf B')$ in the case
of MW5 of Remark 3.4.6 and Classification 5.1.1 (II) in this subsection. A one-sided quotient-sub pair is
determined according to the position of the first arrow $a_1'$ in
the formal product $H'+\Theta'$; then the non-homogeneity of the pair $(\mf A,\mf B)$ is proved.

\medskip

Let $\mf A$ be a bipartite matrix bimodule problem satisfying RDCC condition.
The sequence
$$
(\mf{A},\mf B),(\mf{A}^1,\mf B^1),\cdots, (\mf{A}^{\varsigma},\mf{A}^{\varsigma}),
(\mf{A}^{\varsigma+1},\mf{B}^{\varsigma+1}),\cdots,
(\mf{A}^{\tau},\mf B^\tau)=(\mf A',\mf B')\eqno{(5.5\mbox{-}1)}
$$
satisfies the following conditions:
$R^i$ is trivial for $i\leqslant\varsigma$, the reduction from $\mf{A}^{i}$ to $\mf{A}^{i+1}$ is in the
sense of Lemma 2.3.2 for $i<\varsigma$; $\mf{A}^{\varsigma}$ is
local with $\dz(a^\varsigma_1)=0$, and $R^{\varsigma+1}=k[x]$ in $\mf B^{\varsigma+1}$
after a loop mutation; finally, the reduction from $\mf{A}^{i}$ to $\mf{A}^{i+1}$ is a
regularization for $i>\varsigma$, and $\mf B^\tau=\mf B'$ is in the case of MW5 of Remark 3.4.6 and
Classification 5.1.1 (II). Suppose the index of the first arrow  $a_1^\tau$ of
$\mf B^\tau$ is $(p^\tau,q^\tau)$ in the formal product $\Theta^\tau$, which is sitting at the
$(\textsf p,\textsf q)$-th block partitioned under $\T$.
According to Formula (2.3-7):
$$\begin{array}{c}H^\tau=\sum_{i=1}^{\varsigma} G_\tau^{i}\ast A^{i-1}_1
+(x)\ast A^{\varsigma}_1.
\end{array}\eqno{(5.5\mbox{-}2)}$$
Following discussion will be focused on the reduction blocks $G^{i}_\tau$ of $H^\tau$.

Let $i<\varsigma$, $\vartheta^{i\tau}: R(\mf A^\tau)\rightarrow R(\mf A^i)$
be the induced functor, and $\n^{i\tau}=\vartheta^{i\tau}(1,\cdots,1)$.
There is a simple fact, that any row (column)
index $\rho$ of $H^i+\Theta^i$ in the pair $(\mf A^i,\mf B^i)$ determines a row (column) index
$n^{i\tau}_{1}+\cdots+n^{i\tau}_{\rho}$ of $H^\tau+\Theta^\tau$ in the pair $(\mf A^\tau,\mf B^\tau)$.
Consequently, if the upper (resp. lower, left or right) boundaries of two reduction blocks
$G_i^{j_1},G_i^{j_2}$ in $H^i$ are collinear, if and only if
the corresponding boundaries of two splitting blocks $G_\tau^{j_1},G^{j_2}_\tau$
in $H^\tau$ are collinear. The two blocks $G_i^j$ and $G^j_i(k)$ may not be
distinguished for the sake of convenience in the following statements.

\medskip

{\bf Remark 5.5.1}\, Consider the reduction blocks inside the
$(\textsf p,\textsf q)$-th block partitioned under $\T$. The relative position of the
upper boundaries of $G^i_\tau$ and $G^{i+1}_\tau$ in this block
has three possibilities according to Formulae (2.3-3)--(2.3-5).

(i) The upper boundaries of $G^i_\tau$ and $G^{i+1}_\tau$ are collinear, and this occurs
if and only if the reduction from $\mf A^{i-1}$ to $\mf A^{i}$ is
given by one of the following reduction blocks: $G^{i}=\emptyset$ in a regularization;
$G^i=(\lambda)$ in a loop reduction; $G^{i}=(0), (1)$ or $(0\,1)$ in an edge reduction,
moreover the right boundary of $G^i_\tau$ is not that of
the $(\textsf p,\textsf q)$-th block. In this case their lower boundaries are also collinear.

(ii)\, The upper boundary of $G^{i+1}_\tau$ is strictly lower than that of $G^i_\tau$, and this occurs
if and only if $G^{i}=W$ of size being strictly bigger than $1$ in a loop reduction, or
$G^{i}={1\choose 0}$ or ${{0\,1}\choose{0\,0}}$ in an edge reduction, and the right boundary of $G^i_\tau$ is not that of
the $(\textsf p,\textsf q)$-th block. In this case,
the lower boundaries of $G^i_\tau$ and $G^{i+1}_\tau$ are also collinear.

(iii)\, The lower boundary of $G^{i+1}_\tau$ is the upper boundary of $G^i_\tau$, and this occurs
if and only if the right boundary of $G^i_\tau$ coincides with that
of the $(\textsf p,\textsf q)$-th block.

\medskip

Collect all the reduction blocks of $H^\tau$ inside the $(\textsf p,\textsf q)$-th
block, such that their upper boundaries are above or at that of $a^\tau_1$:
$$G_\tau^{q_1},G_\tau^{q_2},\ldots, G_\tau^{q_u},\quad \mbox {with}\quad
1\leqslant q_1<q_2<\cdots<q_u\leqslant\varsigma.\eqno{(5.5\mbox{-}3)}$$
The set of reduction blocks $\{G_\tau^{q_i}\mid 1\leqslant i\leqslant u\}$ in (5.5-3) is divided
into $h$ groups according to whether the upper boundaries of blocks
are collinear or not, and denoted by $\rho_j$ the common upper boundary of the blocks in the $j$-th group
for $j=1,\cdots,h$, where $\rho_{j+1}$ is
strictly lower than $\rho_{j}$:
$$\{G_\tau^{q_{1,1}},\dots,G_\tau^{q_{1,u_1}}\},
\cdots\cdots, \{G_\tau^{q_{h,1}},\cdots, G_\tau^{q_{h,u_{h}}}\}, \quad u_1+\cdots+u_{h}=u.
\eqno{(5.5\mbox{-}4)}$$
The adjacent blocks $G_\tau^{q_{j,l}}$ and $G_\tau^{q_{j,l+1}}$ in
the $j$-th group have two possibilities:
\ding{172} $G_\tau^{q_{j,l}}$ is in case (i) of Remark 5.5.1, \ding{173} $G_\tau^{q_{j,l}}$ is in case (ii) of
Remark 5.5.1. Then $G_\tau^{q_{j,l+1}}$ comes from
the next reduction with $q_{j,l+1}=q_{j,l}+1$ in \ding{172}. But in \ding{173}, $G_\tau^{q_{j,l+1}}$ follows by
a sequence of reductions with the upper boundaries of the reduction blocks lower
than that of $a_1^\tau$, and the sequence includes at least one reduction in case (iii) of Remark 5.5.1. Finally,
the sequence reaches $G_\tau^{q_{j,l+1}}$ with the upper boundary
$\rho_j$ as a neighbor of $G_\tau^{q_{j,l}}$. Thus $q_{j,l+1}>q_{j,l}+1$.

\medskip

{\bf Lemma 5.5.2}\, The last block $G_\tau^{q_{j,u_j}}$ of the $j$-th group must be
as in case (ii) of Remark 5.5.1 for $j=1,\cdots,h$.

\smallskip

{\bf Proof}\, If $G_\tau^{q_{j,u_j}}$ is in case (iii) of 5.5.1, then $\rho_j$
is lower than $\rho_{j+1}$ for $j<h$, which is a contradiction to the grouping of Formula (5.5-4);
and $a_1^\tau$ is sitting above $\rho_{h}$ for $j=h$, which is a contradiction to the choice
of the sequence (5.5-3).

Suppose $G_\tau^{q_{j,u_j}}$ is in case (i) of 5.5.1.  Then for $j<h$,
the upper boundaries of $G_\tau^{q_{j,u_j}}$ and $G_\tau^{q_{j,u_j}+1}$ coincide, which is
a contradiction to the grouping of (5.5-4). For $j=h$, $G^{q_{h,u_h}}_\tau=0, I,(0\,I)$, or
$\lambda I$, or $\emptyset$ with the height $d\geqslant 1$.
Suppose the next reduction is still in the sense of Lemma 2.3.2 given by $G^{q_{h,u_h}+1}_{\tau}$,
which is denoted by $G'_\tau$ for simplicity.
If $G'_\tau$ is in the case (i) or (ii) of 5.5.1, then $G'_{\tau}$ and $G^{q_{h,u_h}}_\tau$
have the same upper boundary, a contradiction to the grouping of (5.5-4); if $G'_\tau$ is in case of 5.5.1 (iii),
then $a_1^\tau$ locates above $\rho_h$, which is a contradiction to the choice of (5.5-3).
Therefore the reduction in the sense of Lemma 2.3.2 should not be able to continue, and $r+s=q_{h,u_h}$
in the sequence $(\bar\ast')$ before Remark 5.4.4. Thus the hight $d=1$, the parameter $x$
appears by a loop mutation. Since $a^\tau_1$
locates between the upper and lower boundaries of $G^{q_{h,u_h}}_\tau$, $x$ and $a_1^\tau$
are sitting at the $p^\tau$-th row simultaneously, a contradiction to Lemma 5.1.2.
So $G^{q_{j,u_j}}_\tau$ is in the case of 5.5.1 (ii) as desired.\hfill$\Box$

\medskip

{\bf Definition 5.5.3}\, We define $h$ rectangles in $\Theta^\tau$: for $j<h$,
the $j$-th rectangle has the upper boundary $\rho_j$, the lower one
$\rho_{j+1}$, and the left one is the right boundary of $G^{q_{j,u_j}}$,
the right one is that of the $(\textsf p,\textsf q)$-th block.
While the upper boundary of the $h$-th rectangle is $\rho_h$, lower boundary
is that of $G^{q_{h,u_h}}$. The rectangle with upper boundary $\rho_j$ is
said to be the $j$-th {\it ladder}, there are altogether $h$ ladders.

\medskip

The picture below shows an example for $h=3$. Three groups of
Reduction blocks given in sawtooth patterns with some dots,
but the last block in each group is given by a rectangle without dots.
The upper boundaries of the three ladders are shown by dotted lines.

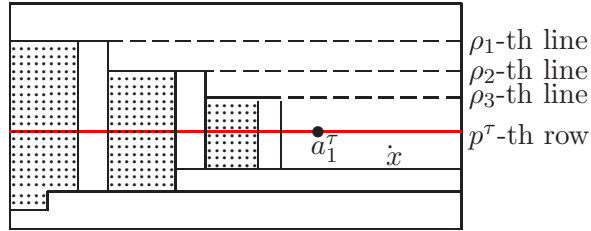
\begin{center}
\setlength{\unitlength}{1mm}
\begin{picture}(80,30)
\put(0,0){\line(1,0){60}}  \put(0,0){\line(0,1){30}}
\put(0,30){\line(1,0){60}} \put(60,0){\line(0,1){30}}

\put(0,2.5){\line(1,0){5}}\put(5,2.5){\line(0,1){2.5}}
\put(5,5){\line(1,0){55}} {\color{red}\put(0,13){\line(1,0){60}}}
\put(9,5){\line(0,1){20}} \put(13,5){\line(0,1){20}}
\put(0,25){\line(1,0){13}}
\multiput(13,25)(3,0){16}{\line(1,0){2}}
\put(22,5){\line(0,1){16}}
\put(26,8){\line(0,1){13}} \put(33,8){\line(0,1){9}}
\put(36,8){\line(0,1){9}} \put(13,21){\line(1,0){13}}
\multiput(25,21)(3,0){12}{\line(1,0){2}}
\put(26,17.5){\line(1,0){10}}\multiput(37,17.5)(3,0){8}{\line(1,0){2}}

\put(22, 8){\line(1,0){38}} \put(40,12){$\bullet$}\put(40, 10){$a^\tau_1$}
\put(50, 10){$\cdot$}\put(50,8.5){$x$}

\put(61, 11.3){$p^\tau$-th row}
\put(61, 16.8){$\rho_3$-th line}
\put(61, 20){$\rho_2$-th line}
\put(61, 24){$\rho_1$-th line}

\multiput(0,-2.5)(1,0){5}{\begin{picture}(0,0)\multiput(0,5)(0,1){2}{$\cdot$}\end{picture}}
\multiput(0,0)(1,0){9}{\begin{picture}(0,0)\multiput(0,5)(0,1){19}{$\cdot$}\end{picture}}
\multiput(13,0)(1,0){9}{\begin{picture}(0,0)\multiput(0,5)(0,1){15}{$\cdot$}\end{picture}}
\multiput(26,2.5)(1,0){7}{\begin{picture}(0,0)\multiput(0,5)(0,1){9}{$\cdot$}\end{picture}}
\end{picture}
\end{center}

{\bf Lemma 5.5.4}\, Let index $r=q_{h,u_h}-1$ in the sequence (5.5-1).
We define a one-sided quotient-sub pair $(\bar{\mf A},\bar{\mf B})=((\mf A^r)^{[m]},(\mf B^r)^{(m)})$
of the pair $(\mf A^r,\mf B^r)$ consisting of the solid arrows
$d_1,\cdots,d_m$ sitting at the $p^r$-row as shown in Picture (4.1-1). Then

(i)\, $m>1$;

(ii)\, all the reduction blocks in $H^\tau$
yielded from some split of $d_2,\cdots,d_m$ locate below the $p^\tau$-th row;

(iii)\, $a_1^\tau$ is split from $d_l$ with $l>1$. If $(\bar{\mf A},\bar{\mf B})$ satisfies
the hypothesis of Theorem 4.6.1 or Corollary 4.6.2, then $d_l$ is a solid edge.

(iv)\, $\varsigma=r+s$ and $\tau=r+t$. Therefore the sequence (5.5-1) coincides with the first
sequence of $(\bar\ast')$ given before Remark 5.4.4.

\smallskip

{\bf Proof}\, (i)\, follows from Lemma 5.5.2.

(ii)\, comes from the choice of the reduction blocks of Formula (5.5-3).

(iii) and (iv) are obvious. \hfill$\Box$

\medskip

{\bf Remark 5.5.5}\, (i) $H^{r}+\Theta^{r}$ of the pair $(\mf A^{r},\mf B^{r})$
has also $h$ ladders in the $(\textsf p,\textsf q)$-th block.
The boundaries of the $j$-th ladder of $H^\tau+\Theta^\tau$ is derived from
that of the $j$-th ladder of $H^{r}+\Theta^{r}$ for $j=1,\cdots,h$
according to the simple fact stated before Remark 5.5.1.

(ii)\, Let $(\mf A_{X}^r,\mf B^r_{X})$ be the induced local pair at $X$
of $(\mf A^r,\mf B^r)$, denote by $h_{_{X}}$ the number
of the inheriting ladders of $H^r_{X}+\Theta^r_{X}$ from $H^r+\Theta^r$, then $h_{_{X}}\leqslant h$.

(iii)\, Return to sequence $(\bar\ast')$ and $(\tilde{\bar\ast}')$ in Subsection 5.4.
It is easy to see that $\tilde H^{r+t}+\tilde\Theta^{r+t}$ has also $h$ ladders,
and the number of rows in the $h$-th ladder in $\tilde H^{r+t}+\tilde\Theta^{r+t}$
is the same as that in $H^{r+t}+\Theta^{r+t}$.

\medskip

{\bf Proposition 5.5.6}\, Let $(\mf A,\mf B)$ be a pair with $\T$ trivial, such that
$\mf{A}=(R,\K,\M,H=0)$ is a bipartite
matrix bimodule problem satisfying RDCC condition. If there exists an induced
pair $(\mf A',\mf B')$ of $(\mf A,\mf B)$ in the case of MW5 defined by Remark 3.4.6, and the sum
$H'+\Theta'$ of $(\mf A',\mf B')$ satisfies Classification 5.1.1 (II), then $\mf{B}$ is not
homogeneous.

\smallskip

{\bf Proof}\, Suppose the induced pair $(\mf A',\mf B')$ is the last term
$(\mf A^{r+t},\mf B^{r+t})$ in the first sequence of Formula $(\bar\ast')$
given before Remark 5.4.4. Keep the notations in the two sequences of $(\bar\ast')$.

We assume in addition that the number of the ladders in $H^{r+t}+\Theta^{r+t}$
is minimal with respect to the property of Classification 5.1.1 (II).

\smallskip

{\bf (I)}\, Let $X$ be given by Definition 4.1.2. If the local pair $(\mf A^r_{X},\mf B^r_{X})$
is wild, using the triangular Formulae of Subsection 3.3, an induced minimal wild
local pair $((\mf A^r_{X})',(\mf B^r_{X})')$
with the parameter $x'$ and the first arrow $a_1'$ is obtained.

(I-1)\, If $(\mf B^r_{X})'$ is in the case of MW3, or MW4, or MW5 with $H_{X}'+\Theta_{X}'$ being
in the case of Classification 5.1.1 (I), then $(\mf A,\mf B)$ is not homogeneous by Proposition
3.4.3, or 3.4.4, or 5.3.4, it is done.

(I-2)\, If $(\mf B^r_{X})'$ is in the case of MW5 and Classification 5.1.1 (II),
then the number of the inheriting ladders $h_{_X}$ in $H^r_{X}+\Theta^r_{X}$ does not exceed
$h$ by Remark 5.5.5 (ii). Suppose $a_1'$ locates at the $h'$-ladder. If $h'<h_X\leqslant h$,
or $h'=h_X<h$, then it contradicts to the minimality assumption on the number of ladders.
If $h'=h_{_X}=h$, since this ladder contains only one row by Lemma 5.5.4, $x'$and $a_1'$ must
locate at the same row, a contradiction to Lemma 5.1.2.

\smallskip

{\bf (II)}\, Suppose $(\mf A^r_{X},\mf B^r_{X})$ is tame infinite, and its quotient-sub-pair
$(\bar {\mf A}_{X},\bar{\mf B}_{X})$ is in case (ii) of Classification 4.2.1.

(II-1) If the one-sided pair $(\bar{\mf A},\bar{\mf B})$ satisfies the hypothesis of Lemma 4.2.3
or 4.4.1, then $(\mf A,\mf B)$ is not homogeneous. In fact the loop $\bar b$ is the unique effective
loop of both $\bar{\mf B}_{X}$ and $\mf B^r_{X}$.

(II-2)\, If $(\bar{\mf A},\bar{\mf B})$ satisfies the hypothesis of Theorem 4.4.2, then the
triangular formulae given in Subsection 3.3 can be used for the local wild pair $(\mf A^{r+2l},\mf B^{r+2l})$
given in Proof 4) of 4.4.2.
If the cases of MW3, MW4, or MW5 and Classification 5.1.1 (I) are reached, it is done. If
MW5 and Classification 5.1.1 (II) is med again, the first arrow must be outside of the $h$-th ladder, which is
a contradiction to the minimal number assumption of the ladders.

\smallskip

{\bf (III)}\, Now the following two cases are considered. First, $\mf B^r_{X}$ is tame infinite,
$\bar{\mf B}_{X}$ is in the case (ii) of Classification 4.2.1,
the pair $(\bar {\mf A},\bar{\mf B})$ is major and the $c$-class arrows satisfy Formula (4.2-6). Second,
$\mf B^r_{X}$ is tame infinite or finite, and $\bar{\mf B}_{X}$ is finite.

Then in both cases $d_l$ of $\bar{\mf B}$, from which $a_1^t$ is split, is a solid edge by Lemma 5.5.4 (iii).
Consequently, Formula $(\tilde{\bar\ast})$ of Theorem 5.4.3
can be used with respect to $d_l$. Keep the notations in the two sequences of $({\tilde{\bar\ast}}')$
of Theorem 5.4.5.

Since $\dz(\tilde{a}^\epsilon_0)=0$ in the pair
$({\sideset{^\cdot}{^\epsilon}{\mathop{\tilde{\mf A}}}},{\sideset{^\cdot}{^\epsilon}{\mathop{\tilde{\mf B}}}})$ of $({\tilde{\bar\ast}}')$
by 5.4.5 (iv), set the edge $\tilde{a}^\epsilon_0\mapsto (1)$ by Proposition 2.2.7. Then all the other arrows splitting from
$d_l$ at the same row are regularized in the further reductions by 5.4.3 (iv).
The induced pair is obviously local and tame infinite or wild type. Then it is possible to
use the triangular formulae of Subsection 3.3 once again, and an induced pair $(\hat {\mf{A}}^1,\hat {\mf{B}}^1)$
in the cases (ii)-(iv) of Classification 3.3.5 is obtained.

\smallskip

(III-1)\, If the induced local pair $(\hat {\mf{A}}^1,\hat {\mf{B}}^1)$ is tame
infinite, then the two-point pair $(\tilde{\mf A}^{r+\epsilon},\tilde{\mf B}^{r+\epsilon})$
satisfies the hypothesis of Proposition 3.4.5, it is done.

(III-2)\, If $(\hat {\mf{A}}^{1},\hat {\mf{B}}^{1})$ is in the case of MW3,
or MW4, or MW5 of Remark 3.4.6 and Classification 5.1.1 (I), it is done.

(III-3)\, If $(\hat {\mf{A}}^{1},\hat {\mf{B}}^{1})$ is in the case of MW5 of Remark 3.4.6
and classification 5.1.1 (II), and suppose in addition, the first arrow of $\hat{\mf{B}}^1$ locates at the
$h_1$-th ladder with $h_1<h$, then there is a contradiction to the
minimality number assumption of the ladders.

(III-4)\, If $(\hat {\mf{A}}^{1},\hat {\mf{B}}^{1})$ is in the case of MW5 of Remark 3.4.6
and classification 5.1.1 (II), and suppose in addition, the first arrow of
$\hat{\mf B}^1$ locates still at the $h$-th ladder, it is needed to do
induction on some pairs of integers.
Denote by $\sigma$ the number of the rows in the $h$-th ladder of $H^{r+t}+\Theta^{r+t}$,
which is a constant after making some bordered matrices by Remark 5.5.5 (iii). And denote by
$m$ the number of the solid arrows in the pair $(\bar{\mf A},\bar{\mf B})$ in Formula (4.1-1),
which is also a constant. Define a finite set with $\sigma m$ pairs:
$$\begin{array}{c}\mathcal S=\{(\varrho, \zeta)\,|\,1\leqslant \varrho\leqslant \sigma,
\zeta=1,\cdots, m\},\\
\mbox{ordered by}\quad(\varrho^1, \zeta^1)\prec (\varrho^2, \zeta^2)
\Longleftrightarrow \varrho_1>\varrho_2, \quad \mbox{or}\quad
\varrho^1=\varrho^2,\,\,\zeta^1<\zeta^2.\end{array}$$
In order to unify notations, the induced minimal wild local
pair $(\mf A^{r+t},\mf B^{r+t})$ in $(\bar\ast')$ is denoted by
$(\hat {\mf{A}},\hat {\mf{B}})$.
Let $(\varrho,\zeta)=(\bar p,l)\in \mathcal S$, where $\bar p$ is
the row-index of the first arrow $\hat a_1=a_1^t$ in $\bar\Theta^t$; $l$ is the subscript of the edge $d_l$, from which $a_1^t$
is split, since $F^t+\bar\Theta^t$ is contained in the $h$-th ladder by Lemma 5.5.4. Similarly, let
$(\varrho^1,\zeta^1)\in \mathcal S$ be determined by the first arrow
$\hat{a}^1_1$ of $\hat{\mf B}^1$. Theorem 5.4.3 (iv) ensures taht
$(\varrho, \zeta)\prec (\varrho^1, \zeta^1)$.

Now the procedure (III) is started once again from the pair $(\hat{\mf A}^{1},\hat{\mf B}^{1})$
instead of $(\hat{\mf A},\hat{\mf B})$. If (III-4) appears repeatedly, then after a finite
number of steps, an induced pair of (III-1)--(III-3) is reached
by induction on $\mathcal S$. \hfill$\Box$

\bigskip
\bigskip

{\large\bf 5.6 Proof of the main theorem}

\bigskip

It is ready to prove Theorem 3 given in the introduction.

\medskip

{\bf  Theorem 5.6.1} Let $(\mf A,\mf B)$ be a pair with $\T$ trivial, such that
$\mf{A}=(R,\K,\M,H=0)$ is a bipartite matrix bimodule problem satisfying RDCC
condition.  If $\mf B$ is of wild type, then $\mf{B}$ is not homogeneous.

\smallskip

{\bf Proof}\, There exists an induced bocs $\mf{B}'$ of the wild bocs $\mf B$,
which is in one of the cases of MW1--MW5 according to
Classification 3.3.2. Proposition 3.4.1--3.4.4 proved that if $\mf B'$ is
in the case of MW1--MW4, then $\mf B$ is not homogeneous. When $\mf A$ is bipartite
and satisfies RDCC condition, Proposition 5.3.4 and 5.5.6 proved that if the induced
bocs $\mf B'$ is in the case of MW5 of Remark 3.4.6, then $\mf B$ is not homogeneous. \hfill$\Box$

\medskip

{\bf Proof of Theorem 3}\,  Let $\Lambda$ be a finite-dimensional basic
algebra over an algebraically closed field $k$. We claim that if $\Lambda$ is of
wild representation type, then mod-$\Lambda$ is not homogeneous.

In fact, let $\mf{A}$ be the matrix bimodule problem
associated with $\Lambda$. Then $\mf{A}$ is bipartite and satisfies RDCC condition by Remark
1.4.4, and it is representation wild type. Therefore the associated bocs $\mf{B}$ is not
homogeneous by Theorem 5.6.1. Note that there is a
one-to-one correspondence between the set of equivalent classes of almost split sequences in
\mod-$\Lambda$ and that of almost split conflations in $R(\mf{B})$ except finitely many equivalent classes
of such sequences, see
\cite{B2} and \cite{ZZ}. Therefore mod-$\Lambda$ is not
homogeneous. \hfill$\Box$

\bigskip
\bigskip
\bigskip

{\large\bf Acknowledgement}

\bigskip

We would like to express our sincere thanks to S.Liu for his
proposal of $\Delta$-algebra, to Y.Han for his
suggestion of the concept on co-bimodule problem. We thank K.Wang,
F.Dai and X.Tang for correcting the English.
Y.Zhang is indebted to R.Bautista for giving the open problem in 1991; and
to W.W.Crawley-Boevey, D.Simson, C.M.Ringel, T.Lei, S.Wang, K.Li, C.Zhao for some discussions
during the long procedure of solving the problem. In particular,
she is grateful to Y.A.Drozd for his invitation to visit
the Institute of Mathematics, National Academy of Sciences of
Ukraine, and helpful conversations.

\end{document}